\documentclass{amsart}
\usepackage{amssymb, amsmath, latexsym, amscd, graphicx, tabularx, color}
\usepackage{rotating}
\usepackage[all]{xy}
\usepackage[dvipdfm]{hyperref}
\usepackage[all]{hypcap}
\usepackage[toc,page]{appendix}
\newtheorem{theorem}{Theorem}[section]
\newtheorem{lemma}[theorem]{Lemma}
\newtheorem{cor}[theorem]{Corollary}
\newtheorem{definition}[theorem]{Definition}
\newtheorem{proposition}[theorem]{Proposition}
\newtheorem{remark}[theorem]{Remark}
\newtheorem{example}[theorem]{Example}

\def\pagenumber{1}
\begin{document}
\setcounter{page}{\pagenumber}
\newcommand{\T}{\mathbb{T}}
\newcommand{\R}{\mathbb{R}}
\newcommand{\Q}{\mathbb{Q}}
\newcommand{\N}{\mathbb{N}}
\newcommand{\Z}{\mathbb{Z}}
\newcommand{\tx}[1]{\quad\mbox{#1}\quad}
\parindent=0pt
\def\SRA{\hskip 2pt\hbox{$\joinrel\mathrel\circ\joinrel\to$}}
\def\tbox{\hskip 1pt\frame{\vbox{\vbox{\hbox{\boldmath$\scriptstyle\times$}}}}\hskip 2pt}
\def\circvert{\vbox{\hbox to 8.9pt{$\mid$\hskip -3.6pt $\circ$}}}
\def\IM{\hbox{\rm im}\hskip 2pt}
\def\COIM{\hbox{\rm coim}\hskip 2pt}
\def\COKER{\hbox{\rm coker}\hskip 2pt}
\def\TR{\hbox{\rm tr}\hskip 2pt}
\def\GRAD{\hbox{\rm grad}\hskip 2pt}
\def\RANK{\hbox{\rm rank}\hskip 2pt}
\def\MOD{\hbox{\rm mod}\hskip 2pt}
\def\DEN{\hbox{\rm den}\hskip 2pt}
\def\DEG{\hbox{\rm deg}\hskip 2pt}

\title[Strong reactions in quantum super PDE's.I]{STRONG REACTIONS IN QUANTUM SUPER PDE's.I:\\ QUANTUM HYPERCOMPLEX EXOTIC SUPER PDE's}

\author{Agostino Pr\'astaro}
\maketitle
\vspace{-.5cm}
{\footnotesize
\begin{center}
Department SBAI - Mathematics, University of Rome La Sapienza, Via A.Scarpa 16,
00161 Rome, Italy. \\
E-mail: {\tt agostino.prastaro@uniroma1.it}
\end{center}
}

\vskip 0.5cm
\centerline{\em This work in three parts is dedicated to Albert Einstein and Max Planck.}
\vskip 0.5cm

\begin{abstract}
In order to encode strong reactions of the high energy physics, by means of nonlinear quantum propagators in the Pr\'astaro's geometric theory of quantum super PDE's, some related geometric structures are further developed and characterized. In particular {\em super-bundles of geometric objects} in the category $\boldsymbol{\mathfrak{Q}_S}$ of quantum supermanifolds are considered and quantum Lie derivative of sections of super bundle of geometric objects are calculated. Quantum supermanifolds with classic limit are classified with respect to the holonomy groups of these last commutative manifolds. A theorem characterizing quantum super manifolds with structured classic limit as super bundles of geometric objects is obtained. A theorem on the characterization of $chi$-flow on suitable quantum manifolds is proved. This solves a previous conjecture too. {\em Quantum instantons} and {\em quantum solitons} are defined as useful generalizations of the previous ones, well-known in the literature. Quantum conservation laws for quantum super PDEs are characterized. Quantum conservation laws are proved work for evaporating quantum black holes too. Characterization of observed nonlinear quantum propagators, in the observed quantum super Yang-Mills PDE, by means of conservation laws and observed energy is obtained. Some previous results by A. Pr\'astaro about generalized Poincar\'e conjecture and quantum exotic spheres, are extended to the category $\mathfrak{Q}_{hyper,S}$ of hypercomplex quantum supermanifolds. (This is the first part of a work divided in three parts. For part II and III see \cite{PRAS28, PRAS29}.)
\end{abstract}

\vskip 0.5cm


\vspace{.08in} \noindent \textbf{Keywords}: Integral (co)bordism groups in quantum (super) PDE's;
Existence of local and global solutions in hypercomplex quantum super PDE's; Conservation laws;
Crystallographic groups; Quantum exotic superspheres; Quantum exotic super PDEs.

\section[Introduction]{\bf Introduction}

\vskip 0.5cm
The algebraic topologic theory of quantum (super) PDE's formulated by A. Pr\'astaro, allows to directly encode quantum phenomena in a category of noncommutative manifolds (quantum (super)manifolds) and to finally solve the problem of unification, at quantum level, of gravity with the fundamental forces \cite{PRAS5,PRAS8,PRAS9,PRAS10,PRAS11,PRAS12,PRAS14,PRAS17,PRAS19,PRAS21,PRAS22}. In particular, this theory allowed to recognize the mechanism of mass creation/distruction, as a natural geometric phenomenon related to the algebraic topologic structure of quantum (super) PDEs encoding the quantum system under study \cite{PRAS22}.

Let us emphasize that with the algebraic topology of quantum (super) PDE's, formulated by A. Pr\'astaro, we can also go beyond the {\em paradox of the wave-particle duality}, and the probabilistic interpretation ({\em Copenhagen interpretation}) of the wave function, introduced in quantum mechanics by De Broglie \cite{BROGLIE1, BROGLIE2}, that even if accepted for practical reasons, was conceptually non well convincing. (Let us recall the Einstein's slogan, {\em ``God does not play dice ..."}, claimed in the famous Fifth Solvay International Conference, 1927.) Really in the geometric formulation of quantum PDEs, we solve the wave-particle duality and obtain a purely geometric formulation of quantum phenomena. In a sense the situation is similar to the paradoxical {\em ether}, as the medium of propagation of electromagnetic radiation, wrongly assumed as a necessary support in the Maxwell's dynamics of electromagnetism. In fact, nowadays we can encode quantum particles, and interactions between them, uniquely as geometric objects ($p$-chains), solutions of suitable quantum (super) PDEs, i.e., PDEs in the category $\mathfrak{Q}_S$ of quantum supermanifolds, as introduced by A. Pr\'astaro. Really it is the same concept of ``fundamental particles", related to the belief that there are some {\em fundamental bricks} building all Universe, in a LEGO blocks game sort, that the noncommutative and nonlinear theory of quantum super PDEs proves to be completely unjustified.\footnote{Nowadays we have also experimental evidences that such fundamental bricks do not exist. What it is considered ``fundamental" can be broken, when enough energy is available and suitable apparatus are arranged for. (See, e.g., electron decay into spinon and orbiton \cite{SCHLAPPA-SCHMITT}, or electrons with fractional electric charges \cite{FRADKIN}, and LHC experiments producing new quasiparticles by breaking proton.)} From this point of view it is not difficult to understand that also the concept of {\em confinement} reserved to quarks can be ``broken" under suitable energy conditions, by generating, e.g., meson decays into quarks, antiquarks and gluons, or baryon decays into quarks and gluons.\footnote{A quark-gluon plasma state, where quarks and gluons are not more confined, has been recently observed in some experiments. (\href{http://press.web.cern.ch/press/pressreleases/releases2010/PR23.10E.html}{"LHC experiments bring new insight into primordial universe" (Press release). CERN. 26 November 2010.})} With this respect, it is also clear that the Gell-Mann's {\em standard model} for hadrons, cannot be considered the last frontier in High Energy Physics !

The geometric theory of quantum (super) PDE's allows us to recover also the positive aspects of classical instantons. In fact, we can define {\em quantum super-instantons}, singular solutions of the quantum super Yang-Mills equations, that being localized in the quantum Minkowskian spacetime $M$, justify their names (similarly to classical instantons). Furthermore, singular solutions can encode tunneling effects, that, as it is well known, is a phenomenon that can be associated to classical instantons. Furthermore quantum particles can be seen as {\em quantum super-solitons} in the framework of solutions of quantum super Yang-Mills equations.

Aim of this work, divided in three parts, is to show how nuclear and subnuclear reactions can be encoded as boundary value problems in the algebraic topology of quantum (super) PDE's, as formulated by A. Pr\'astaro. With this respect, {\em nonlinear quantum propagators} play a fundamental role, beside to their representations by means of elementary reactions. This aspect is related to a decomposition theorem for nonlinear quantum propagators that we prove for solutions of boundary value problems in quantum super PDEs, and that will be developed in the second and third parts. (Part II and Part III are quoted in \cite{PRAS28, PRAS29}.)\footnote{For complementary information on High Energy Physics, related to the subject considered in this work, see e.g., Refs. \cite{ARKANI-FINK-SLAT-WEI,BLANKE-GOLD,CHEW-FRAUTSCHI,CUI,DIAKONOV-PETROV,DIRAC,EDEN,FEYNMAN,FRADKIN,GIAMMARCHI,GRIBOV,
GRIBOV-PONTECORVO,KAISER,MAJORANA,MANDELS-TAMM,POLYAKOV,PONTECORVO,PRADHAN,REGGE,VENEZIANO}.} This first part is instead devoted to some important preparatory subjects on the geometry of quantum super PDEs that will be utilized in part II and part III to encode quantum strong reaction dynamics. In particular, we will generalize to the category $\mathfrak{Q}_S$ of quantum supermanifolds, the {\em structure of super-bundle of geometric objects}, first introduced by A. Pr\'astaro, in the framework of category of commutative manifolds in \cite{PRAS1,PRAS2}. This structure gives a general framework where characterize physical fields in a fully covariant way. In particular it applies to PDEs defined on some higher order $G$-structures.

Another subject that is developed in this first part and that is very important to encode strong reactions, is the characterization of conservation laws for PDE's built in the category ${\frak Q}_S$. These are functions defined on the
integral bordism groups of such equations and belonging to suitable Hopf algebras ({\em full quantum Hopf algebras}). In
particular, we specialize our calculations on the quantum super Yang-Mills equations, quantum black holes as particular solutions there, and the characterization of observed nonlinear quantum propagators by means of conservation laws and observed quantum energy.

The last subject considered in this first part, is devoted to PDE's in the category $\mathfrak{Q}_{hyper,S}$ of quantum hypercomplex supermanifolds, as defined in \cite{PRAS27}, and will focus our attention on {\em quantum exotic super PDE's}, i.e., quantum super PDE's where we can embed {\em quantum exotic super-spheres}. For such Cauchy data we will generalize our previous results on quantum exotic PDE's \cite{PRAS27}. Exotic boundary value problems are of particular interest in strong reactions encoding quantum processes occurring in high energy physics, as we will show in the second and third parts of this work.

In the following we show as it is organized the paper and list the main results. 2. The concept of {\em super bundle of geometric objects} (see \cite{PRAS1,PRAS2}) is generalized to the category $\mathfrak{Q}_S$ for quantum super manifolds. This is necessary in order to obtain fully covariant structures, as requested in generalized mathematical structures that aim encode physical structures. Theorem \ref{holonomic-superbundle-of-geometric-objects-theorem} that characterizes PDEs as superbundles of geometric objects. Theorem \ref{classification-of-quantum-supermanifolds-with-classic-limit} classifies quantum supermanifolds, with classic limit, with respect to the holonomy groups of these last commutative manifolds. Theorem \ref{quantum super-manifolds-with-structured-classic-limit-as-super-bundles-of-geometric-objects} characterizes quantum super manifolds with structured classic limit as super bundles of geometric objects.
Theorem \ref{theorem-quantum-chi-flow}. Here we define {\em quantum $\hat Q^7$} any quantum (super)manifold $M$ having as classic limit $M_C=Q^7$. (This is the spinor bundle $Q^7\to S^3$, with $\mathbb{ C}^2$ fibers.) We say {\em quantum $\chi$-flow} on $M$, any quantum flow that projects on the classic limit in a $\chi$-flow \cite{AKBULUT-SALUR}. Then we have that a quantum $\chi$-flow represents diffeomorphically any homotopy $3$-sphere, $\Sigma^3\subset Q^7$, onto $S^3\subset Q^7$. This also solves a previous conjecture formulated in \cite{AKBULUT-SALUR} for commutative manifolds. Definition \ref{quantum-super-instantons-solitons} gives useful generalizations, in the category $\mathfrak{Q}_S$, of instantons and solitons. 3. Previous Pr\'astaro's results on the quantum conservation laws of quantum super PDE's are resumed. These are necessary in order to obtain the new results on quantum strong reactions given in part II and part III. Theorem \ref{quantum-tunnel-effects-and-quantum-black-holes} characterizes quantum black-hole dynamics by means of quantum integral characteristic supernumbers: these are conserved through a non-weak quantum evaporating black-hole. Theorem \ref{quantum-observed-hamiltonian-as-quantum-conservation-law-of-observed-quantum-yang-mills-pde} and Corollary \ref{conserved-quantum-observed-hamiltonian-in-quantum observed-nonlinear-propagator} give a precise meaning to the phenomenological concept of energy conservation during an observed quantum process.
4. Here we explicitly extend to the category $\mathfrak{Q}_{hyper,S}$ of hypercomplex quantum supermanifolds, previous results, in the category $\mathfrak{Q}_{hyper}$ of hypercomplex quantum manifolds, given in \cite{PRAS27} and in the category of quantum supermanifold \cite{PRAS14,PRAS19}. In particular, Theorem \ref{generalized-poincare-conjecture-in-category-hypercomplex-quantum-super-manifolds} explicitly extends to the category $\mathfrak{Q}_{hyper,S}$ a previous result about the generalized Poincar\'e conjecture in the category $\mathfrak{Q}_{S}$. (See \cite{PRAS19}.) Theorem \ref{theorem-quantum-superspheres-classic-limits-b} classifies diffeomorphic classes of hypercomplex quantum homotopy superspheres. Theorem \ref{bordism-groups-in-quantum-hypercomplex-exotic-super-pdes-stability} classifies integral bordism groups in quantum hypercomplex exotic super PDEs. Theorem \ref{integral-h-cobordism-in-Ricci-flow-super-pdes} characterizes integral h-cobordism in quantum hypercomplex Ricci flow super PDEs.

\section{\bf Super-bundles of geometric objects in the category $\boldsymbol{\mathfrak{Q}_S}$}

The quantum fundamental fields of physics, i.e., electromagnetic, gravitational and
nuclear fields, all must obey the full covariance requirement that can be codified generalizing to the category $\mathfrak{Q}_S$ the {\em structure of super-bundle of geometric objects}, first introduced by A. Pr\'astaro, in the framework of category of commutative manifolds in \cite{PRAS1,PRAS2}.\footnote{See also Pr\'astaro's algebraic topology of PDEs \cite{PRAS3,PRAS4,PRAS5,PRAS6,PRAS7,PRAS11,PRAS13,PRAS15,PRAS18,PRAS20,PRAS23,PRAS24,PRAS25,PRAS26} in order to better understand its generalization to the category of quantum supermanifolds. Interesting applications can be found in \cite{AG-PRA1,AG-PRA2,LYCH-PRAS,PRA-RAS,PRA-REGGE}. For basic information on the geometry of PDE's see \cite{B-C-G-G-G, DUB-FOM-NOV,GOLDSCHMIDT,GROMOV,KRAS-LYCH-VIN}. For basic information on algebraic topology see \cite{HIRSCH,MCCLEARY,MILNOR1,MILNOR2,MILNOR-MOORE,NASH,PONTRJAGIN,SCHOEN-YAU,SMALE,STONG,SWITZER,WALL1,WALL2,WARNER,WHITNEY1,WHITNEY2,YAU}.}

\begin{remark}[The concept of full covariance in $\mathfrak{Q}_S$]\label{full-covariance}
Similarly to what happens in the classical field theory, the concept of {\em full covariance} is fundamental in any quantum field theory. A quantum field is considered a section $s$ of a suitable
fiber bundle $\pi :W\rightarrow M$ in the category $\mathfrak{Q}_S$. To say that $s$ is fully covariant it means
that for any local diffeomorphism $\phi$ of the base $M$
we can calculate the pull-back $\phi ^*s$ of $s$ by means of $\phi$:
$$\phi ^*s\equiv \mathbb{B}(\phi )\circ s\circ\phi ^{-1} $$
where $\mathbb{B}(\phi )$ is a local application on $W$, canonically
associated to $\phi$ on $M$, such that the following diagram is commutative:
$$\xymatrix@C=3cm{W\vert _{U}\ar[d]_{\pi}\ar[r]^{\mathbb{B}(\phi )}&
                 W\vert _{\overline{U}}\ar[d]^{\pi}\\
                 U\ar[r]_{\phi }&\overline{U}}$$
If this circumstance is verified, we say that $\pi :W\rightarrow
M$ is a {\em natural fiber bundle}\label{natural-fiber-bundle} (or {\em fiber bundle of
geometric objects})\label{fiber-bundle-of-geometric-objects}, in the category $\mathfrak{Q}_S$.  We write: $( \pi :W\rightarrow M , {\mathbb{B}} )$.

\begin{table}[t]
\caption{Examples of super fiber bundles of geometric objects in the category $\mathfrak{Q}_S$}
\label{examples-super-fiber-bundles-geometric-objects}
\scalebox{0.8}{$\begin{tabular}{|l|l|}\hline
\multicolumn{2}{|c|}{\footnotesize{\rm Fiber bundle of geometric objects}}\\
\multicolumn{2}{|c|}{\footnotesize{\rm derivated from a fiber bundle of geometric objects $(\pi:W\to M; \mathbb{B})$}}\\
\hline
\hfil{\rm{\footnotesize Name}}\hfil&\hfil{\rm{\footnotesize Definition}}\hfil\\
\hline\hline
{\rm{\footnotesize Tangent bundle}}&{\rm{\footnotesize $( \pi :TM\rightarrow M , {\mathbb{B}}\equiv T )$}}\\
{\rm{\footnotesize Cotangent bundle}}&{\rm{\footnotesize $ ( \pi :T^*M\rightarrow M , {\mathbb{B}}\equiv T^* )$}}\\
{\rm{\footnotesize Bundle of tensors}}&{\rm{\footnotesize $( \pi :T^r_sM\rightarrow M ,
                  {\mathbb{B}}\equiv T^r_s )$}}\\
{\rm{\footnotesize Full quantum tangent bundle}}&{\rm{\footnotesize $( \pi :\widehat{T}M\rightarrow M , {\mathbb{B}}\equiv \widehat{T} )$}}\\
{\rm{\footnotesize Full quantum cotangent bundle}}&{\rm{\footnotesize $ ( \pi :T^+M\rightarrow M , {\mathbb{B}}\equiv T^+ )$}}\\
{\rm{\footnotesize Dot bundle of tensors}}&{\rm{\footnotesize $( \pi :\dot T^r_sM\rightarrow M ,
                  {\mathbb{B}}\equiv \dot T^r_s )$}}\\
{\rm{\footnotesize Bundle of derivative of sections}}&{\rm{\footnotesize $( \pi _1:J{\it D}(W)\rightarrow M;\, \mathbb{B}^{(1)})$}}\\
&{\rm{\footnotesize $\mathbb{B}^{(1)}\equiv J{\it D}(-)\equiv
                    Hom_{\mathbb{R}}(T(-);T\mathbb{B}(-) )|_{J{\it D}(W)}$}}\\
{\rm{\footnotesize pull-back of $Ds,\, \pi_1\circ Ds=id_M$}}&{\rm{\footnotesize $\phi^*Ds=Hom_{\mathbb{R}}(T(\phi);T\mathbb{B}(\phi^{-1}))\circ Ds\circ\phi$}}\\
\hline
\end{tabular}$}

\scalebox{0.8}{$\begin{tabular}{|l|l|}\hline
\multicolumn{2}{|c|}{\footnotesize{\rm Super fiber bundle of geometric objects}}\\
\multicolumn{2}{|c|}{\footnotesize{\rm derived from a super bundle of geometric objects $(\pi_B:B\to M\leftarrow W:\pi_W; \mathbb{B})$}}\\
\hline
\hfil{\rm{\footnotesize Name}}\hfil&\hfil{\rm{\footnotesize Definition}}\hfil\\
\hline\hline
{\rm{\footnotesize Quantum bundle of derivative of sections}}&{\rm{\footnotesize $( \pi_B:B\to M\leftarrow J{\it \hat D}(W);\pi _1;\, \mathbb{B}^{(1)})$}}\\
&{\rm{\footnotesize ${\mathbb{B}}^{(1)}\equiv J{\it \hat D}(-)\equiv
                    Hom_{Z}(T(-);T\mathbb{B}(-) )|_{J\hat{\it D}(W)}$}}\\
{\rm{\footnotesize pull-back of $Ds,\, \pi_1\circ Ds=id_M$}}&{\rm{\footnotesize $\phi^*Ds=Hom_Z(T(\phi);T\mathbb{B}(\phi^{-1}))\circ Ds\circ\phi$}}\\
\hline
\end{tabular}$}
\end{table}

A generalization of this concept is that of {\em super-fiber bundle of
geometric objects}\index{super-fiber bundle of geometric objects} introduced in \cite{PRAS1, PRAS2}.
A quantum physical field is {\em fully covariant}
iff it is a {\em  quantum superfield}\label{quantum superfield}, that is a section of such
a structure.
\end{remark}

\begin{definition}\label{structure of superbundle of geometric objects}
A {\em structure of superbundle of geometric objects} in the category $\mathfrak{Q}_S$, is given by
two fiber bundles $W\mathop{\to}\limits^{\pi _{W}}M\mathop{\leftarrow}
\limits^{  \pi _{B}}B$
over the same base $M$ and a covariant functor
${\mathbb{B}}:\mathcal{C}(B)\rightarrow \mathcal{C}(W) $, where
$\mathcal{C}(B) $, (resp. $\mathcal{C}(W)$) is the category whose objects
are open subbundles of $B$ (respectively,
open subbundles of $W$) and whose morphisms are the
local fiber bundle authomorphisms between those objects such that:

{\em i)} if $B\vert_ U\in Ob(\mathcal{C}(B)) \Rightarrow  {\mathbb{B}}(B\vert_ U)
={\pi ^{-1}}_{W}(U)\in  Ob(\mathcal{C}(W)) $;

{\em ii)} if $f\in  Hom(\mathcal{C}(B))$ with
$f\equiv (f_{B} ,f_{M}):B\vert U\rightarrow  B\vert U'$, then
${\mathbb{B}}(f) \in Hom(\mathcal{C}(W))$ and satisfies:

{\em iii)} $\pi _{W} \circ {\mathbb{B}}(f) =f_{M} \circ \pi _{W}$;

{\em iv)} if $B\vert U \in  Ob(\mathcal{C}(B)) $, $\overline{U}\subset
 U \Rightarrow  {\mathbb{B}}(f)\vert_ {\pi_{W}}^{-1}(\overline{U}) =
{\mathbb{B}}(f\vert_ B\vert_{ \overline{U}})$.

$W$ is called the {\em total bundle} and $B$
the {\em base bundle}. A section of $\pi _{W}$  is called a
{\em quantum superfield of geometric objects}.

The situation is resumed in the commutative diagram in {\em(\ref{commutative-diagram-super-bundle-geometric-objects})}.

\begin{equation}\label{commutative-diagram-super-bundle-geometric-objects}
 \xymatrix@C=60pt{B\vert _{U}\ar[dr]_{\pi_B}\ar[r]^{f_B}&
                 B\vert _{\overline{U}}\ar[dr]^(0.2){\pi_B}& W\vert _{U}\ar[dl]_(0.2){\pi_W}\ar[r]^{\mathbb{B}(f)}&
                 W\vert _{\overline{U}}\ar[dl]^{\pi_W}\\
                & U\ar[r]_{f_M }&\overline{U}&}
\end{equation}

\end{definition}

\begin{example}[$G$-structures in the category $\mathfrak{Q}_S$]\label{g-structures}
$G$-structures in the category $\mathfrak{Q}_S$ are reductions of principal bundles
of $(r|s)$-frames, on a $(m|n)$-dimensional quantum supermanifold, natural generalizations of analogous structures on commutative manifolds \cite{PRAS3}.
\end{example}

\begin{example}
In Tab. \ref{examples-super-fiber-bundles-geometric-objects} are reported some further distinguished examples of (super)-bundles of geometric objects in the category $\mathfrak{Q}_S$.
\end{example}
In the physical applications it is important to consider the following.

\begin{theorem}
Let $(P,M,\pi ;G)$ be a principal
fiber bundle with structure group $G$ in the category $\mathfrak{Q}_S$.
Let $W\equiv P\times F/G$ be a fiber bundle associated to $P$ with fibre $F$.
Then, there exists a canonical covariant functor ${\mathbb{B}}$
such that $(P,W;{\mathbb{B}})$ is a superbundle of geometric objects if we
restrict the category $\mathcal{C}(P)$, where ${\mathbb{B}}$ is defined, to the
subcategory $\mathcal{C}(P)_\bullet$, where $Hom_{\mathcal{C}(P)_\bullet}(P\vert U, P\vert U')$ is the
set of fibered quantum diffeomorphisms $P\vert U\rightarrow  P\vert U'$, such that
the isomorphism on the structure group is the identity.
\end{theorem}

\begin{proof}  The proof can be copied by the classical case \cite{PRAS5}.\footnote{Warn ! Do not consider the concept of super bundle of geometric objects as synonym of fiber bundle in the category $\mathfrak{C}_S$ of supermanifolds. (For a geometric theory of PDEs in $\mathfrak{C}_S$, see \cite{PRAS4}.)} \end{proof}

\begin{definition}[Quantum Lie derivative of section of super bundle of geometric objects]\label{quantum-lie-derivative}
Let $(\phi_{M,\lambda},\phi_{B,\lambda})_{\lambda\in \mathbb{R}}$, be a (local) $1$-parameter group of (local) fiber bundle diffeomorphisms in the category $\mathfrak{Q}_S$. Let $(\xi,\zeta)$ the corresponding couple of quantum vector fields such that the fiber bundle {\em(\ref{commutative-diagram-infinitesimal-generators})} is commutative.

\begin{equation}\label{commutative-diagram-infinitesimal-generators}
  \xymatrix@C=3cm{B|_U\ar[d]_{\pi_B}\ar[r]^{\zeta}&TB|_U\ar[d]^{T(\pi_B)}\\
  U\ar[r]_{\xi}&TU\\}
  \end{equation}
  We call {\em quantum Lie derivative} of a section $s$ of $\pi_W$, with respect to $(\phi_{M,\lambda},\phi_{B,\lambda})_{\lambda\in \mathbb{R}}$, or $(\xi,\zeta)$, the infinitesimal variation , $\partial\widetilde{s}$, of the pull-back $\widetilde{s}_\lambda\equiv\phi_\lambda^*s=\mathbb{B}(\phi_\lambda^{-1})\circ s\circ\phi_{M,\lambda}$, with $\widetilde{s}:\mathbb{R}\times M\to W$. (Here, for sake of simplicity we have denoted $\widetilde{s}$ globally define, but in general it is only locally defined.) Therefore we get
\begin{equation}\label{quantum-lie-derivative}
  \mathcal{L}_{\zeta}s\equiv\partial\widetilde{s}=\frac{d}{d\lambda}(\widetilde{s}_\lambda)|_{\lambda=0}:M\to s^*vTW.
\end{equation}
\end{definition}

\begin{definition}[Full quantum Lie derivative of section of super bundle of geometric objects]\label{full-quantum-lie-derivative}
Let $(\phi_{M,\lambda},\phi_{B,\lambda})_{\lambda\in A}$, be a (local) $1$-parameter group of (local) fiber bundle diffeomorphisms in the category $\mathfrak{Q}_S$. Let $(\xi,\zeta)$ the corresponding couple of vector fields such that the fiber bundle {\em(\ref{commutative-diagram-infinitesimal-generators-a})} is commutative.

\begin{equation}\label{commutative-diagram-infinitesimal-generators-a}
  \xymatrix@C=3cm{B|_U\ar[d]_{\pi_B}\ar[r]^{\zeta}&\widehat{T}B|_U\equiv Hom_Z(A;TB|_U)\ar[d]^{\widehat{T}(\pi_B)\equiv Hom_Z(1_A;T(\pi_B))}\\
  U\ar[r]_{\xi}&\widehat{T}U\equiv Hom_Z(A;TU)\\}
\end{equation}
We call {\em full quantum Lie derivative} of a section $s$ of $\pi_W$, with respect to $(\phi_{M,\lambda},\phi_{B,\lambda})_{\lambda\in A}$, or $(\xi,\zeta)$, the infinitesimal variation , $\partial\widetilde{s}$, of the pull-back $\widetilde{s}_\lambda\equiv\phi_\lambda^*s=\mathbb{B}(\phi_\lambda^{-1})\circ s\circ\phi_{M,\lambda}$, with $\widetilde{s}:A\times M\to W$. (Here, for sake of simplicity we have denoted $\widetilde{s}$ globally define, but in general it is only locally defined.) Therefore we get
\begin{equation}\label{quantum-lie-derivative}
  \mathcal{L}_{\zeta}s\equiv\partial\widetilde{s}=\frac{d}{d\lambda}(\widetilde{s}_\lambda)|_{\lambda=0}:M\to s^*v\widehat{T}W\equiv Hom_Z(A;s^*vTW).
\end{equation}
\end{definition}

\begin{proposition}
Under the same hypotheses of Definition \ref{quantum-lie-derivative} (resp. Definition \ref{full-quantum-lie-derivative}) we get that whether $\pi_W:W\to M$ is a vector bundle, then we get
$\mathcal{L}_{\zeta}s\equiv\partial\widetilde{s}:M\to W$, (resp. $\mathcal{L}_{\zeta}s\equiv\partial\widetilde{s}:M\to \widehat{W}\equiv Hom_Z(A;W)$).
\end{proposition}

\begin{theorem}\label{holonomic-superbundle-of-geometric-objects-theorem}
Let $\hat E_k\subset J{\it\hat D}^k(W)$, be a quantum super PDE, with $(\pi_B:B\to M\leftarrow W:\pi_W;\mathbb{B})$ a super bundle of geometric objects. Let us assume that the map $\pi_{k,0}:\hat E_k\to W$ be surjective, with respect to the natural mapping $\pi_{k,0}:J{\it\hat D}^k(W)\to W$, and let assume also that $\hat E_k$ is formally integrable and completely integrable. Then for any $r$-prolongation $(\hat E_k)_{+r}$, $r\ge 0$, we get the structure of super bundle of geometric objects reported in {\em(\ref{holonomic-superbundle-of-geometric-objects})}.
\begin{equation}\label{holonomic-superbundle-of-geometric-objects}
   (\pi_B:B\to M\leftarrow (\hat E_k)_{+r}:\pi_{k+r};\mathbb{B}^{(k+r)}=J{\it\hat D}^{k+r}(-)_{\bullet}).
\end{equation}
where $J{\it\hat D}^{k+r}(-)_{\bullet}\equiv J{\it\hat D}^{k+r}(\mathbb{B}_{\bullet}(-))|_{\hat E_{k+r}}$, with $\mathbb{B}_{\bullet}(-)$ the restriction of $\mathbb{C}(B)(-)$ to the sub-category $\mathcal{C}_\bullet(B)\subset \mathcal{C}(B)$, such that the corresponding morphisms $f\in Hom(\mathcal{C}_\bullet(B))$ are such that $J{\it\hat D}^{k+r}(\mathbb{B}_\bullet(f))_{\bullet}$ are symmetries of $\hat E_{k+r}$.
We call {\em(\ref{holonomic-superbundle-of-geometric-objects})} a {\em $(k+r)$-holonomic super bundle of geometric objects}.
\end{theorem}

\begin{proof}
This means that for any section $u:M\to (\hat E_k)_{+r}$, we can calculate its quantum Lie derivative, with respect to a $1$-group of (local) fiber bundle diffeomorphisms of $\pi_B:B\to M$, that should be also symmetries of $\hat E_k$, hence of $(\hat E_k)_{+r}$ too. In particular for holonomic sections $u=D^{k+r}s$, namely $s$ is a solution of $\hat E_k$, we get $\partial\widetilde{D^{k+r}s}:M\to (D^{k+r}s)^*vTJ{\it\hat D}^{k+r}(W)\cong J{\it\hat D}^{k+r}(s^*vTW)$, such that the diagram (\ref{commutative-diagram-b}) is commutative.
\begin{equation}\label{commutative-diagram-b}
 \xymatrix@C=3cm{(D^{k+r}s)^*vTJ{\it\hat D}^{k+r}(W)\ar@{=}^(0.57){\thicksim}[r]  &J{\it\hat D}^{k+r}(s^*vTW)\ar@{=}[d]\\
 \hat E_{k+r}[s]\equiv(D^{k+r}s)^*vT(\hat E_k)_{+r}\ar@{^{(}->}[u]\ar@{^{(}->}[r]&J{\it\hat D}^{k+r}(E[s])\\}
\end{equation}

Therefore the quantum Lie derivative of a solution of $\hat E_k$ can be identified with a solution of the linearized equation $\hat E_k[s]\subset J{\it\hat D}^{k}(E[s])$, with $E[s]\equiv s^*vTW$.\footnote{Warn ! The set of all solutions of the linearized equation $\hat E_k[s]\equiv (D^{k}s)^*vT\hat E_k$ is larger than ones obtained as Lie derivative of $s$, with respect to some $1$-parameter symmetry group  of some super bundle of geometric objects $(\pi_B:B\to M\leftarrow W:\pi_W;\mathbb{B})$.}

\end{proof}

\begin{theorem}[Classification of quantum supermanifolds with classic limit]\label{classification-of-quantum-supermanifolds-with-classic-limit}
We can classify quantum supermanifolds, having a classic limit, $\pi_C:M\to M_C$, on the ground of its classic limit $M_C$. In other words, we can say $M$ has some {\em property} whether its classic limit has this property. (For example we can say that $M$ is {\em orientable} when $M_C$ is so.) In Tab. \ref{classification-quantum-supermanifolds-holonomy-groups-classic limits} are reported some quantum supermanifolds classified with respect to the holonomic group of their classic limit.
\end{theorem}

\begin{proof}
This follows from standard properties of classifications of fiber bundles in algebraic topology, (see e.g., \cite{SWITZER}), and taking into account that each fiber over $p\in M_C$ is contractible to a point.
\end{proof}
\begin{example}
Quantum (super)manifolds with classic limits having nontrivial structure can be considered of particular interest whether from the mathematical point of view or physical one. In particular, recent works investigate on Calaby-Yau mirror-pairs of submanifolds of $G_2$-manifolds and $Spin(7)$-manifolds. (See e.g., \cite{AKBULUT-SALUR}.)\footnote{{\em Calabi-Yau manifolds} are compact, complex K\"ahler manifolds that have trivial first Chern classes (over $\mathbb{R}$). Yau proved (1977-1979) a conjecture by Calabi (1957) that there exists on every CY-manifold a K\"ahler metric with vanishing Ricci curvature \cite{YAU,CALABI}. One important class of $G_2$ manifold are the ones from $CY$-manifolds. Let $(X,\omega,\Omega)$ be a complex $3$-dimensional CY-manifold with K\"ahler form $\omega$ and a nowhere vanishing holomorphic $3$-form $\Omega$, then $X^6\times S^1$ has holonomy group $SU(3)\subset G_2$, hence is a $G_2$ manifold. In this case $\varphi=\Re\Omega+\omega\wedge dt$. Similarly, $X^6\times \mathbb{R}$ gives a noncompact $G_2$ manifold.} These are Riemannian manifolds classified on the ground of their holonomy groups. Let us recall in Tab. \ref{classification-quantum-supermanifolds-holonomy-groups-classic limits} the classification of possible holonomy groups $Hol(M_C,g)$ for simply connected Riemannian manifolds $(M_C,g)$ which are not locally a product space and not locally a Riemannian symmetric space (shortly denote them by INS-SCRM here). (This classification was first obtained by M. Berger (1955) and further improved by some other mathematicians (D. Alekseevski, Brown-Gray, Hitchin ...).  R.L. Bryant (1987) first proved existence of metrics with holonomy $G_2$ (resp. $Spin(7)$) on $7$-dimensional manifolds (resp. $8$-dimensional manifolds) \cite{BERGER, BRYANT, GRAY, HITCHIN}.

The group $G_2$ can be identified with the group of authomorphisms of octonions $\mathbb{O}$ or equivalently the subgroup of $GL(7,\mathbb{ R})$ that preserves a suitable $3$-form $\varphi_0\in\Omega^3(\mathbb{ R}^7)$, given in {\em(\ref{fundamental-form})}.
\begin{equation}\label{fundamental-form}
\left\{
\begin{array}{ll}
  \varphi_0&=\frac{1}{6} \epsilon_{ijk}dx^i\wedge dx^j\wedge dx^k\in\Omega^3(\mathbb{R}^7)\\
  &=dx^1\wedge dx^2\wedge dx^3+dx^1\wedge dx^4\wedge dx^5+dx^1\wedge dx^6\wedge dx^7+dx^2\wedge dx^4\wedge dx^6\\
  &-dx^2\wedge dx^5\wedge dx^7-dx^3\wedge dx^4\wedge dx^7-dx^3\wedge dx^5\wedge dx^6.\\
\end{array}
\right.
\end{equation}

 A smooth $7$-dimensional manifold $M_C$ has a {\it $G_2$-structure} if its tangent bundle reduces to a $G_2$ bundle, or equivalently if there is a $3$-form $\varphi\in\Omega^3(M_C)$ such that at each $p\in M_C$ the pair $(T_pM_C,\varphi(p))$ is isomorphic to $(T_0\mathbb{ R}^7,\varphi_0)$. Then the $3$-form $\varphi$ identifies an orientation $\mu\in\Omega^7(M_C)$, and $\mu$ determines a metric $g_{\varphi}$. A manifold with $G_2$-structure $(M_C,\varphi)$, is called a {\it $G_2$ manifold} if $Hol(M_C,g_{\varphi})\subset G_2$, or equivalently $\nabla_{g_{\varphi}}\varphi=0\, \Leftrightarrow\, \{d\varphi=0,\, d(\star_{g_{\varphi}}\varphi)=0\}$. (Because $G_2$ is a connected, simply connected group, a connected $7$-dimensional manifold with a $G_2$-structure is orientable and admits a spin structure, i.e., its first two Stiefel-Whitney classes vanish. (Gray (1969) and R. Bryant (2005).) $G_2$ manifolds can be also characterized as critical points of a suitable functional on the $3$-forms (N. Hitchin (2000).).

In \cite{AKBULUT-SALUR} it is defined {\it mirror pair} a couple of {\em Calabi-Yau manifolds} if their complex structures are induced from the same calibration $3$-form, $\varphi$, in a $G_2$ manifold. It assigns to a $G_2$ manifold $(M_C,\varphi,\Lambda)$, with the calibration $3$-form $\varphi$ and a oriented $2$-plane field $\Lambda$, a pair of parametrized tangent bundle valued $2$- and $3$-forms of $M_C$. These forms can be used to define different complex and symplectic structures on certain $6$-dimensional subbundles of $TM_C$. More precisely, let $(M_C,\varphi)$ be a $G_2$ manifold. A $4$-dimensional submanifold $X\subset M_C$ is called {coassociative} if $\varphi|_X=0$. A $3$-dimensional submanifold $Y\subset M_C$ is called {\em associative} if $\varphi|_Y=vol(Y)$. This condition is equivalent to the condition $\chi|_Y=0$, where $\chi\in\Omega^3(M,TM)\equiv C^\infty(TM\bigotimes\Lambda^0_3M_C)$ is the $TM_C$-valued $3$-form on $M_C$ given by $<\chi(u,v,w),z>=\star\varphi(u,v,w,z)$. We can also define a $TM$-valued $2$-form $\psi\in C^\infty(TM_C\bigotimes\Lambda^0_2M_C)$, given by $<\psi(u,v),w>=\varphi(u,v,w)=<u\times v,w>$. To any $3$-dimensional submanifold $Y\subset(M_C,\varphi)$, $\chi$ assigns a normal vector field, which vanishes when $Y$ is associative. For any associative manifold $Y\subset (M_C,\varphi)$ with a non-vanishing oriented $2$-plane field, $\chi$ defines a complex structure on its normal bundle. In particular, any coassociative submanifold $X\subset M_C$ has an almost complex structure if its normal bundle has a non vanishing section. Two CY-manifolds are {\em mirror pairs} of each other if their complex structures are induced from the same calibration $3$-form in a $G_2$ manifold. Furthermore, we call them {\em strong mirror pairs} if their normal vector fields $\xi$ and $\xi'$ are homotopic to each other through non-vanishing vector fields. One can assign to a $G_2$ manifold $(M_C,\varphi,\Lambda)$, with the calibration $3$-form $\varphi$ and oriented $2$-plane field $\Lambda$, a pair of parametrized tangent bundle valued $2$- and $5$-forms of $M_C$. These forms can be used to define differential complex and symplectic structures on certain $6$-dimensional subbundles of $TM_C$. When these bundles are integrated they give mirror CY-manifolds. In a similar way one can recognize mirror dual $G_2$ manifolds inside of a $Spin(7)$ manifold $(M_C^8,\psi)$. In case $M_C$ admits an oriented $3$-plane field, by iterating this process one can obtain Calabi-Yau submanifolds pairs in $M_C$ whose complex and symplectic structures determine each other via the calibration form of the ambient $G_2$ (or $Spin(7)$) manifold.
\begin{table}[h]
\caption{Classification of quantum supermanifolds on the ground of the Holonomy groups of their classic limit $\pi_C:M\to M_C$.}
\label{classification-quantum-supermanifolds-holonomy-groups-classic limits}
\begin{tabular}{|l|l|l|l|}
\hline
{\rm{\footnotesize $\dim_\mathbb{ R} M_C$}}&{\rm{\footnotesize $Hol(M_C,g)$}}&{\rm{\footnotesize Manifold's name}}&{\rm{\footnotesize Further characterizations}}\\
\hline\hline
{\rm{\footnotesize $n$}}&{\rm{\footnotesize $SO(n)$}}&{\rm{\footnotesize Orientable}}&{\rm{\footnotesize }}\\
\hline
{\rm{\footnotesize $2n$}}&{\rm{\footnotesize $U(n)$}}&{\rm{\footnotesize K\"ahler}}&{\rm{\footnotesize orientable:\hskip 2pt[$U(n)\subset SO(2n)$]}}\\
\hline
{\rm{\footnotesize $2n$}}&{\rm{\footnotesize $SU(n)$}}&{\rm{\footnotesize Calabi-Yau}}&{\rm{\footnotesize Ricci-flat, K\"ahler, orientable:}}\\
&&&{\rm{\footnotesize [$SU(n)\subset U(n)\subset SO(2n)$]}}\\
\hline
{\rm{\footnotesize $4n$}}&{\rm{\footnotesize $Sp(n).Sp(1)$}}&{\rm{\footnotesize Quaternionic-K\"ahler}}&{\rm{\footnotesize Einstein}}\\
\hline
{\rm{\footnotesize $4n$}}&{\rm{\footnotesize $Sp(n)$}}&{\rm{\footnotesize Hyperk\"ahler}}&{\rm{\footnotesize  Ricci-flat, Calabi-Yau, K\"ahler, orientable:}}\\
&&&{\rm{\footnotesize [$Sp(n)\subset SU(2n)\subset U(2n)\subset SO(4n)$]}}\\
\hline
{\rm{\footnotesize $7$}}&{\rm{\footnotesize $G_2$}}&{\rm{\footnotesize $G_2$-manifold}}&{\rm{\footnotesize Ricci-flat, spin, orientable:\hskip 2pt[$G_2\subset SO(7)$]}}\\
\hline
{\rm{\footnotesize $8$}}&{\rm{\footnotesize $Spin(7)$}}&{\rm{\footnotesize $Spin(7)$-manifold}}&{\rm{\footnotesize Ricci-flat}}\\
\hline
\multicolumn{4}{l}{\rm{\footnotesize A manifold with $G_2$ structure $(M_C,\varphi)$ is called a $G_2$ manifold if $Hol(M_C,g_{\varphi})\subset G_2$,}}\\
\multicolumn{4}{l}{\rm{\footnotesize or equivalently $\mathop{\nabla}\limits^{(g_{\varphi})}\varphi=0$.}}\\
\end{tabular}
\end{table}

\end{example}

\begin{theorem}[Quantum super manifolds with structured classic limit as super bundles of geometric objects]\label{quantum super-manifolds-with-structured-classic-limit-as-super-bundles-of-geometric-objects}
Let $M$ be a quantum super manifolds with structured classic limits, $\pi_C:M\to M_C$. A super bundle of geometric objects on $M$, $(\pi_C:B\to M\leftarrow W:\pi_W)$ identifies a super bundle of geometric objects on $M_C$ iff $Hom(\mathcal{C}(B))$ and $Hom(\mathcal{C}(W))$ are restricted to (local) fiber bundles morphisms and the functor $\mathbb{B}$ admits such a restriction, such that the diagrams in {\em(\ref{fiber-bundles-commutative-diagrams})} are commutative.

\begin{equation}\label{fiber-bundles-commutative-diagrams}
 \scalebox{0.8}{$\xymatrix@C=2cm{B\ar[d]_{\pi_B}\ar[r]^{f_B}&B\ar[d]^{\pi_B}\\
 M\ar[d]_{\pi_C}\ar[r]^{f_B}&M\ar[d]^{\pi_C}\\
 M_C\ar[r]_{f_C}&M_C\\}\hskip 0.5cm
\xymatrix@C=2cm{W\ar[d]_{\pi_W}\ar[r]^{f_W}&W\ar[d]^{\pi_W}\\
 M\ar[d]_{\pi_C}\ar[r]^{f_B}&M\ar[d]^{\pi_C}\\
 M_C\ar[r]_{f_C}&M_C\\}$}
 \end{equation}
Then, when above conditions are satisfied, we say that the super bundle of geometric objects on $M$, $(\pi_C:B\to M\leftarrow W:\pi_W)$, is a {\em classic-regular super bundle of geometric objects} on the quantum supermanifold $M$.
\end{theorem}

\begin{proof}
In fact any super bundle of geometric object on $M$, i.e., $(\pi_C:B\to M\leftarrow W:\pi_W;\mathbb{B})$, one has the natural fiber bundle structures on $M_C$ reported in (\ref{fiber-bundles-commutative-diagrams-a}).
\begin{equation}\label{fiber-bundles-commutative-diagrams-a}
 \xymatrix{B\ar@/_2pc/[ddr]_{\pi_{B,C}}\ar[dr]^{\pi_B}&&W\ar[ld]_{\pi_{W}}\ar@/^2pc/[ldd]^{\pi_{W,C}}\\
 &M\ar[d]_{\pi_C}&\\
 &M_C&\\}
 \end{equation}
 Therefore for any morphism $(f)=(f_B,f_M,f_{M_C})$ must be the diagram in (\ref{fiber-bundles-commutative-diagrams-b}) commutative.
 \begin{equation}\label{fiber-bundles-commutative-diagrams-b}
\xymatrix@C=60pt{B\ar@/_2pc/[ddr]_{\pi_{B,C}}\ar[dr]_{\pi_B}\ar[r]^{f_B}&
                 B\ar@/_1pc/[ddr]_{\pi_{B,C}}\ar[dr]^(0.2){\pi_B}& W\ar@/^1pc/[ldd]^{\pi_{W,C}}\ar[dl]_(0.2){\pi_W}\ar[r]^{\mathbb{B}(f)}&
                 W\ar[dl]^{\pi_W}\ar@/^2pc/[ldd]^{\pi_{W,C}}\\
                & M\ar[d]_{\pi_C}\ar[r]^{f_M }&M\ar[d]^{\pi_C}&\\
                & M_C\ar[r]_{f_C }&M_C&\\}
 \end{equation}
\end{proof}

\begin{example}
The {\em $(k+r)$-holonomic super bundle of geometric objects} $(\pi_B:B\to M\leftarrow (\hat E_k)_{+r}:\pi_{k+r};\mathbb{B}^{(k+r)}=J{\it\hat D}^{k+r}(-))$ is a classic-regular super bundle of geometric objects on the quantum supermanifold $M$.
\end{example}

\begin{remark}[$\chi$-flow] In \cite{AKBULUT-SALUR} it is considered the spinor bundle $Q^7\to S^3$ (with $\mathbb{ C}^2$ fibers). $Q$ is just a $G_2$ manifold where the two $6$-dimensional submanifolds $S^2\times\mathbb{ R}^4\subset Q$ and $S^3\times\mathbb{ R}^3\subset Q$ constitute a mirror pair. The zero section $S^3\subset Q$ is an {\em associative submanifold}, i.e., $\varphi|_{S^3}\equiv vol(S^3)$. This condition is equivalent to the condition $\chi|_{S^3}\equiv 0$, where $\chi\in C^\infty(TQ\bigotimes\Lambda^0_3(Q))$, is defined by $<\chi(u,v,w),z>=(\star_{g_{\varphi}}\varphi)(u,v,w,z)$. The equivalence of these conditions follows from $\varphi(u,v,w)^2+|\chi(u,v,w)|^2/4=|u\, \wedge\, v\, \wedge\, w|^2$. Then $\chi$ identifies a flow on a $3$-dimensional submanifold $f:Y\to (Q,\varphi)$, called {\em $\chi$-flow}, described by $(\partial t.f)=\chi(f_*vol(Y))$. (See also some previous works by R. L. Bryant \& M-S. Salamon (1989) and N. Hitchin (2000).) Furthermore, in \cite{AKBULUT-SALUR} the following conjecture is made.

\textbf{Conjecture.} {\em Since one can imbed any homotopy $3$-sphere $\Sigma^3$ into $Q$ (homotopic to the zero-section), one can conjecture that the $\chi$-flow on $\Sigma^3\subset Q$, takes $\Sigma^3$ diffeomorphically onto the zero section $S^3$.}

This conjecture is justified since $\Sigma^3\cong S^3$, as it is nowadays well proved. (See \cite{PRAS23, PRAS24, PRAS26}.) In fact, the PDE's algebraic topology as introduced by A. Pr\'astaro, to characterize global solutions of PDE's, can be used also to solve this conjecture when applied to the $\chi$-flow equation.
\end{remark}

More explicitly we have the following lemma.

\begin{lemma}\label{lemma-chi-flow}
The chi-flow equation admits a solution that diffeomorphically relates any $3$-dimensional homotopy sphere $\Sigma^3$, with $S^3$.
\end{lemma}

\begin{proof}
Let us consider the following fiber bundle $\pi:W\equiv\mathbb{R}\times Q^7\to\mathbb{R}\times S^3$, $(t,x^k)_{1\le k\le 7}\mapsto(t,x^k)_{1\le k\le 3}$. Here the fiber is $\mathbb{C}^2\cong\mathbb{R}^4$. One has the canonical embedding of $S^3$ into $W$, for any $t\in\mathbb{R}$: $(t,x^k)_{1\le k\le 3}\to(t,x^1,x^2,x^3,0,\cdots,0)$, identified by the zero section of $\pi$. In (\ref{local-representation-jeometric-objects}) are given the local representations of the geometric objects above introduced.

\begin{equation}\label{local-representation-jeometric-objects}
\scalebox{0.8}{$\left\{
\begin{array}{ll}
  \varphi&=\frac{1}{6}\epsilon_{ijk}dx^i\wedge dx^j\wedge dx^k \\
  &= dx^1\wedge dx^2\wedge dx^3+dx^1\wedge dx^4\wedge dx^5+dx^1\wedge dx^6\wedge dx^7\\
  &+dx^2\wedge dx^4\wedge dx^6-dx^2\wedge dx^5\wedge dx^7-dx^3\wedge dx^4\wedge dx^7-dx^3\wedge dx^5\wedge dx^6\\
 g_{\varphi}&=\delta_{ij}dx^i\otimes dx^j\\
 \mu&=dx^1\wedge\cdots\wedge dx^7\\
 \star_{\varphi}\varphi=&\frac{1}{24}\epsilon_{ijkl}dx^i\wedge dx^j\wedge dx^k\wedge dx^l\\
 &= dx^4\wedge dx^5\wedge dx^6\wedge dx^7+dx^2\wedge dx^3\wedge dx^6\wedge dx^7+dx^2\wedge dx^3\wedge dx^4\wedge dx^5\\
 &+dx^1\wedge dx^3\wedge dx^5\wedge dx^7-dx^1\wedge dx^3\wedge dx^4\wedge dx^6-dx^1\wedge dx^2\wedge dx^5\wedge dx^6\\
 &-dx^1\wedge dx^2\wedge dx^4\wedge dx^7\\
 \partial x_i\wedge\partial x_j&=\epsilon_{ijk}\partial x_k,\, \hbox{\rm(cross product)}. \\
 \chi&=\chi^{p}_{ijk}\partial x_p\otimes dx^i\wedge dx^j\wedge dx^k\\
 &=\frac{1}{24}\epsilon_{ijk}{}^p\partial x_p\otimes dx^i\wedge dx^j\wedge dx^k\\
\end{array}
\right.$
}\end{equation}

In (\ref{chi-equation-coordinates}) is given the $\chi$-flow equation in local coordinates.
\begin{equation}\label{chi-equation-coordinates}
 \left\{
   \begin{array}{ll}
 E_1\subset J{\it D}(W)&   \left\{(\partial t.f^p)=\sum_{1\le p\le 7}\frac{1}{24}\det(\partial x_i.f^j)\epsilon^p_{123}\right\}_{1\le i,j\le 3}:\\
&(\partial t.f^1)= 0 \\
      &(\partial t.f^2)= 0 \\
       &(\partial t.f^3)= 0 \\
      &(\partial t.f^s)=\frac{1}{24}(\det(\partial x_i.f^j))_{1\le i,j\le 3}\, 4\le s\le 7. \\
   \end{array}
   \right.
\end{equation}

Since $\Sigma^3$ is diffeomorphic to $S^3$, let us consider $\phi:\Sigma^3\to S^3$ this diffeomorphism. Then a solution of equation (\ref{chi-equation-coordinates}) is written in (\ref{solution-chi-equation-coordinates}).
\begin{equation}\label{solution-chi-equation-coordinates}
\left\{
   \begin{array}{l}
     f^k(t,x^1,\cdots,x^7)= \phi^k(x^1,x^2,x^3)\, 1\le k\le 3. \\
      f^s(t,x^1,\cdots,x^7)=\frac{1}{24}j(\phi)\, t\\
   \end{array}
   \right.
\end{equation}
where $j(\phi)$ is the determinant of the jacobian of diffeomorphism $\phi$. Therefore, the proof of the lemma is down.
\end{proof}
\begin{theorem}[Quantum $\chi$-flow]\label{theorem-quantum-chi-flow}
We define {\em quantum $\hat Q^7$} any quantum (super)manifold $M$ having as classic limit $M_C=Q^7$. We say {\em quantum $\chi$-flow} on $M$, any quantum flow that projects on the classic limit in a $chi$-flow. Then we have that a quantum $\chi$-flow represents diffeomorphically any homotopy $3$-sphere, $\Sigma^3\subset Q^7$, onto $S^3\subset Q^7$.
\end{theorem}

\begin{proof}
This follows directly from above definitions and Lemma \ref{lemma-chi-flow}.
\end{proof}

\begin{theorem}
A gauge theory in the category $\mathfrak{Q}_S$ is fully covariant.
\end{theorem}

\begin{proof}
In fact, any gauge theory in the category $\mathfrak{Q}_S$, can be identified by the following super bundle of geometric objects on a quantum supermanifold $M$: $$\left(\pi_P:P\to M\leftarrow Hom_Z(TM;\mathfrak{g}):\pi;\mathbb{B}(-)=Hom_Z(-;1_{\mathfrak{g}})\right),$$ where $\pi_P:P\to M$ is a $G$-principal bundle bundle over $M$, in the category $\mathfrak{Q}_S$, and $\mathfrak{g}$ is the Lie superalgebra corresponding to the quantum Lie supergroup $G$ of the $P$. Therefore one has the commutative diagram reported in (\ref{fiber-bundles-commutative-diagrams-gauge-theory}).

\begin{equation}\label{fiber-bundles-commutative-diagrams-gauge-theory}
 \xymatrix{P\ar@{=}[r] &B\ar@/_2pc/[ddr]_{\pi_{B,C}}\ar[dr]^{\pi_B}&&W\ar[ld]_{\pi_W}\ar@/^2pc/[ldd]^{\pi_{W,C}}\ar@{=}[r]& Hom_Z(TM;\mathfrak{g})\\
 &&M\ar[d]_{\pi_C}&&\\
 &&M_C&&\\}
 \end{equation}
This is a classic-regular super bundle of geometric objects on $M$. (For details on gauge theory in the category $\mathfrak{Q}_S$ see Refs. \cite{PRAS14,PRAS19,PRAS21,PRAS22}.)
\end{proof}

 In the following we generalize the definition of instanton and soliton of the classical field theory.

\begin{definition}[Quantum super-instantons and quantum super-solitons]\label{quantum-super-instantons-solitons}
A {\em quantum super-instanton} is a solution of the quantum super Yang-Mills equation $\widehat{(YM)}$ with
non-trivial topology.

A {\em quantum super-soliton} is a sectional compact quantum super-instanton.
\end{definition}

\begin{remark}
Let us emphasize that the classic limits of quantum supermanifolds considered in Definition \ref{quantum-super-instantons-solitons} do not necessarily coincide with usual instantons and solitons respectively. (Compare with the definitions usually adopted in commutative differential geometry.\footnote{These represent solutions of the Yang-Mills equation on four-dimensional Euclidean space, considered as the Wick rotation of Minkowski spacetime. (Let us recall the pioneering result by A. M. Polyakov \cite{POLYAKOV} proving that instantons effects in $3$-dimensional QED, coupled to a scalar field (i.e., Higgs field) lead to a massive photon.) In a paper, with the title ``{\em Super bundle of geometric objects, Yang-Mills gauge field and $Spin^G$-instantons}", announced in \cite{PRA-REGGE}, and appeared in Section 4.8 of the book \cite{PRAS5}, were developed some relations between super bundle of geometric objects on classical instantons and gravitational instantons. In particular were there applied the fully covariant theory of classical spinor fields, first previously developed in \cite{PRAS1,PRAS2}.} Compare also with their non-commutative extension given in \cite{BRADEN-NEKRASOV,NEKRASOV}.)
\end{remark}
\section{\bf Conservation Laws in Quantum Super PDEs}\label{sec-conservation-laws-quantum-super-pdes}

Conservation laws are considered for PDE's built in the category ${\frak Q}_S$ of
quantum supermanifolds. These are functions defined on the
integral bordism groups of such equations and belonging to
suitable Hopf algebras ({\em full quantum Hopf algebras}). In
particular, we specialize our calculations on the quantum super
Yang-Mills equations and quantum black holes.

In this section we shall resume some our fundamental definition
and result for PDE's in the {\em category of quantum
supermanifolds}, ${\frak Q}_S$, where the objects are just quantum
supermanifolds, and the morphisms are maps of class $Q^k_w$,
$k\in\{0,1,2,\cdots,\infty,\omega\}$ \cite{PRAS11,PRAS12,PRAS14,PRAS17,PRAS20,PRAS21}. A small subcategory
is $\mathfrak{C}_S\subset\mathfrak{Q}_S$ of supermanifolds as defined in
\cite{PRAS4}.

Let $\pi :W\rightarrow M$ be a fiber
bundle, in the category ${\frak Q}_S$, such that $ \dim
W=(m|n,r|s)$, over the quantum superalgebra $ B\equiv A\times E$
and $\dim M=(m|n)$ over $A$ and such that $E$ is a $Z$-module,
with $Z\equiv Z(A)\subset A$, the center of $A$. The {\em quantum
$k$-jet-derivative space} $J\hat D^k(W)$ of $\pi:W\to M$, is the
$k$-jet-derivative space of sections of $\pi$, belonging to the
class $Q^k_w$. The $ k$-jet-derivative $J{\hat D}^k(W)$ is a
quantum supermanifold modeled on the quantum superalgebra, ({\em quantum $k$-holonomic superalgebra}), reported in (\ref{quantum-k-holonomic-superalgebra}).
\begin{equation}\label{quantum-k-holonomic-superalgebra}
\left\{
\begin{array}{l}
 B_k\equiv\mathop{\prod}\limits_{0\le s\le
k}(\mathop{\prod}\limits_{i_1+\cdots+i_s\in{\Bbb Z}_2,i_r\in{\Bbb
Z}_2}\mathop{\widehat{A}}\limits^{s}{}_{i_1\cdots i_s}(E))\\
 \mathop{\widehat{A}}\limits^{s}{}_{i_1\cdots i_s}(E)\equiv
Hom_Z(A_{i_1}\otimes_Z\cdots\otimes_ZA_{i_s};E)\\
 \mathop{\widehat{A}}\limits^{0}(E)\equiv A\times E\\
 \mathop{\widehat{A}}\limits^{1}{}_i(E)\equiv\widehat{A}_{0}(E)\times
\widehat{A}_{1}(E)\equiv Hom_Z(A_0;E)\times Hom_Z(A_1;E).\\
\end{array}
\right.
\end{equation}

Each
$\mathop{\widehat{A}}\limits^{s}{}_{i_1\cdots i_s}(E)$ is a
quantum superalgebra with $\mathbb{Z}_2$-gradiation induced by $ E$.
Hence $ \mathop{\widehat{A}}\limits^{s}{}_{i_1\cdots
i_s}(E)_q\equiv
Hom_Z(A_{i_1}\otimes_Z\cdots\otimes_ZA_{i_s};E_p)$, $ i_r,p,
q\in\mathbb{Z}_2$, $ q\equiv i_1+\cdots+i_s+p$. If $(x^A,y^B)_{1\le
A\le m+n, 1\le B\le r+s}$ are fibered quantum coordinates on the
quantum supermanifold $W$ over $M$, then $(x^A,y^B,y^B_A,\cdots,
y^B_{A_1\cdots A_k})$ are fibered quantum coordinates on $ J\hat
D^k(W)$ over $ M$, with the following gradiations: $|x^A|=|A|$, $
|y^B|=|B|$, $ |y^B_{A_1\cdots A_s}|=|B|+|A_1|+\cdots|A_s|$. Note,
also, that there is not symmetry in the indexes $ A_i$. $ J{\hat
D}^k(W)$ is an affine bundle over $ J\hat D^{k-1}(W)$ with
associated vector bundle $ \pi^*_{k,0}Hom_Z(\dot S^k_0M;vTW)$,
where $\dot S^k_0M$ is the $k$-times symmetric tensor product of
$TM$, considered as a bundle of $Z$-modules over $M$, and
$\pi_{k,0}:J{\hat D}^k(W)\to W$ is the canonical surjection.
Another important example is $\hat J^k_{m|n}(W)$, that is the {\em
$k$-jet space for quantum supermanifolds} of dimension $(m|n)$
(over $A$) contained in the quantum supermanifold $W$. This
quantum supermanifold locally looks like $J\hat D^k(W)$. Set $
J\hat D^{\infty}(W)\equiv\mathop{\rm lim}\limits_{\leftarrow\hskip
2pt k }J\hat D^{k}(W)$, $\hat J^\infty_{m|n}(W)\equiv\mathop{\rm
lim}\limits_{\leftarrow\hskip 2pt  k}\hat J^{k}_{m|n}(W)$. These
are quantum supermanifolds modeled on $B\equiv\prod_k B_k$.

A {\em quantum super PDE} of order $k$ on the fibre bundle
$\pi:W\to M$, defined in the category of quantum supermanifolds,
${\frak Q}_S$, is a subset $\hat E_k\subset J\hat D^k(W)$, or
$\hat E_k\subset \hat J^k_{m|n}(W)$. A geometric theory of quantum
(super) PDE's can be formulated introducing suitable hypotheses of
regularity on $\hat E_k$. (See \cite{PRAS14}.)

The characterization of global solutions of a PDE $\hat
E_k\subseteq \hat J^k_{m|n}(W)$, in the category ${\frak Q}_S$,
can be made by means of its integral bordism groups
$\Omega_{p|q}^{\hat E_k}$, $p\in\{0,1,\dots,m-1\}$,
$q\in\{0,1,\dots,n-1\}$. Let us shortly recall some fundamental
definitions and results about. Let $f_i:X_i\to \hat E_k$,
$f_i(X_i)\equiv N_i\subset \hat E_k$, $i=1,2$, be
$(p|q)$-dimensional admissible compact closed smooth integral
quantum supermanifolds of $\hat E_k$. The admissibility requires
that $N_i$ should be contained into some solution $V\subset \hat
E_k$, identified with a $(m,n)$-chain, with coefficients in $A$. Then, we say that they are {\em $\hat E_k$-bordant} if
there exists a $(p+1|q+1)$-dimensional smooth quantum
supermanifolds $f:Y\to \hat E_k$, such that $\partial Y=X_1\sqcup
X_2$, $f|_{X_i}=f_i$, $i=1,2$, and $V\equiv f(Y)\subset \hat E_k$
is an admissible integral quantum supermanifold of $\hat E_k$ of
dimension $(p+1|q+1)$. We say that $N_i$, $i=1,2$, are {\em $\hat
E_k$-quantum-bordant} if there exists a $(p+1|q+1)$-dimensional
smooth quantum supermanifolds $f:Y\to \hat J^k_{m|n}(W)$, such
that $\partial Y=X_1\sqcup X_2$, $f|_{X_i}=f_i$, $i=1,2$, and
$V\equiv f(Y)\subset \hat J^k_{m|n}(W)$ is an admissible integral
manifold of $\hat J^k_{m|n}(W)$ of dimension $(p+1|q+1)$. Let us
denote the corresponding bordism groups by $\Omega_{p|q}^{\hat
E_k}$ and $\Omega_{p|q}(\hat E_k)$, $p\in\{0,1,\dots,m-1\}$,
$q\in\{0,1,\dots,n-1\}$, called respectively
{\em$(p|q)$-dimensional integral bordism group} of $\hat E_k$ and
{\em$(p|q)$-dimensional quantum bordism group} of $\hat E_k$.
Therefore these bordism groups work, for $(p,q)=(m-1,n-1)$, in the
category of quantum supermanifolds that are solutions of $\hat
E_k$. Let us emphasize that singular solutions of $\hat E_k$ are,
in general, (piecewise) smooth quantum supermanifolds into some
prolongation $(\hat E_k)_{+s}\subset \hat J^{k+s}_{m|n}(W)$, where
the set, $\Sigma(V)$, of {\em singular points} of a solution $V$
is a non-where dense subset of $V$. Here we consider {\em
Thom-Boardman singularities}, i.e., $q\in\Sigma(V)$, if
$(\pi_{k,0})_*(T_qV)\not\cong T_qV$. However, in the case where
$\hat E_k$ is a differential equation of finite type, i.e., the
symbols $\hat g_{k+s}=0$, $s\ge 0$, then it is useful to include
also in $\Sigma(V)$, discontinuity points, $q,q'\in V$, with
$\pi_{k,0}(q)=\pi_{k,0}(q')=a\in W$, or with
$\pi_{k}(q)=\pi_{k}(q')=p\in M$, where
$\pi_k=\pi\circ\pi_(k,0):\hat J^k_{m|n}(W)\to M $. We denote such
a set by $\Sigma(V)_S$, and, in such cases we shall talk more
precisely of {\em singular boundary} of $V$, like $(\partial
V)_S=\partial V\setminus\Sigma(V)_S$. Such singular solutions are
also called {\em weak solutions}.

Let us define some notation to distinguish between some integral
bordisms.

\begin{definition}
Let $\Omega_{m-1|n-1}^{\hat
E_k}$, (resp. $\Omega_{m-1|n-1,s}^{\hat E_k}$, resp.
$\Omega_{m-1|n-1,w}^{\hat E_k}$), be the integral bordism group
for $(m-1|n-1)$-dimensional smooth admissible regular integral
quantum supermanifolds contained in $\hat E_k$, borded by smooth
regular integral quantum supermanifold-solutions, (resp.
piecewise-smooth or singular solutions, resp. singular-weak
solutions), of $\hat E_k$.
\end{definition}

\begin{theorem}\cite{PRAS14}
One has the exact commutative diagram {\em(\ref{ddd})}. Therefore, one has the canonical isomorphisms:

$K^{\hat E_k}_{m-1|n-1,w/(s,w)}\cong K^{\hat
E_k}_{m-1|n-1,s}$; $\Omega^{\hat E_k}_{m-1|n-1}/K^{\hat
E_k}_{m-1|n-1,s}\cong \Omega^{\hat E_k}_{m-1|n-1,s}$;

$\Omega^{\hat E_k}_{m-1|n-1,s}/K^{\hat
E_k}_{m-1|n-1,s,w}\cong\Omega^{\hat E_k}_{m-1|n-1,w}$;
$\Omega^{\hat E_k}_{m-1|n-1}/K^{\hat
E_k}_{m-1|n-1,w}\cong\Omega^{\hat E_k}_{m-1|n-1,w}$.

If $\hat E_k$ is formally quantum superintegrable, then one has the following isomorphisms:

$\Omega^{\hat E_k}_{m-1|n-1}\cong\Omega^{\hat
E_\infty}_{m-1|n-1}\cong\Omega^{\hat E_\infty}_{m-1|n-1,s}$;
$\Omega^{\hat E_k}_{m-1|n-1,w}\cong\Omega^{\hat
E_\infty}_{m-1|n-1,w}$.

\begin{equation}\label{ddd}
\scalebox{0.8}{$\xymatrix{&0\ar[d]&0\ar[d]&0\ar[d]&\\
0\ar[r]&K^{\hat E_k}_{m-1|n-1,w/(s,w)}\ar[d]\ar[r]&K^{\hat E_k}_{m-1|n-1,w}\ar[r]\ar[d]&K^{\hat E_k}_{m-1|n-1,s,w}\ar[d]\ar[r]& 0\\
0\ar[r]&K^{\hat E_k}_{m-1|n-1,s}\ar[r]\ar[d]&\Omega^{\hat E_k}_{m-1|n-1}\ar[d]\ar[r]&\Omega^{\hat E_k}_{m-1|n-1,s}\ar[d]\ar[r]&0\\
&0\ar[r]&\Omega^{\hat E_k}_{m-1|n-1,w}\ar[d]\ar[r]&\Omega^{\hat E_k}_{m-1|n-1,w}\ar[d]\ar[r]&0\\
&&0&0&\\}$}
\end{equation}
\end{theorem}

\begin{theorem}\label{regularity-and-integral-bordism-groups}
Let $\hat
E_k\subset\hat J^k_{m|n}(W)$ be a quantum super PDE that is
formally quantum superintegrable, and completely superintegrable.
We shall assume that the symbols $\hat g_{k+s}\not=0$, $s=0,1$.
(This escludes the case $k=\infty$.) Then one has the following
isomorphisms: $\Omega_{p|q,s}^{\hat E_k}\cong\Omega_{p|q,w}^{\hat
E_k}\cong\Omega_{p|q}(\hat E_k)$, with $p\in\{0,\dots,m-1\}$ and
$q\in\{0,\dots,n-1\}$.\end{theorem}

\begin{proof}
In fact, in these cases any
weak solution identifies a singular solution, by connecting its
branches by means of suitable pieces of fibres. Furthermore, since
$\hat E_{k+1}$ is a strong retract of $\hat J^{k+1}_{m|n}(W)$, we
can deform any quantum bording $V\subset\hat J^{k+1}_{m|n}(W)$,
$\dim V=(m|n)$, with $\partial V\subset\hat E_{k+1}$, into a
(singular) solution of $\hat E_{k+1}$, hence into a solution of
$\hat E_{k}$. (For details see Refs. \cite{PRAS14}.)
\end{proof}

\begin{cor}
Let $\hat E_k\subset \hat J^k_{m|n}(W)$ be a quantum super PDE, that is
formally superintegrable and completely superintegrable. One has
the following isomorphisms: $\Omega_{m-1|n-1,w}^{\hat
E_k}\cong\Omega_{m-1|n-1}(\hat E_k)\cong\Omega_{m-1|n-1,w}^{\hat
E_{k+h}} \cong\Omega_{m-1|n-1,w}^{\hat
E_\infty}\cong\Omega_{m-1|n-1,w} (\hat
E_{k+h})\cong\Omega_{m-1|n-1}(\hat E_\infty)$.
\end{cor}

In order to distinguish between quantum integral supermanifolds
$V$ representing singular solutions, where $\Sigma(V)$ has no
discontinuities, and quantum integral supermanifolds where
$\Sigma(V)$ contains discontinuities, we can also consider
``conservation laws" valued on quantum integral supermanifolds
$N$ representing the integral bordism classes $[N]_{\hat
E_k}\in\Omega_{p|q}^{E_k}$.

\begin{definition}
Let us define the {\em space of quantum integral conservation laws} of $\hat E_k\subset \hat J^k_{m|n}(W)$ the $Z$-module given in {\em(\ref{space-quantum-integral-conservation-laws})}.
\begin{equation}\label{space-quantum-integral-conservation-laws}
  \left\{
\begin{array}{ll}
 \hat{\frak
I}(\hat E_k) &\equiv\oplus_{p,q\ge 0 }{{\widehat{\Omega}^{p|q}(\hat
E_k)\cap d^{-1}(C\widehat{\Omega}^{p+1|q+1}(\hat
E_k))}\over{d\widehat{\Omega}^{p-1|q-1}(E_k)\oplus\{C\widehat{\Omega}^{p|q}(\hat
E_k)\cap d^{-1}(C\widehat{\Omega}^{p+1|q+1}(\hat
E_k)))\}}}  \\
 & \equiv\oplus_{p,q\ge 0 }\hat{\frak I}(\hat
E_k)^{p|q}.\\
\end{array}
  \right.
\end{equation}

Here $C\widehat{\Omega}^{p|q}(\hat E_k)$ denotes the
space of all quantum $(p|q)$-forms on $\hat E_k$. Then we define
{\em quantum integral characteristic supernumbers} of $N$, with $[N]_{\hat
E_k}\in\widehat{\Omega}_{p|q}^{\hat E_k}$, the numbers $\hat
i[N]\equiv<[N]_{\hat E_k},[\alpha]>\in B$, for all
$[\alpha]\in\hat{\frak I}(\hat E_k)^{p|q}$.
\end{definition}
Then, one has the following theorems.

\begin{theorem}\label{conservation-laws-and-integral-characteristic-supernumbers}\cite{PRAS14}
Let us assume that
$\hat{\frak I}(\hat E_k)^{p|q}\not=0$. One has a natural
homomorphism: ${\underline{j}}_{p|q}:\Omega_{p|q}^{\hat E_k}\to
Hom_A(\hat{\frak I}(\hat E_k)^{p|q};A)$, $[N]_{\hat E_k}\mapsto
{\underline{j}}_{p|q}([N]_{\hat E_k})$,
${\underline{j}}_{p|q}([N]_{\hat
E_k})([\alpha])=\int_N\alpha\equiv<[N]_{\hat E_k},[\alpha]>$.
Then, a necessary condition that $N'\in[N]_{\hat E_k}$ is the
following: $\hat i[N]=\hat i[N']$,  $\forall[\alpha]\in\hat{\frak
I}(\hat E_k)^{p|q}$. Furthermore, if the {\em classic limit},
$N_C$, of $N$ is orientable then above condition is
sufficient also in order to say that $N'\in[N]_{\hat E_k}$.
\end{theorem}

\begin{cor}
Let $\hat
E_k\subseteq \hat J^k_{m|n}(W)$ be a quantum super PDE. Let us
consider admissible $(p|q)$-dimensional, $0\le p\le m-1$, $0\le
q\le n-1$, integral quantum supermanifolds, with orientable
classic limits. Let $N_1\in[N_2]_{E_k}\in\Omega_{p|q}^{\hat E_k}$,
then there exists a $(p+1|q+1)$-dimensional admissible integral
quantum supermanifold $V\subset \hat E_k$, such that $\partial
V=N_1\sqcup N_2$, where $V$ is without discontinuities iff the
integral supernumbers of $N_1$ and $N_2$ coincide.
\end{cor}

Above considerations can be generalized to include more
sophisticated quantum solutions of quantum super PDEs.

\begin{definition}
Let $\hat
E_k\subset \hat J^k_{m|n}(W)$ be a quantum super PDE and let $B$
be a quantum superalgebra. Let us consider the following chain
complex {\em(bigraded bar quantum chain complex of $\hat E_k$)}:
$\{\bar C_{\bullet|\bullet}(\hat E_k;B),\partial\}$, induced by
the $\mathbb{Z}_2$-gradiation of $B$ on the corresponding bar
quantum chain complex of $\hat E_k$, i.e., $\{\bar
C_{\bullet}(\hat E_k;B),\partial\}$. (See Refs.{\em\cite{PRAS14}}.) More
precisely $\bar C_{p}(\hat E_k;B)$ is the free two-sided
$B$-module of formal linear combinations with coefficients in $B$,
$\sum\lambda_i c_i$, where $c_i$ is a singular $p$-chain
$f:\triangle^p\to\hat E_k$, that extends on a neighborhood
$U\subset{\Bbb R}^{p+1}$, such that $f$ on $U$ is differentiable
and $Tf(\triangle^p)\subset\hat{\bf E}^k_{m|n}$, where $\hat{\bf
E}^k_{m|n}$ is the Cartan distribution of $\hat E_k$.
\end{definition}

\begin{theorem}\label{theorem-to-recall}\cite{PRAS14}
The homology $\bar
H_{\bullet|\bullet} (\hat E_k;B)$ of the bigraded bar quantum
chain complex of $\hat E_k$ is isomorphic to {\em(closed) bar
integral singular $(p|q)$-bordism groups}, with coefficients in
$B$, of $\hat E_k$: ${}^B\bar\Omega_{{p|q} ,s}^{\hat E_k}\cong
\bar H_{q|q} (\hat E_k;B) \cong (\bar\Omega_{p,s}^{\hat
E_k}\otimes_{\Bbb K}B_0)\bigoplus (\bar\Omega_{q,s}^{\hat
E_k}\otimes_{\Bbb K}B_1)$, $p\in\{0,1,\dots,m-1\}$,
$q\in\{0,1,\dots,n-1\}$. (If $B={\Bbb K}$ we omit the apex $B$).
If $\hat E_k\subset \hat J_{m|n}^k(W)$ is formally quantum
superintegrable and completely superintegrable, and the symbols
$\hat g_{k+s}\not=0$, then one has the following canonical
isomorphisms: $ {}^A\bar\Omega_{{p|q} ,s}^{\hat
E_k}\cong\Omega_{{p|q} ,w}^{\hat E_k}\cong\Omega_{{p|q} ,s}^{\hat
E_k}\cong\Omega_{{p|q}}(\hat E_k)$. Furthermore, the quantum
$(p|q)$-bordism groups $\Omega_{p|q}(\hat E_k)$ is an extension of
a subgroup of ${}^A{\underline{\Omega}}_{p|q,s}(W)\cong
H_{p|q}(W;A)$, and the integral $(p|q)$-bordism group
$\Omega_{p|q}^{\hat E_k}$ is an extension of the quantum
$(p|q)$-bordism group.
\end{theorem}

\begin{cor}\label{corollary-to-recall}
Let $\hat
E_k\subset \hat J^k_{m|n}(W)$ be a quantum super PDE, that is
formally superintegrable and completely superintegrable. One has
the following isomorphisms: $\Omega_{m-1|n-1,w}^{\hat
E_k}\cong\Omega_{m-1|n-1}(\hat E_k)\cong\Omega_{m-1|n-1,w}^{\hat
E_{k+h}} \cong\Omega_{m-1|n-1,w}^{\hat
E_\infty}\cong\Omega_{m-1|n-1,w} (\hat
E_{k+h})\cong\Omega_{m-1|n-1}(\hat E_\infty)
\cong{}^A{\underline{\Omega}}_{{m-1|n-1},s}(W)\cong
H_{{m-1|n-1}}(W;A)$.
\end{cor}

\begin{definition}
The {\em full space of $(p|q)$-conservation laws}, (or {\em full $(p|q)$-Hopf
superalgebra}), of $\hat E_k$ is the following one: ${\bf
H}_{p|q}(E_k)\equiv B^{\Omega_{p|q}^{\hat E_k}}$, where $B\equiv\prod_{k} B_k$. We call {\em full Hopf superalgebra}, of
$\hat E_k$, the following: ${\bf H}_{m-1|n-1}(\hat E_\infty)\equiv
B ^{\Omega_{m-1|n-1}^{\hat E_\infty}}$.
\end{definition}

\begin{definition}\label{space-differential-conservation-laws}
The {\em space of (differential) conservation laws} of $\hat E_k\subset \hat
J^k_{m|n}(W)$, is ${\frak C}ons(\hat E_k)=\hat{\frak I}(\hat
E_\infty)^{m-1|n-1}$.
\end{definition}

\begin{theorem}\cite{PRAS14}
The full
$(p|q)$-Hopf superalgebra of a quantum super PDE $\hat E_k\subset
\hat J^k_{m|n}(W)$ has a natural structure of quantum Hopf
superalgebra. Quantum Hopf
algebras are generalizations of such algebras.
\end{theorem}

\begin{proposition}
The space of
conservation laws of $\hat E_k$ has a canonical representation in
${\bf H}_{m-1|n-1}(\hat E_\infty)$.
\end{proposition}

\begin{proof}
In fact, one has the following
homomorphism $j:{\frak C}ons(E_k)\to{\bf H}_{m-1|n-1}(\hat
E_\infty)$, $j[\alpha]([N]_{\hat E_\infty})=<[\alpha],[N]_{\hat
E_\infty}>=\int_{N_C}i^*\alpha\in B$, where $i:N_C\to N$ is the
canonical injection.
\end{proof}

\begin{theorem}
Set: ${\bf K}^{\hat
E_k}_{m-1|n-1,w/(s,w)}\equiv B ^{K^{E_k}_{n-1,w/(s,w)}}$, ${\bf
K}^{\hat E_k}_{m-1|n-1,w}\equiv B ^{K^{\hat E_k}_{m-1|n-1,w}}$,
${\bf K}^{\hat E_k}_{m-1|n-1,s,w}\equiv B^{K^{\hat
E_k}_{m-1|n-1,(s,w)}}$, ${\bf K}^{\hat E_k}_{m-1|n-1,s}\equiv
B^{K^{\hat E_k}_{m-1|n-1,s}}$, ${\bf H}_{m-1|n-1}(\hat E_k)\equiv
B^{\Omega^{\hat E_k}_{m-1|n-1}}$, ${\bf H}_{m-1|n-1,s}(\hat
E_k)\equiv B ^{\Omega^{\hat E_k}_{m-1|n-1,s}}$, ${\bf
H}_{m-1|n-1,w}(\hat E_k)\equiv B^{\Omega^{\hat E_k}_{m-1|n-1,w}}$.
One has the following canonical isomorphisms:
$$\begin{array}{l}
{\bf K}^{\hat E_k}_{m-1|n-1,w/(s,w)}\cong{\bf K}^{K^{\hat
E_k}_{m-1|n-1,s}};\\
{\bf K}^{\hat E_k}_{m-1|n-1,w}/{\bf K}^{\hat
E_k}_{n-1,s,w}\cong{\bf K}^{K^{\hat E_k}_{m-1|n-1,w/(s,w)}};\\
{\bf H}_{m-1|n-1}(\hat E_k)/{\bf H}_{m-1|n-1,s}(\hat E_k)\cong{\bf
K}^{\hat E_k}_{m-1|n-1,s};\\
{\bf H}_{m-1|n-1}(\hat E_k)/{\bf H}_{m-1|n-1,w}(\hat E_k)\cong{\bf
K}^{\hat E_k}_{m-1|n-1,w} \\
\cong{\bf H}_{m-1|n-1,s}(\hat E_k)/{\bf H}_{m-1|n-1,w}(\hat
E_k)\cong{\bf K}^{\hat E_k}_{m-1|n-1,s,w}.\\
\end{array}$$
\end{theorem}

\begin{proof}
The proof is obtained directly
by duality of the exact commutative diagram (\ref{ddd}).
\end{proof}

\begin{theorem}
Under the same
hypotheses of Theorem \ref{regularity-and-integral-bordism-groups}, one has the following canonical
isomorphism: ${\bf H}_{m-1|n-1,s}(E_k)\cong{\bf
H}_{m-1|n-1,w}(\hat E_k)$. Furthermore, we can represent
differential conservation laws of $E_k$ in ${\bf
H}_{m-1|n-1,w}(\hat E_k)$.
\end{theorem}

\begin{proof}
Let us note that $\hat{\frak
I}(\hat E_k)^{m-1|n-1}\subset\hat{\frak I}(\hat
E_\infty)^{m-1|n-1}$. If $j:{\frak C}ons(\hat E_k)\to{\bf
H}_{m-1|n-1}(\hat E_\infty)$, is the canonical representation of
the space of the differential conservation laws in the full Hopf
superalgebra of $\hat E_k$, (corresponding to the integral bordism
groups for regular smooth solutions), it follows that one has also
the following canonical representation $j|_{\hat{\frak I}(\hat
E_k)^{m-1|n-1}}:\hat{\frak I}(\hat E_k)^{m-1|n-1}\to{\bf
H}_{m-1|n-1,s}(\hat E_k)\cong {\bf H}_{m-1|n-1,w}(\hat E_k)$. In
fact, for any $N'\in[N]_{\hat E_k,s}\in\Omega_{m-1|n-1,s}^{\hat
E_k}\cong\Omega_{m-1|n-1,w}^{\hat E_k}$, one has
$\int_{N'}\beta=\int_{N}\beta$, for any $[\beta]\in \hat{\frak
I}(\hat E_k)^{m-1|n-1}$.
\end{proof}

\begin{theorem}[Quantum tunnel effects and quantum black holes]\label{quantum-tunnel-effects-and-quantum-black-holes}
The quantum supergravity equation $\hat E_2\subset
J\hat D^2(i^*\bar C)$ admits global solutions having a change of
sectional topology ({\em quantum tunnel effects}). In general
these solutions are not globally representable as second
derivative of sections of the fiber bundle $i^*\bar C\to
N$. $\hat E_2$ admits solutions that represent {\em evaporations of quantum black holes}.
\end{theorem}

\begin{proof}

\begin{table}[h]
\caption{Dynamic Equation on macroscopic shell: $\hat E_2[i]\subset J\hat D^2(i^*\bar C)$ and Bianchi
identity.}
\label{observed-dynamic-equation-and-bianchi-identity}
\begin{tabular}{|c|l|}
\hline{\footnotesize{\rm Fields equations}}&{\footnotesize{\rm$(\partial\omega^\gamma_{ab}.L)-\partial_\mu(\partial\omega^{\gamma\mu}_{ab}.L)=0$
                  (curvature equation)}}\\
&{\footnotesize{\rm$(\partial\theta^\gamma_\alpha.L)-\partial_\mu(\partial\theta^{\gamma\mu}_\alpha.L)=0$
                   (torsion equation)}}\\
{\footnotesize{\rm $(\hat E_2[i])$}}&{\footnotesize{\rm $(\partial\psi^\gamma_{\beta i}.L)
                  -\partial_\mu(\partial\psi^{\gamma\mu}_{\beta i}.L)=0$
                  (gravitino equation)}}\\
&{\footnotesize{\rm $(\partial A^\gamma.L)-\partial_\mu(\partial A^{\gamma\mu}.L)=0$
                  (Maxwell's equation)}}\\
\hline
{\footnotesize{\rm Bianchi identity}}&{\footnotesize{\rm $(\partial x_{[\gamma }.R^{ab}_{\beta\alpha ]})
                            +2\omega ^a_{e[\gamma}R^{eb}_{\beta\alpha]}=0$}}\\
      &{\footnotesize{\rm$(\partial x_{[\gamma }.R^{\alpha }_{\beta \omega ]})+
                            \omega ^{\alpha b}_{[\gamma}R_{\beta\omega] b}+\frac{1}{2}
                            (C\gamma ^{\alpha})_{\delta\mu}\psi^\delta_{j[\gamma}
                            \rho _{\beta\omega]}^{\mu j}=0$}}\\
{\footnotesize{\rm$(B[i])$}}&{\footnotesize{\rm $(\partial x_{[\gamma } .\rho ^{\beta i}_{\omega \alpha ]})
                           +\frac{1}{2}(\sigma_{ab})^{\beta i}_{\delta j}\omega ^{ab}_{[\gamma }\rho _{\omega\alpha]}^{\delta j}=0$}}\\
&{\footnotesize{\rm$(\partial x_{[\gamma }F_{\beta \alpha]})+\frac{1}{2}C_{\delta\mu}
                            \epsilon_{ij}\psi _{[\gamma}^{\delta i}
                            \rho ^{\mu j}_{\beta\alpha ]}=0$}}\\
\hline
{\footnotesize{\rm Fields}}&{\footnotesize{\rm $R^{ab}_{\mu \nu }=(\partial x_{[\mu} .\omega ^{ab}_{\nu]} )
                   +2\omega ^a _{e[\mu}\omega ^{eb}_{\nu ]}$ (curvature)}}\\
       &{\footnotesize{\rm$R^{\alpha }_{\mu \nu }=
                   (\partial x_{[\mu }.{\theta} ^{\alpha }_{\nu]} )+
                   {\omega} ^{\alpha }_{\beta[\mu},{\theta} ^{\beta}_{\nu]}
                   +\frac{1}{2}(C\gamma^\alpha)_{\beta\delta}
                   \psi^\beta_{j[\mu}\psi^{\delta j}_{\nu]}$
                   (torsion)}}\\
                   &{\footnotesize{\rm ${\rho} ^{\beta i}_{\mu \nu }
                   =(\partial x_{[\mu} .{\psi }^{\beta i}_{\nu] })
                   +\frac{1}{2}(\sigma_{ab})^{\beta i}_{\gamma j}\omega^{ab}_{[\mu}
                   \psi^{\gamma j}_{\nu]}\hskip 3pt\hbox{\rm({\rm gravitino})}$}}\\
                   &{\footnotesize{\rm$F_{\mu \nu }=(\partial x_{[\mu}.A_{\nu]} )+\frac{1}{2}
                   C_{\beta\gamma}\epsilon_{ij}\psi^{\beta i}_{[\mu}\psi^{\gamma j}_{\nu]}$
                   (electromagnetic field)}}\\
                   \hline
                   \end{tabular}
                   \end{table}

We shall consider, now, the {\em quantum $N=2$ superPoincar\'e
group} over a quantum superalgebra $A=A_0\oplus A_1$, that is a
quantum Lie supergroup $G$ having as quantum Lie superalgebra
$\widehat{\frak g}$ one identified by the following infinitesimal
generators: $\{Z_K\}_{1\le K\le
19}\equiv\{J_{\alpha\beta},P_\alpha,\overline{Z},Q_{\beta
i}\}_{0\le\alpha,\beta\le 3;1\le a\le 2}$, such that
$J_{\alpha\beta}=-J_{\beta\alpha},P_\alpha,\overline{Z}\in
Hom_Z(A_0;{\frak g})$, $Q_{\beta i}\in Hom_Z(A_1;{\frak g})$. The
corresponding nonzero ${\Bbb Z}_2$-graded brackets are the
following: $[J_{\alpha\beta},J_{\gamma\delta}]=\eta _{\beta\gamma}
J_{\alpha\delta}+\eta _{\alpha\delta}J_{\beta\gamma}-\eta
_{\alpha\gamma}J_{\beta\delta}-\eta_{\beta\delta}J_{\alpha\gamma}$,
 $[P_\alpha,P_\beta]=-8e^2J_{\alpha\beta},\quad
[J_{\alpha\beta},P_\gamma]=\eta _{\beta\gamma}P_\beta-\eta
_{\alpha\gamma}P_\beta$, $[J_{\alpha\beta},Q_{\gamma
i}]=(\sigma_{\alpha\beta})^{\mu j}_\gamma Q_{\mu j},\quad
[Q_{\beta i},Q_{\mu
j}]=(C\gamma^\alpha)_{\beta\mu}\delta_{ij}P_\alpha+C
_{\beta\mu}\epsilon_{ij}\overline{Z}$. Here $ C_{\alpha\beta}$ is
the antisymmetric charge conjugation matrix, $
\sigma_{\beta\mu}={1\over 4}[\gamma_\beta,\gamma_\mu]$, with
$\gamma^\mu$ the Dirac matrices. $\overline{Z}$ commutes with all
the other ones. One
has $\dim G=(d|N_2)=(11|8)$, and we will consider the following
principal bundle in the category of quantum supermanifolds: {}$P$
is a quantum supermanifold of dimension $(15|8)$; $M$ is a quantum
supermanifold of dimension $(4|N_1)=(4|0)$, identified, for the sake of simplicity,
with a quantum Minkowski space-time. Then a
pseudoconnection can be written by means of the following
fullquantum differential $1$-forms on $P$:
${}_\rceil\mu^K=\mu^K_HdY^H$, $(\mu^K_H)
=(\frac{1}{2}\omega^{\alpha\beta}_H,\theta^\mu_H,A_H,\psi^{aj}_H)
$. With respect to a section $s:M\to P$ we get:
$(s^*{}_\rceil\mu)^K=\bar\mu^K_\gamma dX^\gamma$, $
(\bar\mu^K_\gamma) =(\frac{1}{2}\bar\omega^{\alpha\beta}_\gamma
,\bar\theta^\mu_\gamma,\bar A_\gamma,\bar\psi^{\alpha j}_\gamma
)$, where $\bar\omega^{\alpha\beta}_\gamma$ is the usual
Levi-Civita connection, $\bar\theta^\mu_\gamma$ is the vierbein,
$\bar A_\gamma$  is the electromagnetic field and
$\bar\psi^{aj}_\gamma$ is the usual spin ${3\over 2}$ field. The
blow up structure: $ \pi ^*\hat C(P)\hookrightarrow
Hom_Z(TP;{\frak g})$ implies that we can identify our fields with
sections $ {}_\rceil\mu$ of the fiber bundle $\bar\pi:\bar C\equiv
Hom_Z(TM;{\frak g})\to M$. ($ \hat C(P)\cong J\hat D(P)/G$ is the
fiber bundle, over $ M$, of principal quantum connections on the $
G$-principal fiber bundle $\pi:P\to M$.) The corresponding
curvatures can be written in the form: $ {}_\rceil
R^K_{\beta\alpha}=(\partial x_{\beta}
\mu^K_{\alpha})+C^K_{IJ}[\mu^I_{\beta},\mu^J_{\alpha}]_+$. The
local expression of the {\em dynamic equation}, $\hat
E_2[i]\subset J\hat D^2(i^*\bar C)$, evalued on a macroscopic
shell, i.e., an embedding $i:N\to M$, of a globally hyperbolic,
$p$-connected manifold $N$, $0\le p\le 3$, is given by the quantum
super PDE reported in Tab.3.1, where $ L:J\hat
D(\underline{E})\to\widehat{A}$ is a quantum Lagrangian function.
Possible Lagrangian densities are polynomial in the curvature,
(see example below), hence we can assume that they give formally
quantum superintegrable, and completely quantum superintegrable,
quantum super PDE's. Then, assuming that $\hat E_2[i]$ is formally
integrable and completely superintegrable, the integral bordism
groups of $ \hat E_2[i]$ and its fullquantum $ p$-Hopf
superalgebras, can be calculated. More precisely, we use the fact
that $\hat C(P)\to M$ is a contractible fiber bundle of dimension
$(4|0,44|32)$ over the quantum superalgebra $A\times\widehat
A=(A_0\times
A_1)\times\mathop{\widehat{A}}\limits^1{}_0(A)\times\mathop{\widehat{A}}\limits^1{}_1(A)$,
and that $N$ is topologically trivial. In fact, we can apply
Theorem \ref{theorem-to-recall} and Corollary \ref{corollary-to-recall}, to obtain the quantum and integral
bordism groups of $ {\hat E}_2[i]$: $\Omega_{p,s}^{\hat
E_2[i]}\cong \Omega_{p,w}^{\hat E_2[i]_{+\infty}}\cong 0$, for
$p=1,2,3$ and $\Omega_{0,s}^{\hat E_2[i]}\cong \Omega_{0,w}^{\hat
E_2[i]_{+\infty}}\cong A$. Therefore, we have that $
1$-dimensional admissible integral closed quantum submanifolds
contained into $ {\hat E}_2[i]$, ({\em admissible quantum closed
strings}), can propagate and interact between them by means of $
2$-dimensional admissible integral quantum manifolds contained
into $ \hat J_4^2(i^*\bar C)$, or by means of $ 2$-dimensional
admissible integral quantum manifolds contained into $ {\hat
E}_2[i]$, in such a way to generate (quantum) tunnel effects.
Finally, as a consequence of the triviality of the $
3$-dimensional integral bordism grooups, we get the existence of
global quantum solutions of such equations.

Let us now see that theorem can be proved by
using surgery techniques and taking into account that for the
$3$-dimensional integral bordism group of $\hat E_k$ one has
$\Omega_{3,s}^{\hat E_k}=0=\Omega_{3,w}^{({\hat E}_2)_{+\infty}}$.
In fact a boundary value problem for $\hat E_k[i]$ can be directly
implemented in the manifold $\hat E_k[i]\subset J\hat D^2(i^*\bar
C)\subset \hat J^2_4(i^*\bar C)$ by requering that a
$3$-dimensional compact space-like (for some $t=t_0$), admissible
integral manifold $ B\subset\hat E_k[i]$ propagates in $\hat
E_k[i]$ in such a way that the boundary $
\partial B$ describes a fixed $3$-dimensional
time-like integral manifold $Y\subset\hat E_k[i]$. (We shall
require that the boundary $\partial B$ of $B$ is orientable.) $Y$
is not, in general, a closed (smooth) manifold. However, we can
solder $Y$ with two other compact $3$-dimensional integral
manifolds $ X_i$, $ i=1,2$, in such a way that the result is a
closed $ 3$-dimensional (smooth) integral manifold $ Z\subset\hat
E_k[i]$. More precisely, we can take $X_1=B$ so that $
\widetilde{Z}\equiv X_1\bigcup_{\partial B}Y$ is a $
3$-dimensional compact integral manifold such that $
\partial\widetilde{Z}\equiv C$ is a $2$-dimensional
space-like integral manifold. We can assume that $ C$ is an
orientable manifold. Then, from the triviality of the integral
bordism group, it follows that $\partial X_2=C$, for some
space-like compact $3$-dimensional integral manifold $
X_2\subset\hat E_k[i]$. Set $ Z\equiv\widetilde{Z}\bigcup_CX_2$.
Therefore, one has $Z=X_1\bigcup_{\partial B}Y\bigcup_CX_2$. Then,
again from the triviality of the integral bordism group, it
follows also that there exists a $4$-dimensional integral (smooth)
manifold $V\subset\hat E_k[i]$ such that $
\partial V=Z$. Hence the integral manifold $ V$ is a
solution of our boundary value problem between the times $t_0$ and
$ t_1$, where $t_0$ and $t_1$ are the times corresponding to the
boundaries where are soldered $X_i$, $i=1,2$ to $Y$. Now, this
process can be extended for any $t_2>t_1$. So we are able to find
(smooth) solutions for any $ t>t_0$, hence (smooth) solutions for
any $ t>t_0$, therefore, global (smooth) solutions. Remark that in
order to assure the smoothness of the global solution so built it
is enough to develop such construction in the infinity
prolongation $\hat E_k[i]_{+\infty}$ of $\hat E_k[i]$. Finally
note that in the set of solutions of $ \hat E_k[i]$ there are ones
that have change of sectional topology. In fact the $
3$-dimensional integral bordism groups are trivial:
$\Omega_{3,s}^{\hat E_k[i]}=0=\Omega_{3,w}^{\hat
E_k[i]_{+\infty}}$.

Let us, now, consider the dynamics of a quantum black-hole.
In order to obtain such
solutions we must have a Cauchy integral data with a geometric
black hole $B$ embedded in a compact $ 3$-dimensional integral
manifold $N$, $B\subset N$, such that its boundary $\partial N$
propagates with a fixed flow. Then a solution, with quantum tunnel
effect of such boundary problem, can describe an evaporation process
of such black hole. Above results assure the existence of such
solutions and a way to build them.

In order to represent such results with respect to a quantum relativistic observer, let us consider
the space of observed quantum integral conservation laws.
\begin{definition}\label{observed-quantum-integral-conservation-law}
$$\left\{
\begin{array}{ll}
  \hat{\frak
I}(\hat E_k[i])& \equiv\oplus_{q\ge 0 }{{\widehat{\Omega}^{q}(\hat
E_k[i])\cap d^{-1}(C\widehat{\Omega}^{q+1}(\hat
E_k[i]))}\over{d\widehat{\Omega}^{q-1}(\hat
E_k[i])\oplus\{C\widehat{\Omega}^{q}(\hat E_k[i])\cap
d^{-1}(C\widehat{\Omega}^{q+1}(\hat
E_k[i])))\}}}\\
  & \equiv\oplus_{q\ge 0 }\hat{\frak I} (\hat
E_k[i])^{q}.\\
\end{array} \right.$$
Here $C\hat\Omega^q(\hat E_k[i])$ denotes the space
of all Cartan quantum $q$-forms on $\hat E_k[i]$. (See also
{\em \cite{PRAS14}}.) Then we define {\em quantum integral
characteristic supernumbers} of $N$, with
$[N]\in\widehat{\Omega}_{q}^{\hat E_k[i]}$, the numbers $\hat
i[N]\equiv<[N],[\alpha]>\in B$, for all $[\alpha]\in\hat{\frak
I}(\hat E_k[i])^{q}$.
\end{definition}
One has the following lemma.
\begin{lemma}\cite{PRAS14}
Let us assume that
$\hat{\frak I}(\hat E_k[i])^{q}\not=0$. One has a natural
homomorphism: ${\underline{j}}_{q}:\Omega_{q}^{\hat E_k[i]}\to
Hom_A(\hat{\frak I}(\hat E_k[i])^{q};A)$, $[N]\mapsto
{\underline{j}}_{q}([N])$,
${\underline{j}}_{q}([N])([\alpha])=\int_N\alpha\equiv<[N],[\alpha]>$.
Then, a necessary condition that $N'\in[N]\in\Omega_q^{\hat
E_k[i]}$ is the following: $\hat i[N]=\hat i[N']$,
$\forall[\alpha]\in\hat{\frak I}(\hat E_k[i])^{q}$. Furthermore,
if $N$ is orientable then above condition is sufficient also in
order to say that $N'\in[N]$.\end{lemma}

Therefore a quantum evaporation black-hole process can be described by means
of quantum smooth integral manifolds, and therefore for such a
process ``conservation laws" are not destroyed. By the way, as we
can have also weak solutions around a quantum black-hole, we can
assume also that interactions with such objects could be described
by means of weak-solutions, like shock-waves. Therefore we shall
more precisely talk of {\em weak quantum black-holes} and {\em
non-weak quantum black-holes}, according if they are described
respectively by means of weak solutions, or non-weak solutions. As a by-product we get that all the quantum integral
characteristic supernumbers are conserved through a non-weak quantum evaporating black-hole.
\end{proof}

Let us conclude this section by considering the following theorem that gives some important quantum conservation laws for what we shall develop in part II and in part III. (Some other quantum conservation laws will be considered in Appendix B and in part III too.)

\begin{theorem}[Observed quantum Hamiltonian and quantum $r$-momentum as quantum conservation laws of ${\widehat{(YM)}[i]}$]\label{quantum-observed-hamiltonian-as-quantum-conservation-law-of-observed-quantum-yang-mills-pde}
The observed quantum super Yang-Mills PDE $\widehat{(YM)}[i]$ admits the following important quantum conservation laws:
 \begin{equation}\label{quantum-conservation-laws-a}
   \left\{
   \begin{array}{l}
     \omega_H=(-1)^{\mu+1}[ (\partial y^\mu_\beta.L)y^\beta_4-\delta^\mu_4 L]\otimes  dx^1\wedge\cdots\wedge \widetilde{dx^\mu}\wedge\cdots\wedge dx^4\\
     p_r=(-1)^{\mu+1}[ (\partial y^\mu_\beta.L)y^\beta_r-\delta^\mu_r L]\otimes  dx^1\wedge\cdots\wedge \widetilde{dx^\mu}\wedge\cdots\wedge dx^4,\hskip 0.5cm r=1,2,3.\\
   \end{array}\right.
 \end{equation}
 The corresponding charges $H[i|t]=\int_{\sigma_t}\omega_H\in A$ and $P_r[i|t]=\int_{\sigma_t}p_r\in A$, $r=1,2,3$, are called respectively {\em observed phenomenological quantum energy} and {\em observed quantum $r$-momentum} of the nonlinear quantum propagator $V\subset  {\widehat{(YM)}[i]}$, $\partial V=N_0\bigcup P\bigcup N_1$, $\partial P=\partial N_0\bigcup\partial N_1$, or of the quantum particle encoded by the transverse section $\sigma_t\subset V$.

 In general, namely for any admissible nonlinear quantum propagator $V$ of ${\widehat{(YM)}[i]}$, $H[i|t]$ and $P_r[i|t]$ do not necessitate to be constant.
\end{theorem}
\begin{proof}
Let $\xi=\bar\xi^\mu\partial x_\mu+\xi^\beta\partial y_\beta$ be an infinitesimal symmetry of the fiber bundle $E[i]\to N$, (quantum fibered coordinates $(x^\alpha,y^\beta$), that encodes an infinitesimal symmetry of ${\widehat{(YM)}[i]}$. Then considered the quantum Lagrangian $L:J\hat D(E[i])\to A$ of ${\widehat{(YM)}[i]}$, we get the following quantum conservation law of ${\widehat{(YM)}[i]}$:
\begin{equation}\label{quantum-conservation-laws-b}
  \left\{
  \begin{array}{l}
    \beta=Z\rfloor\eta:N\to A\otimes\Lambda^0_3N,\, \eta=dx^1\wedge dx^2\wedge dx^3\wedge dx^4:N\to\Lambda^0_4N \\
   Z=Z^\mu\otimes\partial x_\mu:N\to A\otimes TN\\
   Z^\mu=[(\partial y^\mu_\beta.L)y^\beta_\sigma-\delta^\mu_\sigma L]\bar\xi^\sigma-(\partial y^\mu_\beta.L)\xi^\beta.\\
  \end{array}
  \right.
\end{equation}

More precisely one has:
\begin{equation}\label{conservation-laws-c}
  \beta=(-1)^{\mu+1}\left\{[ (\partial y^\mu_\beta.L)y^\beta_\sigma-\delta^\mu_\sigma L]\bar\xi^\sigma-(\partial y^\mu_\beta.L)\xi^\beta)]\right\}\otimes  dx^1\wedge\cdots\wedge \widetilde{dx^\mu}\wedge\cdots\wedge dx^4.
    \end{equation}
    This means that $d\beta|_V=0$. In particular for time-translations, namely with $\xi=\partial x_4$ and for $r$-space-translationsm namely with $\xi=\partial x_r$, we get the conservation laws reported in (\ref{quantum-conservation-laws-a}). Really these vector fields are infinitesimal symmetries of ${\widehat{(YM)}[i]}$, since these equations do not explicitly depend on the coordinate $x^\alpha$, $\alpha=1,2,3,4$. One can also directly verify that the conservations laws in (\ref{quantum-conservation-laws-a}) satisfy the conditions $\omega_H|_V=0$ and  $p_r|_V=0$. (For details see Appendix A.)

\begin{definition}[Defect phenomenological quantum energy and defect quantum $r$-momentum of an observed nonlinear quantum propagator of ${\widehat{(YM)}[i]}$]\label{lost-quantum-energy-a}
Given an observed nonlinear quantum propagator $V$ of $\widehat{(YM)}[i]$, such that $\partial V=N_0\sqcup P\sqcup N_1$, where $N_i$, $i=0,1$, are $3$-dimensional space-like admissible Cauchy data of $\widehat{(YM)}[i]$, and $P$ is a suitable time-like $3$-dimensional integral manifold with $\partial P=\partial N_0\sqcup\partial N_1$, we call respectively {\em defect phenomenological quantum energy}

\begin{equation}\label{defect-phenomenological-quantum-energy}
  \mathfrak{H}[V]_{\partial}=\int_P\omega_H
\end{equation}

and {\em defect quantum $r$-momentum}

\begin{equation}\label{defect-phenomenological-quantum-energy}
  \mathfrak{P}_r[V]=\int_P p_r,\, r=1,2,3.
\end{equation}
\end{definition}

 Despite $\omega_H$ and $p_r$ are conservation laws, it is not assured that their charges $H[i|t]$ and $P_r[i|t_r]$ are constant for any nonlinear quantum propagator. In fact, we get
$$0=\int_{V}d\omega_{H}|V=\int_{\partial V}\omega_{H}|_{\partial V}=\int_{N_0}\omega_{H}|_{N_0}-\int_{N_1}\omega_{H}|_{N_1}+\int_{P}\omega_{H}|_{P}.$$
We get
\begin{equation}\label{conserved-quantum-observed-hamiltonian-a}
 H[i|t_0]-H[i|t_1]=-\int_{P}\omega_{H}|_{P}\in A.
\end{equation}
Therefore, $ H[i|t_0]=H[i|t_1]$ iff $\int_{P}\omega_{H}|_{P}=0\in A$. In general the condition $\mathfrak{H}[V]_\partial=0$ is not verified. In fact, by considering that on $P$ the coordinates $(x^\alpha)$ are not more independent, and taking independent only $(x^1,x^2,x^4)$, in some neighbourhood of regular points of $P$, we can write
\begin{equation}\label{boundary-restrictions-conservation-laws}
  \int_P\omega_H|_P=\int_P[-x^3_4T^4_4-x^3_1T^1_4-x^3_2T^2_4+T^3_4]\otimes dx^1\wedge dx^2\wedge dx^4
\end{equation}
where $T^\alpha_\beta$ is the quantum energy-momentum tensor. (See in Appendix.) Thus if the $3$-dimensional manifold $P$ is orientable and smooth one can write
\begin{equation}\label{energy-momentum-condition-regularity}
x^3_4T^4_4+x^3_1T^1_4+x^3_2T^2_4=T^3_4,
\end{equation}
hence $\omega_H|_P=0$. However, since $V$ is in general a singular manifold, and also $P$ is so, it follows that we can not use (\ref{energy-momentum-condition-regularity}) for a generic nonlinear quantum propagator. Therefore, $\mathfrak{H}[V]_\partial $ does not necessitate to be zero. (For more details see Appendix C.)

A similar conclusion holds for $\mathfrak{P}_r[V]$. In fact, we get

 \begin{equation}\label{boundary-restrictions-conservation-laws-a}
  \int_Pp_r|_P=\int_P[-x^3_1T^1_r-x^3_2T^2_r-x^3_4T^4_r+T^3_r]\otimes dx^1\wedge dx^2\wedge dx^4.
\end{equation}
 Therefore, if $P$ is a $3$-dimensional orientable, smooth manifold one can write $x^3_1T^1_r+x^3_2T^2_r+x^3_4T^4_r=T^3_r$, hence $p_r|P=$ and $\mathfrak{P}_r=0$, $r=1,2,3$. In such a case $H[i|t]\in A$ and $P_r[i|t]\in A$ are constants. However, since $V$ and $P$ are in general singular manifolds in almost all the interesting interactions in particle physics, we cannot state that the quantum charges $H[i|t]\in A$ and $P_r[i|t]\in A$ are conserved in all quantum interactions.
\end{proof}
\begin{cor}[Criterion zero defect quantum phenomenological energy]\label{criterion-zero-defect-quantum-phenomenological-energy}
$H[i|t]$ is constant iff $\mathfrak{H}[V]_\partial=0$. This is surely the case when $P$ is a $3$-dimensional orientable smooth manifold.
\end{cor}

\begin{cor}[Criterion zero defect quantum $r$-momentum]\label{criterion-zero-defect-quantum-r-momentum}
$P_r[i|t]$, $r=1,2,3$, is constant iff $\mathfrak{P}_r[V]=0$. This is surely the case when $P$ is a $3$-dimensional orientable smooth manifold.
\end{cor}

\begin{example}[Steady-state nonlinear quantum propagators]\label{steady-state-nonlinear-quantum-propagators}
Let $V\subset {\widehat{(YM)}[i]}$ a nonlinear quantum propagator encoding a steady-state (or, in particular, a stationary state). These are solutions where $y^\beta_0=const$ (or $y^\beta_0=0$). From a geometrical point of view one has the following structure $V\cong B\times \Delta T$, namely $V$ it results into a cylinder where $B$ is a space-like $3$-dimensional manifold. Let us assume that $B$ is a $3$-dimensional orientable smooth manifold. Then $V$ is also a $4$-dimensional orientable smooth manifold and $\partial V=\partial B\times I=P$. Therefore $P$ is a $3$-dimensional orientable smooth manifold, and as a by product of Corollary \ref{criterion-zero-defect-quantum-phenomenological-energy} and Corollary \ref{criterion-zero-defect-quantum-r-momentum} we get that $H[i|t]_\partial\in A$ and $P_r[i|t]\in A$ are constants.
\end{example}

\begin{remark}[Quantum Hamiltonian vs phenomenological quantum hamiltonian]
Let us emphasize that $H[i|t]=\int_{\sigma_t}\omega_H$ does not represent the full observed quantum energy  content of the space-like particle encoded by $\sigma_t$. In fact $H[i|t]\not=\int_{\sigma_t}H\otimes dx^1\wedge dx^2\wedge dx^3\in A$, where
\begin{equation}\label{observed-quantum-hamiltonian}
H=(\partial y^\alpha_\beta.L)y^\beta_\alpha-L:{\widehat{(YM)}[i]}\to A.
\end{equation}
On the other hand in classical field theory the generalization of the classical Hamiltonian for particles, namely $H=(\partial \dot x_i.L)\dot x^i-L$ is interpreted just by $H=(\partial y^\alpha_b.L)y^\beta_\alpha-L$. Therefore also in quantum field theory we shall consider the full quantum energy encoded by (\ref{observed-quantum-hamiltonian}), instead of the observed phenomenological quantum Hamiltonian, $\omega_H$ only. Really $H$ allows to consider all the collective contributions, caused by the fact that particles are not point-like objects.
\end{remark}

\begin{theorem}[Observed quantum energy and observed phenomenological quantum energy]\label{observed-quantum-energy-and-observed-phenomenological-quantum-energy}
The relation between observed quantum energy and observed phenomenological quantum energy is given by the following
\begin{equation}\label{relation-observed-quantum-energy-and-observed-phenomenological-quantum-energy}
  H[i|t]=\int_{\sigma_t}\omega_H|_{\sigma_t}=H[i|t]_0-H[i|t]_{00}
  \end{equation}
  where
  $H[i|t]_0=\int_{\sigma_t} H\otimes dx^1\wedge dx^2\wedge dx^3\in A$ and $H[i|t]_{00}=\int_{\sigma_t}[(\partial y^j_\beta.L)y^\beta_j]\otimes dx^1\wedge dx^2\wedge dx^3\in A$.

$H[i|t]_0$ does not necessitate to be constant. In fact one has
\begin{equation}\label{relation-observed-quantum-energy-and-observed-phenomenological-quantum-energy-a}
  H[i|t_0]_0-H[i|t_1]_0=H[i|t_0]_{00}-H[i|t_1]_{00}-\int_{P}\omega_H|_P.
  \end{equation}
  Therefore, $H[i|t_0]_0=H[i|t_1]_0$ iff $H[i|t_0]_{00}-H[i|t_1]_{00}-\int_{P}\omega_H|_P=0$.

\end{theorem}
\begin{proof}
The proof follows directly from above definitions and results.
\end{proof}

\begin{definition}[Defect observed quantum energy of an observed nonlinear quantum propagator of ${\widehat{(YM)}[i]}$]\label{lost-quantum-energy-b}
Given an observed nonlinear quantum propagator $V$ of $\widehat{(YM)}[i]$, such that $\partial V=N_0\sqcup P\sqcup N_1$, where $N_i$, $i=0,1$, are $3$-dimensional space-like admissible Cauchy data of $\widehat{(YM)}[i]$, and $P$ is a suitable time-like $3$-dimensional integral manifold with $\partial P=\partial N_0\sqcup\partial N_1$, we call {\em defect observed quantum energy}

\begin{equation}\label{defect-phenomenological-quantum-energy}
  \mathfrak{H}[V]=H[i|t_0]_{00}-H[i|t_1]_{00}-\mathfrak{H}[V]_{\partial}.
\end{equation}
\end{definition}
\begin{cor}[Criterion zero defect observed quantum energy]\label{criterion-zero-defect-observed-quantum-energy}
$H[i|t]_0$ is constant iff $\mathfrak{H}[V]=0$.

We can also write
\begin{equation}\label{conserved-quantum-observed-hamiltonian-in-quantum observed-nonlinear-propagator}
 H[i|t_0]_0=H[i|t_1]_0 \hskip 3pt {\rm mod}\hskip 3pt \mathfrak{H}[V]\in A
\end{equation}
\end{cor}

\begin{remark}
Corollary \ref{criterion-zero-defect-observed-quantum-energy} gives a precise meaning to the phenomenological statement that with respect to an observer a quantum process must conserve the total mass-energy.
\end{remark}

\begin{example}\label{stationary-states-conserved-observed-quantum-energy}
If the nonlinear quantum propagator $V\subset \widehat{(YM)}[i]$, is such that it encodes a stationary-state, then $H[i|t]_0$ is constant. In fact in such a case one has that $\mathfrak{H}[V]_{\partial}=0$ and $H[i|t_0]_{00}-H[i|t_1]_{00}=0$, hence $\mathfrak{H}[V]=0$, and from Corollary \ref{criterion-zero-defect-observed-quantum-energy} we can conclude that $H[i|t]_0$ is constant.
\end{example}

\begin{example}[Strong reactions with jet quenching in  ${\widehat{(YM)}[i]}$]\label{jet-quenching-in-observed-quantum-yang-mills-pde}
Nonlinear quantum propagators with defect quantum energy can encode strong reactions between ultra-relativistic heavy-ion collisions, where interactions between the high-momentum parton and the hot, dense medium produced in the collisions, lead to loss of energy. This phenomenon is called {\em jet quenching}.\footnote{Parton model was proposed by R. Feynmann in 1969 for high-energy hadron collisions and actually usually referred as quark-gluon model. Nowadays there exist experimental evidences for the quenching phenomenon. (See, e.g., \cite{AAD-ET-AL} and CERN-Press Releases reported therein.)}
\end{example}

\section{\bf Quantum hypercomplex exotic super PDE's}\label{quantum-exotic-super-pdes-section}

In this section we consider PDE's in the category $\mathfrak{Q}_{hyper,S}$ of quantum hypercomplex supermanifolds, as defined in \cite{PRAS27}, and focus our attention on ``quantum exotic super PDE's", i.e., quantum super PDE's where we can embed ``quantum exotic supersheres". For such Cauchy data we will generalize our previous results on quantum exotic PDE's \cite{PRAS27}. Such ``exotic" boundary value problems are of particular interest in strong reactions encoding quantum processes occurring in high energy physics, as we will prove in part II and in part III.
\begin{remark}
Let us remark that in the other sections of this paper we refer to the category  $\mathfrak{Q}_{S}$ instead that  $\mathfrak{Q}_{hyper,S}$. This is made for three reasons. The first is for convenience, since quantum micro-worlds can be encoded in $\mathfrak{Q}_{S}$. The second is that one can directly generalize intrinsic results obtained in the geometry of PDEs in the category $\mathfrak{Q}_{S}$ to similar ones in the category $\mathfrak{Q}_{hyper,S}$, as we have proved in \cite{PRAS27}, and we will also see in this section. The third reason is that nonassociative algebras, as arise for example in some quantum hypercomplex algebras, can be considered as subsets of their enveloping algebras. These last being necessarily associative, give a general criterion to encode nonassociative algebras in larger associative frameworks.
\end{remark}

\begin{definition}[Quantum homotopy $(m|n)$-supersphere]\label{quantum-homotopy-m-n-supersphere}
We call quantum homotopy $(m|n)$-supersphere (with respect to a quantum hypercomplex superalgebra $A$) a smooth, compact, closed $(m|n)$-dimensional quantum supermanifold $M\equiv\hat\Sigma^{m|n}$, that is homotopy equivalent to the ($m|n)$-dimensional quantum supersphere $\hat S^{m|n}$, with classic regular structure $\bar\pi_C:M\to M_C$, where $M_C$ is a homotopy $m$-sphere, and such that the homotopy equivalence between $M$ and $\hat S^{m|n}$ is realized by a commutative diagram {\em(\ref{commutative-diagram-homotopy-equivalence-quantum-m-n-supersphere})}.
\end{definition}

\begin{equation}\label{commutative-diagram-homotopy-equivalence-quantum-m-n-supersphere}
\xymatrix@C=1cm{M\ar[d]_{\bar\pi_C}\ar[r]^{f}&\hat S^{m|n}\ar[d]^{\pi_C}\ar@{=}[r]&A^{m|n}\bigcup\{\infty\}\ar[d]\ar@{=}[r]&(A_0^m\times A_1^n)\bigcup\{\infty\}\ar[d]\\
M_C\ar[r]_{f_C}&S^m\ar@{=}[r]&\mathbb{K}^m\bigcup\{\infty_C\}\ar@{=}[r]&\mathbb{K}^m\bigcup\{\infty_C\}\\}
\end{equation}
$\mathbb{K}=\mathbb{R},\, \mathbb{C}$ and $\pi_C$ is induced by the canonical mappings $c:A\to\mathbb{K}$ and $\infty\mapsto\infty_C$.
\begin{remark}\label{remark-quantum-superspheres-classic-limits}
Let $\hat\Sigma^{m|n}_1$ and $\hat\Sigma^{m|n}_2$ be two quantum diffeomorphic, quantum homotopy $(m|n)$-superspheres: $\hat\Sigma^{m|n}_1\cong\hat\Sigma^{m|n}_2$. Then the corresponding classic limits $\hat\Sigma^{m|n}_{1,C}$ and $\hat\Sigma^{m|n}_{2,C}$ are diffeomorphic too: $\hat\Sigma^{m|n}_{1,C}\cong\hat\Sigma^{m|n}_{2,C}$. This remark is the natural consequence of the fact that quantum diffeomorphisms here considered respect the fiber bundle structures of quantum homotopy $(m|n)$-superspheres with respect their classic limits: $\pi_C:\hat\Sigma^{m|n}\to\hat\Sigma^{m|n}_C$. Therefore quantum diffeomorphisms between quantum homotopy $(m|n)$-superspheres are characterized by a couple $(f,f_C):(\hat\Sigma_1^{m|n},\hat\Sigma_{1,C}^{m|n})\to(\hat\Sigma_2^{m|n},\hat\Sigma_{2,C}^{m|n})$ of mappings related by the commutative diagram in {\em(\ref{commutative-diagram-a-lemma-quantum-spheres-classic-limits-a})}.
\begin{equation}\label{commutative-diagram-a-lemma-quantum-spheres-classic-limits-a}
 \xymatrix@C=3cm{\hat\Sigma^{m|n}_1\ar[d]_{\pi_{1,C}}\ar[r]^{f}&\hat\Sigma^{m|n}_2\ar[d]^{\pi_{2,C}}\\
 \hat\Sigma^{m|n}_{1,C}\ar[r]_{f_C}&\hat\Sigma^{m|n}_{2,C}\\}
\end{equation}
There $f$ is a quantum diffeomorphism between quantum supermanifolds and $f_C$ is a diffeomorphism between manifolds. Note that such diffeomorphisms of quantum homotopy $(m|n)$-superspheres allow to recognize that $\hat\Sigma^{m|n}_1$ has also $\hat\Sigma^{m|n}_{2,C}$ as classic limit, other than $\hat\Sigma^{m|n}_{1,C}$. (See commutative diagram in {\em(\ref{commutative-diagram-b-lemma-quantum-spheres-classic-limits-a})}.)
\begin{equation}\label{commutative-diagram-b-lemma-quantum-spheres-classic-limits-a}
 \xymatrix@R=1cm@C=2cm{&\hat\Sigma^{m|n}_1\ar[ddl]^(0.7){\pi_{1,C}}\ar[ddr]_(0.7){\pi'_{1,C}}\ar[dr]^{f}&\\
 \hat\Sigma^{m|n}_2\ar[ur]^{f^{-1}}\ar[d]_{\pi'_{2,C}}\ar[rr]^{1_{\hat\Sigma^{m|n}_2}}&&\hat\Sigma^{m|n}_2\ar[d]^{\pi_{2,C}}\\
 \hat\Sigma^{m|n}_{1,C}\ar[rr]_{f_C}&&\hat\Sigma^{m|n}_{2,C}\\}
\end{equation}

This clarifies that the classic limit of a quantum homotopy $(m|n)$-supersphere is unique up to diffeomorphisms.
\end{remark}
Let us also emphasize that (co)homology properties of quantum homotopy $(m|n)$-superspheres are related to the ones of $m$-spheres, since we here consider classic regular objects only.

\begin{lemma}\label{cohomology-quantum-homotopy-m-n-superspheres}
In {\em(\ref{cohomology-properties-quantum-homotopy-m-n-superspheres})} are reported the cohomology spaces for quantum homotopy $(m|n)$-superspheres.

\begin{equation}\label{cohomology-properties-quantum-homotopy-m-n-superspheres}
  H^p(\hat\Sigma^{m|n};\mathbb{Z})\cong  H^p(\hat S^{m|n};\mathbb{Z})\cong  H^p(S^{m};\mathbb{Z})=\left\{\begin{array}{ll}
                                                                                                       0 & p\not=0,\, m\\
                                                                                                       \mathbb{Z}& p=0,\, m.\\
                                                                                                     \end{array}
  \right.
\end{equation}
\end{lemma}

\begin{proof}
Let us first calculate the homology groups in integer coefficients $\mathbb{Z}$, of quantum $(m|n)$-superspheres. In (\ref{homology-properties-quantum-m-n-superspheres}) are reported the homology spaces for quantum $m\not=0$.
\begin{equation}\label{homology-properties-quantum-m-n-superspheres}
  H_p(\hat S^{m|n};\mathbb{Z})
  \cong  H_p(S^{m};\mathbb{Z})=\left\{\begin{array}{ll}
  0 & p\not=0,\, m\\
  \mathbb{Z}& p=0,\, m.\\
  \end{array}
  \right.
\end{equation}
Furthermore, for $m=0$ we get
$$ H_p(\hat S^{0|n};\mathbb{Z})
  \cong  H_p(S^{0};\mathbb{Z})=\left\{\begin{array}{ll}
  0 & p\not=0\\
  \mathbb{Z}\bigoplus\mathbb{Z}& p=0.\\
  \end{array}\right.$$
  Above formulas can be obtained by the reduced Mayer-Vietoris sequence applied to the triad $(\hat S^{m|n},\hat D^{m|n}_+,\hat D^{m|n}_-)$ since we can write $\hat S^{m|n}=\hat D^{m|n}_+\bigcup\hat D^{m|n}_-$, where $\hat D^{m|n}_+$ and $\hat D^{m|n}_-$ are respectively the north quantum $(m|n)$-superdisk and south quantum $(m|n)$-superdisk that cover $\hat S^{m|n}$. Taking into account that $\hat D^{m|n}_+\bigcap\hat D^{m|n}_-=\hat S^{m-1|n-1}$, we get the long exact sequence (\ref{reduced-quantum-m-n-supersphere-mayer-vietoris-exact-sequence}).

  \begin{equation}\label{reduced-quantum-m-n-supersphere-mayer-vietoris-exact-sequence}
   \scalebox{0.7}{$ \xymatrix{\cdots\ar[r]&\widetilde{ H}_{p}(\hat S^{m-1|n-1};\mathbb{Z})\ar[r]&\widetilde{ H}_{p}(\hat D_+^{m|n};\mathbb{Z})\bigoplus \widetilde{ H}_{p}(\hat D_-^{m|n};\mathbb{Z})\ar[r]&\widetilde{H}_{p}(\hat S^{m|n};\mathbb{Z})\ar[d]_{\partial}\\
    &\widetilde{ H}_{p-1}(\hat S^{m|n};\mathbb{Z})\ar[d]_{\partial}&\ar[l]\widetilde{ H}_{p-1}(\hat D_+^{m|n};\mathbb{Z})\bigoplus \widetilde{\hat H}_{p-1}(\hat D_-^{m|n};\mathbb{Z})&\ar[l]\widetilde{ H}_{p-1}(\hat S^{m-1|n-1};\mathbb{Z})\\
    &\widetilde{ H}_{p-2}(\hat S^{m-1|n-1};\mathbb{Z})\ar[r]&\widetilde{ H}_{p-2}(\hat D_+^{m|n};\mathbb{Z})\bigoplus \widetilde{ H}_{p-2}(\hat D_-^{m|n};\mathbb{Z})\ar[r]&\widetilde{H}_{p-2}(\hat S^{m|n};\mathbb{Z})\ar[d]_{\partial}\\
    &\vdots\ar[d]_{\partial}&&\\
    &\widetilde{ H}_{0}(\hat S^{m-1|n-1};\mathbb{Z})\ar[r]&\widetilde{ H}_{0}(\hat D_+^{m|n};\mathbb{Z})\bigoplus \widetilde{ H}_{0}(\hat D_-^{m|n};\mathbb{Z})\ar[r]&\widetilde{ H}_{0}(\hat S^{m|n};\mathbb{Z})\ar[r]&0.\\}$}
  \end{equation}
Taking into account that $\widetilde{ H}_{0}(\hat D_-^{m|n};\mathbb{Z})=0$ e get $\widetilde{ H}_{p}(\hat S^{m|n};\mathbb{Z})\cong\widetilde{H}_{p-1}(\hat S^{n-1|n-1};\mathbb{Z})$ and $\widetilde{ H}_{0}(\hat S^{m|n};\mathbb{Z})\cong 0$.
Therefore, we get
\begin{equation}\label{reduced-homology-properties-quantum-m-n-superspheres}
  \widetilde{ H}_p(\hat S^{m|n};\mathbb{Z})
  \cong  \left\{\begin{array}{ll}
  \mathbb{Z} & \hbox{\rm if $p=m$}\\
  0& \hbox{\rm if $p\not=m$}\\
  \end{array}
  \right\}\Rightarrow\, \left\{\begin{array}{l}
   H_p(\hat S^{0|n};\mathbb{Z}) =\left\{ \begin{array}{ll}
  \mathbb{Z}\bigoplus\mathbb{Z} & \hbox{\rm if $m=0$}\\
  0& \hbox{\rm if $p\not=0$}\\
  \end{array}\right. \\
   H_p(\hat S^{m|n};\mathbb{Z}) =\left\{ \begin{array}{ll}
  \mathbb{Z} & \hbox{\rm if $p=0,\, m$}\\
  0& \hbox{\rm if $p\not=0,\, m$.}\\
  \end{array}\right. \\
  \end{array}
  \right.
\end{equation}

Therefore we get formulas (\ref{homology-properties-quantum-m-n-superspheres}). To conclude the proof we shall consider that $ H^p(\hat S^{m|n};\mathbb{Z})\cong Hom_{\mathbb{Z}}( H_p(\hat S^{m|n};\mathbb{Z});\mathbb{Z})$. Furthermore, quantum homotopy $(m|n)$-superspheres have same (co)homology of quantum superspheres since are homotopy equivalent to these last ones.

\end{proof}

\begin{lemma}\label{quantum-euler-characteristic-quantum-superspheres-classic-limits-b}
The {\em quantum Euler characteristic numbers} for quantum homotopy $(m|n)$-superspheres are reported in {\em(\ref{euler-characteristic-properties-quantum-homotopy-superspheres})}. These coincide with the corresponding quantum Euler characteristic numbers of quantum $(m|n)$-superspheres and with the Euler characteristic numbers of usual $m$-spheres. Furthermore they are the same of the corresponding {\em total quantum Euler characteristic numbers}. See in {\em(\ref{euler-characteristic-properties-quantum-homotopy-superspheres})}.

\begin{equation}\label{euler-characteristic-properties-quantum-homotopy-superspheres}
\left\{\begin{array}{ll}
\hat \chi(\hat\Sigma^{m|n})&= \hat \chi(\hat S^{m|n})=\chi(S^{m})=(-1)^0\beta_0+(-1)^m\beta_m=1+(-1)^m\\
&=\left\{\begin{array}{ll}
0 & \hbox{\rm $m=$ odd}\\
2& \hbox{\rm $m=$ even.}\\
\end{array}
  \right\}\\
  &={}^{Tot}\hat \chi(\hat\Sigma^{m|n})={}^{Tot}\hat \chi(\hat S^{m|n}).\\
  \end{array}
  \right.
\end{equation}
\end{lemma}

\begin{proof}
We have considered that $\hat S^{m|n}$ admits the following quantum-supercell decomposition: $\hat S^{m|n}=\hat e^{m|n}\bigcup \hat e^{0|0}$, where $\hat e^{m|n}=\hat D^{m|n}$ is a $(m|n)$-dimensional quantum supercell, with respect to the quantum superalgebra $A$, and $\hat e^{0|0}=\hat D^{0|0}$ is the $(0|0)$-dimensional quantum supercell with respect to $A$. Therefore we can consider the {\em quantum homological Euler characteristic} $\hat\chi(\hat S^{m|n})$ of $\hat S^{m|n}$, given by formulas (\ref{quantum-homology-euler-characteristic-properties-quantum-m-n-superspheres}).
\begin{equation}\label{quantum-homology-euler-characteristic-properties-quantum-m-n-superspheres}
\begin{array}{ll}
\hat\chi(\hat S^{m|n})&=(-1)^0\dim_A H_0(\hat S^{m|n};A)+(-1)^m\dim_A H_m(\hat S^{m|n};A)\\
&=(-1)^0\dim_A A+(-1)^m\dim_A A\\
&=1+(-1)^m\\
&=\left\{\begin{array}{ll}
0 & \hbox{\rm $m=$ odd}\\
2& \hbox{\rm $m=$ even.}\\
\end{array}\right.\\
\end{array}
\end{equation}
 So the homological quantum Euler characteristic of the quantum $(m|n)$-supersphere is the same of the homological Euler characteristic of the usual $m$-sphere. Furthermore, since quantum homotopy $(m|n)$-superspheres are homotopy equivalent to quantum $(m|n)$-superspheres, it follows that the quantum Euler characteristic of a quantum homotopy $(m|n)$-supersphere is equal to the one of $\hat S^{m|n}$. Moreover, also the total quantum Euler characteristic numbers for quantum (homotopy) $(m|n)$-superspheres coincide with the ones of $S^m$. In fact, we can consider the {\em quantum total-homological Euler characteristic} ${}^{Tot}\hat\chi(\hat S^{m|n})$ of $\hat S^{m|n}$, given by formulas (\ref{total-quantum-homology-euler-characteristic-properties-quantum-m-n-superspheres}).
\begin{equation}\label{total-quantum-homology-euler-characteristic-properties-quantum-m-n-superspheres}
\scalebox{0.8}{$\begin{array}{ll}
{}^{Tot}\hat\chi(\hat S^{m|n})&=\sum_{p\ge 0}(-1)^p\dim_A {}^{Tot}H_p(\hat S^{m|n};A)\\
&=\sum_{p\ge 0}(-1)^p\dim_A[ \bigoplus_{r+s=p}H_{r|s}(\hat S^{m|n};A)]\\
&=\sum_{p\ge 0}(-1)^p\dim_A[ \bigoplus_{r+s=p}(H_{r}(\hat S^{m|n};A_0)\oplus H_{s}(\hat S^{m|n};A_1)])\\
&=\sum_{p\ge 0}(-1)^p\dim_A[ \bigoplus_{r+s=p}(H_{r}(\hat S^{m|n};\mathbb{Z})\otimes_{\mathbb{Z}}A_0)\oplus H_{s}(\hat S^{m|n};\mathbb{Z})\otimes_{\mathbb{Z}}A_1)])\\
&=\sum_{p\ge 0}(-1)^p\dim_A[ \bigoplus_{r+s=p}(H_{r}(S^{m};\mathbb{Z})\otimes_{\mathbb{Z}}A_0)\oplus H_{s}( S^{m};\mathbb{Z})\otimes_{\mathbb{Z}}A_1)])\\
&=1+(-1)^m.\\
\end{array}$}
\end{equation}
 So also the total homological quantum Euler characteristic of the quantum $(m|n)$-supersphere is the same of the homological Euler characteristic of the usual $m$-sphere. Furthermore, since quantum homotopy $(m|n)$-superspheres are homotopy equivalent to quantum $(m|n)$-superspheres, it follows that the total quantum Euler characteristic of a quantum homotopy $(m|n)$-supersphere is equal to the one of $\hat S^{m|n}$.
 \end{proof}

\begin{theorem}[Generalized Poincar\'e conjecture in the category $\mathfrak{Q}_{hyper,S}$]\label{generalized-poincare-conjecture-in-category-hypercomplex-quantum-super-manifolds}
Let $A$ be a hypercomplex quantum superalgebra with center $Z(A)$ a Noetherian
$\mathbb{K}$-algebra, $\mathbb{K}=\mathbb{R}$ or
$\mathbb{K}=\mathbb{C}$. Let $M$ be a classic regular, closed
compact quantum supermanifold of dimension $(m|n)$, in the category $\mathfrak{Q}_{hyper,S}$, homotopy
equivalent to $\hat S^{m|n}$, (hence this last is the quantum super
CW-substitute of $M$). Then we get that $M\approx \hat S^{m|n}$, i.e.,
$M$ is also homeomorphic to $\hat S^{m|n}$, and $M_C\approx S^m$,
i.e., the classic limit $M_C$ of $M$ is homeomorphic to the classic
limit $S^m$ of $\hat S^{m|n}$.
\end{theorem}

\begin{proof}
In \cite{PRAS19} we have proved the generalized quantum Poincar\'e conjecture in the category $\mathfrak{Q}_S$ of quantum supermanifolds, by considering the quantum Ricci flow PDE just in the category $\mathfrak{Q}_S$. Then, by using similar arguments to prove Theorem 3.10 and Theorem 3.11 in \cite{PRAS27}, we can state that generalized quantum Poincar\'e conjecture works also in the category $\mathfrak{Q}_{hyper,S}$ of quantum hypercomplex supermanifolds. Therefore the quantum $(m|n)$-supersphere $\hat S^{m|n}$, considered the quantum super CW-substitute of any quantum homotopy $(m|n)$-supersphere $\hat\Sigma^{m|n}$, is just homeomorphic to this last one: $\hat\Sigma^{m|n}\thickapprox \hat S^{m|n}$.
\end{proof}

\begin{theorem}\label{theorem-quantum-superspheres-classic-limits-b}
Let $\hat\Theta_{m|n}$ be the set of equivalence classes of quantum diffeomorphic quantum homotopy $(m|n)$-superspheres over a quantum (hypercomplex) superalgebra $A$ (and with Noetherian center $Z(A)$).\footnote{Quantum diffeomorphisms are meant in the sense specified in Remark \ref{remark-quantum-superspheres-classic-limits}.} In $\hat\Theta_{m|n}$ it is defined an additive commutative and associative composition map such that $[\hat S^{m|n}]$ is the zero of the composition. Then one has the exact commutative diagram reported in {\em(\ref{short-exact-sequence-lemma-quantum-superspheres-classic-limits-b})}.
\begin{equation}\label{short-exact-sequence-lemma-quantum-superspheres-classic-limits-b}
 \xymatrix{&&&0\ar[d]&\\
 0\ar[r]&\hat\Upsilon_{m|n}\ar@{^{(}->}[r]&\hat\Theta_{m|n}\ar[r]^{j_C}&\Theta_m\ar[d]\ar[r]&0\\
 &&&\hat\Theta_{m|n}/\hat\Upsilon_{m|n}\ar[d]&\\
 &&&0&\\}
\end{equation}
where $\Theta_m$ is the set of equivalence classes for diffeomorphic homotopy $m$-spheres and $j_C$ is the canonical mapping $j_C:[\hat\Sigma^{m|n}]\mapsto[\hat\Sigma^{m|n}_C]$. One has the canonical isomorphisms:
\begin{equation}\label{canonical-isomorphisms-in-short-exact-sequence-lemma-quantum-superspheres-classic-limits-b}
   \mathbb{Z}\bigotimes_{\hat\Upsilon_{m|n}} \mathbb{Z} \hat\Theta_{m|n}\cong\mathbb{Z} \Theta_m,\, \hbox{\rm as right $\hat\Theta_{m|n}$-modules}.
\end{equation}
\end{theorem}

\begin{proof}
After above Remark \ref{remark-quantum-superspheres-classic-limits} we can state that the mapping $j_C$ is surjective. In other words we can write
$$\hat\Theta_{m|n}=\bigcup_{[\hat\Sigma^{m|n}_C]\in\Theta_m}(\hat\Theta_{m|n})_{[\hat\Sigma^{m|n}_C]}.$$
The fiber $(\hat\Theta_{m|n})_{[\hat\Sigma^{m|n}_C]}$ is given by all classes $[\hat\Sigma^{m|n}]$ such their classic limits are diffeomorphic, hence belong to the same class in $\Theta_m$. Furthermore one has $\ker(j_C)=j_C^{-1}([S^m])\equiv\hat\Upsilon_{m|n}\subset\hat\Theta_{m|n}$. Therefore, we can state that $\hat\Theta_{m|n}$ is an extension of $\Theta_m$ by $\hat\Upsilon_{m|n}$. Such extensions are classified by $H^2(\Theta_m;\hat\Upsilon_{m|n})$.\footnote{In Tab. \ref{homologies-ciclic-group} are reported useful formulas to explicitly calculate these groups.}

\begin{table}[h]
\caption{Homology of finite cyclic group $\mathbb{Z}_i$ of order $i$.}
\label{homologies-ciclic-group}
\begin{tabular}{|c|c|}
  \hline
  \hfil{\rm{\footnotesize $ r$}}\hfil& \hfil{\rm{\footnotesize $H_r(\mathbb{Z}_i;\mathbb{Z})$}}\hfil\\
  \hline
 \hfil{\rm{\footnotesize $0$}}\hfil& \hfil{\rm{\footnotesize $\mathbb{Z}$}}\hfil\\
 \hline
  \hfil{\rm{\footnotesize $r$ odd}}\hfil& \hfil{\rm{\footnotesize $\mathbb{Z}_i$}}\hfil\\
\hline
  \hfil{\rm{\footnotesize $r>0$ even}}\hfil& \hfil{\rm{\footnotesize $0$}}\hfil\\
\hline
\end{tabular}
\end{table}
The composition map in $\hat\Theta_{m|n}$ is defined by {\em quantum fibered connected sum}, i.e., a connected sum on quantum supermanifolds that respects the connected sum on their corresponding classic limits. More precisely let $M\to M_C$ and $N\to N_C$ be connected $(m|n)$-dimensional classic regular quantum supermanifolds. We define quantum fibered connected sum of $M$ and $N$ the classic regular $(m|n)$-dimensional quantum supermanifold $M\sharp N\to M_C\sharp N_C$, where
\begin{equation}\label{quantum-fibered-connected-sum}
 \left\{
 \begin{array}{ll}
   M\sharp N&=(M\setminus\hat D^{m|n})\bigcup(\hat S^{m-1|n-1}\times\hat D^{1|1})\bigcup(N\setminus\hat D^{m|n})\\
   &\\
   M_C\sharp N_C&=(M_C\setminus D^m)\bigcup( S^{m-1}\times D^1)\bigcup(N_C\setminus D^m).\\
 \end{array}
 \right.
\end{equation}
Then the additive composition law is $+:\hat\Theta_{m|n}\times\hat\Theta_{m|n}\to\hat\Theta_{m|n}$, $[M]+[N]=[M\sharp N]$. $[\hat S^{m|n}]$ is the zero of this addition. In fact, since $\overline{\hat S^{m|n}\setminus\hat D^{m|n}}\bigcup_{\hat S^{m-1|n-1}}(\hat S^{m-1|n-1}\times\hat D^{1|1})\cong\hat D^{m-1|n-1}$, we get
\begin{equation*}
    \begin{array}{ll}
      M\sharp\hat S^{m|n}&\cong\overline{M\setminus\hat D^{m|n}}\bigcup_{\hat S^{m-1|n-1}}(\overline{\hat S^{m|n}\setminus\hat D^{m|n}}\bigcup_{\hat S^{m-1|n-1}}(\hat S^{m-1|n-1}\times\hat D^{1|1}))\\
      &\\
      &\cong\overline{M\setminus\hat D^{m|n}}\bigcup_{\hat S^{m-1|n-1}}\hat D^{m|n}\cong M.\\
    \end{array}
\end{equation*}
Analogous calculus for $M_C$ completes the proof.
\end{proof}

In the following remark we will consider some examples and further results to better understand some relations between quantum homotopy superspheres and their classic limits.

\begin{example}[Quantum homotopy $(7|n)$-supersphere]\label{quantum-homotopy-7-n-superspheres-with-hypercomplex-quantum-superalgebra-a} Let us calculate the extension classes for quantum homotopy $(7|n)$-superspheres in the category $\mathfrak{Q}_{hyper,S}$, assumed classic regular, i.e., having the fiber bundle structure $\pi_C:\hat\Sigma^{7|n}\to \Sigma^{7}$.
\begin{equation}\label{short-exact-sequence-lemma-quantum-7-n-superspheres-classic-limits-c}
 \xymatrix{0\ar[r]&\hat\Upsilon_{7|n}\ar@{^{(}->}[r]&\hat\Theta_{7|n}\ar[r]^{j_C}&\Theta_7\ar[r]&0},\hskip 3pt n\ge 0.
\end{equation}
 These are given by $H^2(\Theta_7;\hat\Upsilon_{7|n})\cong H^2(\mathbb{Z}_{28};\hat\Upsilon_{7|n})$. We get
 \begin{equation}\label{calculation-short-exact-sequence-lemma-quantum-7-n-superspheres-classic-limits-c}
    H^2(\mathbb{Z}_{28};\hat\Upsilon_{7|n}) =Hom_{\mathbb{Z}}(H_2(\mathbb{Z}_{28};\mathbb{Z});\hat\Upsilon_{7|n})
    =Hom_{\mathbb{Z}}(0;\hat\Upsilon_{7|n})\bigoplus Ext_{\mathbb{Z}}(H_1(\mathbb{Z}_{28};\mathbb{Z});\hat\Upsilon_{7|n}).
 \end{equation}
We shall prove that $Ext_{\mathbb{Z}}(H_1(\mathbb{Z}_{28};\mathbb{Z});\hat\Upsilon_{7|n})=\hat\Upsilon_{7|n}/28\cdot \hat\Upsilon_{7|n}$. (We have used the fact that $H_2(\mathbb{Z}_{28};\mathbb{Z})=0$.) Let us look in some detail to this $\mathbb{Z}$-module. By using
the projective resolution of $\mathbb{Z}_{28}$ given in {\em(\ref{projective-resolution-z28})},
\begin{equation}\label{projective-resolution-z28}
    \xymatrix{0\ar[r]&\mathbb{Z}\ar[r]^{\mu=.28}&\mathbb{Z}\ar[r]^{\epsilon}&\mathbb{Z}_{28}\ar[r]&0}
\end{equation}
we get the exact sequence {\em(\ref{projective-resolution-z28-derived})}.

\begin{equation}\label{projective-resolution-z28-derived}
    \xymatrix{0\ar[r]&Hom_{\mathbb{Z}}(\mathbb{Z}_{28};\hat\Upsilon_{7|n})\ar@{=}[d]\ar[r]^{\epsilon_*}&
    Hom_{\mathbb{Z}}(\mathbb{Z};\hat\Upsilon_{7|n})\ar@{=}[d]^{\wr}\ar[r]^{\mu_*}&Hom_{\mathbb{Z}}(\mathbb{Z};\hat\Upsilon_{7|n})\ar@{=}[d]^{\wr}\\
    0\ar[r]&Hom_{\mathbb{Z}}(\mathbb{Z}_{28};\hat\Upsilon_{7|n})\ar[r]^(0.6){\epsilon_*}&
    \hat\Upsilon_{7|n}\ar[r]^{\mu_*}&\hat\Upsilon_{7|n}\\}
\end{equation}
Therefore we get
$$Ext_{\mathbb{Z}}(\mathbb{Z}_{28};\hat\Upsilon_{7|n})=\hat\Upsilon_{7|n}/\IM(\mu_*).$$
In order to see what is $\IM(\mu_*)$ we can use analogous considerations made in Example 5.6 in \cite{PRAS27}. We get $\IM(\mu^*)=28.\hat\Upsilon_{7|n}$, hence $Ext_{\mathbb{Z}}(\mathbb{Z}_{28};\hat\Upsilon_{7|n})=\hat\Upsilon_{7|n}/28\cdot\hat\Upsilon_{7|n}$. The particular structure of this module, depends on the particular hypercomplex quantum superalgebra $A$ considered.
For example, take $A=\mathbb{C}$. One has $\hat S^{7|n}=\hat S^7$, hence, since $\hat S^7\to S^7$, is just the fiber bundle $S^{14}\to S^7$, we can easily copy the result in Example 5.6 in \cite{PRAS27}, to conclude that $\hat\Theta_{7|n}=\hat\Theta_{7}=\mathbb{Z}_2\bigoplus\mathbb{Z}_{28} $, $\forall n\ge 0$.
\end{example}

\begin{example}[Quantum homotopy $(m|n)$-superspheres for the limit case $A=\mathbb{R}$]\label{quantum-homotopy-m-n-superspheres-with-real-numbers-quantum-algebra}
In the limit case where the quantum algebra is $A=\mathbb{R}$, for a quantum homotopy $(m|n)$-supersphere $\hat\Sigma^{m|n}$ one has just $\hat\Sigma^{m|n}=\hat\Sigma^{m|n}_C=\Sigma^m$, hence $\pi_C=id_{\Sigma^m}$. Furthermore $\hat\Theta_{m|n}=\Theta_m$ and $\hat\Upsilon_{m|n}=0=[S^m]\in\Theta_m$. In particular if $m=\{1,2,3,4,5,6\}$, we get $\hat\Theta_{m|n}=\Theta_m=\hat\Upsilon_{m|n}=0$, $\forall n\ge 0$. (For the smooth case $m=4$ see \cite{PRAS26}.)
\end{example}

\begin{theorem}[Homotopy groups of quantum $(m|n)$-supersphere]\label{homotopy-groups-quantum-m-n-supersphere}
Quantum homotopy $(m|n)$-superspheres cannot have, in general, the same homotopy groups of $m$-spheres:\footnote{In other words, quantum homotopy $(m|n)$-superspheres are not homotopy equivalent to the $m$-sphere.}
\begin{equation}\label{non-isomorphism-homotopy-groups-quantum-m-n-supersphere-sphere}
    \pi_{k}(\hat\Sigma^{m|n})\cong\pi_{k}(\hat S^{m|n})\not=\pi_{k}(S^m).
\end{equation}

Furthermore, $S^m$ can be identified with a contractible subspace, yet denoted $S^m$, of $\hat S^{m|n}$. There exists a mapping $\hat S^{m|n}\to S^m$, but this is not a retraction, and the inclusion $S^m\hookrightarrow\hat S^{m|n}$, cannot be a homotopy equivalence.
\end{theorem}

\begin{proof}
Since must necessarily be $\pi_k(\hat \Sigma^{m|n})\cong \pi_k(\hat S^{m|n})$, $k\ge 0$, it is enough prove theorem for $\hat S^{m|n}$. We shall first recall some useful definitions and results of Algebraic Topology, here codified as lemmas.
\begin{definition}
A pair $(X,A)$ has the {\em homotopy extension property} if a homotopy $f_t:A\to Y$, $t\in I$, can be extended to  homotopy $f_t:X\to Y$ such that $f_0:X\to Y$ is a given map.
\end{definition}
\begin{lemma}
If $(X,A)$ is a CW pair, then it has the homotopy extension property.
\end{lemma}
\begin{lemma}\label{contractibility-and-homotopy-equivalence}
If the pair $(X,A)$ satisfies the homotopy extension property and $A$ is contractible, then the quotient map $q:X\to X/A$ is a homotopy equivalence.
\end{lemma}
Let us consider that we can represent $S^m$ into $\hat S^{m|n}$ by a continuous mapping $s:S^m\to\hat S^{m|n}$, defined by means of the commutative diagram in (\ref{commutative-diagram-section-quantum-m-n-supersphere}).
\begin{equation}\label{commutative-diagram-section-quantum-m-n-supersphere}
  \xymatrix@C=2cm{\hat S^{m|n}\ar[d]_{\pi_C}\ar@{=}[r]&A^{m|n}\bigcup\{\infty\}\\
  S^m\ar@{=}[r]&\mathbb{R}^m\bigcup\{\infty\}\ar[u]_{s\equiv(\epsilon^m,0,,id_{\infty})}\\}
\end{equation}
where $\epsilon^m:\mathbb{R}^m\to A_0^m\subset A^m$ is induced by the canonical ring homomorphism $\epsilon:\mathbb{R}\to A$. $s$ is a section of $\pi$: $\pi\circ s=id_{S^m}$. Let us yet denote by $S^m$ the image of $s$. So we can consider the canonical couple $(\hat S^{m|n},S^m)$ as a CW pair, hence it has the homotopy extension property. $S^m$ is not a contractible subcomplex of $\hat S^{m|n}$, so in general the quotient map $\hat q:\hat S^{m|n}\to \hat S_{m|n}/S^m$ is not a homotopy equivalence.
We have the following lemma.

\begin{lemma}
The couple $(S^m,\infty)$ can be deformed into $(\hat S^{m|n},\infty)$ to the base point $\{\infty\}$.
\end{lemma}
\begin{proof}
In fact, let $p\in\hat S^{m|n}\setminus S^m$. Then the inclusion $i:S^m\hookrightarrow\hat S^{m|n}$ is nullhomotopic since $\hat S^{m|n}\setminus\{\infty\}\thickapprox A^{m|n}$ (homeomorphism).
\end{proof}
Since $S^m$ is contractible into $\hat S^{m|n}$, to the point $\infty\in\hat S^{m|n}$, the quotient map $\hat q:\hat S^{m|n}\to \hat S^{m|n}/S^m$ can be deformed into quotient mapping $\hat q_t$ over deformed quotient spaces $X_t\equiv\hat S^{m|n}/S^m_t$, with $S^m_t\equiv f_t(S^m)\subset\hat S^{m|n}$, for some homotopy $f:I\times S^m\to \hat S^{m|n}$, such that $X_0=\hat S^{m|n}/S^m$, $X_1=\hat S^{m|n}$ and $\hat q_1=id_{\hat S^{m|n}}$. (See diagram (\ref{diagram-deformation-quotient-mappings}).)

\begin{equation}\label{diagram-deformation-quotient-mappings}
\xymatrix{\hat S^{m|n}\ar[rdd]_(0.5){\hat q_1}\ar[rd]^(0.4){\hat q_t}\ar[r]^(0.4){\hat q_0=\hat q}&\hat S^{m|n}/S^m\equiv X_0\ar@{.}[d]\\
&\hat S^{m|n}/S^m_t\equiv X_t\ar@{.}[d]\\
&\hat S^{m|n}/\{\infty\}=\hat S^{m|n}\equiv X_1}
\end{equation}
But this does not assure that $\hat q$ is a homotopy equivalence.\footnote{Rally $S^m$ is contractible in $\hat S^{m|n}$, but is not a contractible sub-complex of $\hat S^{m|n}$. This clarifies the meaning of Lemma \ref{contractibility-and-homotopy-equivalence}. For example, in the case $A=\mathbb{C}$, one has that $\hat S^{1|n}/S^1$ is not homotopy equivalent to $S^2\cong_A \hat S^{1|n}$. In fact $\pi_2(S^2)=\mathbb{Z}$ and $\pi_2(S^2/S^1)\cong\pi_2(S^2\vee S^2)\cong H_2(S^2\vee S^2;\mathbb{Z})=\mathbb{Z}\bigoplus\mathbb{Z}$.} Let us, now, consider also some further lemmas.
\begin{lemma}
If $(X,A)$ is a CW pair and we have attaching maps $f,\, g:A\to X_0$ that are homotopic, then $X_0\bigcup_fX_1\backsimeq X_0\bigcup_gX_1$ $\hbox{\rm rel $X_0$}$ (homotopy equivalence).
\end{lemma}

\begin{lemma}
If $(X,A)$ satisfies the homotopy extension property and the inclusion $A\hookrightarrow X$ is a homotopy equivalence, then $A$ is a deformation retract of $X$.
\end{lemma}

\begin{lemma}
A map $f:X\to Y$ is a homotopy equivalence iff $X$ is a deformation retract of the mapping cylinder $M_f$.
\end{lemma}

Let us emphasize that we have a natural continuous mapping $\pi_C:\hat Sm|n\to S^m$, i.e., the surjection between the quantum $(m|n)$-supersphere and its classic limit, identified by the commutative diagram (\ref{commutative-diagram-section-quantum-m-n-supersphere}). The inclusion $i:S^m\hookrightarrow\hat S^{m|n}$ cannot be a deformation retract (and neither a strong deformation retract), otherwise $i$ should be a homotopy equivalence.\footnote{It is enough to consider the counterexample when $A=\mathbb{C}$ and $\hat S^{1|n}=\mathbb{C}\bigcup\{\infty\}=\mathbb{R}^2\bigcup\{\infty\}=S^2$. Then $S^1$ cannot be homotopy equivalent to $\hat S^{1|n}=S^2$, since $\pi_1(S^1)=\mathbb{Z}$ and $\pi_1(S^2)=0$.}
However, $\pi_C:\hat S^{m|n}\to S^m$, cannot be neither a retraction, otherwise their homotopy groups should be related by the split short exact sequence (\ref{split-short-exact-sequence-quantum-homotopy-groups-m-n-supersphere}),
\begin{equation}\label{split-short-exact-sequence-quantum-homotopy-groups-m-n-supersphere}
  \xymatrix{0\ar[r]&\pi_k(S^m,\infty)\ar@<1ex>[r]^{i_*}&\ar@<1ex>[l]^{r_*}\pi_k(\hat S^{m|n},\infty)\ar[r]&\pi_k(\hat S^{m|n},S^m,\infty)\ar[r]&\infty \\}
\end{equation}
hence we should have the splitting given in (\ref{splittings-homotopy-groups}). (For details on relations between homotopy groups and retractions see, e.g. \cite{PRAS11}.)

\begin{equation}\label{splittings-homotopy-groups}
   \pi_k(\hat S^{m|n},\infty)\cong \IM(i_*)\bigoplus\ker(r_*)\cong\pi_k(S^m,\infty)\bigoplus\ker(r_*).
\end{equation}
But this cannot work. In fact, in the case $A=\mathbb{C}$, we should have the commutative diagram (\ref{commutative-diagram-section-quantum-m-n-supersphere-complex-case}) with exact horizontal lines.
\begin{equation}\label{commutative-diagram-section-quantum-m-n-supersphere-complex-case}
  \xymatrix{0\ar[r]&\pi_1(S^1,\infty)\ar@{=}[d]\ar@<1ex>[r]^{i_*}&\ar@<1ex>[l]^{r_*}\pi_1(\hat S^1,\infty)\ar@{=}[d]\ar[r]&\pi_1(\hat S^1,S^1,\infty)\ar@{=}[d]\ar[r]&\infty \\
  0\ar[r]&\mathbb{Z}\ar[r]&0\ar[r]&\pi_1(\hat S^1,S^1,\infty)\ar[r]&0\\}
\end{equation}
This should imply that $\pi_1(S^1,\infty)=0$, instead that $\mathbb{Z}$, hence the bottom horizontal line in (\ref{commutative-diagram-section-quantum-m-n-supersphere-complex-case}) cannot be an exact sequence, hence $\pi_C:\hat S^1\cong S^2\to S^1$ cannot be a retraction !
\end{proof}

\begin{cor}
Quantum homotopy superspheres cannot be homotopy equivalent to $S^m$, except in the case that the quantum algebra $A$ reduces to $\mathbb{R}$.
\end{cor}

Quantum homotopy groups for quantum supermanifolds are introduced in \cite{PRAS19}.
\begin{theorem}[Quantum homotopy groups of quantum $(m|n)$-supersphere]\label{quantum-homotopy-groups-quantum-m-n-supersphere}
Quantum homotopy $(m|n)$-superspheres have quantum homotopy groups isomorphic to homotopy groups of $m$-spheres:
\begin{equation}\label{non-isomorphism-homotopy-groups-quantum-m-n-supersphere-sphere}
    \hat\pi_{k}(\hat\Sigma^{m|n})=\hat\pi_k(\hat S^{m|n})\cong\pi_{k}(S^m).
\end{equation}
\end{theorem}
\begin{proof}
In fact, we can prove for examples that $\hat\pi_{k}(\hat S^{m|n})\cong 0$, for $k<m$, and $\hat\pi_{m}(\hat S^{m|n})\cong \mathbb{Z}$. For this it is enough to reproduce the anoalogous proofs for the commutative spheres, by substituting cells with quantum supercells.
For example we can have the following quantum versions of analogous propositions for commutative CW complexes.

\begin{lemma}
Let $X$ be a quantum CW-complex admitting a decomposition in two quantum subcomplexes $X=A\bigcup B$, such that $A\bigcap B=C\not=\varnothing$. If $(A,C)$ is $m$-connected and $(B,C)$ is $n$-connected, $m,\, n\ge0$, then the mappings $\hat\pi_k(A,C)\to\hat\pi_k(X,B)$ induced by inclusion is an isomorphism for $k<m+n$, and a surjection for $k=m+n$.
\end{lemma}

\begin{lemma}[Quantum Freudenthal suspension theorem]\label{quantum-freudenthal-suspension-theorem}
The quantum suspension map $\hat\pi_k(\hat S^{m|n})\to\hat\pi_{k+1}(\hat S^{m+1|n+1})$ is an isomorphism for $k<2m-1$, and a surjection for $k=2m-1$.\footnote{This holds also for quantum suspension $\hat\pi_k(X)\to\hat\pi_{k+1}(\hat S X)$, for an $(m-1)$-connected quantum CW-complex $X$.}
\end{lemma}
As a by-product we get the isomorphism $\hat\pi_m(\hat S^{m|n})\cong\mathbb{Z}$.

\end{proof}

\begin{remark}
Let us emphasize that Theorem \ref{quantum-homotopy-groups-quantum-m-n-supersphere} does not allow to state that $S^m$ is a deformation retract of $\hat S^{m|n}$, as one could conclude by a wrong application of the Whitehead's theorem, reported in the following lemma.

\begin{lemma}[Whitehead's theorem]\label{whitehead-theorem}
If a map $f:X\to Y$ between connected CW complexes induces isomorphisms $f_*:\pi_m(X)\to\pi_*(Y)$ for all $m$, then $f$ is a homotopy equivalence. Furthermore, if $f$ is the inclusion of a subcomplex $f:X\hookrightarrow Y$, then $X$ is a deformation retract of $Y$.
\end{lemma}

In fact, in the case $i:S^m\hookrightarrow\hat S^{m|n}$ we are talking about different CW structures. One for $S^m$ is the usual one, the other, for $\hat S^{m|n}$ is the quantum super-CW structure. In order to easily understand the difference let us refer again to the case $A=\mathbb{C}$. Here one has $\pi_1(S^1)=\mathbb{Z}=\hat\pi_1(\hat S^{1|n}=\hat S^{1|0}=\hat S^1)$, but $\hat\pi_1(\hat S^1)=[\hat S^1,\hat S^1]=[S^2,S^2]=\pi_2(S^2)$. Furthermore, $\pi_1(\hat S^1)=[S^1,S^2]=0\not=\pi_1(S^1)$. Therefore, $\hat S^{1|n}$, with respect to the usual CW complex structure, has its first homotopy group zero, hence different from the first homotopy group of its classic limit $S^1$. In fact $S^1$ is not a deformation retract of $S^2=\hat S^1$. (Therefore there is not contradiction with the Whitehead's theorem.)
\end{remark}
Moreover, it is useful to formulate the quantum version of the Whitehead's theorem and some related lemmas. These can be proved by reproducing analogous proofs by substituting CW complex structure with quantum CW complex structure in quantum supermanifolds.

\begin{theorem}[Quantum Whitehead theorem]\label{quantum-whitehead-theorem}
If a map $f:X\to Y$ between connected quantum CW complexes induces isomorphisms $f_*:\hat\pi_m(X)\to\hat\pi_*(Y)$ for all $m$, then $f$ is a homotopy equivalence. Furthermore, if $f$ is the inclusion of a quantum subcomplex $f:X\hookrightarrow Y$, then $X$ is a quantum deformation retract of $Y$.
\end{theorem}

\begin{lemma}[Quantum compression lemma]\label{quantum-compression-lemma}
Let $(X,A)$ be a quantum CW pair and let $(Y,B)$ be any quantum pair with $B\not=\varnothing$. Let us assume that for each $m$ $\hat\pi_m(Y,B,y_{0})=0$, for all $y_0\in B$, and $X\setminus A$ has quantum supercells of dimension $(m|n)$. Then, every map $f:(X,A)\to(Y,B)$ is homotopic ${\rm rel}_A$ to a map $X\to B$.\footnote{When $m=0$, the condition $\hat\pi_m(Y,B,y_{0})=0$, for all $y_0\in B$, means that $(Y,B)$ is $0$-connected. Let us emphasize that there is not difference between $0$-connected and quantum $0$-connected. In fact $[\hat S^{0|n},Y]=\hat\pi_0(Y)=\pi_0(Y)=[S^0,Y]$, since $\hat S^{0|n}\backsimeq S^0=(\{a\},\{b\})$, i.e., homotopy equivalent to a set of two points. However, after Theorem \ref{quantum-homotopy-groups-quantum-m-n-supersphere}, there is not difference between the notion of {\em quantum $p$-connected} (i.e., $\hat\pi_k=0$, $k\le p$), quantum (homotopy) $(m|n)$-supersphere, and {\em $p$-connected} (i.e., $\pi_k=0$, $k\le p$), (homotopy) $(m|n)$-supersphere. In other words, a quantum homotopy $(m|n)$-supersphere is quantum $(m-1)$-connected as well as its classic limit is $(m-1)$-connected.}
\end{lemma}

\begin{lemma}[Quantum extension lemma]\label{quantum-extension-lemma}
Let $(X,A)$ be a quantum CW pair and let $f:A\to Y$ be a mapping with $Y$ a path-connected quantum supermanifold. Let us assume that $\hat\pi_{m-1}(Y)=0$, for all $m$, such that $X\setminus A$ has quantum cells of dimension $m$. Then, $f$ can be extended to a map $f:X\to Y$.
\end{lemma}
\begin{proof}
The proof can be done inductively. Let us assume that $f$ has been extended over the quantum $(m-1|n-1)$-superskeleton. Then, an extension over quantum $(m|n)$-supercells exists iff the composition of the quantum supercell's attaching map $\hat S^{m-1|n-1} \to \hat X^{m-1|n-1}$ with $f:\hat X^{m-1|n-1}\to Y$ is null homotopic.
\end{proof}
As a by-product of above results we get also the following theorems that relate quantum homotopy groups and quantum relative homotopy groups.

\begin{theorem}[Quantum exact long homotopy sequence]\label{quantum-exact-long-homotopy-sequence}
One has the exact sequence {\em(\ref{quantum-exact-long-homotopy-sequence-a})}.
\begin{equation}\label{quantum-exact-long-homotopy-sequence-a}
   \xymatrix{\cdots&\hat\pi_m(A,x_0)\ar[r]^{\hat i_*}&\hat\pi_m(X,x_0)\ar[r]^{\hat j_*}&\hat\pi_m(X,A,x_0)\ar[d]^{\hat \partial}\\
   &\hat\pi_{0}(X,x_0)&\ar[l]\cdots&\ar[l]\hat\pi_{n-1}(A,x_0)\\ }
\end{equation}
where $\hat i_*$ and $\hat j_*$ are induced by the inclusions $\hat i:(A,x_0)\hookrightarrow(X,x_0)$ and $\hat j:(X,x_0,x_0)\hookrightarrow(X,A,x_0)$ respectively. Furthermore, $\hat\partial$ comes from the following composition $(\hat S^{m-1|n-1},s_0)\hookrightarrow(\hat D^{m|n},\hat S^{m-1|n-1},s_0)\to(X,A,x_0)$, hence $\hat\partial[f]=[f_{|\hat S^{m-1|n-1}}]$.
\end{theorem}
\begin{theorem}[Quantum Hurewicz theorem]\label{quantum-hurewicz-theorem}
The exact commutative diagram in {\em(\ref{commutative-diagram-quantum-hurewicz-theorem})} relates (quantum) homotopy groups and (quantum) homology groups for (quantum) homotopy $(m|n)$-spheres, $m\ge 2$. The morphisms $a$ and $b$ are isomorphisms for $p\le m$ and epimorphisms for $p=m+1$.\footnote{Compare with analogous theorem in \cite{PRAS19} for quantum supermanifolds.}
\begin{equation}\label{commutative-diagram-quantum-hurewicz-theorem}
   \xymatrix{0\ar[r]&\pi_p(S^m)\ar[d]_{a}\ar@<0.5ex>[r]&\ar@<0.5ex>[l]\hat\pi_p(\hat S^{m|n})\ar[d]_{b}\ar[r]&0\\
   0\ar[r]&H_p(S^m;\mathbb{Z})\ar[d]\ar@<0.5ex>[r]&\ar@<0.5ex>[l] H_p(\hat S^m;\mathbb{Z})\ar[d]\ar[r]&0\\
   &0&0&\\}
\end{equation}
\end{theorem}

The following propositions are also stated as direct results coming from Theorem \ref{quantum-homotopy-groups-quantum-m-n-supersphere} and analogous propositions for topologic spaces.

\begin{proposition}\label{quantum-homotopic-properties-spheres}
The following propositions are equivalent for $i\le m-1$.

{\rm 1)} $S^i\to S^n$ is homotopic to a constant map.

{\rm 2)} $S^i\to S^n$ extends to a map $D^{i+1}\to S^n$.

{\rm 3)} $\hat S^{i|j}\to \hat S^{m|n}$ is homotopic to a constant map.

{\rm 4)} $\hat S^{i|j}\to \hat S^{m|n}$ extends to a map $\hat D^{i+1|j+1}\to \hat S^{m|n}$.
\end{proposition}

\begin{proof}
1) and 2) follow from the fact that $S^m$ is $(m-1)$-connected, and 3) and 4) from the fact that $\hat S^{m|n}$ is quantum $(m-1)$-connected. Furthermore, let us recall the following related result of Algebraic Topology.\footnote{There exists also a relative version of Lemma \ref{equivalent-conditions-for-m-connected-space}, saying that $\pi_i(X,A,x_0)=0$, for all $x_0\in A$,  is equivalent to one of the following propositions. (a1) Every map $(D^i,\partial D^i)\to(X,A)$ is homotopic ${\rm rel}\, \partial D^i$, to a map $D^i\to A$. (a2) Every map $(D^i,\partial D^i)\to(X,A)$ is homotopic through such maps to a map $D^i\to A$. (a3) Every map $(D^i,\partial D^i)\to(X,A)$ is homotopic through such maps to a constant map $D^i\to A$.}
\begin{lemma}\label{equivalent-conditions-for-m-connected-space}
The following propositions are equivalent.

{\rm (i)} The space $X$ is $m$-connected.

{\rm (ii)} Every map $f:S^i\to X$ is homotopic to a constant map.

{\rm (ii)} Every map $f:S^i\to X$ extends to a map $D^{i+1}\to X$.
\end{lemma}
\end{proof}

Let us, now, consider quantum super PDE's, in the category $\mathfrak{Q}_{hyper,S}$, with respect to quantum homotopy $(m|n)$-superspheres.

\begin{definition}[Quantum hypercomplex exotic super PDE's]\label{quantum-hypercomplex-exotic-super-pdes}
Let $\hat E_k\subset \hat J^k_{m|n}(W)$ be a $k$-order PDE on the fiber bundle $\pi:W\to M$ in the category $\mathfrak{Q}_{hyper,S}$, with $\dim_AM=(m|n)$ and $\dim_BW=(n|m,r|s)$, where $B=A\times E$ and $E$ is also a $Z(A)$-module. We say that $\hat E_k$ is a {\em quantum exotic PDE} if it admits Cauchy integral manifolds $N\subset \hat E_k$, $\dim N=(m-1|n-1)$, such that one of the following two conditions is verified.

{\em(i)} $\hat\Sigma^{m-2|n-2}\equiv\partial N$ is a quantum exotic supersphere of dimension $(m-2|n-2)$, i.e. $\hat\Sigma^{m-2|n-2}$ is homeomorphic to $\hat S^{m-2|n-2}$, ($\Sigma^{m-2|n-2}\thickapprox \hat S^{m-2|n-2}$) but not diffeomorphic to $\hat S^{m-2|n-2}$, ($\hat \Sigma^{m-2|n-2}\not\cong \hat S^{m-2|n-2}$).

{\em(ii)} $\varnothing=\partial N$ and $N\thickapprox \hat S^{m-1|n-1}$, but $N\not\cong \hat S^{m-1|n-1}$.\footnote{For complementary information see \cite{PRAS27}.}
\end{definition}

\begin{definition}[Quantum hypercomplex exotic-classic super PDE's]\label{quantum-hypercomplex-exotic-classsic-super-pdes}
Let $\hat E_k\subset \hat J^k_{m|n}(W)$ be a $k$-order super PDE as in Definition \ref{quantum-hypercomplex-exotic-super-pdes}. We say that $\hat E_k$ is a {\em quantum exotic-classic super PDE} if it is a quantum exotic super PDE, and the classic limit of the corresponding Cauchy quantum exotic supermanifolds are also exotic homotopy spheres.
\end{definition}

 From above results we get also the following one.

\begin{lemma}\label{lemma-quantum-hypercomplex-exotic-classic-super-pdes}
A quantum super PDE $\hat E_k\subset \hat J^k_{m|n}(W)$, where $m$ is such that $\Theta_{m-1}=0$, cannot be a quantum exotic-classic PDE, in the sense of Definition \ref{quantum-hypercomplex-exotic-classsic-super-pdes}.
\end{lemma}

\begin{lemma}\label{lemma-quantum-super-spheres-classic-limits-d}
For $m\in\{1,2,3,4,5,6\}$, one has the isomorphism reported in {\em(\ref{isomorphisms-lemma-quantum-super-spheres-classic-limits-d})}.

\begin{equation}\label{isomorphisms-lemma-quantum-super-spheres-classic-limits-d}
   \hat\Theta_{m|n}\cong\hat\Upsilon_{m|n}.
\end{equation}
In correspondence of such dimensions on $m$ we cannot have quantum exotic-classic super PDE's.
\end{lemma}

\begin{proof}
Isomorphisms in (\ref{isomorphisms-lemma-quantum-super-spheres-classic-limits-d}), follow directly from above lemmas, and the fact that $\Theta_{m}=0$ for $m\in\{1,2,3,4,5,6\}$. (See Refs. \cite{PRAS24, PRAS26}.)
\end{proof}

\begin{example}[The quantum hypercomplex Ricci flow super-equation]\label{the-quantum-hypercomplex-ricci-flow-super-equation}
As a by-product of Theorem \ref{generalized-poincare-conjecture-in-category-hypercomplex-quantum-super-manifolds} it follows that under the same hypotheses there adopted on the quantum algebra $A$, it follows that the quantum Ricci flow equation in the category $\mathfrak{Q}_{hyper,S}$, is a quantum exotic super PDE. On the other hand, such a PDE cannot be quantum exotic-classic for $m<7$. (See \cite{PRAS23, PRAS24, PRAS26}.) (For complementary information on the Ricci flow equation see also the following Refs. \cite{AG-PRA1,DUB-FOM-NOV,HAMIL1,HAMIL2,HAMIL3,HAMIL4,HAMIL5,KERVAIRE-MILNOR,MILNOR1,MILNOR2,MILNOR2a,MOISE1,MOISE2,NASH,PER1, PER2,PONTRJAGIN,SMALE,WALL1,WALL2}.)
\end{example}

\begin{example}[The quantum hypercomplex Navier-Stokes super-equation]\label{the-quantum-hypercomplex-navier-stokes-super-equation}
The quantum Navier-Stokes equation can be encoded on the quantum super-extension of the affine fiber bundle $\pi:W\equiv M\times\mathbf{I}\times\mathbb{R}^2\to M$, $(x^\alpha,\dot x^i,p,\theta)_{0\le\alpha\le 3,1\le i\le 3}\mapsto(x^\alpha)$. (See Refs. \cite{PRAS5,PRAS7,PRAS26} for the Navier-Stokes equation in the category of commutative manifolds and \cite{PRAS9,PRAS11,PRAS27} for its quantum extension on quantum manifolds.) Therefore, Cauchy manifolds are $(3|3)$-dimensional quantum supermanifolds. For such dimension do not exist exotic spheres. Therefore, the Navier-Stokes equation cannot be a quantum exotic-classic super PDE. Similar considerations hold for PDE's of the quantum super-extensions of continuum mechanics PDE's.
\end{example}

\begin{example}[The quantum hypercomplex $(m|n)$-d'Alembert super-equation]\label{the-quantum-hypercomplex-m-n-d'alembert-super-equation}
The quantum $(m|n)$-d'Alembert super-equation on $A^{m|n}$ cannot be a quantum exotic-classic super PDE for quantum $(m|n)$-dimensional Riemannian manifolds, with $m<7$, in the category $\mathfrak{Q}_{hyper,S}$. (For complementary information on the geometric structure of the d'Alembert PDE in the category of commutative manifolds and quantum manifolds see Refs. \cite{PRAS7,PRAS8,PRAS11,PRAS13,PRAS25,PRAS26,PRAS27}.)
\end{example}

\begin{example}[The quantum hypercomplex Einstein super-equation]\label{the-quantum-hypercomplex-einstein-super-equation}
Considerations similar to ones made in Example \ref{the-quantum-hypercomplex-m-n-d'alembert-super-equation}, hold for the quantum Einstein super-equation in the category $\mathfrak{Q}_{hyper,S}$.
\end{example}

\begin{theorem}[Integral bordism groups in quantum hypercomplex exotic super PDE's in the category $\mathfrak{Q}_{hyper,S}$ and stability]\label{bordism-groups-in-quantum-hypercomplex-exotic-super-pdes-stability}
Let $\hat E_k\subset \hat J^k_{m|n}(W)$ be a quantum exotic formally integrable and completely integrable super PDE on the fiber bundle $\pi:W\to M$, in the category $\mathfrak{Q}_{hyper,S}$, such that $\hat g_k\not=0$ and $\hat g_{k+1}\not=0$.\footnote{The fiber bundle $\pi:W\to M$ is as in Definition \ref{quantum-hypercomplex-exotic-super-pdes}, hence $\dim_AM=(m|n)$, $\dim_BW=(n|m,r|s)$, with $E$ endowed with a $Z(A)$-module structure too.}
Then there exists a bi-graded topologic spectrum $\Xi_{i|j}$ such that for the singular integral $(p|q)$-(co)bordism groups can be expressed by means of suitable bigraded homotopy groups as reported in {\em(\ref{singular-integral-p-q-co-bordism-groups-quantum-exotic-super-pdes})}.
\begin{equation}\label{singular-integral-p-q-co-bordism-groups-quantum-exotic-super-pdes}
    \left\{
    \begin{array}{l}
      \Omega_{p|q,s}^{\hat E_k}=\mathop{\lim}\limits_{(i|j)\to(\infty|\infty)}\hat\pi_{p+i|q+j}(\hat E_k^+\wedge\Xi_{i|j})\\
      \Omega^{p|q,s}_{\hat E_k}=\mathop{\lim}\limits_{(i|j)\to(\infty|\infty)}[\hat S^{i|j}\hat E_k^+,\Xi_{p+i|q+j}]\\
      \end{array}
    \right\}_{p\in\{0,1,\cdots,m-1\},q\in\{0,1,\cdots,n-1\}}.
\end{equation}

Furthermore, the singular integral bordism group for admissible smooth closed compact Cauchy manifolds, $N\subset \hat E_k$, is given in {\em(\ref{singular-integral-bordism-group-m-1-n-1})}.
\begin{equation}\label{singular-integral-bordism-group-m-1-n-1}
    \Omega_{m-1|n-1,s}^{\hat E_k}\cong H_{m-1|n-1}(W;A).
\end{equation}
In the {\em quantum homotopy equivalence full admissibility hypothesis}, i.e., by considering admissible only $(m-1|n-1)$-dimensional smooth Cauchy integral supermanifolds identified with quantum homotopy superspheres, and assuming that the space of conservation laws is not trivial,
one has $\Omega^{\hat E_k}_{m-1|n-1,s}=0$. Then $\hat E_k$ becomes a quantum extended $0$-crystal super PDE. Therefore, there exists a global singular attractor, in the sense that all Cauchy supermanifolds, identified with quantum homotopy $(m-1|n-1)$-superspheres, bound singular manifolds.

Furthermore, if in $W$ we can embed all the quantum homotopy $(m-1|n-1)$-superspheres, and all such supermanifolds identify admissible smooth $(m-1|n-1)$-dimensional Cauchy supermanifolds of $\hat E_k$), then two of such Cauchy supermanifolds bound a smooth solution iff they are diffeomorphic and one has the following bijective mapping: $\Omega^{\hat E_k}_{m-1|n-1}\leftrightarrow\hat\Theta_{m-1|n-1}$.

Moreover, if in $W$ we cannot embed all quantum homotopy $(m-1|n-1)$-superspheres, but only $\hat S^{m-1|n-1}$, then in the {\em quantum supersphere full admissible hypothesis}, i.e., by considering admissible only quantum $(m-1|n-1)$-dimensional smooth Cauchy integral supermanifolds identified with $\hat S^{m-1|n-1}$, then $\Omega^{\hat E_k}_{m-1|n-1}=0$. Therefore $\hat E_k$ becomes a quantum $0$-crystal super PDE and there exists a global smooth attractor, in the sense that two of such smooth Cauchy supermanifolds, identified with $\hat S^{m-1|n-1}$ bound quantum smooth supermanifolds. Instead, two Cauchy supermanifolds identified with quantum exotic $(m-1|n-1)$-superspheres bound by means of quantum singular solutions only.

All above quantum smooth or quantum singular solutions are unstable. Quantum smooth solutions can be stabilized.
\end{theorem}

\begin{proof}
The relations (\ref{singular-integral-p-q-co-bordism-groups-quantum-exotic-super-pdes}) and (\ref{singular-integral-bordism-group-m-1-n-1}) can be proved by a direct extension of analogous characterizations of integral bordism groups of PDE's in the category of commutative manifolds and quantum PDEs. (See \cite{PRAS7,PRAS8,PRAS26,PRAS27}.) Then the rest of the proof follows directly by using above results in this section, and following a road similar to the proof of Theorem 5.38 given in \cite{PRAS27}.\end{proof}

Similarly one can prove the following theorem that extends in the category $\mathfrak{Q}_{hyper,S}$ an analogous theorem in the category of commutative manifolds and in the category $\mathfrak{Q}_{hyper}$. (See \cite{PRAS26,PRAS27}.)

\begin{theorem}[Integral h-cobordism in quantum hypercomplex Ricci flow super PDE's]\label{integral-h-cobordism-in-Ricci-flow-super-pdes}
The quantum Ricci flow equation for quantum $(m|n)$-dimensional Riemannian supermanifolds, admits that starting from a quantum $(m|n)$-dimensional supersphere $\hat S^{m|n}$, we can dynamically arrive, into a finite time, to any quantum $(m|n)$-dimensional homotopy supersphere $M$. When this is realized with a smooth solution, i.e., solution with characteristic flow without singular points, then $\hat S^{m|n}\cong M$. The other quantum homotopy spheres $\hat\Sigma^{m|n}$, that are homeomorphic to $\hat S^{m|n}$ only, are reached by means of singular solutions.

For $1\le m\le 6$,  quantum hypercomplex Ricci flow super PDE's cannot be quantum exotic-classic ones. In particular, the case $m=4$, is related to the proof that the smooth Poincar\'e  conjecture is true.
\end{theorem}

\begin{appendices}

\appendix{\bf Appendix A: Quantum energy-momentum tensor and quantum Euler-Lagrange equation.}\label{appendixa}
\renewcommand{\theequation}{A.\arabic{equation}}
\setcounter{equation}{0}  

We can explicitly see that $d\omega_H|V=0$ and $d p_r|_V=0$ for solutions of $\widehat{(YM)}[i]$, by considering that
\begin{equation}\label{calculation-a}
  \left\{
  \begin{array}{l}
    d\omega_H=(\partial x_\mu.[(\partial y^\mu_\beta.L)y^\beta_4-\delta^\mu_4 L])\otimes dx^1\wedge dx^2\wedge dx^3\wedge dx^4\\
   dp_r=(\partial x_\mu.[(\partial y^\mu_\beta.L)y^\beta_r-\delta^\mu_r L])\otimes dx^1\wedge dx^2\wedge dx^3\wedge dx^4\\
   \end{array}
  \right.
\end{equation}
Set $T^\mu_\alpha=(\partial y^\mu_\gamma.L)y^\gamma_\alpha-\delta^\mu_\alpha L$. This is the quantum energy-momentum tensor. One can see that $(\partial x_\mu.T^\mu_\alpha)=0$ for solutions of the Euler-Lagrange equations, namely $(\partial y_\beta.L)-\partial x_\mu.(\partial y^\mu_\beta.L)=0$. In fact, taking into account that
\begin{equation}\label{calculus-b}
\left\{
\begin{array}{ll}
(\partial x_\alpha.L)&=(\partial y_\beta.L)(\partial x_\alpha.y^\beta)+(\partial y^\gamma_\beta.L)(\partial x_\alpha.y^\beta_\gamma)\\
&=\partial x_\gamma(\partial y^\beta_\alpha(\partial y^\gamma_\beta.L))\\
\end{array}
\right.
\end{equation}
we get
$$\delta^\gamma_\alpha(\partial x_\gamma.L)=\partial x_\gamma(\partial y^\beta_\alpha(\partial y^\gamma_\beta.L)).$$

Therefore, we have
$$\partial x_\gamma.[(\partial y^\gamma_\beta.L)y^\beta_\alpha-\delta^\gamma_\alpha L]=0=(\partial x_\gamma.T^\gamma_\alpha).$$

\vskip 0.5cm
\appendix{\bf Appendix B: Quantum angular-momentum tensor and quantum Euler-Lagrange equation.}\label{appendixb}
\renewcommand{\theequation}{B.\arabic{equation}}
\setcounter{equation}{0}  

The observed quantum Yang-Mills super PDEs are also invariant under infinitesimal transformations of the Lorentz group. These generate the following quantum conservation laws:
\begin{equation}\label{quantum-conserved-laws-angular-momentum-tensor }
\left\{\begin{array}{ll}
          \beta_{\mu\nu}&=<(M^\lambda_{\mu\nu}\otimes\partial x_\lambda), dx^1\wedge dx^2\wedge dx^3\wedge dx^4> \\
         =& \sum_{1\le \lambda\le 4}(-1)^{1+1} (M^\lambda_{\mu\nu}\otimes dx^1\wedge dx^2\wedge\cdots \widetilde{dx^\lambda}\cdots\wedge dx^3\wedge dx^4 \\
       \end{array}
  \right.
\end{equation}
where $(M^\lambda_{\mu\nu}=x^\alpha\delta_{\alpha\mu}T^\lambda_\mu-x^\alpha\delta_{\alpha\nu}T^\lambda_\mu-S^\lambda_{\mu\nu}$ with
$S^\lambda_{\mu\nu}=(s_{\mu\nu})^A_B y^B(\partial y^\lambda_A.L)$. One has $M^\lambda_{\mu\nu}=-M^\lambda_{\nu\mu}$ and on any nonlinear quantum propagator $V$, it results $d\beta_{\mu\nu}|V=0$, namely $(\partial x_\lambda.M^\lambda_{\nu\mu})=0$. Therefore we get also
\begin{equation}\label{fundamental-relation-between-spin-torsion-energy-momentum-tensor}
  (\partial x_\lambda.S^\lambda_{\mu\nu})=T_{\mu\nu}-T_{\nu\mu}.
\end{equation}
$M=\partial x_\lambda M^\lambda_{\mu\nu}\otimes dx^\mu\wedge dx^\nu: N\to A\otimes \Lambda^0_2N\otimes TN$ is called the {\em observed quantum angular momentum tensor}. The quantum charge corresponding to the quantum conservation law $\beta_{\mu\nu}$ is
\begin{equation}\label{quantum-angular-momentum-charge}
  M_{\mu\nu}[i|t]=\int_{\sigma_t}\beta_{\mu\nu}=-\int_{\sigma_t}M^4_{\mu\nu}\otimes dx^1\wedge dx^2\wedge dx^3.
\end{equation}
Therefore we get for any nonlinear quantum propagator $V$, $\partial V=N_0\bigcup P\bigcup N_1$, $\partial P=\partial N_0\bigcup \partial N_1$,
\begin{equation}\label{quantum-angular-momentum-charge-a}
  0=\int_Vd\beta=M_{\mu\nu}[i|t_0]-M_{\mu\nu}[i|t_1]+\int_P\beta_{\mu\nu}].
\end{equation}
Let us define {\em defect quantum $(\mu\nu)$-angular-momentum} of $V$
\begin{equation}\label{defect-quantum-angular-momentum}
  \mathfrak{M}_{\mu\nu}[V]=\int_P\beta_{\mu\nu}|_P\in A.
\end{equation}
Then we see that $M_{\mu\nu}[i|t]$ is constant iff $\mathfrak{M}_{\mu\nu}[V]=0$. If $P$ is an orientable smooth $3$-dimensional manifold, we can see that $\mathfrak{M}_{\mu\nu}[V]=0$. In fact, we can write
\begin{equation}\label{defect-quantum-angular-momentum-a}
  \mathfrak{M}_{\mu\nu}[V]=\int_P[M^3_{\mu\nu}-M^1_{\mu\nu}x^3_1-M^2_{\mu\nu}x^3_2-M^4_{\mu\nu}x^3_4]\otimes dx^1\wedge dx^2\wedge dx^4.
\end{equation}
Since we assumed $P$ oriented and smooth, in some neighbourhood of any point of $P$ one has a relation $x^3=x^3(x^1,x^2,x^4)$, hence the following relation $M^3_{\mu\nu}=M^1_{\mu\nu}x^3_1+M^2_{\mu\nu}x^3_2+M^4_{\mu\nu}x^3_4$ holds everywhere on $P$: As a by product we get also $\mathfrak{M}_{\mu\nu}[V]=0$. However, we can also consider nonlinear quantum propagators, where the implicit function theorem does not work on $P$, namely we can consider singular nonlinear quantum propagators. We call such propagators {\em exotic nonlinear quantum propagators}. For such propagators $ \mathfrak{M}_{\mu\nu}[V]\not=0$, and $ M_{\mu\nu}[i|t]$ is not more constant.\footnote{Exotic nonlinear quantum propagators have also non-zero $\mathfrak{H}[V]_{\partial}$, $\mathfrak{H}[V]$ and $\mathfrak{P}_r[V]$, $r=1,2,3$. (For details see in section \ref{sec-conservation-laws-quantum-super-pdes}.)} Instead $ M_{\mu\nu}[i|t]$ is surely constant in nonlinear quantum propagators encoding steady-states of orientable smooth $3$-dimensional particles. In fact in such cases $y^\beta_0=const$ and $V\cong B\times I$, where $B$ is a $3$-dimensional space-like manifold. Thus if $B$ is a smooth orientable manifold, one has that also $P\cong \partial B\times I$ is an orientable smooth $3$-dimensional manifold, hence $ \mathfrak{M}_{\mu\nu}[V]=0$.

\appendix{\bf Appendix C: Proof that for singular solutions defect quantum energy can be non-zero.}\label{appendixc}
\renewcommand{\theequation}{C.\arabic{equation}}
\setcounter{equation}{0}  

In this appendix we give an explicit proof that for singular solutions of the observed quantum Yang-Mills PDEs $\widehat{(YM)}[i]$ the defect quantum energy $\mathfrak{H}[V]_{\partial}$, of a nonlinear quantum propagator $V\subset \widehat{(YM)}[i]$, $\partial V=N_0\bigcup P\bigcup N_1$, $\partial P=\partial N_0\bigcup N_1$, does not necessitate to be zero. In other words we shall assume that $P$ is an integral singular $3$-chain $P=\sum_ia^iu_i$, with
$a_i\in A$ and $u_i:\Delta^3\to \widehat{(YM)}[i]\subset \hat J^2_4(E[i])$
an integral $3$-simplex, defined in a neighborhood $U\subset\mathbb{R}^4$ of $\bigtriangleup^3$
such that $T(u_i)(\bigtriangleup^3)\subset \hat{\mathbf{E}}^2_4$,
where $\hat{\mathbf{E}}^2_4\subset T\widehat{(YM)}[i]$ is the Cartan distribution
of $\widehat{(YM)}[i]$. Then the evaluation of $\omega_H:\widehat{(YM)}[i]\to A\otimes \Lambda^0_3\widehat{(YM)}[i]$ on $P$ is given in (\ref{evaluation-quantum-conservation-law-on-chain}).
\begin{equation}\label{evaluation-quantum-conservation-law-on-chain}
\scalebox{0.8}{$\begin{array}{ll}
<\omega_H,P>&=<\omega_H,\sum_ia^iu_i>=\sum_ia^i<\omega_H,u_i>=\sum_ia_i\int_{\Delta^3}u_i^*\omega_H\\
&\\
&=\sum_ia_i\int_{\Delta^3}[T^1_4 j(u_i)^{234}_{123}-T^2_4j(u_i)^{134}_{123}+T^3_4j(u_i)^{124}_{123}-T^4_4j(u_i)^{123}_{123}]\otimes d\xi^1\wedge d\xi^2\wedge d\xi^3.
\end{array}$}
\end{equation}
There $\{\xi^k\}_{1\le k\le 3}$ denotes a coordinate system on $\Delta^3$, and $u_i:\Delta^3\to \hat J^2_4(E[i])$ is locally represented by functions resumed in (\ref{evaluation-quantum-conservation-law-on-chain-a}).
\begin{equation}\label{evaluation-quantum-conservation-law-on-chain-a}
\left\{\begin{array}{ll}
x^\alpha\circ u_i&=u_i^\alpha(\xi^k)\\
&\\
y^j_\beta\circ u_i&=(u_i)^j_\beta(\xi^k),\, 0\le|\beta|\le 2.\\
\end{array}\right.
\end{equation}
Furthermore the jacobian matrix $j(u_i)$ is reported in (\ref{evaluation-quantum-conservation-law-on-chain-b}) to clarify our notation.
\begin{equation}\label{evaluation-quantum-conservation-law-on-chain-b}
j(u_i)=\left(
  \begin{array}{ccccc}
    (\partial \xi_1.u_i^1) & (\partial \xi_1.u_i^2) & (\partial \xi_1.u_i^3) & (\partial \xi_1.u_i^4) & (\partial \xi_1.(u_i)^j_\beta) \\
    (\partial \xi_2.u_i^1) & (\partial \xi_2.u_i^2) & (\partial \xi_2.u_i^3) & (\partial \xi_2.u_i^4) & (\partial \xi_2.(u_i)^j_\beta) \\
    (\partial \xi_3.u_i^1) & (\partial \xi_3.u_i^2) & (\partial \xi_3.u_i^3) & (\partial \xi_3.u_i^4) & (\partial \xi_3.(u_i)^j_\beta) \\
    \end{array}
\right)
\end{equation}
Therefore one has $$u_i^*dx^1\wedge\cdots\wedge \widetilde{dx^\alpha}\wedge \cdots\wedge dx^4=j(u_i)^{1\cdots\widetilde{\alpha}\cdots 4}_{123}d\xi^1\wedge d\xi^3\wedge  d\xi^3,$$
where $j(u_i)^{1\cdots\widetilde{\alpha}\cdots 4}_{123}$ is the minor of the jacobian matrix of $u_i$ obtained taking the rows $1$, $2$, and $3$ and the first $1$, $\cdots$, $\widetilde{\alpha}$, $\cdots$, $4$, columns. Since for singular $P$, $j(u_i)^{1\cdots\widetilde{\alpha}\cdots 4}_{1\hskip 2pt 2\hskip 2pt 3}$ are not all different from zero, it follows that the right-term, between square brackets, in (\ref{evaluation-quantum-conservation-law-on-chain-c})
\begin{equation}\label{evaluation-quantum-conservation-law-on-chain-c}
  u_i^*\omega_H= [T^1_4 j(u_i)^{234}_{123}-T^2_4j(u_i)^{134}_{123}+T^3_4j(u_i)^{124}_{123}-T^4_4j(u_i)^{123}_{123}]\otimes d\xi^1\wedge d\xi^2\wedge d\xi^3,
\end{equation}
can be non-zero.\footnote{Of course $u_i^*\omega_H=0$ in the {\em completely degenerate points}, namely where $j(u_i)^{1\cdots\widetilde{\alpha}\cdots 4}_{1\hskip 2pt 2\hskip 2pt 3}=0$, $\forall 1\le\widetilde{\alpha}\le 4$.} On the other hand, when $P$ is a smooth oriented $3$-dimensional manifold, then we can write, for example,
\begin{equation}\label{evaluation-quantum-conservation-law-on-chain-d}
\left\{\begin{array}{ll}
u_i^1&=\xi^1\\
u_i^2&=\xi^2\\
u_i^3&=u^3(\xi^1,\xi^2,\xi^3)\\
u_i^4&=\xi^3.\\
\end{array}\right.
\end{equation}
Thus we get $j(u_i)^{234}_{123}=-(u_i)^3_1$, $j(u_i)^{134}_{123}=(u_i)^3_2$, $j(u_i)^{124}_{123}=1$, $j(u_i)^{123}_{123}=(u_i)^3_3=u^3_4$.  Hence we have

\begin{equation}\label{evaluation-quantum-conservation-law-on-chain-e}
  u_i^*\omega_H= [-T^1_4(u_i)^3_1 -T^2_4(u_i)^3_2-T^4_4(u_i)^3_4+T^3_4]\otimes d\xi^1\wedge d\xi^2\wedge d\xi^3=0
\end{equation}
for the tensor properties of the quantum energy-momentum tensor $T^\alpha_\beta$.

\end{appendices}


\begin{thebibliography}{2020}



\bibitem{AAD-ET-AL} G. Aad et al., \textit{Observation on centrality dependent dijet asymmetry in lead-lead collisions at $\sqrt{s_{NN}}=276$ Tev with the ATLAS detector at the LHC}, Phys. Rev. Lett. \textbf{105(17)}(2010), 252303--17.
{[CERN Press Releases,  \href{http://press.web.cern.ch/press-releases/2010/11/lhc-experiments-bring-new-insight-primordial-universe}{\textit{LHC experiments bring new insight into primordial universe}}, November 26, 2010. Retrieved December 2, 2010.]}

\bibitem{AG-PRA1} R. P. Agarwal and A. Pr\'astaro, \textit{Geometry of PDE's.III(I):
Webs on PDE's and integral bordism groups. The general theory}. Adv.
Math. Sci. Appl. \textbf{17(1)}(2007), 239-266; \textit{Geometry
of PDE's.III(II): Webs on PDE's and integral bordism groups.
Applications to Riemannian geometry PDE's}, Adv. Math. Sci. Appl.
\textbf{17(1)}(2007), 267-281.

\bibitem{AG-PRA2} R. P. Agarwal and A. Pr\'astaro, \textit{Singular PDE's geometry and boundary value problems}.
J. Nonlinear Conv. Anal. \textbf{9(3)}(2008), 417-460; \textit{On singular PDE's geometry and boundary value problems}. Appl. Anal. \textbf{88(8)}(2009), 1115-1131.

\bibitem{AKBULUT-SALUR} S. Akbulut and S. Salur, \textit{Mirror duality via $G_2$ and $Spin(7)$ manifolds},
Arithmetic and Geometry Around Quantization, \"O. Ceyhan et. al. (eds.) Progress in Mathematics, Springer Science + Business Media LLC (2010), 279. DOI: 10.1007/978-0-8176-4831-2-1.


\bibitem{ARKANI-FINK-SLAT-WEI} N. Arkani-Hamed, D. P. Finkbeiner, T. S. Slatyer and N. Weiner, \textit{A theory of dark matter}. Phys. Rev. D \textbf{79}(2009), 015014--015020.

 \bibitem{BERGER} M. Berger, \textit{Classification des espaces homog\'en\'es symm\'etriques irr\'educibles}. C. R. acad. Sci., Paris \textbf{240}(1955), 2370--2372; \textit{Sur les groupes d'holonomie homog\'enes des vari\'et\'es riemanniennes}. Bull. Soc. Math. Fr. \textbf{83}(1955), 279--330.


\bibitem{BLANKE-GOLD} R. Blankenbecler and M. L. Goldberger, \textit{Behavior of scattering amplitudes at high energies, bound states, and resonances}. Phys. Rev. \textbf{126(2)}(1962), 766--786.

\bibitem{BRADEN-NEKRASOV} H. W. Braden and N. A. Nekrasov, \textit{Space-time foam from non-commutative instantons}. Commun. Math. Phys. \textbf{249(3)}(2004), 431-448.


\bibitem{BROGLIE1} L. Broglie de, \textit{Recherches sur la th\'eorie des quanta}. Annales de Physique \textbf{10(3)}(1925), 22-128.

\bibitem{BROGLIE2} L. Broglie de, \textit{The wave nature of electron}. Nobel lecture, 12 December 1929. In Nobel Lectures in Physics (1901-1995). CD-Rom edn. Singapore: World Scientific.

\bibitem{BRYANT} R. L. Bryant \textit{A survey of Riemannian metrics with special holonomy groups},  Prog. Int. Cong. Mth., Berkeley/Calif. 1986; \textit{Holonomy and special geometries}, Bourginon, J-P. (ed.), Dirac operators: yesterday and today. Proceedings of the summer school and workshop, Beirut, Lebanon, August 27-September 7, 2001. Someville, MA: International Press, 71-90(2005).


\bibitem{B-C-G-G-G} R. L. Bryant, S. S. Chern, R. B. Gardner, H. L. Goldshmidt and P.
A. Griffiths, \emph{Exterior Differential Systems}, Springer-Verlag, New York, 1991.

\bibitem{CALABI} E. Calabi, \textit{On K\"ahler manifolds with vanishing canonical class}. Princeton Math. Ser. \textbf{12}(1957), 78-89.

\bibitem{CHEW-FRAUTSCHI} G. F. Chew and S. C. Frautschi, \textit{Principle of equivalence for all strongly interacting particles within the $S$-matrix framework}, Phys. Rev. Lett. \textbf{7}(1961), 394-397; \textit{Regge trajectories and the principle of maximum strenth for strong interactions}. Phys. Rev. Lett. \textbf{8}(1962), 41-44.

\bibitem{CUI} H-Y. Cui, \textit{Derivation of Gell-Mann-Nishijima formula from the electromagnetic field modes of a hadron}, \href{http://arxiv.org/pdf/1001.0226.pdf}{\tt arXiv:1001.0226v2[physics-gen-ph]}.

\bibitem{DIAKONOV-PETROV} D. Diakonov and V. Petrov, \textit{A heretical view on linear Regge trajectories}, arXiv: hep-ph/0312144.

\bibitem{DIRAC} P. A. M. Dirac, \textit{Relativistic wave equation}, Proc. Royal Soc. London, Serie A, Math. Phys. Sci. \textbf{155(886)}(1936), 447-459.


\bibitem{DUB-FOM-NOV} B. A. Dubrovin, A. T. Fomenko and S. P. Novikov, \textit{Modern Geometry-Methods and Applications. Part I; Part II; Part III.}, Springer-Verlag, New York 1990. (Original Russian edition: Sovremennaja Geometrie: Metody i Prilo\v{z}enia. Moskva: Nauka, 1979.)

\bibitem{EDEN} R. J. Eden, \textit{Regge poles and elementary particles}, Rep. Prog. Phys. \textbf{34}(1971), 995--1053.


\bibitem{FEYNMAN} R. P. Feynman, \textit{The theory of positrons}, Phys. Rev. \textbf{76}(1949), 749--759; \textit{Space-time approach to quantum electrodynamics}, Phys. Rev. \textbf{76(6)}(1949), 769--789; \textit{QED: Strange Theory of Light and Matter}. Princeton Univ. Press, Princeton, NJ, 1985.


\bibitem{FRADKIN} E. Fradkin, \textit{Quantum physics: Debut of the quarter electron}, Nature. \textbf{452}(2008), 822--823.

\bibitem{GIAMMARCHI} M. G. Giammarchi et al., \textit{Search for electron decay mode $e\to\gamma+\nu$ with prototype of Borexino detector}, Physics Letters B. \textbf{525}(2002), 29--40.

\bibitem{GOLDSCHMIDT} H. Goldshmidt, \textit{Integrability criteria for systems of non-linear
partial differential equations}. J. Differ. Geom. \textbf{1}(1967), 269-307.

\bibitem{GRAY} A. Gray, \textit{A note on manifolds whose holonomy group is a subgroup of $Spin(n)Sp(1)$}. Mich. Math. J. \textbf{16}(1969), 125--128.


\bibitem{GRIBOV} V. N. Gribov, \textit{The theory of complex angular momenta}, Gribov lectures on theoretical physics. With a foreword by Yuri Dokshitzer and an introduction by Yuri Dokshitzer and Leonid Frankfurt. Cambridge University Press. Xii, 297 pp. ISBN D-521-81834-6/pbk.

\bibitem{GRIBOV-PONTECORVO} V. N. Gribov and B. M. Pontecorvo, \textit{Neutrino astronomy and lepton charge}, Phys. Lett. \textbf{B 28}(1969), 493--496.

\bibitem{GROMOV} M. Gromov, \textit{Partial Differential Relations}. Springer-Verlag,
Berlin 1986.
\bibitem{HAMIL1} R. S. Hamilton, \textit{Three-manifolds with positive Ricci curvature}. J. Differ. Geom.
 \textbf{17}(1982), 255-306.

\bibitem{HAMIL2} R. S. Hamilton, \textit{Four-manifolds with positive Ricci curvature operator}.
J. Differ. Geom. \textbf{24}(1986), 153-179.

\bibitem{HAMIL3} R. S. Hamilton, \textit{Eternal solutions to the Ricci flow}. J. Differ. Geom.
\textbf{38}(1993), 1-11.

\bibitem{HAMIL4} R. S. Hamilton, \textit{The formation of singularities in the Ricci flow}. Surveys in
Differential Geometry, International Press, 1995, \textbf{2}(1995), 7--136.

\bibitem{HAMIL5} R. S. Hamilton, \textit{A compactness property for solutions of the Ricci flow on
three-manifolds}.  Comm. Anal. Geom. \textbf{7}(1999), 695--729.



\bibitem{HITCHIN} N. Hitchin \textit{The moduli space of complex Lagrangian manifolds},  Suppl. J. Differential Geom. \textbf{7}(2000), 327--345.


\bibitem{HIRSCH} M. Hirsch \textit{Differential Topology},  Springer-Verlag, New York, 1976.


\bibitem{KAISER} D. Kaiser \textit{Physics and Feynman Diagrams},  American Scientists \textbf{93}(2005), 156--165.

\bibitem{KERVAIRE-MILNOR} M. A. Kervaire and J. W. Milnor, \textit{Groups of homotopy spheres: I}, Ann. of Math. \textbf{77(3)}(1963), 504--537.

\bibitem{KRANTZ} S. G. Krantz \textit{Complex Analysis: The Geometric Viewpoint},  The Carus Mathematical Monographs, Second Edition, USA, \textbf{23}(2004).


\bibitem{KRAS-LYCH-VIN} I. S. Krasilshchik, V. V. Lychagin and A. M.
Vinogradov, \textit{Jet Spaces and Nonlinear Partial Differential
Equations}, Gordon \& Breach, N. Y. 1986.

\bibitem{LYCH-PRAS} V. Lychagin and A. Pr\'astaro, \textit{Singularities of Cauchy data,
characteristics, cocharacteristics and integral cobordism},
Diff. Geom. Appls. \textbf{4}(1994), 283--300.

 \bibitem{MAJORANA} E. Majorana, \textit{Teoria simmetrica dell'elettrone e del positrone}. Nuovo Cimento \textbf{14}(1937), 171--184.

\bibitem{MANDELS-TAMM} L. I. Mandelshtam and I. E. Tamm, \textit{The uncertainty relation between energy and time in nonrelativistic quantum mechanics}. J. Physics. \textbf{9(4)}(1945), 249--254.

\bibitem{MCCLEARY} J. McCleary, \textit{User's guide to spectral sequences}. Publish or Perish Inc., Delaware, 1985.

\bibitem{MILNOR1} J. Milnor, \textit{On manifolds homeomorphic to the $7$-sphere}. Ann. of Math. \textbf{64(2)}(1956), 399--405.

\bibitem{MILNOR2} J. Milnor, \textit{The Steenrod algebra and its dual}. Ann. of Math. \textbf{67(2)}(1958), 150--171.

\bibitem{MILNOR2a} J. Milnor, \textit{Morse theory}. Ann. of Math. Studies. Princeton University Press, Princeton N.J, 1963.

\bibitem{MILNOR-MOORE} J. Milnor and J. Moore, \textit{On the structure of Hopf algebras}. Ann. of Math. \textbf{81(2)}(1965), 211--264.

\bibitem{MOISE1} E. Moise, \textit{Affine structures in $3$-manifolds. V. The triangulation theorem and Hauptvermuntung}. Ann. of Math. Sec. Ser. \textbf{56}(1952), 96--114.

\bibitem{MOISE2} E. Moise, \textit{Geometric topology in dimension $2$ and $3$}. Springer-Verlag, Berlin, 1977.

\bibitem{NASH} J. Nash, \textit{Real algebraic manifolds}. Ann. of Math.
\textbf{56(2)}(1952), 405--421.

\bibitem{NEKRASOV} N. A. Nekrasov, \textit{Instantons and the 11th dimension}. Phil. Trans. R. Soc. London A \textbf{359}(2001), 1405--1412.


\bibitem{PER1} G. Perelman, \textit{The entropy formula for the Ricci flow and its geometry applications}, \href{http://arxiv.org/pdf/math/0211159.pdf}{\tt arXiv:math/0211159}.

\bibitem{PER2} G. Perelman, \textit{Ricci flow with surgery on three-manifolds}, \href{http://arxiv.org/pdf/math/0303109.pdf}{\tt arXiv:math/0303109}.


\bibitem{POLYAKOV} A. M. Polyakov, \textit{Quark confinement and topology of gauge groups}. Nucl. Phys. \textbf{B 120}(1977), 429--458.


\bibitem{PONTRJAGIN} L. S. Pontrjagin, \textit{Smooth manifolds and their applications homotopy theory}. Amer. Math. Soc. Transl. \textbf{11}(1959), 1--114.

 \bibitem{PONTECORVO} B. M. Pontecorvo, \textit{Neutrino experiments and the problem of conservation of leptonic charge}. Zh. Exp. Teor. Fiz. \textbf{53}(1967), 1717--1725.

\bibitem{PRADHAN} T. Pradhan, \textit{Electron decay}, \href{http://arxiv.org/pdf/hep-th/0312325.pdf}{\tt arXiv:hep-th/0312325v1}.


\bibitem{PRAS1} A. Pr\'astaro, \textit{Spinor super
bundles of geometric objects on $spin^{G}$ space-time structures},
Boll. Unione Mat. Ital. \textbf{(6)1-B}(1982), 1015--1028.

\bibitem{PRAS2} A. Pr\'astaro, \textit{Gauge geometrodynamics}, Riv. Nuovo Cimento \textbf{5(4)} (1982),
1--122.

\bibitem{PRAS3} A. Pr\'astaro, \textit{Cobordism of PDE's}, Boll. Unione Mat. Ital. \textbf{(7)5-B}(1991), 977--1001.

\bibitem{PRAS4} A. Pr\'astaro, \textit{Quantum geometry of PDE's}, Rep. Math. Phys. \textbf{30(3)}(1991), 273--354;
\textit{Geometry of super PDE's}, in: Geometry of Partial Differential Equations, A. Pr\'astaro \& Th. M. Rassias (eds.), World Scientific Publishing, River Edge, NJ, (1994), 259--315; \textit{Geometry of quantized super PDE's}, in: The Interplay Between Differential Geometry and Differential Equations, V. Lychagin (ed.), Amer. Math. Soc. Transl. \textbf{2/167}(1995), 165--192; \textit{Quantum geometry of super PDE's}, Rep. Math. Phys. \textbf{37(1)}(1996), 23--140;
 \textit{(Co)bordism in PDEs and quantum PDEs}, Rep. Math. Phys. \textbf{38(3)}(1996), 443--455.


\bibitem{PRAS5} A. Pr\'astaro, \textit{Geometry of PDE's and Mechanics}. World Scientific Publ., Denvers, USA, 1996.

\bibitem{PRAS6} A. Pr\'astaro, \textit{Quantum and integral (co)bordisms in partial differential equations}. Acta Appl. Math. \textbf{51}(1998), 243--302.

\bibitem{PRAS7} A. Pr\'astaro, \textit{(Co)bordism groups in PDE's}. Acta Appl. Math. \textbf{59(2)}(1999), 111--202.

\bibitem{PRAS8} A. Pr\'astaro, \textit{(Co)bordism groups in quantum PDE's}. Acta Appl. Math. \textbf{64(2/3)}(2000), 111--217.
\bibitem{PRAS9} A. Pr\'astaro, \textit{Quantum manifolds and integral (co)bordism groups in quantum partial differential equations}, Nonlinear Anal. Theory Methods Appl.
\textbf{47/4}(2001), 2609--2620.


\bibitem{PRAS10} A. Pr\'astaro, \href{http://atlas-conferences.com/c/a/k/r/14.htm}{\textit{Quantum super Yang-Mills equations: Global existence and mass-gap}}, Dynamic
Syst. Appl. \textbf{4}(2004), 227--232. (Eds. G. S. Ladde, N. G.
Madhin and M. Sambandham), Dynamic Publishers, Inc., Atlanta, USA.
ISBN:1-890888-00-1.

\bibitem{PRAS11} {A. Pr\'astaro}, \textit{Quantized Partial Differential Equations}, World Scientific Publ., Singapore, 2004.

\bibitem{PRAS12} A. Pr\'astaro, \textit{Conservation laws in quantum super PDE's}, \href{Melbourne-Florida-2005.pdf}{Proceedings of the Conference on Differential
\& Difference Equations and Applications (eds. R. P. Agarwal \& K.
Perera)}, Hindawi Publishing Corporation, New York (2006), 943--952.

\bibitem{PRAS13} A. Pr\'astaro, \textit{Geometry of PDE's. I: Integral bordism groups in PDE's}. J. Math. Anal. Appl. \textbf{319}(2006), 547--566; \textit{Geometry of PDE's. II: Variational PDE's and integral bordism groups}.
J. Math. Anal. Appl. \textbf{321}(2006), 930--948; \textit{(Co)bordism groups in quantum super PDE's.I: Quantum supermanifolds}, Nonlinear Anal. Real World Appl. \textbf{8(2)}(2007), 505--538.

\bibitem{PRAS14} A. Pr\'astaro, \textit{(Co)bordism groups in quantum super PDE's.I: Quantum supermanifolds},
Nonlinear Anal. Real World Appl. \textbf{8(2)}(2007), 505--538; \textit{(Co)bordism groups in quantum super PDE's.II: Quantum super PDE's}, Nonlinear Anal. Real World Appl. \textbf{8(2)}(2007), 480--504; \textit{(Co)bordism groups in quantum super PDE's.III: Quantum super Yang-Mills equations}, Nonlinear Anal. Real World Appl. \textbf{8(2)}(2007), 447--479.

\bibitem{PRAS15} A. Pr\'astaro, \textit{(Un)stability and bordism groups in PDE's}.
Banach J. Math. Anal. \textbf{1(1)}(2007), 139--147.


\bibitem{PRAS16} A. Pr\'astaro, \textit{Geometry of PDE's. IV: Navier-Stokes equation and integral bordism groups}. J. Math. Anal. Appl. \textbf{338(2)}(2008), 1140--1151.

\bibitem{PRAS17} A. Pr\'astaro, \href{http://atlas-conferences.com/cgi-bin/abstract/catb-14}{\textit{On quantum black-hole solutions of quantum super Yang-Mills equations}},
Dynamic Syst. Appl. \textbf{5}(2008), 407--414. (Eds. G. S. Ladde,
N. G. Madhin C. Peng \& M. Sambandham), Dynamic Publishers, Inc.,
Atlanta, USA. ISBN: 1-890888-01-6.

\bibitem{PRAS18} A. Pr\'astaro, \textit{Extended crystal PDE's stability.I: The general theory}.
Math. Comput. Modelling \textbf{49(9-10)}(2009), 1759--1780; \textit{Extended crystal PDE's stability.II: The extended crystal MHD-PDE's}. Math. Comput. Modelling \textbf{49(9-10)}(2009), 1781--1801; \textit{On the extended crystal PDE's stability.I: The $n$-d'Alembert extended crystal PDE's}. Appl. Math. Comput.
\textbf{204(1)}(2008), 63--69; \textit{On the extended crystal PDE's stability.II:
Entropy-regular-solutions in MHD-PDE's}. Appl. Math. Comput. \textbf{204(1)}(2008), 82--89.

\bibitem{PRAS19} A. Pr\'astaro, \textit{Surgery and bordism groups in quantum partial differential equations.I: The quantum Poincar\'e conjecture}. Nonlinear Anal. Theory Methods Appl. \textbf{71(12)}(2009), 502--525; \textit{Surgery and bordism groups in quantum partial differential equations.II: Variational quantum PDE's}. Nonlinear Anal. Theory Methods Appl. \textbf{71(12)}(2009), 526--549.


\bibitem{PRAS20} A. Pr\'astaro, \href{http://link.springer.com/chapter/10.1007/978-1-4939-1106-6_18}{\textit{Extended crystal PDE's}}.  \href{http://www.springer.com/mathematics/analysis/book/978-1-4939-1105-9}{Mathematics Without Boundaries: Surveys in Pure Mathematics. P. M. Pardalos and Th. M. Rassias (Eds.) Springer-Heidelberg New York Dordrecht London, (2014), 415--481.} ISBN 978-1-4939-1106-6 (Online) 978-1-4939-1105-9 (Print). \href{http://dx.doi.org/10.1007/978-1-4939-1106-6}{DOI: 10.1007/978-1-4939-1106-6}.
    \href{http://arxiv.org/abs/0811.3693}{\tt arXiv:0811.3693[math.AT]}.


\bibitem{PRAS21} A. Pr\'astaro, \textit{Quantum extended crystal PDE's}, \href{http://nonlinearstudies.com/index.php/nonlinear}{Nonlinear Studies \textbf{18(3)}(2011), 447--485.} \href{http://arxiv.org/abs/1105.0166}{\tt arXiv:1105.0166[math.AT]}.


\bibitem{PRAS22} A. Pr\'astaro, \textit{Quantum extended crystal super PDE's}, Nonlinear Analysis. Real World Appl. \textbf{13(6)}(2012), 2491--2529. \href{http://dx.doi.org/10.1016/j.nonrwa.2012.02.014}{DOI: 10.1016/j.nonrwa.2012.02.014.}. \href{http://arxiv.org/abs/0906.1363}{\tt arXiv:0906.1363[math.AT].}

\bibitem{PRAS23} A. Pr\'astaro, \textit{Exotic heat PDE's}, Commun. Math. Anal. \textbf{10(1)}(2011), 64--81. \href{http://arxiv.org/abs/1006.4483}{\tt arXiv: 1006.4483[math.GT].}

\bibitem{PRAS24} A. Pr\'astaro, \href{http://link.springer.com/chapter/10.1007/978-3-642-28821-0_15}{\textit{Exotic heat PDE's.II}}. \href{http://www.springer.com/mathematics/applications/book/978-3-642-28820-3}{Essays in Mathematics and its Applications. In Honor of Stephen Smale's 80th Birthday.  P. M. Pardalos and Th. M. Rassias (Eds.) Springer-Heidelberg New York Dordrecht London (2012), 369--419.} ISBN 978-3-642-28820-3 (Print) 978-3-28821-0 (Online).
\href{http://dx.doi.org/10.1007/978-3-642-28821-0}{DOI: 10.1007/978-3-642-28821-0}.
    \href{http://arxiv.org/abs/1009.1176}{\tt arXiv: 1009.1176[math.AT]}.

\bibitem{PRAS25} A. Pr\'astaro, \href{http://link.springer.com/chapter/10.1007/978-1-4614-3498-6_36}{\textit{Exotic $n$-d'Alembert PDE's and stability}}. \href{http://www.springer.com/mathematics/book/978-1-4614-3497-9?cm_mmc=NBA-_-Jun-12_EAST_10792128-_-product-_-978-1-4614-3497-9}{Nonlinear Analysis: Stability, Approximation and Inequalities. (Dedicated to Themistocles M. Rassias for his 60th birthday.) G. Georgiev (USA), P. Pardalos (USA) and H. M. Srivastava (Canada) (eds.), Springer Optimization and its Applications Volume 68(2012), 571-586}. ISBN 978-1-4614-3498-6.
\href{http://dx.doi.org/10.1007/978-1-4614-3498-6}{DOI: 10.1007/978-1-4614-3498-6}.
  \href{http://arxiv.org/abs/1011.0081}{\tt arXiv:1011.0081[math.AT]}.

\bibitem{PRAS26} A. Pr\'astaro, \href{http://link.springer.com/chapter/10.1007/978-1-4939-1106-6_18}{\textit{Exotic PDE's}}. \href{http://www.springer.com/mathematics/analysis/book/978-1-4939-1123-3}{Mathematics Without Boundaries: Surveys in Interdisciplinary Research. P. M. Pardalos and Th. M. Rassias (Eds.) Springer-Heidelberg New York Dordrecht London, (2014), 471--531.} ISBN 978-1-4939-1123-3 (print) 978-1-4939-1124-0 (eBook).
    \href{http://dx.doi.org/10.1007/978-1-4939-1124-0}{DOI: 10.1007/978-1-4939-1124-0}.
\href{http://arxiv.org/abs/1101.0283}{\tt arXiv:1101.0283[math.AT]}.

\bibitem{PRAS27} A. Pr\'astaro, \textit{Quantum exotic PDE's}. Nonlinear Anal. Real World Appl. \textbf{14(2)}(2013), 893--928. \href{http://dx.doi.org/10.1016/j.nonrwa.2012.04.001}{DOI: 10.1016/j.nonrwa.2012.04.001.}. \href{http://arxiv.org/abs/1106.0862}{\tt arXiv:1106.0862[math.AT]}.

\bibitem{PRAS28} A. Pr\'astaro, \textit{Strong reactions in quantum super PDE's. II: Nonlinear quantum propagators}.
\href{http://arxiv.org/abs/1205.2894}{\tt arXiv:1205.2894[math.AT]}.

\bibitem{PRAS29} A. Pr\'astaro, \textit{Strong reactions in quantum super PDE's. III: Exotic quantum supergravity}.
\href{http://arxiv.org/abs/1206.4856}{\tt arXiv:1206.4856[math.AT]}.

\bibitem{PRA-RAS} A. Pr\'astaro and Th. M. Rassias, \textit{Ulam stability in geometry of PDE's}.
Nonlinear Funct. Anal. Appl. \textbf{8(2)}(2003), 259--278.

\bibitem{PRA-REGGE} A. Pr\'astaro \& T. Regge, \textit{The group structure of supergravity}, Ann. Inst. H.
Poincar\'e Phys. Th\'eor. \textbf{44(1)}(1986), 39--89.

\bibitem{REGGE} T. Regge, \textit{Introduction to complex orbital moments}. Nuovo Cimento \textbf{14}(1959), 951--976.

\bibitem{SCHAFER} R. D. Schafer, \textit{An Introduction to Nonassociative Algebras}.
Academic Press, New York (1966). New edition, Dover Publications, New York (1995).


\bibitem{SCHLAPPA-SCHMITT} J. Schlappa, T. Schmitt et al., \textit{Spin-orbital separation in the quasi-one-dimensional Mott insulator $Sr_2CuO_3$}, Nature \textbf{18.04}(2012). doi: 101038/nature10974.

\bibitem{SCHOEN-YAU} R. S. Schoen and S. T. Yau, \textit{Conformally flat manifolds, Kleinian groups and scalar curvature}, Invent. Math. \textbf{92(1)}(1988), 47--71.

\bibitem{SMALE} S. Smale, \textit{Generalized Poincar\'e conjecture in dimension greater than four}. Ann. of Math. \textbf{74(2)}(1961), 391--406.

 \bibitem{STONG} R. E. Stong, \textit{Notes on Bordism Theories}. Amer. Math. Studies. Princeton Univ. Press, Princeton, 1968.

\bibitem{SWITZER} A. S. Switzer,  \textit{Algebraic Topology-Homotopy and Homology}, Springer-Verlag, Berlin,
1976.

\bibitem{THOM} R. Thom, \textit{Quelques propri\'et\'e globales des vari\'et\'es
diff\'erentieles}, Comm. Math. Helv. \textbf{28}(1954), 17--86.

\bibitem{VENEZIANO} G. Veneziano, \textit{Construction of a crossing-symmetric, Regge-behaved amplitude for linearly trajectories}. Nuovo Cim. \textbf{57A}(1968), 190--197.

\bibitem{WALL1} C. T. C. Wall, \textit{Determination of the cobordism ring}. Ann. of Math. \textbf{72}(1960), 292--311.

\bibitem{WALL2} C. T. C. Wall,  \textit{Surgery on Compact Manifolds}, London Math. Soc. Monographs \textbf{1},
Academic Press, New York, 1970; 2nd edition (ed. A. A. Ranicki),
Amer. Math. Soc. Surveys and Monographs \textbf{69}, Amer. Math.
Soc., 1999.

\bibitem{WARNER} F. W. Warner,  \textit{Foundations of Differentiable Manifolds and Lie
Groups}, Scott, Foresman and C., Glenview, Illinois, USA, 1971.


\bibitem{WHITNEY1} H. Whitney, \textit{Differentiable manifolds}. Ann. of Maths. \textbf{37}(1936), 647--680.

\bibitem{WHITNEY2} H. Whitney, \textit{The general type of singularity of a set of $2n-1$ smooth functions of $n$ variables}. Duke Math. J. \textbf{10}(1943), 161--173.

\bibitem{YAU} S. T. Yau, \textit{Calabi's conjecture and some new results in algebraic geometry}. Proc. Natl. Acad. Sci. USA \textbf{74}(1977), 1798--1799.

\end{thebibliography}
\end{document}


\setcounter{page}{\pagenumber}
\newcommand{\T}{\mathbb{T}}
\newcommand{\R}{\mathbb{R}}
\newcommand{\Q}{\mathbb{Q}}
\newcommand{\N}{\mathbb{N}}
\newcommand{\Z}{\mathbb{Z}}
\newcommand{\tx}[1]{\quad\mbox{#1}\quad}
\parindent=0pt
\def\SRA{\hskip 2pt\hbox{$\joinrel\mathrel\circ\joinrel\to$}}
\def\tbox{\hskip 1pt\frame{\vbox{\vbox{\hbox{\boldmath$\scriptstyle\times$}}}}\hskip 2pt}
\def\circvert{\vbox{\hbox to 8.9pt{$\mid$\hskip -3.6pt $\circ$}}}
\def\IM{\hbox{\rm im}\hskip 2pt}
\def\COIM{\hbox{\rm coim}\hskip 2pt}
\def\COKER{\hbox{\rm coker}\hskip 2pt}
\def\TR{\hbox{\rm tr}\hskip 2pt}
\def\GRAD{\hbox{\rm grad}\hskip 2pt}
\def\RANK{\hbox{\rm rank}\hskip 2pt}
\def\MOD{\hbox{\rm mod}\hskip 2pt}
\def\DEN{\hbox{\rm den}\hskip 2pt}
\def\DEG{\hbox{\rm deg}\hskip 2pt}

\title[Strong reactions in quantum super PDE's.II]{STRONG REACTIONS IN QUANTUM SUPER PDE's.II:\\ NONLINEAR QUANTUM PROPAGATORS}

\author{Agostino Pr\'astaro}
\maketitle
\vspace{-.5cm}
{\footnotesize
\begin{center}
Department SBAI - Mathematics, University of Rome La Sapienza, Via A.Scarpa 16,
00161 Rome, Italy. \\
E-mail: {\tt agostino.prastaro@uniroma1.it}
\end{center}
}

\vskip 0.5cm
\centerline{\em This work in three parts is dedicated to Albert Einstein and Max Planck.}
\vskip 0.5cm

\begin{abstract}
In this second part, of a work in three parts, devoted to encode strong reactions of the high energy physics, in the algebraic topologic theory of quantum super PDE's, (previously formulated by A. Pr\'astaro), decomposition theorems of integral bordisms in quantum super PDEs are obtained. (For part I and part III see \cite{PRAS28, PRAS29}.) In particular such theorems allow us to obtain representations of nonlinear quantum propagators in quantum super PDE's, by means of elementary ones ({\em quantum handle decompositions} of nonlinear quantum propagators). These are useful to encode nuclear and subnuclear reactions in quantum physics. Pr\'astaro's geometric theory of quantum PDE's allows us to obtain constructive and dynamically justified answers to some important open problems in high energy physics. In fact a {\em Regge-type} relation between {\em reduced quantum mass} and {\em quantum phenomenological spin} is obtained. A dynamical {\em quantum Gell-Mann-Nishijima formula} is given. An existence theorem of observed local and global solutions with electric-charge-gap, is obtained for quantum super Yang-Mills PDE's, $\widehat{(YM)}[i]$, by identifying a suitable constraint, $\widehat{(YM)}[i]_w\subset \widehat{(YM)}[i]$, {\em quantum electromagnetic-Higgs PDE}, bounded by a quantum super partial differential relation $\widehat{(Goldstone)}[i]_w\subset \widehat{(YM)}[i]$, {\em quantum electromagnetic Goldstone-boundary}. An electric neutral, connected, simply connected observed quantum particle, identified with a Cauchy data of $\widehat{(YM)}[i]$, it is proved do not belong to $\widehat{(YM)}[i]_w$.
Existence of {\em $Q$-exotic nonlinear quantum propagators} of $\widehat{(YM)}[i]$,  i.e., nonlinear quantum propagators that do not respect the quantum electric-charge conservation is obtained.

By using integral bordism groups of quantum super PDE's, a {\em quantum crossing symmetry theorem} is proved. As a by-product existence of {\em massive photons}  and {\em massive neutrinos} are obtained. A dynamical proof that quarks can be broken-down is given too. A {\em quantum time}, related to the observation of any nonlinear quantum propagator, is calculated. Then an {\em apparent quantum time} estimate for any reaction is recognized. A criterion to identify solutions of the quantum super Yang-Mills PDE encoding {\em(de)confined quantum systems} is given. {\em Supersymmetric particles} and {\em supersymmetric reactions} are classified on the ground of integral bordism groups of the quantum super Yang-Mills PDE $\widehat{(YM)}$. Finally, existence of the {\em quantum Majorana neutrino} is proved. As a by-product, the existence of a new quasi-particle, that we call {\em quantum Majorana neutralino}, is recognized made by means of two quantum Majorana neutrinos, a couple $(\widetilde{\nu_e},\widetilde{\bar\nu_e})$, supersymmetric partner of $(\nu_e,\bar\nu_e)$, and two Higgsinos.
\end{abstract}

\vskip 0.5cm

\noindent {\bf AMS Subject Classification:} 55N22, 58J32, 57R20; 58C50; 58J42; 20H15; 32Q55; 32S20.

\vspace{.08in} \noindent \textbf{Keywords}: Integral (co)bordism groups in quantum (super) PDE's;
Existence of local and global solutions in hypercomplex quantum super PDE's; Conservation laws;
Crystallographic groups; Quantum singular PDE's; Quantum exotic superspheres; Quantum reactions; Quantum Regge-type trajectories; Quantum Gell-Mann-Nishijima formula; $Q$-exotic nonlinear quantum propagators; Quantum crossing symmetries in quantum super PDEs; quantum massive photons; quantum massive neutrinos; Dark matter; Apparent quantum time; Quantum (de)confinement; Quantum supersymmetric partners; Quantum Majorana neutrino.

\section[Introduction]{\bf Introduction}

\vskip 0.5cm
The algebraic topologic theory of quantum (super) PDE's formulated by A. Pr\'astaro, allows to directly encode quantum phenomena in a category of noncommutative manifolds (quantum (super)manifolds) and to finally solve the problem of unification, at quantum level, of gravity with the fundamental forces \cite{PRAS5,PRAS8,PRAS9,PRAS10,PRAS11,PRAS12,PRAS14,PRAS17,PRAS19,PRAS21,PRAS22}. In particular, this theory allowed to recognize the mechanism of mass creation/distruction, as a natural geometric phenomenon related to the algebraic topologic structure of quantum (super) PDEs encoding the quantum system under study \cite{PRAS22}.

Aim of this second part is to explicitly prove that nuclear and subnuclear reactions can be encoded as boundary value problems in the Pr\'astaro's algebraic topology of quantum super PDEs, and can be represented in terms of elementary reactions. (For part I and part III see \cite{PRAS28, PRAS29}.)

It is important to emphasize that quantum conservation laws do not necessarily produce conservation of quantum charges in quantum reactions. In fact, it is also important consider the topological structure of the corresponding nonlinear quantum propagators encoding these reactions. In \cite{PRAS28} (Part I), we have shown this fact for the observed quantum energy. In this second part we characterize observed nonlinear quantum propagators $V$ of the observed quantum super Yang-Mills PDE, $\widehat{(YM)}[i]$, with respect to the total quantum electric-charge. Then we define {\em $Q$-exotic nonlinear quantum propagators} ones where there is a non-zero defect quantum electric-charge, $\mathfrak{Q}[V]\in A$, in the corresponding encoded reactions. ($A$ is the fundamental quantum superalgebra in $\widehat{(YM)}[i]$.) This important phenomenon, that is related to the gauge invariance of $\widehat{(YM)}[i]$, was non-well previously understood, since the gauge invariance was wrongly interpreted as a condition that necessarily produce the conservation of electric-charge in reactions.
Really just the gauge invariance is the main origin of such phenomenon, but beside the structure of the nonlinear quantum propagator. This fundamental aspect of quantum reactions in $\widehat{(YM)}[i]$, gives strong theoretical support to the guess about existence of quantum reactions where the ``electric-charge" is not conserved. The electric-charge conservation law was quasi a dogma in particle physics. However, there are in the world many heretical experimental efforts to prove existence of decays like the following $e^-\to \gamma+\nu$, i.e. electron decay into a photon and neutrino. In this direction some first weak experimental evidences were recently obtained.\footnote{See \cite{GIAMMARCHI}. Some other exotic decays were also investigated, as for example the exotic neutron's decay: $n\to p+\nu+\bar\nu$. \cite{NORMANN-BACHALL-GOLDHABER}.} With this respect, one cannot remark the singular role played, in the history of the science in these last 120 years, by the electron, a very small and light particle. In fact, at the beginning of the last century was just the electron to cause break-down in the Maxwell and Lorentz physical picture of the micro-world, until to produce a completely new point of view, i.e. the quantum physics. Now, after 120 years the electron appears to continue do not accept the place that physicists have reserved to it in the {\em world-puzzle}.

In the following we show as it is organized the paper and list the main results. 2.
Theorem \ref{quantum-superhandle-decomposition-quantum-nonlinear-propagator}: A representation of nonlinear quantum propagators by means {\em quantum exchangions} and {\em quantum virtual particles} is given. 3. Theorem \ref{quantum-scattering-processes-in-quantum-yang-mills-d-4-quantum-super-minkowski} characterizes nonlinear quantum propagators in quantum super Yang-Mills PDE, on quantum super Minkowski space-time. 4. Theorem \ref{reduced-quantum-mass-torsion}: The quantum mass is represented by means of the quantum torsion of the corresponding solution. Corollary \ref{phenomenological-quantum-regge-type-rule}: A direct relation between square {\em reduced quantum mass} and {\em phenomenological quantum spin} is given.\footnote{In particular this result shows how the phenomenological {\em Regge-type trajectories} emerge from the geometric theory of quantum (super) PDE's formulated by A. Pr\'astaro.}
Theorem \ref{quantum-gell-mann-nishijima-formula} gives a dynamical {\em quantum Gell-Mann-Nishijima formula}. By means of this formula one obtains a dynamical interpretation of {\em quantum hypercharge} and {\em quantum $3^{th}$ isospin component}. An existence theorem of observed local and global solutions with electric-charge-gap, is obtained too for observed quantum super Yang-Mills PDE's, $\widehat{(YM)}[i]$. It is identified a suitable constraint, $\widehat{(YM)}[i]_w\subset \widehat{(YM)}[i]$, {\em observed quantum electromagnetic-Higgs PDE}, bounded by a quantum super partial differential relation $\widehat{(Goldstone)}[i]_w\subset \widehat{(YM)}[i]$, {\em observed quantum electromagnetic Goldstone-boundary}. An observed nonlinear quantum propagator $V\subset\widehat{(YM)}[i]$, crossing the observed quantum electromagnetic Goldstone-boundary loses (or acquires) the property to have an electric-charge gap.  An electric neutral, connected, simply connected observed quantum particle, identified with a Cauchy data of $\widehat{(YM)}[i]$, cannot be contained into $\widehat{(YM)}[i]_w$.
Theorem \ref{non-conservation-of-quantum-electric-charge} proves that the quantum electric charge is not necessarily a conserved quantum law for any nonlinear quantum propagator of $\widehat{(YM)}[i]$. In fact, it is proved the existence of {\em$Q$-exotic nonlinear quantum propagators} in $\widehat{(YM)}[i]$, encoding observed reactions that do not respect the conservation of the quantum electric-charge.
Theorem \ref{quantum-crossing-symmetry}: By using integral bordism groups of quantum super PDE's, a {\em quantum crossing symmetry theorem} is proved.
Theorem \ref{quantum-massive-photons-theorem}: Existence of solutions of quantum super Yang-Mills PDEs, $\widehat{(YM)}$, representing productions of {\em quantum massive photons}, is proved.\footnote{These can be identified with quantum neutral massive vector bosons. The annihilation electron-positron must necessarily produce an intermediate {\em quantum virtual massive photon}.}
Theorem \ref{a-quantum-massive-photons}: Theorem \ref{quantum-massive-photons-theorem} is generalized to other quantum electric-charged particles. Theorem \ref{existence-of-quantum-massive-neutrino} states existence of quantum massive neutrinos, i.e., there exist nonlinear quantum propagators of $\widehat{(YM)}$ encoding decays of massive quasi-particles into a couple (neutrino,antineutrino).\footnote{Quantum massive photons, quantum massive neutrinos and $a$-quantum massive photons could interpret the so-called ``dark matter" that nowadays is a spellbinding object of active research. This exotic matter is related to the geometric structure of $\widehat{(YM)}$ that identifies the subequation $\widehat{(Higgs)}$. A global solution $V\subset \widehat{(YM)}$, crossing the
quantum Goldstone-boundary of $\widehat{(Higgs)}$, acquires (or loses) mass. (See \cite{PRAS22},  Example \ref{massive-virtual-particles-dark-matter} and footnote at page \pageref{dark-matter}, concernig $\pi^+$-photoproduction.)} Theorem \ref{quarks-break-down}: It is proved that quarks cannot be considered fundamental particles, i.e., they can be broken-down. Theorem \ref{quantum-virtual-anomaly-massive-particles}: Existence of massive quasi-particles, with masses that apparently contradict the conservation of mass-energy, is proved. This is related to the new concept of {\em quantum time} (see Definition \ref{quantum-relativistic-observed-time}) that comes from the interaction between quantum relativistic frame (i.e., an observer) and the quantum system, (hence also to the topological structure of nonlinear quantum propagators).
Theorem \ref{de-confinement-criterion}: A criterion to identify solutions of $(\widehat{YM})$, encoding confined quantum systems, is proved. Lemma \ref{supersymmetric-partners-particles} and Lemma \ref{supersymmetric-partners-quantum-nonlinear-propagators} state existence of supersymmetric particles and supersymmetric reactions by means of integral bordism groups of $\widehat{(YM)}$. Theorem \ref{existence-quantum-majorana-neutrino}: Existence of {\em quantum Majorana neutrino} and a new quasi-particle, that we call {\em quantum Majorana neutralino}, for its similarity with the so-called neutralino, are obtained, by using algebraic topologic properties of $\widehat{(YM)}$.

\section{\bf Surgery in Quantum Super PDEs}\label{sec-surgery-quantum-super-pdes}

Let us assume as a prerequisite of this section the knowledge of some previous works by A. Pr\'astaro on quantum (super) PDEs. However, to fix more directly ideas we recall also some fundamental definitions and results that are soon related to the subject here considered.

\begin{figure}[h]
\centerline{\includegraphics[width=6cm]{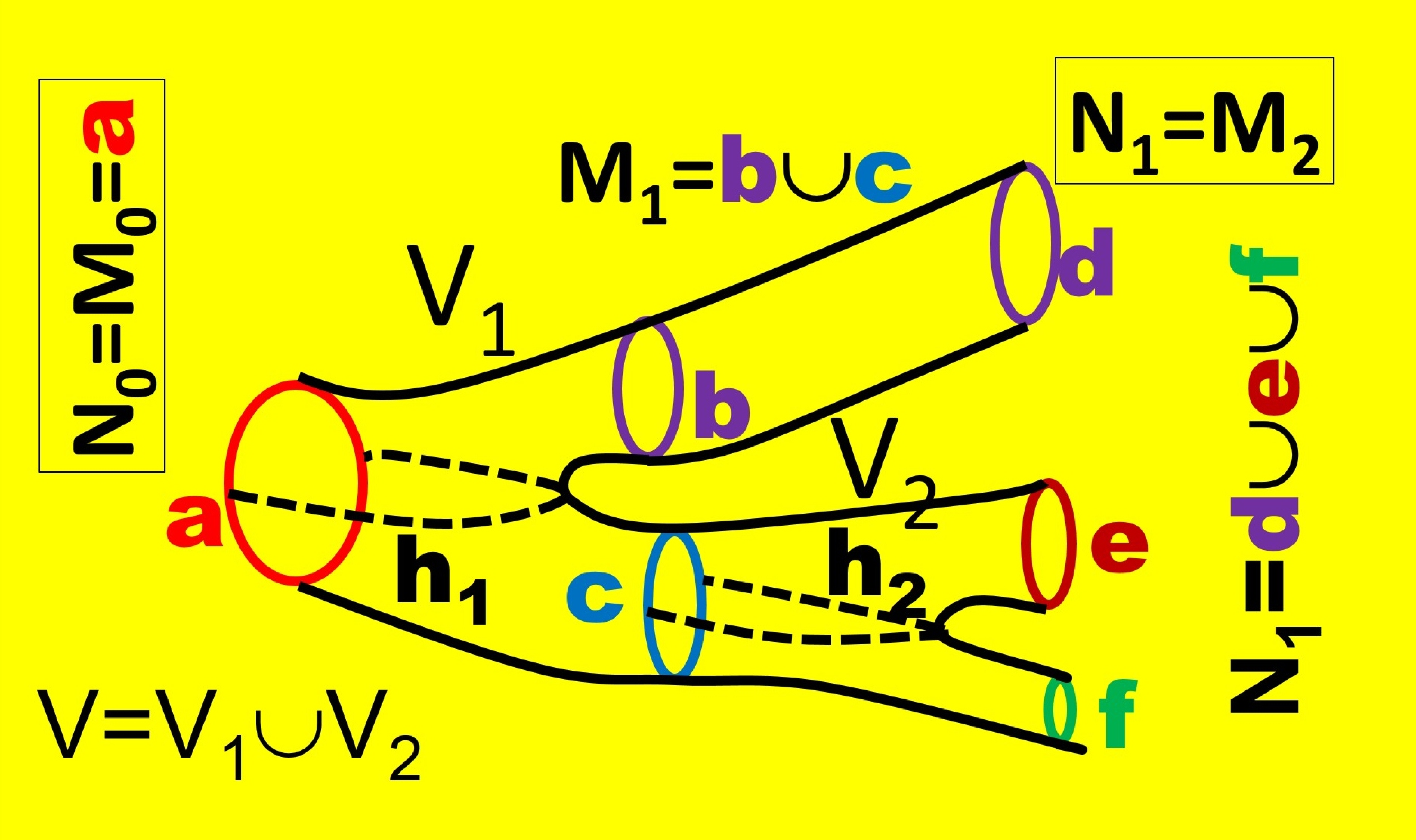}}
\renewcommand{\figurename}{Fig.}
\caption{Representation of quantum superhandle decomposition of nonlinear quantum propagator $V=V_1\bigcup_{M_1} V_2$. $\partial V=N_0\bigcup P\bigcup N_1$, $\partial V_1=M_0\bigcup P_1\bigcup M_1$, $\partial V_2=M_1\bigcup M_2$, with $M_0=N_0=a$, $M_2=N_1$. Furthermore $M_1=b\bigcup c$ and $N_1=d\bigcup e\bigcup f$. $\partial P_1=\partial M_0\bigcup \partial M_1$ and $\partial P_2=\partial M_1\bigcup \partial M_2$. $V_1$ is (up to homeomorphisms) $N_0$ with attached the handle $h_1$, and $V_2$ is obtained by attaching the handle $h_2$ to $M_1$, in such a way that $H_2$ it results attached to $c$ only. (Therefore $d$ is diffeomorphic to $b$.) One has $P=P_1\bigcup_{M_1}P_2$. (For graphic reasons $P$, $P_1$ and $P_2$ do not explicitly appear in the picture.)}
\label{figure-proof-quantum-superhandle-decomposition-quantum-nonlinear-propagator}
\end{figure}

Let $A=A_0\times A_1$ be a quantum superalgebra in the sense of Pr\'astaro.
\begin{definition}
An $(m|n)$-dimensional bordism in the category $\mathfrak{Q}_S$ of quantum supermanifolds, consists of the following $(W;M_0,f_0;M_1,f_1)$, where $W$ is a compact quantum supermanifold of dimension $(m|n)$ and closed $(m-1|n-1)$-dimensional quantum supermanifolds $M_0$ and $M_1$, such that $\partial W=N_0\sqcup N_1$, and quantum diffeomorphisms $f_i:M_o\to N_i$, $i=0,1$. An $(m|n)$-dimensional h-bordism, (resp. s-bordism) is a $(m|n)$-dimensional bordism as above, such that the inclusions $N_i\hookrightarrow W$, $i=0,1$, are homotopy equivalences, (resp. simply homotopy equivalences). We will simply denote also by $(W;M_0,M_1)$ an $(m|n)$-dimensional bordism in the category $\mathfrak{Q}_S$.\footnote{In this paper we will consider also more general bordisms, i.e., $(m|n)$-dimensional compact quantum supermanifolds $V$ such that $\partial V=N_0\bigcup P\bigcup N_1$, where $N_0$, $P$ and $N_1$ are $(m-1|n-1)$-dimensional quantum super manifolds such that $\partial P=\partial N_0\bigcup \partial N_1$. If $\partial N_0=\partial N_1=\varnothing$, then $P=\varnothing$. Then considerations similar to the ones made in this sections can be extended to these more general cases. See \cite{PRAS14,PRAS19}, Fig. \ref{figure-proof-quantum-superhandle-decomposition-quantum-nonlinear-propagator} and the next section.}
\end{definition}

\begin{theorem}[Quantum handle decomposition of nonlinear quantum propagator in quantum super PDEs]\label{quantum-superhandle-decomposition-quantum-nonlinear-propagator}
Let $\hat E_k\subset \hat J^k_{m+1|n+1}(W)$ be a quantum super PDE in the category $\mathfrak{Q}_S$ and let $V\subset \hat E_k$ be a compact solution of a boundary value problem, i.e., $\partial V=N_0\sqcup N_1$.\footnote{Let us recall that $V$ such that $\partial V=N_0\sqcup N_1$, is called {\em quantum non-linear propagator} between the two Cauchy data $N_0$ and $N_1$. This is the non-linear extension of the concepts of quantum propagator, usually used to quantize a classical field theory. For more details see Refs. \cite{PRAS19,PRAS21,PRAS22}.} Let us assume that $\hat E_k$ is formally integrable and completely integrable. Then there exists a quantum-super-handle-presentation {\em(\ref{quantum-superhandle-presentation-solution})} of the nonlinear quantum propagator $V$.
\begin{equation}\label{quantum-superhandle-presentation-solution}
    V_{1}\bigcup  V_{2}\bigcup\cdots\bigcup  V_{s}\thickapprox V
\end{equation}
where $(V_j;M_{j-1},M_j)$ is an adjoint elementary cobordism with index $p_j|q_j$, such that $$0|0\le p_1|q_1\le p_2|q_2\le\cdots\le p_s|q_s\le m+1|n+1$$ and $M_0=N_0$, $M_k=N_1$. In {\em(\ref{quantum-superhandle-presentation-solution})}, the symbol $\thickapprox$ denotes homeomorphism.
\end{theorem}

\begin{proof}
Let $M$ be a quantum supermanifold of dimension $(m|n)$ and let us consider an embedding
$$\phi:\hat S^{p|q}\times\hat D^{m-p|n-q}\to M.$$
Let us consider the following $(m|n)$-dimensional quantum supermanifold
$$\scalebox{0.8}{$M'\equiv\left(M\setminus{\rm int}(\phi(\hat S^{p|q}\times\hat D^{m-p|n-q}))\right)\bigcup_{\phi(\hat S^{p|q}\times\hat S^{m-p-1|n-q-1})}\left(\hat D^{p+1|q+1}\times\hat S^{m-p-1|n-q-1}\right)$.}$$
We say that $M'$ is obtained from $M$ by a {\em $(p|q)$-surgery}, i.e., cutting out ${\rm int}(\hat S^{p|q}\times\hat D^{m-p|n-q})$ and gluing in $\hat D^{p+1|q+1}\times\hat S^{m-p-1|n-q-1}$. The process of surgery is related to cobordism and handle attaching ones. Let $(X,\partial X)$ be a $(m+1|n+1)$-dimensional quantum supermanifold with boundary $\partial X$, and let
$$\phi:\hat S^{p|q}\times\hat D^{m-p|n-q}\to \partial X\equiv M$$
be an embedding. Set
$$X'\equiv X\bigcup_{\phi(\hat S^{p-1|q-1}\times\hat D^{m-p+1|n-q+1})}\left(\hat D^{p|q}\times\hat D^{m-p+1|n-q+1}\right).$$
We call $X'$ obtained from $X$ by {\em attaching a $(m+1|n+1)$-dimensional quantum superhandle of index $p|q$}. One has
$$\scalebox{0.8}{$M'\equiv\partial X'=\left(\partial X\setminus{\rm int}(\phi(\hat S^{p-1|q-1}\times\hat D^{m-p+1|n-q+1}))\right)\bigcup_{\phi(\hat S^{p-1|q-1}\times\hat S^{m-p|n-q})}\left(\hat D^{p|q}\times\hat S^{m-p|n-q}\right)$.}$$
We say that $M'$ is obtained from $M$ by a quantum $(p-1|q-1)$-surgery.
The surgery allows us to obtain a cobordism
$$W\equiv \left(M\times \hat D^{1|1}\right)\bigcup_{\phi(\hat S^{p-1|q-1}\times\hat D^{m-p+1|n-q+1}\times(1,1))}\left(\hat D^{p|q}\times\hat D^{m-p+1|n-q+1}\right)$$
with $\partial W=M\sqcup M'$.\footnote{In order to fix ideas we report in Fig. \ref{figure-handle-attachment-vs-cobordism} some examples in dimension $(m|n)=(1|1),\, (2|2)$, showing the relation between attaching quantum superhandles and cobordism.
In Fig. \ref{figure-handle-attachment-vs-cobordism-on-non-closed-manifold} is represented the case of a $(2|2)$-dimensional quantum supermanifold $X$ with boundary $\partial X$. Let us recall that $\hat S^{m|n}\equiv A^{m|n}\bigcup\{\infty\}$ is the $(m|n)$-dimensional quantum supersphere and $\hat D^{m|n}\subset \hat S^{m|n}$, such that $\partial\hat D^{m|n}\cong \hat S^{m-1|n-1}$ is the quantum supermanifold called quantum $(m|n)$-dimensional superdisk. One has $A^{0|0}=A^0=\{0\}$ and $\partial\hat D^{1|1}\cong \hat S^{0|0}=\{0\}\bigcup\{\infty\}$. (For more information about see \cite{PRAS19}.)}
\begin{figure}[h]
\centerline{\includegraphics[width=4cm]{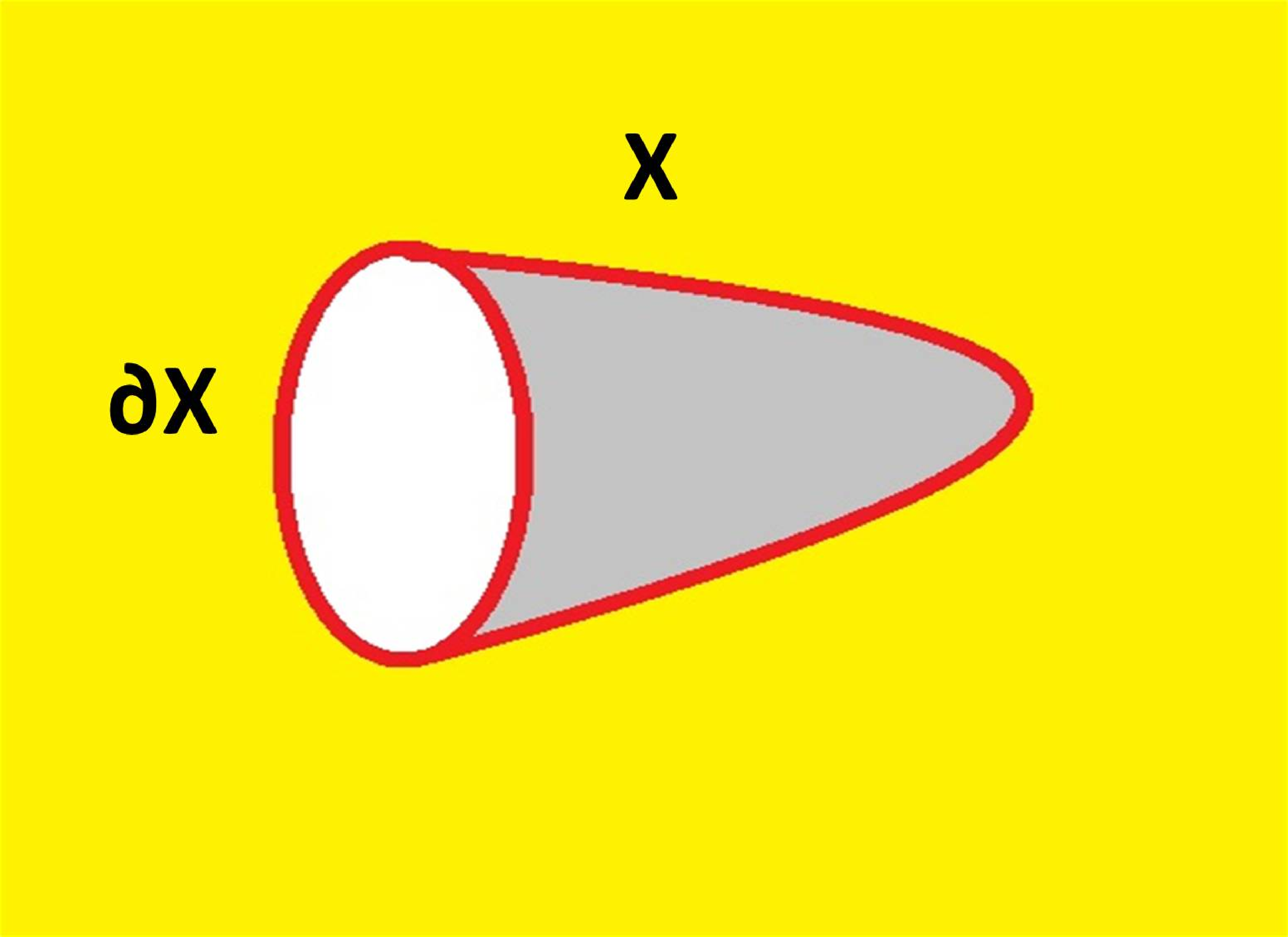} \includegraphics[width=3.85cm]{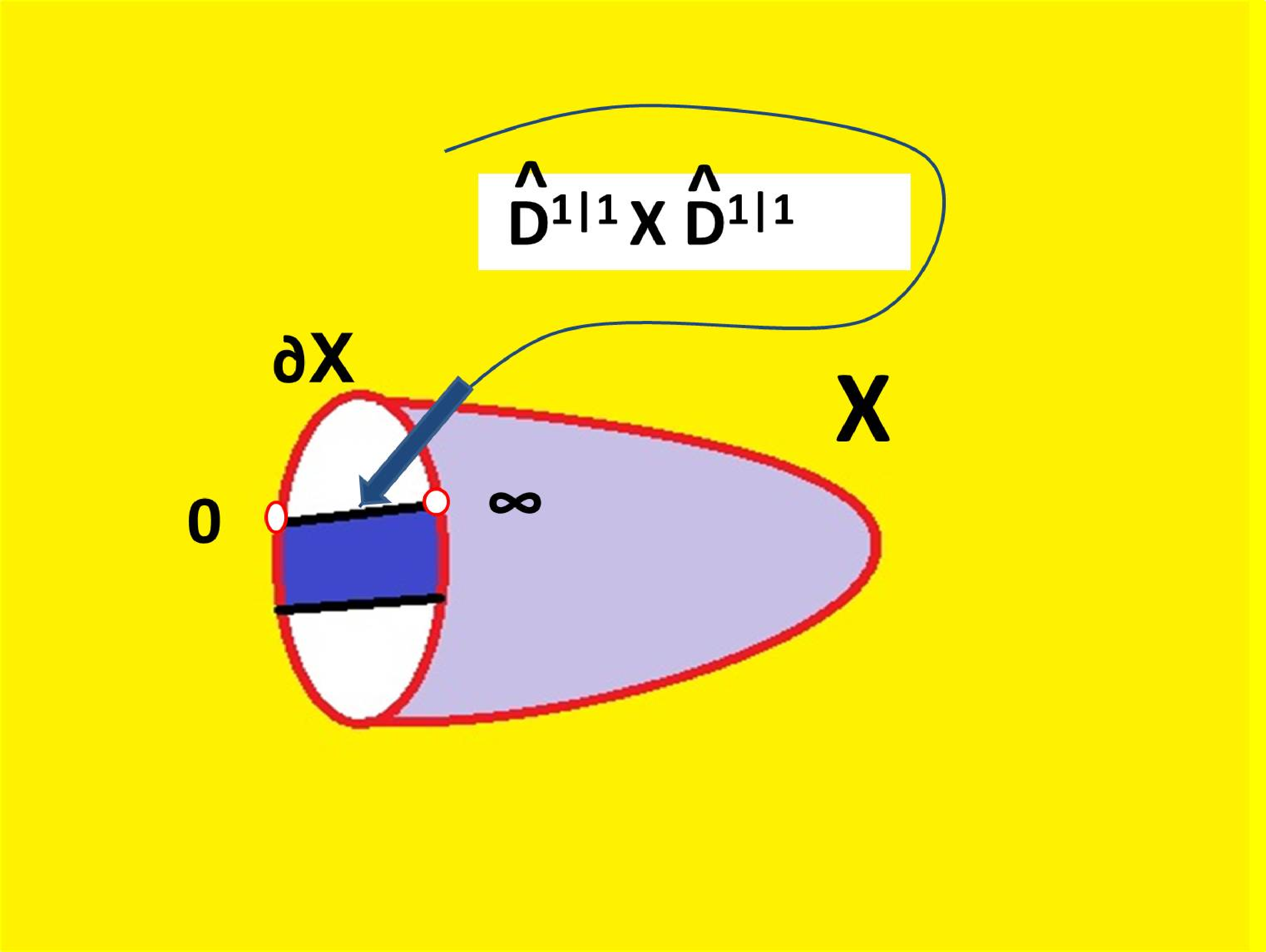}
\includegraphics[width=3.85cm]{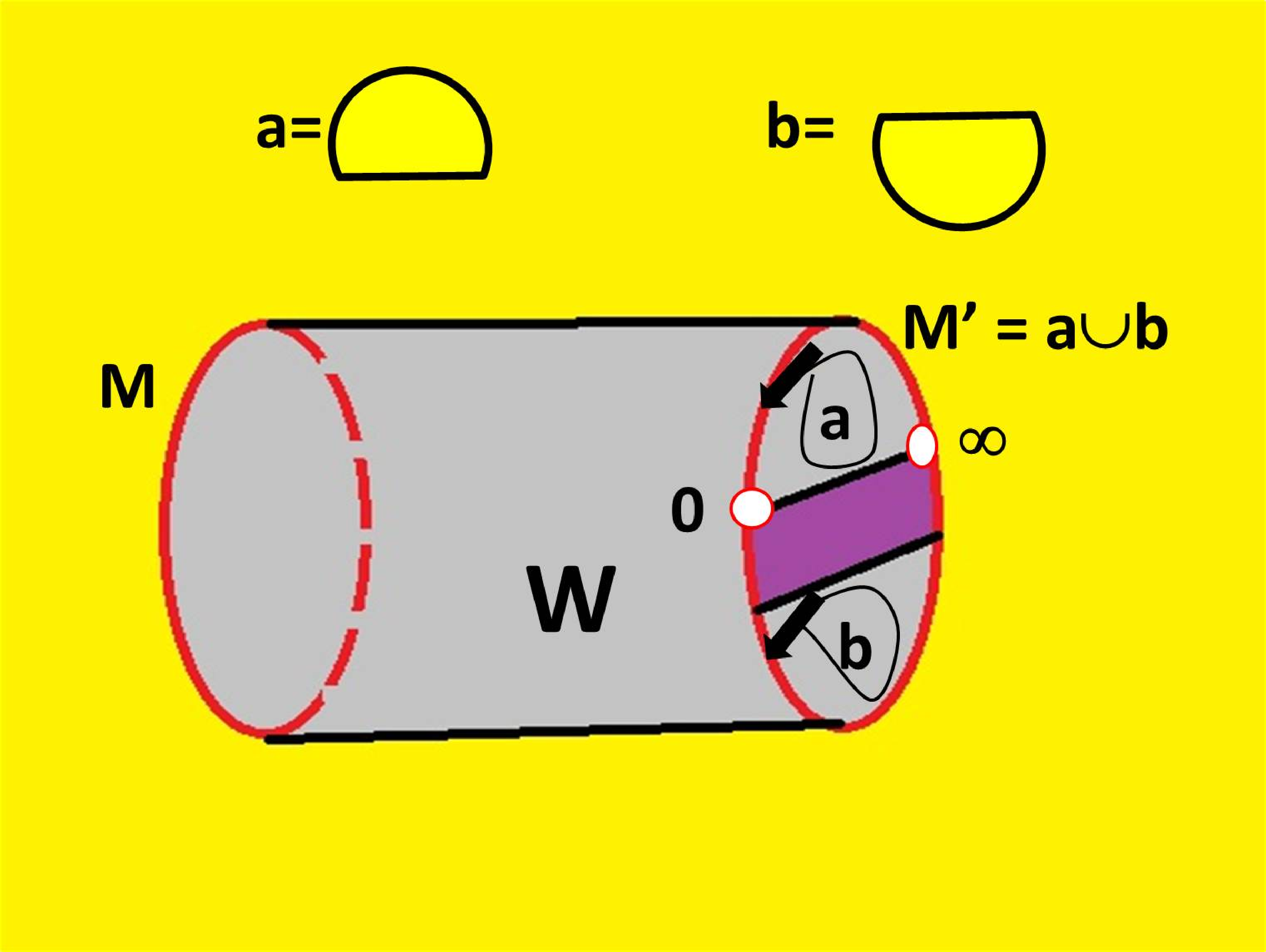}}
\renewcommand{\figurename}{Fig.}
\caption{Relation between attaching quantum superhandles and cobordisms on non-closed quantum supermanifold. The picture is made taking a $(2|2)$-dimensional quantum supermanifold $X$ with boundary $\partial X$. Then $X'\equiv X\bigcup_{\hat S^{0|0}\times \hat D^{1|1}}\hat h^{1|1}$ is obtained from $X$ by attaching on the boundary $M\equiv\partial X$, a $(2|2)$-dimensional quantum superhandle of index $1|1$, $\hat h^{1|1}=\hat D^{1|1}\times \hat D^{1|1}$. For the boundary $M'\equiv\partial X'$ one has $M'=(\partial X\setminus {\rm int}(\hat S^{0|0}\times \hat D^{1|1})\bigcup_{\hat S^{0|0}\times \hat S^{0|0}}\hat D^{1|1}\times \hat S^{0|0}$. This is a disjoint union of two of quantum supercircles. Then the $(2|2)$-dimensional quantum supermanifold $W=(M\times \hat D^{1|1})\bigcup_{\hat S^{0|0}\times \hat D^{1|1}\times\{\infty\}}\hat D^{1|1}\times \hat D^{1|1}$ is the bordism obtained by attaching to $M\times \hat D^{1|1}$ the $(2|2)$-dimensional quantum superhandle of index $1|1$: $\hat h^{1|1}=\hat D^{1|1}\times \hat D^{1|1}$. In fact, one has $\partial W=M\sqcup M'$.}
\label{figure-handle-attachment-vs-cobordism-on-non-closed-manifold}
\end{figure}

\begin{figure}[h]
\centerline{\includegraphics[width=5cm]{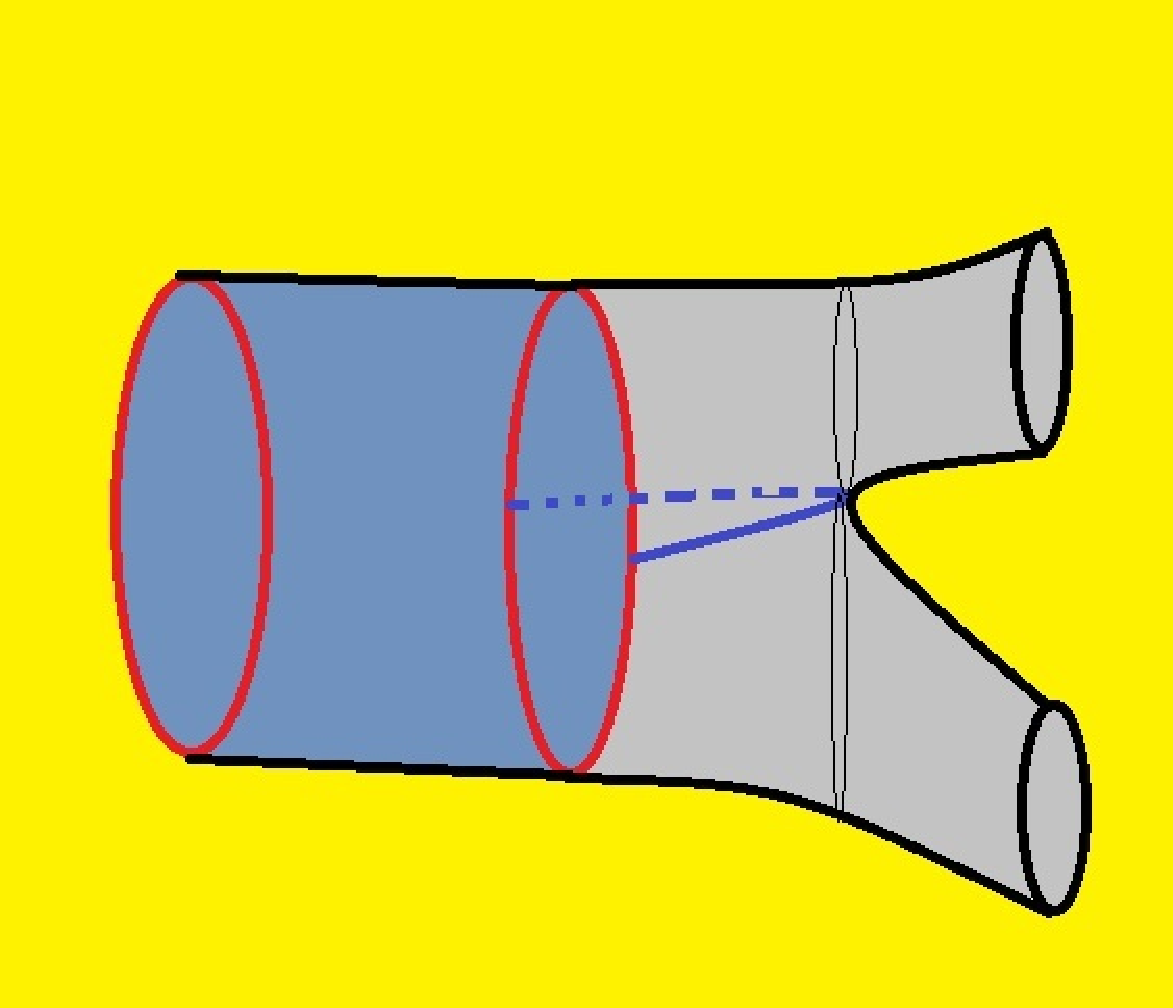} \includegraphics[width=5cm]{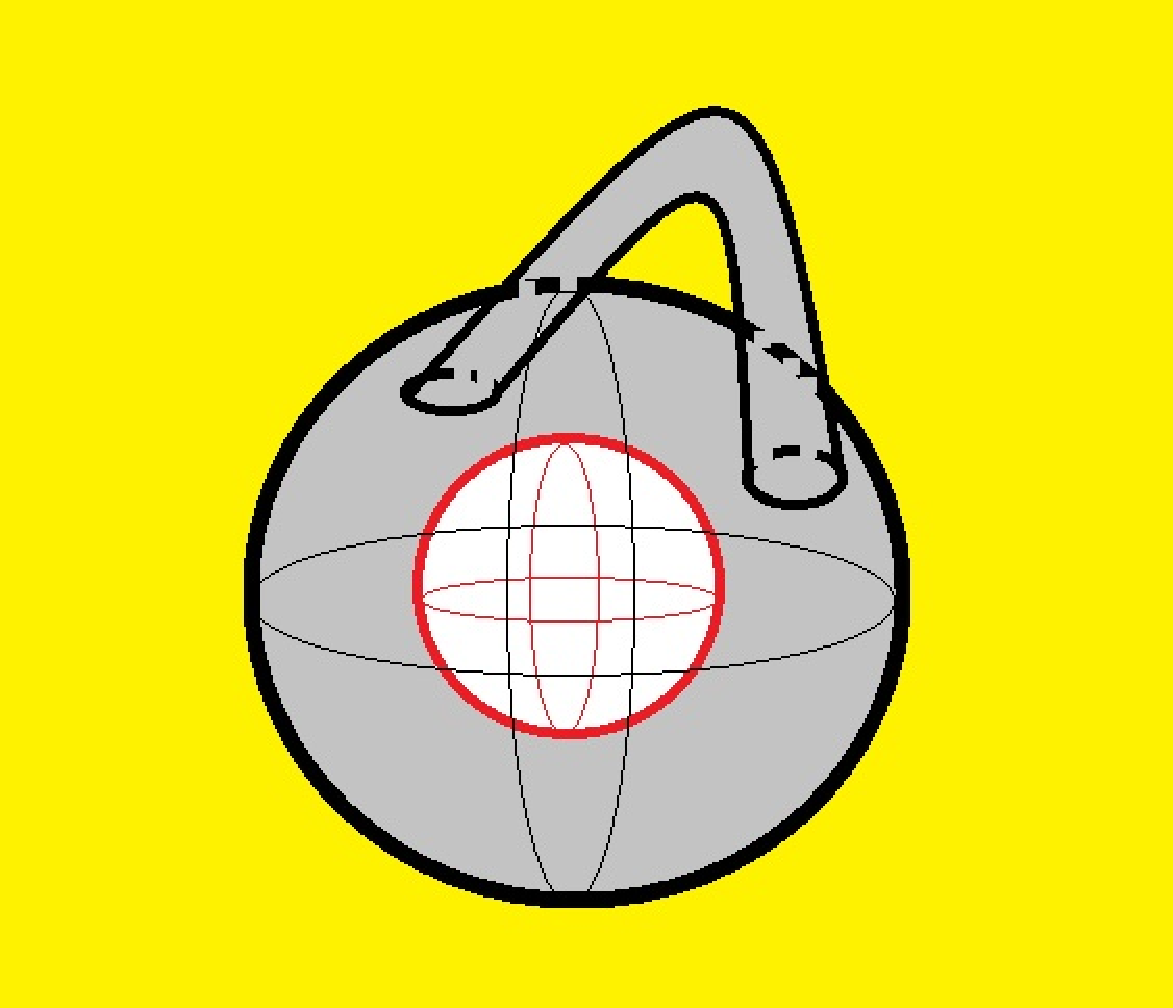}}
\renewcommand{\figurename}{Fig.}
\caption{Relation between attaching quantum superhandles and cobordisms. In the figure on the left, $\dim M=1|1=\dim M'$, $M=\hat S^{1|1}$ and $M'=M\bigcup_{\hat S^{0|0}=\{0,\infty\}}(\hat D^{1|1}\times \hat S^{0|0})$. The $(2|2)$-dimensional quantum supermanifold blue colored is $M\times \hat D^{1|1}$. $W$ is the $(2|2)$-dimensional quantum supermanifold (blue-grey colored) bording $M$ and $M'$, $\partial W=M\sqcup M'$. In the figure on the right, $\dim M=(2|2)=\dim M'$, $M=\hat S^{2|2}$ and $M'=M\bigcup_{\hat S^{0|0}\times \hat S^{1|1}}(\hat D^{1|1}\times \hat S^{1|1})=M\bigcup (\hat D^{1|1}\times \hat S^{1|1})$. $W=(M\times \hat D^{1|1})\bigcup_{\hat S^{0|0}\times \hat D^{2|2}\times\{\infty\}}\hat D^{1|1}\times \hat D^{2|2}=(M\times \hat D^{1|1})\bigcup \hat h^{1|1}$ is the $(3|3)$-dimensional quantum manifold, (grey colored), bording $M$ and $M'$, $\partial W=M\sqcup M'$, where $\hat h^{1|1}$ is the $(3|3)$-dimensional quantum superhandle of index $1|1$.}
\label{figure-handle-attachment-vs-cobordism}
\end{figure}

We shall use the following lemma.
\begin{lemma}\label{quantum superhandle-decomposition-bordism}

{\em 1)} Every bordism $(W;M,M')$ in the category $\mathfrak{Q}_S$, with $\dim W=(m+1|n+1)$, $\dim M=\dim M'=m|n$, has a quantum superhandle decomposition as the union of a finite sequence

$$(W;M,M')=(W_1;M_0,M_1)\bigcup(W_2;M_1,M_2)\bigcup\cdots(W_k;M_{k-1},M_k)$$
of adjoint elementary conditions $(W_j;M_{j-1},M_j)$ with index $p_j|q_j$, such that $$0|0\le p_1|q_1\le p_2|q_2\le\cdots\le p_s|q_s\le m+1|n+1$$ and $M_0=M$, $M_k=M'$.

{\em 2)} Closed $(m|n)$-dimensional quantum super manifolds $M$, $M'$ are cobordant iff $M'$ can be obtained from $M$ by a sequence of surgeries.

{\em 3)} Every closed $(m|n)$-dimensional quantum supermanifold can be obtained from $\varnothing$ by attaching quantum superhandles.
\end{lemma}

\begin{proof}
See Theorem 2.19 in \cite{PRAS19}.\footnote{Handle decomposition of bordism for manifolds has been introduced by Stephen Smale \cite{SMALE}. Lemma \ref{quantum superhandle-decomposition-bordism} generalizes to quantum supermanifolds an analogous result for commutative manifolds. However, the proof of Lemma \ref{quantum superhandle-decomposition-bordism} is not found on the Morse function, (see, e.g., \cite{MILNOR2a}), but on the quantum CW-substitutes for quantum supermanifolds. }
\end{proof}
Let us remark that the proof of Lemma \ref{quantum superhandle-decomposition-bordism} is related to the quantum CW substitute structure of nonlinear quantum propagators, hence the relations considered in this lemma are to consider homeomorphisms. On the other hand the handle decomposition of the nonlinear quantum propagator $V$, is related also to the quantum cell composition that always holds for a $(m+1|n+1)$-dimensional quantum supermanifold, according to its quantum CW-structure. But the relation between a quantum superdisk $\hat D^{r|s}$ and a $(r|s)$-dimensional quantum superhandle $\hat h^{p|q}$ of index $p|q$, is in general an homeomorphism $\hat h^{p|q}\thickapprox \hat D^{p|q}\times \hat D^{r-p|s-q}$, that cannot be reduced to a diffeomorphism in the category $\mathfrak{Q}_S$.\footnote{For example, let us consider the category of complex manifolds, say $\mathfrak{C}$. There the quantum algebra is the $\mathbb{R}$ algebra of complex numbers $\mathbb{C}$, and the morphisms are holomorphic mappings, hence differentiable mappings, having $\mathbb{C}$-linear derivatives. Diffeomorphisms in $\mathfrak{C}$ are usually called biholomorphic mappings. In fact, it is well known that do not exist biholomorphic mappings between an unit ball in $\mathfrak{C}$, (here identified with a {\em quantum $m$-disk}, $\hat D^m$), and the complex $m$-disk, $D^m$, ({\em generalized Poincar\'e's theorem} \cite{KRANTZ}).} In Tab. \ref{handle-decompositions-quantum-supermanifolds-examples} are reported some useful examples of quantum superhandle decompositions of quantum supermanifolds.
\begin{table}[h]
\renewcommand{\tablename}{Tab.}
\caption{Examples of quantum superhandle decompositions of quantum supermanifolds.}
\label{handle-decompositions-quantum-supermanifolds-examples}
\scalebox{0.8}{$
\begin{tabular}{|l|l|l|l|}
\hline
\hfil{\rm{\footnotesize Name}}\hfil&\hfil{\rm{\footnotesize Symbol}}\hfil&\hfil{\rm{\footnotesize Handle-decomposition}}\hfil&\hfil{\rm{\footnotesize Figures}}\hfil\\
\hline
\hline
\hfil{\rm{\footnotesize Quantum $(m|n)$-supersphere}}\hfil&\hfil{\rm{\footnotesize $\hat S^{m|n}$}}\hfil&\hfil{\rm{\footnotesize $\hat h^{0|0}\bigcup\hat h^{m|n}$}}\hfil&\hfil{\rm{\footnotesize \includegraphics[width=1cm]{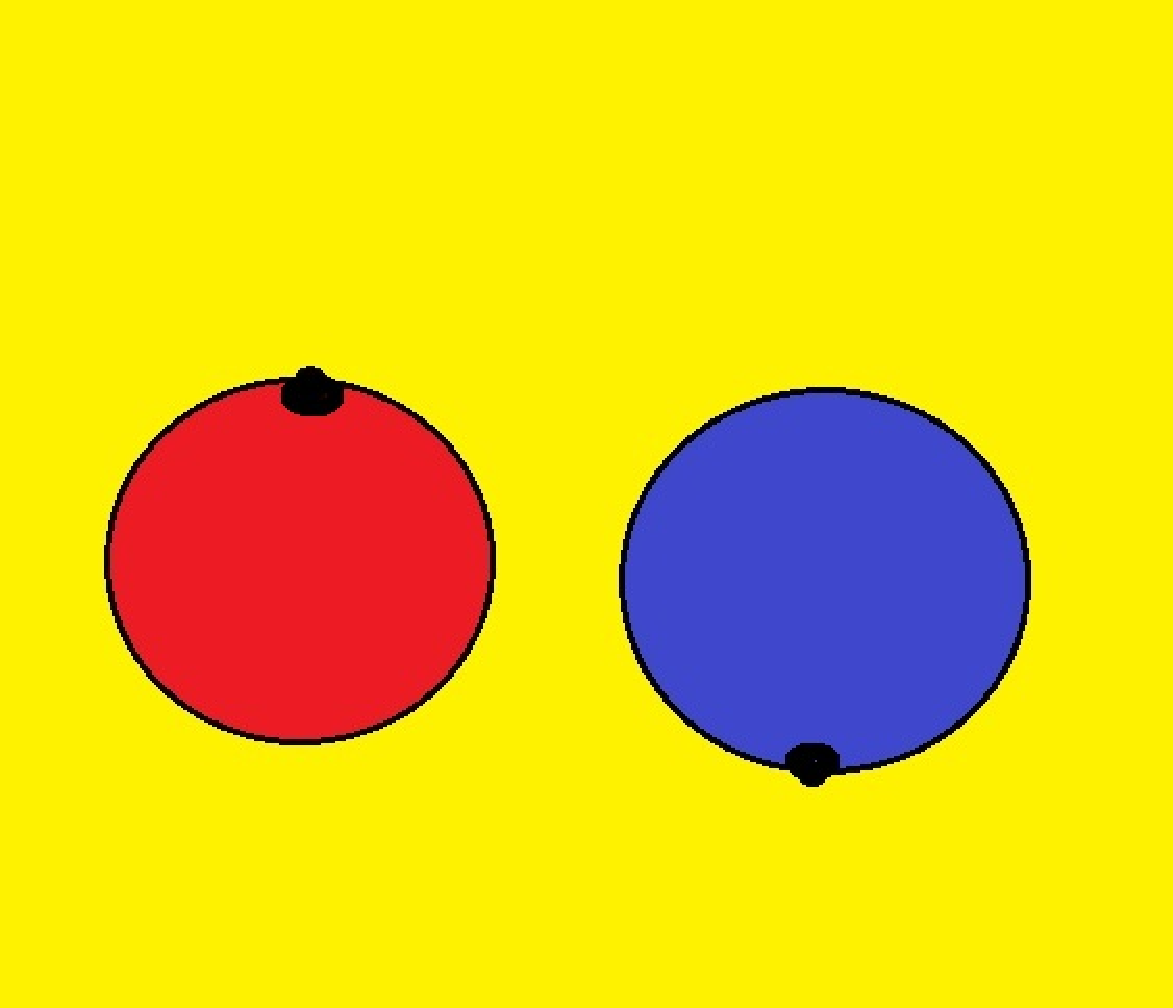}}}\hfil\\
\hline
\hfil{\rm{\footnotesize Quantum $(m+1|n+1)$-supercobordism}}\hfil&\hfil{\rm{\footnotesize $(\hat D^{m+1|n+1};\varnothing,\hat S^{m|n})$}}\hfil&\hfil{\rm{\footnotesize $\hat h^{0|0}$}}\hfil&\hfil{\rm{\footnotesize \includegraphics[width=1cm]{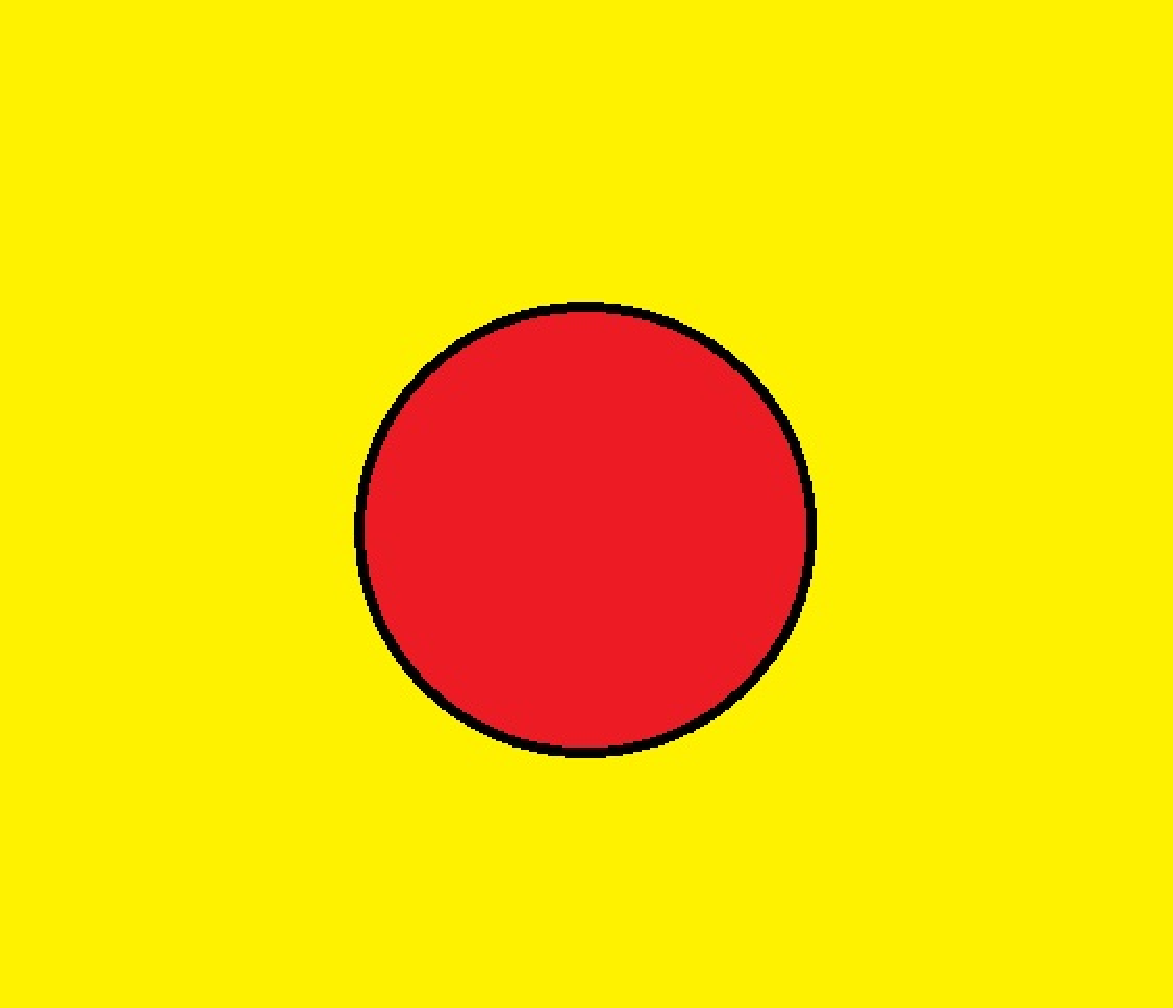}}}\hfil\\
\hline
\hfil{\rm{\footnotesize Quantum $(2|2)$-supertorus}}\hfil&\hfil{\rm{\footnotesize $\hat T^{1|1}=\hat S^{1|1}\times\hat S^{1|1}$}}\hfil&\hfil{\rm{\footnotesize $\hat h^{0|0}\bigcup \hat h^{1|1}\bigcup\hat h^{1|1}\bigcup \hat h^{2|2}$}}\hfil&\hfil{\rm{\footnotesize \includegraphics[width=1cm]{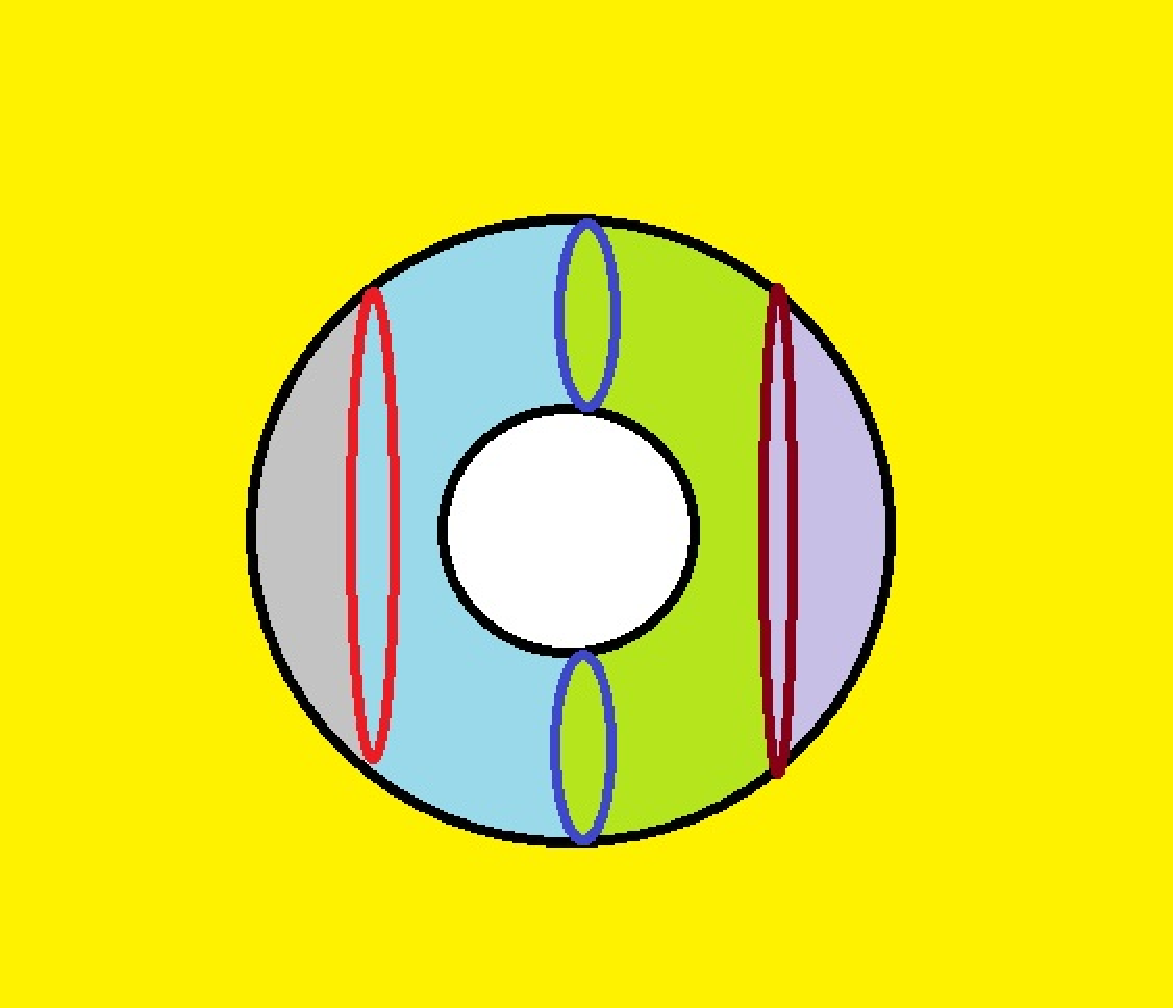}}}\hfil\\
\hline
\hfil{\rm{\footnotesize Quantum punctured M\"obius band}}\hfil&\hfil{\rm{\footnotesize $(\hat M_{ob}\setminus \hat D^{2|2};\hat S^{1|1},\hat S^{1|1})$}}\hfil&\hfil{\rm{\footnotesize $\hat S^{1|1}\times\hat D^{1|1}\bigcup\hat h^{1|1}$}}\hfil&\hfil{\rm{\footnotesize \includegraphics[width=1cm]{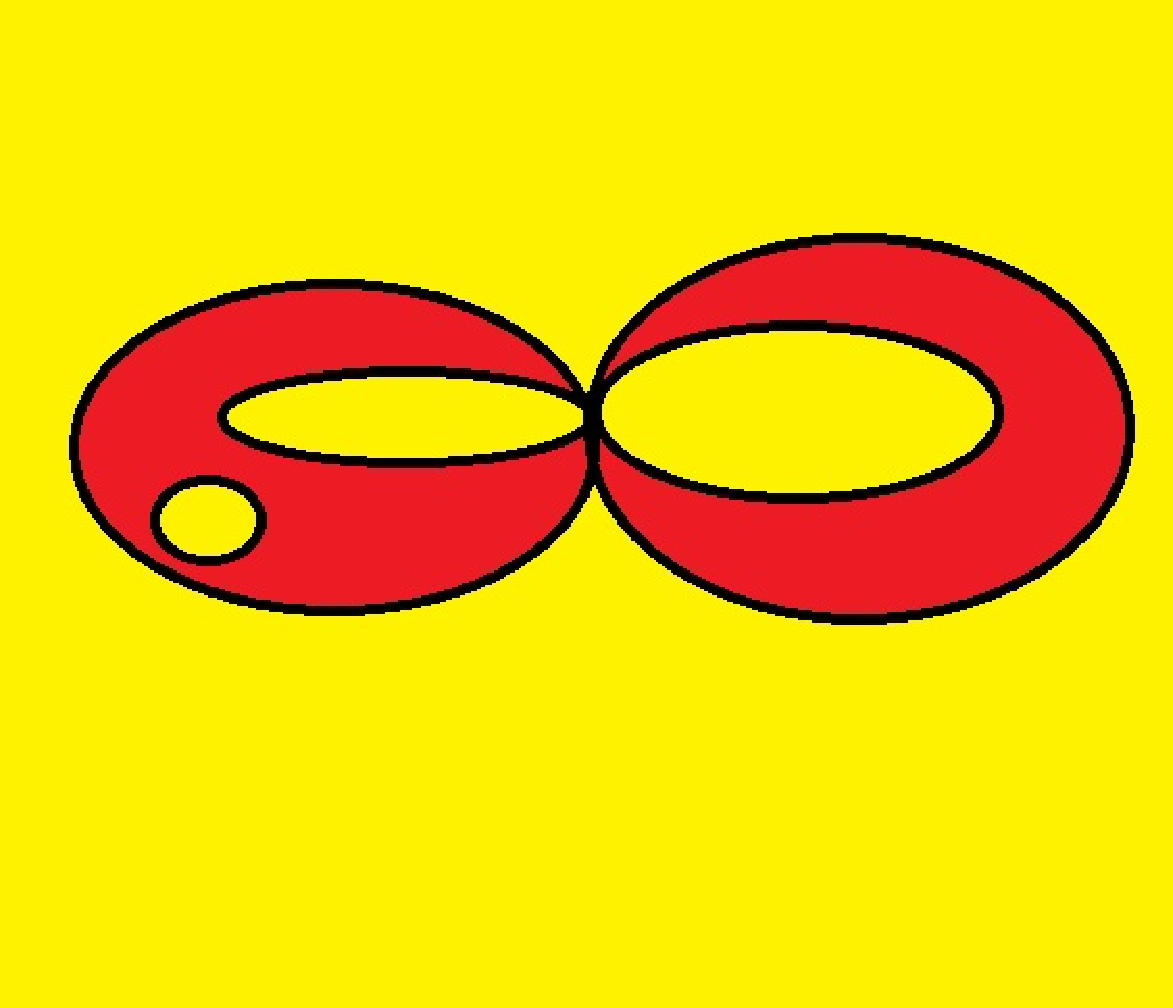}}}\hfil\\
\hline
\end{tabular}$}
\end{table}

Since $V$ is a $(m+1|n+1)$-dimensional bordism between $(m|n)$-dimensional quantum supermanifolds $N_0$ and $N_1$, we can determine its handle decomposition (\ref{quantum-superhandle-presentation-solution}), that identifies intermediate final Cauchy quantum supermanifolds $M_{j}$ that are obtained from the previous one $M_{j-1}$ by means of an integral surgery. Let us emphasize that each $M_{j}$ must necessarily be an integral quantum supermanifold having the same dimension of $N_0$ since it belongs to the same bordism class $[N_0]$. Furthermore, since each quantum supermanifold $M_j$ is contained into $V$, it follows that  $M_{j}$ can be identified with a Cauchy $(m|n)$-dimensional quantum supermanifold of $\hat E_k$, i.e. an admissible $(m|n)$-dimensional quantum integral supermanifold of $\hat E_k\subset \hat J^k_{m+1|n+1}(W)$. Therefore $M_j$ belongs to the same bordism class $[N_0]\in\Omega_{m|n}^{\hat E_k}$ of $N_0$ in $\hat E_k$. The situation is pictured in Fig. \ref{figure-proof-quantum-superhandle-decomposition-quantum-nonlinear-propagator}.
\end{proof}
\begin{definition}[Quantum exchangions and quantum virtual particles]\label{quantum-integral-handle-exchange}
   Let $\hat E_k\subset \hat J^k_{m+1|n+1}(W)$ be a quantum super PDE in the category $\mathfrak{Q}_S$ and let $V\subset \hat E_k$ be a compact solution of a boundary value problem, i.e., $\partial V=N_0\sqcup N_1$. We call {\em quantum exchangions} the quantum integral superhandle that are attached to obtain intermediate Cauchy manifolds  $N_{j|j}$ defined in the proof of Theorem \ref{quantum-superhandle-decomposition-quantum-nonlinear-propagator}. We call {\em quantum virtual particles} the intermediate quantum supermanifolds $M_j$ appearing there. (See, also Fig. \ref{figure-proof-quantum-superhandle-decomposition-quantum-nonlinear-propagator}.)
\end{definition}

\begin{remark}
Note the structural difference between quantum exchangions and virtual particles. In fact, the first are quantum $(m+1|n+1)$-chains, hence having the same dimension of solutions of $\hat E_k$. Instead, the second ones are quantum $(m|n)$-chains, hence having the same dimension of the initial and final Cauchy data.\footnote{It is useful to emphasize that quantum exchangions and quantum virtual particles, well interpret the meaning of ``exchange particles" introduced in Feynman diagrams (1949) to represent particle interactions, like photons (electro-magnetic interactions), $W$ and $Z$ particles (weak interactions), gluons and gravitons (strong interactions). (See, e.g., Refs. \cite{FEYNMAN, KAISER}.) Furthermore, in the framework of Theorem \ref{quantum-superhandle-decomposition-quantum-nonlinear-propagator} we understand also that {\em Regge trajectories}  and {\em Regge resonances} (1959), introduced in the phenomenological theory of strong reactions by the pioneering works of T. Regge \cite{REGGE} (and principally developed also by R. Blanckebecker and M. L. Goldberger \cite{BLANKE-GOLD}, G. F. Chew and S. C. Fraustschi \cite{CHEW-FRAUTSCHI}, V. N. Gribov \cite{GRIBOV} and G. Veneziano \cite{VENEZIANO}), find their interpretation as quantum exchangions and quantum virtual particles.}
\end{remark}

\begin{proposition}\label{quantum-exchangions-and-cobordism-class}
A quantum exchangion does not change the integral bordism class of $N_{j}$ with respect to $N_0$ one.
\end{proposition}
\begin{proof}
In fact $N_0$ and $M_{j}$ must necessarily belong to the same quantum integral bordism class in $\Omega_{m|n}^{\hat E_k}$ since $\partial V_{1}=N_0\sqcup M_{2}$ and $\partial V_{2}=M_2\sqcup M_{3}$, and so on.
\end{proof}
\begin{proposition}
Whether $N_0$ is the disjoint union of two (or more) components, e.g., $N_0=a\sqcup b$, then a quantum exchangion can change the quantum numbers of $a$ and $b$, even if the total quantum number of $a\sqcup b$ does not change.
\end{proposition}
\begin{proof}
This follows directly from Proposition \ref{quantum-exchangions-and-cobordism-class}.
\end{proof}

\begin{definition}[Fundamental quantum particles]\label{fundamental-quantum-particles}
For a quantum system, encoded by a quantum PDE $\hat E_k\subset \hat J^k_{m+1|n+1}(W)$, we define {\em quantum $(m|n)$-particles}, admissible quantum integral $(m|n)$-chains $N\subset \hat E_k\subset \hat J^k_{m+1|n+1}$.

We call {\em fundamental quantum $(m|n)$-particles} for the quantum PDE $\hat E_k$, ones that cannot be decomposed into other quantum particles.
\end{definition}

\begin{proposition}
The fundamental quantum $(m|n)$-particles of quantum PDE $\hat E_k\subset \hat J^k_{m+1|n+1}(W)$, are identified with the bordism classes of the weak integral bordism group $\Omega_{m|n,w}^{\hat E_k}$, (resp. singular integral bordism group $\Omega_{m|n,s}^{\hat E_k}$, resp. integral bordism group $\Omega_{m|n}^{\hat E_k}$). Then we will distinguish into {\em weak-fundamental quantum  $(m|n)$-particles}, {\em singular-fundamental quantum  $(m|n)$-particles}, resp. {\em fundamental quantum  $(m|n)$-particles} respectively, for the quantum PDE $\hat E_k\subset \hat J^k_{m+1|n+1}(W)$. Furthermore the exact commutative diagram {\em(\ref{relations-between-fundamental-quantum-particles-classes})} shows how such fundamental quantum particles are related.

\begin{equation}\label{relations-between-fundamental-quantum-particles-classes}
\begin{array}{ccccccccc}
&&{\scriptstyle 0}&&{\scriptstyle 0}&&{\scriptstyle 0}&&\\ &&\downarrow&&\downarrow&&\downarrow&&\\
{\scriptstyle 0}&\to&{\scriptstyle K^{\hat E_k}_{m-1|n-1,w/(s,w)}}&\to&{\scriptstyle K^{\hat E_k}_{m-1|n-1,w}}&\to&
{\scriptstyle K^{\hat E_k}_{m-1|n-1,s,w}}&\to&{\scriptstyle 0}\\
&&\downarrow&&\downarrow&&\downarrow&&\\
{\scriptstyle 0}&\to&{\scriptstyle K^{\hat E_k}_{m-1|n-1,s}}&\to&{\scriptstyle \Omega^{\hat E_k}_{m-1|n-1}}&\to&
{\scriptstyle \Omega^{\hat E_k}_{m-1|n-1,s}}&\to&{\scriptstyle 0}\\
&&\downarrow&&\downarrow&&\downarrow&&\\
&&{\scriptstyle 0}&\to&{\scriptstyle \Omega^{\hat E_k}_{m-1|n-1,w}}&\to&{\scriptstyle \Omega^{\hat E_k}_{m-1|n-1,w}}&\to
&{\scriptstyle 0}\\
&&&&\downarrow&&\downarrow&&\\
&&&&{\scriptstyle 0}&&{\scriptstyle 0}&&\\
\end{array}
\end{equation}

\end{proposition}
\begin{proof}
The proof follows directly from Definition \ref{fundamental-quantum-particles}, Theorem \ref{quantum-superhandle-decomposition-quantum-nonlinear-propagator} and Proposition 3.3 in \cite{PRAS14}.
\end{proof}

\begin{remark}
Let us emphasize that the concept of {\em fundamental quantum particles}, is strictly related to the quantum (super) PDE $\hat E_k$, or in other words, to the quantum theory encoded by $\hat E_k$. Whether we aim to forget the framework defined by $\hat E_k$, and to discuss about fundamental quantum particles in general, we see that the unique fundamental one is $\varnothing$. In fact, according to Lemma \ref{quantum superhandle-decomposition-bordism}(3), we can build any quantum particles just starting from $\varnothing$, and adding quantum handles or quantum-cells.\footnote{From this geometric structural point of view, we clearly understand that also some particles, that are usually considered fundamental, cannot be considered so. For example electrons and quarks are not quantum fundamental particles whether considered as geometric objects. This agrees with experimental evidences that so-called ``{\em quasiparticles}" with fractional charges, or with separated quantum numbers (e.g., decay $e^-\to {\rm spinon}+{\rm orbiton}$), were detected. (See, e.g. \cite{FRADKIN, SCHLAPPA-SCHMITT}.) Furthermore, it is important to distinguish the concept of quantum fundamental particle, from the one of quantum stable particle. Without any quantum theory, i.e., without any quantum (super) PDE, for the first the unique one is $\varnothing$, and for the second, it is a no sense. On the other hand, with respect to a quantum (super) PDE $\hat E_k$ we can understand that the concept of stability is just related to the integral bordism groups of $\hat E_k$, hence from this point of view the concept of ``fundamental particle" and ``stable particle" become related each other one.}
\end{remark}

\begin{example}[Two-body high energy reactions]\label{two-body-high-energy-reactions-examples}
Let us consider the following typical two-body strong reactions:
$\pi^-+p\to\pi^0+n$, $\pi^++p\to\pi^++p$, $\pi^++p\to p+\pi^+$. The first reaction can be considered obtained with a nonlinear quantum propagator having an intermediate virtual neutral $\rho$ meson, with charge exchange.  In the second reaction the nonlinear quantum propagator has an intermediate virtual double electric charged particle and a neutral quantum pomeron $p_{om}$, quantum exchangion $\hat h^{1|1}$, representing an elastic scattering. In the third reaction can be considered the nonlinear quantum propagator has an intermediate virtual neutron $n$, in backward scattering. In other words, in all these reactions the final Cauchy data can be obtained from the initial one by means of two adjoint elementary bordisms.\footnote{The quantum virtual particles there involved can be considered generalizations in the geometric theory of quantum (super)PDEs, of objects called reggeons in analogous reactions considered in the framework of the phenomenological theory of strong reactions developed in the first two decades of the second half of the last century. (See, e.g., \cite{EDEN}.)  However, let us emphasize that {\em pomeron} is better represented as quantum exchangion than a quantum virtual particle. In fact, pomeron, as usually considered in the phenomenological theory of strong reactions, does not carry charge \cite{EDEN}. On the other hand in the elastic scattering $\pi^++p\to\pi^++p$ it is impossible that the quantum virtual particle should be neutral one, whether the nonlinear quantum propagator is not exotic. (See \cite{PRAS29}.)}
\end{example}

\begin{example}[Proton-proton chain reactions]\label{proton-proton-chain-reactions-examples}
In the following we shall see that solutions, representing proton-proton chain reactions, typical in the Sun for the production of ${}^4_2He$, can be encoded with nonlinear quantum propagators, admitting handle decompositions. In general one can classify such reactions into five basic groups.\footnote{Let us recall the nucleus-notation ${}^A_ZX$, where $X$ denotes the chemical symbol, $A=Z+N$ is the mass number, with $Z$ the atomic number=number of protons and $N$ the number of neutrons.}

{\em (I) Production of ${}^4_2He$ with intermediate particles $\{e^+,\, \nu_e,\, \gamma,\, {}^2_1H\equiv {}^2_1D,\, {}^3_2He\}$. \footnote{In the Sun this chain reaction is dominant in the temperature range $10-14\, MK$. In presence of electrons, on has also the secondary reaction $e^++e^-\to 2\gamma + 1.02 Mv$.}}
\begin{equation}\label{proton-proton-chain-reactions-examples-a}
\left\{
\begin{array}{l}
  {}^1_1H+{}^1_1H\to e^++\nu_e+{}^2_1H+\, 0.42\, MeV\\
   {}^2_1D+{}^1_1H\to {}^3_2He+\gamma+ 5.49\, MeV\\
  {}^3_2H+{}^3_2He\to 2\, {}^1_1H+{}^4_2He+\, 12.86\, MeV\\
\end{array}
\right.
\end{equation}

{\em (II) Production of ${}^4_2He$ with intermediate particles $\{\nu_e,\, \gamma,\, {}^7_4Be,\, {}^7_3Li\}$. \footnote{In the Sun this chain reaction is dominant in the temperature range $14-23\, MK$.}}
\begin{equation}\label{proton-proton-chain-reactions-examples-b}
\left\{
\begin{array}{l}
  {}^3_2He+{}^4_2He\to {}^7_4Be+\gamma\\
   {}^7_4Be+e^-\to{}^7_3Li+\nu_e+\, 0.383-0.861\, MeV\\
  {}^7_3Li+{}^1_1H\to 2\, {}^4_2He\\
\end{array}
\right.
\end{equation}

{\em (III) Production of ${}^4_2He$ with intermediate particles $\{e^+,\, \nu_e,\, \gamma,\, {}^7_4Be,\, {}^8_4Be,\, {}^8_5B\}$. \footnote{In the Sun this chain reaction is dominant in the temperatures that exceeds $23\, MK$.}}
\begin{equation}\label{proton-proton-chain-reactions-examples-c}
\left\{
\begin{array}{l}
  {}^3_2He+{}^4_2He\to+{}^7_4Be+\gamma\\
   {}^7_4Be+{}^1_1H\to {}^8_5B+\gamma\\
  {}^8_5B\to  {}^8_4Be+e^++\nu_e\\
  {}^8_4Be\to  2\, {}^4_2He\\
\end{array}
\right.
\end{equation}

{\em (IV) Production of ${}^4_2He$ from the  ${}^3_2He$ interaction with proton.}
\begin{equation}\label{proton-proton-chain-reactions-examples-d}
{}^3_2He+{}^1_1H\to{}^4_2He+e^++\nu_e+\, 18.8\, MeV.
\end{equation}

{\em (V) Production of deuterium from the electron capture by two protons.}
\begin{equation}\label{proton-proton-chain-reactions-examples-d}
{}^1_1H+e^-+{}^1_1H\to{}^2_1D+\nu_e.
\end{equation}
In Fig. \ref{figure-proton-proton-chain-reactions-examples} we represent some solutions corresponding to reactions in {\em(I)}. By using the definition introduced in \cite{PRAS22}, we can say that in both reactions there represented, the nonlinear quantum propagator $V$, is a {\em quantum matter-solution}. Furthermore, in the second reaction, the handle decomposition of $V$, identifies a {\em Goldstone piece} $V_{\partial}=V\bigcap(\widehat{Goldstone})$. This is the part arriving to $\gamma$.
\begin{figure}[h]
\centerline{\includegraphics[width=5cm]{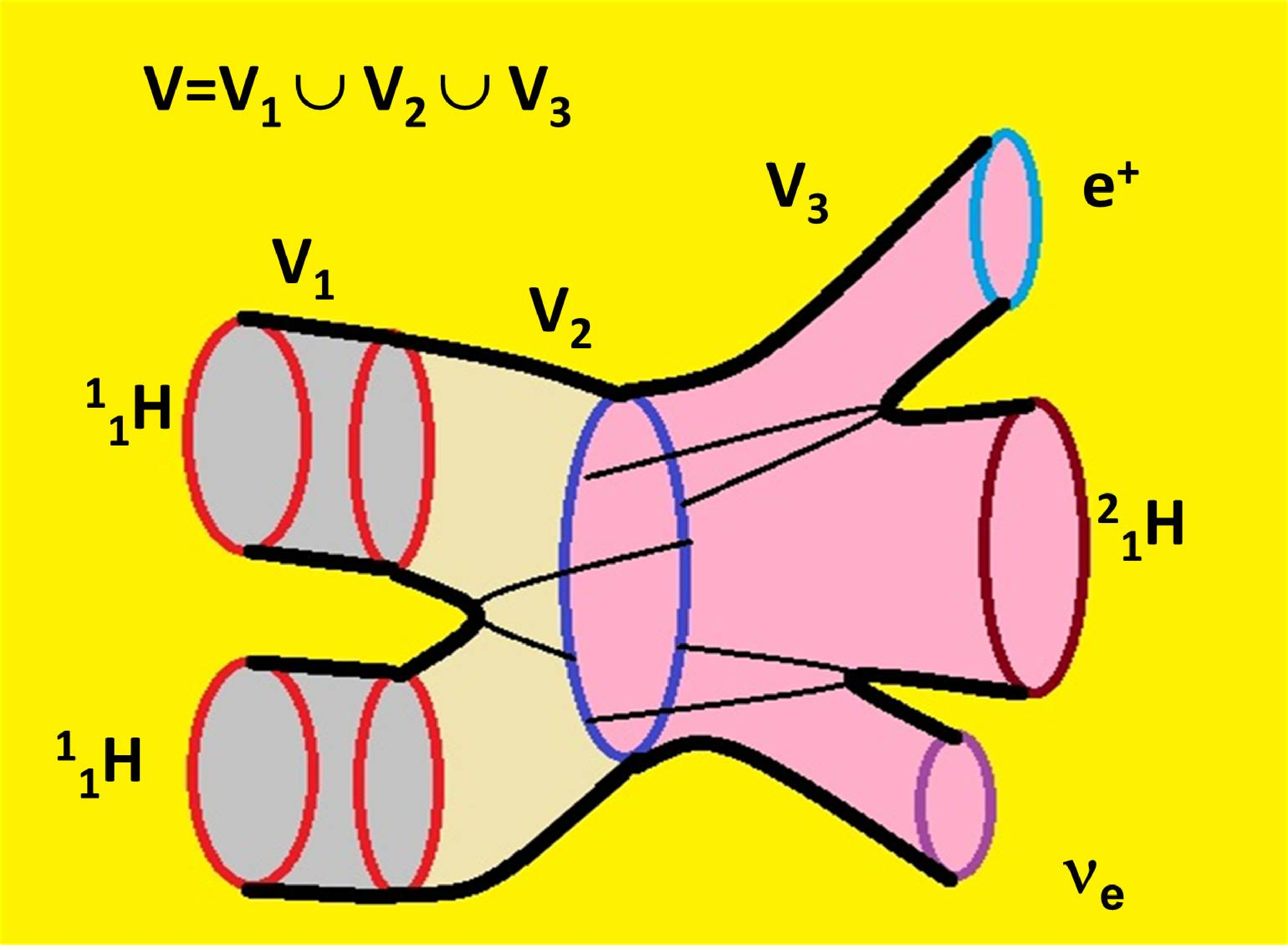} \includegraphics[width=4.9cm]{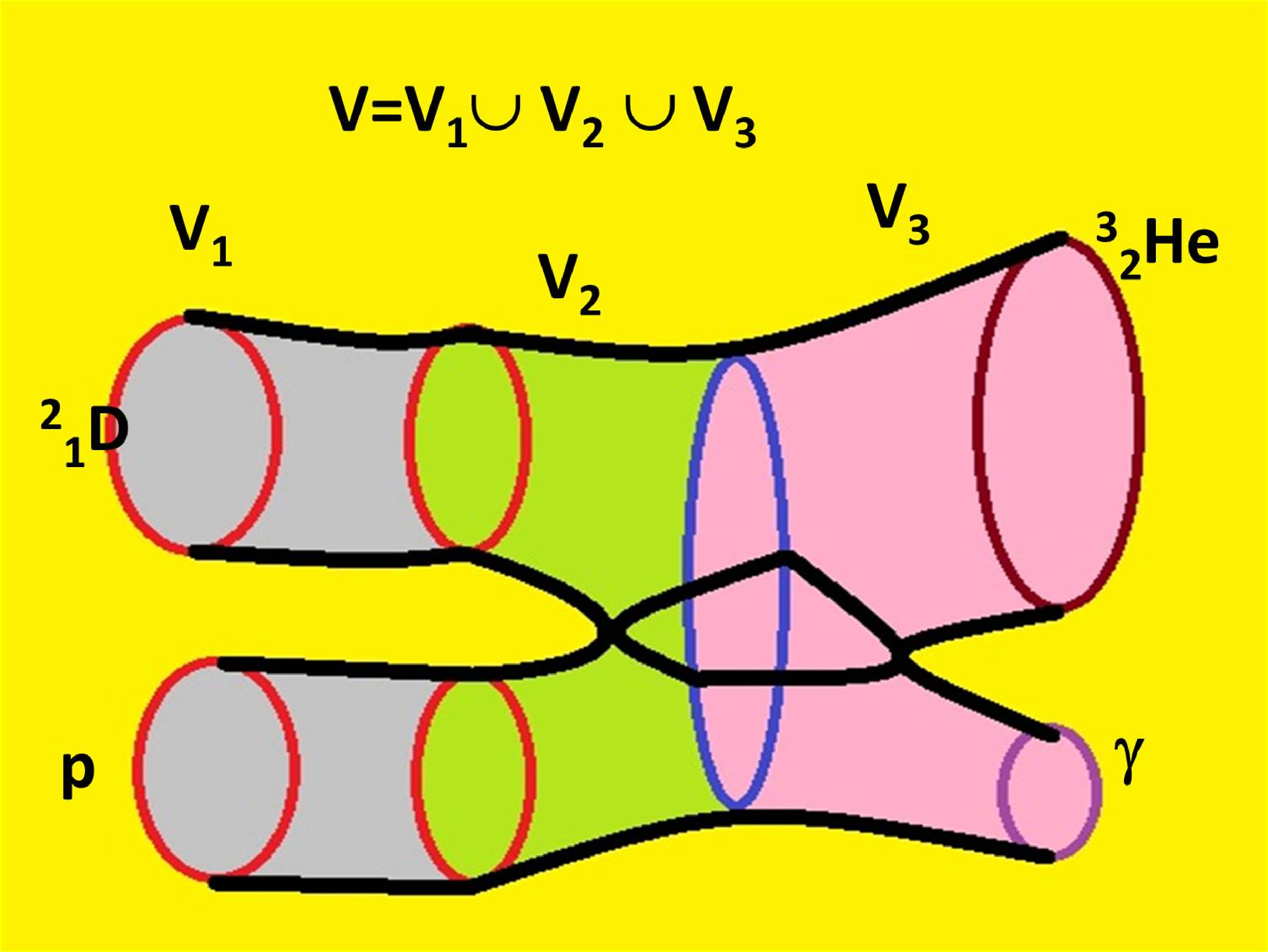}}
\renewcommand{\figurename}{Fig.}
\caption{Representation of some solutions of $p+p$ reaction chain. From left to right. ${}^1_1H+{}^1_1H\to e^++\nu_e+{}^2_1H$. Here the nonlinear quantum propagator $V$ has a three step decomposition:  $V_{1}\bigcup  V_{2}\bigcup  V_{3}$. Furthermore there is a virtual quantum particle $M_3$ with a double electric charge. One has $V_1=M_0\times\hat D^{1|1}$, $V_2=\hat h^{1|1}$ and $V_3=\hat h^{1|1}\bigcup \hat h^{1|1}$. The final Cauchy data is represented by three disjoint supermanifolds. In the reaction ${}^2_1D+{}^1_1H\to {}^3_2He+\gamma$, the nonlinear quantum propagator $V$ has a three step decomposition too:  $V_{1}\bigcup  V_{2}\bigcup  V_{3}$, but  $V_1=M_0\times\hat D^{1|1}$, $V_2=\hat h^{1|1}$ and $V_3=\hat h^{1|1}$. There is a double electric charged virtual quantum particle $M_3$, and the final Cauchy data is represented by two disjoint quantum supermanifolds.}
\label{figure-proton-proton-chain-reactions-examples}
\end{figure}
\end{example}

\section{\bf Surgery in Quantum Super Yang-Mills PDEs}\label{sec-surgery-quantum-super-yang-mills-pdes}

Let us emphasize, now, that in the reactions considered in the above section, we have not really specified the contribution of a specific quantum PDE, i.e., we have not considered some specific integral bordism group. On the other hand, without the introduction of these fundamental structures, what we can do is reproduce some phenomenological theory like dispersion relations and Regge models that cannot be but predictive dynamical theories. Therefore to dynamically encode, for example, proton-proton reactions we shall use the quantum super Yang-Mills PDE, and to solve suitable boundary value problems. We have formulated this general theory in some our previous works \cite{PRAS21,PRAS22}. In the following we shall consider some particular applications.

\begin{definition}[Quantum scattering processes in $\widehat{(YM)}$]\label{quantum-scattering-processes-in-quantum-yang-mills}
We call {\em quantum scattering process} in $\widehat{(YM)}$ any boundary value problem where are fixed two disjoint Cauchy data $N_0,\, N_1\in\widehat{(YM)}$. The first is called {\em initial Cauchy data} and the second the {\em final Cauchy data}. A {\em solution of a quantum scattering process} in $\widehat{(YM)}$ is any nonlinear quantum propagator $V\subset\widehat{(YM)}$, such that $\partial V=N_0\sqcup N_1$, if $\partial N_0=\partial N_1=\varnothing$, otherwise  $\partial V=N_0\sqcup P\sqcup N_1$, with $P\subset\widehat{(YM)}$, integral manifold such that $\partial P=\partial N_0\sqcup\partial N_1$.
\end{definition}

\begin{theorem}[Quantum scattering processes in $\widehat{(YM)}$ on $4$-dimensional quantum super Minkowskian manifold]\label{quantum-scattering-processes-in-quantum-yang-mills-d-4-quantum-super-minkowski}
Let us assume that the base manifold $M$ of $\widehat{(YM)}$ is a $(4|4)$-dimensional quantum super Minkowski manifold.
Then $\widehat{(YM)}$ is a quantum extended crystal PDE. Furthermore, if we consider admissible only integral boundary manifolds, with
orientable classic limit, and with zero characteristic quantum
supernumbers, ({\em full admissibility hypothesis}), one has:
$\Omega_{3|3}^{\widehat{(YM)}}=0$, and $\widehat{(YM)}$ becomes a
quantum $0$-crystal super PDE. Hence we get existence of global
$Q^\infty_w$ solutions for any boundary condition of class
$Q^\infty_w$.

{\em Elementary nonlinear quantum propagators} $V$ are $(4|4)$-dimensional quantum supermanifolds with boundary $\partial V=N_0\bigcup P\bigcup N_1$, such that $V$ is homeomorphic to $\hat D^{4|4}$. One can classify $V$ on the ground of its boundary $\partial V$. This can be diffeomorphic to a quantum exotic supersphere $\hat\Sigma^{3|3}$ or to a quantum supersphere $\hat S^{3|3}$.\footnote{In general nonlinear quantum propagators are not elementary ones, but can be decomposed in elementary ones.}
\end{theorem}

\begin{proof}
In $D=4$, the usual $N$-supersymmetric extension $\mathfrak{g}$ of the
Poincar\'e algebra $\mathfrak{p}=\mathfrak{s}\mathfrak{o}(1,3)\oplus
\mathfrak{t}$, is a $\mathbb{Z}_2$-graded vector space
$\mathfrak{g}=\mathfrak{g}_0\oplus\mathfrak{g}_1$, with a graded Lie
bracket, such that $\mathfrak{g}_0=\mathfrak{p}\oplus \mathfrak{b}$,
where $\mathfrak{b}$ is a reductive Lie algebra, such that its
self-adjoint part is the tangent space to a real compact Lie
group.\footnote{A {\em reductive} Lie algebra is the sum of a
semisimple and an abelian Lie algebra. Since a {\em semisimple} Lie
algebra is the direct sum of simple algebras, i.e., non-abelian Lie
algebras, $\mathfrak{l}_i$, where the only ideals are $\{0\}$ and
$\{\mathfrak{l}_i\}$, it follows that $\mathfrak{b}$ can be
represented in the form
$\mathfrak{b}=\mathfrak{a}\oplus\sum_i\mathfrak{l}_i$.} Furthermore
$\mathfrak{g}_1=(\frac{1}{2},0)\otimes
\mathfrak{s}\oplus(0,\frac{1}{2})\otimes\mathfrak{s}^*$, where
$(\frac{1}{2},0)$ and $(0,\frac{1}{2})$ are specific representations
of the Poincar\'e algebra. Both components are conjugate to each
other under the $*$ conjugation. $\mathfrak{s}$ is a $N$-dimensional
complex representation of $\mathfrak{b}$ and $\mathfrak{s}^*$ its
dual representation.\footnote{If $\rho:\mathfrak{g}\to L(V)$ is a
representation of Lie algebra, its dual $\bar\rho:\mathfrak{g}\to
L(\bar V)$, working on the dual space $\bar V$, is defined by
$\bar\rho(u)=\overline{-\rho(u)}$, $\forall u\in\mathfrak{g}$.} Note
also that the Lie bracket for the odd part is usually denoted by
$\{,\}$ in theoretical physics. Then with such a notation one has
\begin{equation}\label{bracket-odd-part}
\{Q^i_\alpha,Q^j_\beta\}=\delta^{ij}(\gamma^\mu
C)_{\alpha\beta}P_\mu+U^{ij}(C)_{\alpha\beta}+V^{ij}(C\gamma_5)_{\alpha\beta}
\end{equation}
where $U^{ij}=-U^{ji}$, $V^{ij}=-V^{ji}$ are the $(N-1)N$ central
charges, $C$ is the (antisymmetric) charge conjugation matrix,
$(Q^i_\alpha)_{i=1,\dots,N}$, are the $N$ Majorana spinor
supersymmetry charge generators. The dynamical components
$\hat\mu^i$, $i=1,\dots,N$, of the quantum fundamental field,
corresponding to the generators $Q^i$, are called {\em quantum
gravitinos}. So in a quantum $N$-SG-Yang-Mills PDE, one
distinguishes $N$ quantum gravitino types, (and $(N-1)N$ central
charges).

Then a quantum superextension of $\mathfrak{g}$ is
 $A\otimes_{\mathbb{R}}\mathfrak{g}$, where $A$ is a quantum superalgebra.
 This can be taken $A\subseteq L(\mathcal{H})$, where $\mathcal{H}$ is a super-Hilbert space.
 (See also Refs.\cite{PRAS14}.)
\begin{table}[h]
\caption{Supersymmetric semi-simple tensor extension Poincar\'e algebra in $D=4$.}
\label{supersymmetric-semi-simple-extension-poincare-algebra-d-4}
\begin{tabular}{|l|}
\hline
{\footnotesize$ [J_{\alpha\beta},J_{\gamma\delta}]=\eta
_{\beta\gamma} J_{\alpha\delta}+\eta
_{\alpha\delta}J_{\beta\gamma}-\eta
_{\alpha\gamma}J_{\beta\delta}-\eta_{\beta\delta}J_{\alpha\gamma},\quad [P_\alpha,P_\beta]=cZ_{\alpha\beta}$}\\
{\footnotesize$ [J_{\alpha\beta},P_\gamma]=\eta_{\beta\gamma}P_\alpha-\eta_{\alpha\gamma}P_\beta,\quad
[J_{\alpha\beta},Z_{\gamma\delta}]=\eta_{\alpha\delta}Z_{\beta\gamma}+\eta_{\beta\gamma}Z_{\alpha\delta}-
\eta_{\alpha\gamma}Z_{\beta\delta}-\eta_{\beta\delta}Z_{\alpha\gamma}$}\\
{\footnotesize$ [Z_{\alpha\beta},P_\gamma]=\frac{4a^2}{c}(\eta_{\beta\gamma}P_\alpha-\eta_{\alpha\gamma}P_\beta),\quad
[Z_{\alpha\beta},Z_{\gamma\delta}]=\frac{4a^2}{c}(\eta_{\alpha\delta}Z_{\beta\gamma}+\eta_{\beta\gamma}Z_{\alpha\delta}-
\eta_{\alpha\gamma}Z_{\beta\delta}-\eta_{\beta\delta}Z_{\alpha\gamma})$}\\
{\footnotesize$[J_{\alpha\beta},Q_\gamma]=-(\sigma_{\alpha\beta}Q)_\gamma,\quad
[P_\alpha,Q_\gamma]=a(\gamma_\alpha Q)_\gamma,\quad
[Z_{\alpha\beta},Q_\gamma]=-\frac{4a^2}{c}(\sigma_{\alpha\beta}Q_\gamma)$}\\
{\footnotesize$
[Q_\alpha,Q_\beta]=-b[\frac{2a}{c}(\gamma^\delta C)_{\alpha\beta}P_\delta
+(\sigma^{\gamma\delta}C)_{\alpha\beta}Z_{\gamma\delta}]$}\\
\hline
\end{tabular}
\end{table}
In Tab. \ref{supersymmetric-semi-simple-extension-poincare-algebra-d-4} are reported supersymmetric semi-simple tensor extensions of Poincar\'e algebra in $D=4$ too.
There $a$, $b$ and $c$ are constants. This algebra admits the following splitting:
$\mathfrak{s}\mathfrak{o}(3,1)\oplus \mathfrak{o}\mathfrak{s}\mathfrak{p}(1,4)$, where
$\mathfrak{s}\mathfrak{o}(3,1)$ is the $4$-dimensional Lorentz algebra and $\mathfrak{o}\mathfrak{s}\mathfrak{p}(1,4)$
is the orthosymplectic algebra. Then, by considering the quantum superextension
 $A\otimes_{\mathbb{R}}[\mathfrak{s}\mathfrak{o}(3,1)\oplus \mathfrak{o}\mathfrak{s}\mathfrak{p}(1,4)]$,
 where $A$ is a quantum superalgebra, and with respect to the splitting $\hat\mu={}_{\circledR}\hat\mu+{}_{\copyright}\hat\mu+{}_{\maltese}\hat\mu$ of the fundamental field $\hat\mu$, we get:
\begin{equation}\label{splitted-fundamental-quantum-field-example}
    \left\{
\begin{array}{l}
{}_{\circledR}\hat\mu=P_\alpha\hat\theta^\alpha_\gamma dx^\gamma\\
{}_{\copyright}\hat\mu=J_{\alpha\beta}\hat\omega^{\alpha\beta}_\gamma dx^\gamma\\
{}_{\maltese}\hat\mu=[\overline{Z}_{\alpha\beta}\hat A_\gamma^{\alpha\beta}+Q_{\alpha i}\phi^{\alpha i}_\gamma ]dx^\gamma.\\
\end{array}
    \right.
\end{equation}

The dynamic equation are resumed in Tab. \ref{local-expression-yang-mills-and-bianchi-identity}.

\begin{table}[h]
\caption{Local expression of
$\widehat{(YM)}\subset J\hat D^2(W)$ and Bianchi identity $(B)\subset J\hat D^2(W)$.}
\label{local-expression-yang-mills-and-bianchi-identity}
\begin{tabular}{|l|c|}
\hline
{\footnotesize\rm(Field equations)\hskip 1cm $E^A_{K}\equiv-(\partial_{B}.\hat R^{BA}_K)+[\widehat
C^H_{KR}\hat\mu^R_{C},\hat R^{[AC]}_H]_+=0$}&{\footnotesize\rm$\widehat{(YM)}$}\\
\hline \hline
\hfil{\footnotesize\rm $\hat R^K_{A_1A_2}=\left[(\partial X_{
A_1}.\hat\mu^K_{A_2})+\frac{1}{2}\widehat{C}{}^K_{IJ}[\hat\mu^I_{A_1},\hat\mu^J_{A_2}]_+\right]$}\hfil&{\footnotesize\rm(Fields)}\\
\hline
 {\footnotesize\rm(Bianchi identities)\hskip 0.25cm $B^K_{HA_1A_2}\equiv(\partial X_{H}.\hat
R^K_{A_1A_2})+\frac{1}{2} \widehat{C}{}^K_{IJ}[\bar\mu^I_{H},\hat R^J_{A_1A_2}]_+=0$}&{\footnotesize\rm$(B)$}\\
\hline
\multicolumn {2}{l}{\footnotesize\rm$\hat R^K_{A_1A_2}:\Omega_1\subset J\hat D(W)\to\mathop{\widehat{A}}\limits^2;\quad
B^K_{HA_1A_2}:\Omega_2\subset J\hat D^2(W)\to\mathop{\widehat{A}}\limits^3;\quad E^A_{K}:\Omega_2\subset
J\hat D^2(W)\to\mathop{\widehat{A}}\limits^{3}.$}\\ \end{tabular}
\end{table}

We call {\em quantum graviton} a quantum metric $\widehat{g}$
obtained by a solution $\hat\mu$ of $\widehat{(YM)}$, via the
corresponding quantum vierbein.

Since $H_3(M;\mathbb{K})=0$, we get

\begin{equation}\label{triviality-integral-bordism-group}
\left\{
\begin{array}{ll}
 \Omega_{3|3,w}^{\widehat{(YM)}}&\cong \Omega_{3|3,s}^{\widehat{(YM)}}\\
 &\cong A_0\bigotimes_{\mathbb{K}}H_3(W;\mathbb{K})\bigoplus A_1\bigotimes_{\mathbb{K}}H_3(W;\mathbb{K})\\
 &\cong A_0\bigotimes_{\mathbb{K}}H_3(M;\mathbb{K})\bigoplus A_1\bigotimes_{\mathbb{K}}H_3(M;\mathbb{K})=0\\                                                                                                  \end{array}
 \right.
\end{equation}

So $\widehat{(YM)}$ is a quantum extended crystal super PDE. However, in general, $\widehat{(YM)}$ is not a quantum $0$-crystal super PDE. In fat one has the following short exact sequence

$$\xymatrix{0\ar[r]&\ker(j)\ar[r]&\Omega_{3|3}^{\widehat{(YM)}}\ar[r]^{j}&\Omega_{3|3,s}^{\widehat{(YM)}}\ar[r]&0 }$$
hence $\Omega_{3|3}^{\widehat{(YM)}}\cong \ker(j)\not=0$. Note that $\ker(j)$ is made by $[N]\in\Omega_{3|3}^{\widehat{(YM)}}$ such that $N=\partial V$, where $V$ is some $(4|4)$-dimensional quantum supermanifolfd identified with a submanifold of $\hat J^2_{4|4}(W)$. However,
if we consider admissible only integral boundary manifolds, with
orientable classic limit, and with zero characteristic quantum
supernumbers, ({\em full admissibility hypothesis}), one has:
$\Omega_{3|3}^{\widehat{(YM)}}=0$, and $\widehat{(YM)}$ becomes a
quantum $0$-crystal super PDE. Hence we get existence of global
$Q^\infty_w$ solutions for any boundary condition of class
$Q^\infty_w$.

Then we get the exact
commutative diagram (\ref{Yang-Mills-Reinhart-bordism-groups-relation}). (For notation see \cite{PRAS21, PRAS22}.)

\begin{equation}\label{Yang-Mills-Reinhart-bordism-groups-relation}
\xymatrix{0\ar[r]&K^{\widehat{(YM)}}_{3|3;2}\ar[r]&\Omega_{3|3}^{\widehat{(YM)}}\ar[r]&
\mathop{\Omega}\limits_c{}_{6}^{\widehat{(YM)}}\ar[d]\ar[r]\ar[dr]&0&\\
&0\ar[r]&K^\uparrow_{6}\ar[r]&\Omega^\uparrow_{6}\ar[r]&
\Omega_{6}\ar[r]& 0\\}
\end{equation}
Taking into account the result by Thom on the unoriented cobordism
groups \cite{THOM}, we can calculate
$\Omega_6\cong\mathbb{Z}_2\bigoplus\mathbb{Z}_2\bigoplus\mathbb{Z}_2$.
Then, we can represent $\Omega_6$ as a subgroup of a $3$-dimensional
crystallographic group type $[G(3)]$. In fact, we can consider the
amalgamated subgroup $D_2\times\mathbb{Z}_2\star_{D_2}D_4$, and
monomorphism $\Omega_6\to D_2\times\mathbb{Z}_2\star_{D_2}D_4$,
given by $(a,b,c)\mapsto(a,b,b,c)$. Alternatively we can consider
also $\Omega_6\to D_4\star_{D_2}D_4$.  (See Appendix C in
\cite{PRAS20} for amalgamated subgroups of $[G(3)]$.) In any case the
crystallographic dimension of $\widehat{(YM)}$ is $3$ and the
crystallographic space group type are $D_{2d}$ or $D_{4h}$ belonging
to the tetragonal syngony.

$\bullet$\hskip 2pt If the initial and final Cauchy data are non-closed compact bounded $3|3$-dimensional quantum supermanifolds, we can consider the exact commutative diagram reported in (\ref{exact-commutative-diagram-bigraded-bar-quantum-chain-complex-pde}). There $\bar B_{\bullet|\bullet} (\hat E_k;A)=\ker(\partial|_{\bar
C_{\bullet|\bullet} (\hat E_k;A)})$, $\bar Z_{\bullet|\bullet}
(\hat E_k;A)=\IM(\partial|_{\bar C_{\bullet|\bullet} (\hat
E_k;A)})$, $\bar H_{\bullet|\bullet} (\hat E_k;A)=\bar
Z_{\bullet|\bullet} (\hat E_k;A)/\bar B_{\bullet|\bullet} (\hat
E_k;A)$. Furthermore,
$$\left\{\begin{array}{l}
b\in[a]\in \bar Bor_{\bullet|\bullet}(\hat E_k;A)\Rightarrow
a-b=\partial c,\quad
                 c\in \bar C_{\bullet|\bullet}(\hat E_k;A),\\
b\in[a]\in \bar Cyc_{\bullet|\bullet}(\hat E_k;A)\Rightarrow \partial(a-b)=0,\\
b\in[a]\in {}^A{\underline{\Omega}}_{\bullet|\bullet,s}(\hat
E_k)\Rightarrow
                 \left\{\begin{array}{l}
                          \partial a=\partial b=0\\
                          a-b=\partial c,\quad c\in \bar C_{\bullet|\bullet}(\hat E_k;A)\\
                          \end{array}
                                \right\}.\\
\end{array}\right.$$
It follows the canonical isomorphism:
${}^A{\underline{\Omega}}_{\bullet|\bullet,s}(\hat E_k)\cong \bar
H_{\bullet|\bullet}(\hat E_k;A)$. As $\bar C_{\bullet|\bullet}(\hat
E_k;A)$ is a free two-sided projective $A$-module, one has the
unnatural isomorphism: $$\bar Bor_{\bullet|\bullet}(\hat E_k;A)\cong
{}^A{\underline{\Omega}}_{\bullet|\bullet,s}(\hat E_k)\bigoplus \bar
Cyc_{\bullet|\bullet}(\hat E_k;A).$$

Then $a\equiv N_0\bigcup P$ and $b\equiv N_1$ belong to the same bordism group. More precisely one has
$$b\in[a]\in \bar Bor_{3|3}^{\hat E_k},\ \Leftrightarrow a-b=\partial c,\, c\equiv V,$$
such that, $\partial V=N_0\bigcup P\bigcup N_1$.
In general one has the following relation with the (closed) bordism group $\Omega_{3|3}^{\hat E_k}$:
$$\bar Bor_{3|3}^{\hat E_k}\cong \Omega_{3|3}^{\hat E_k}\bigoplus \bar Cyc_{3|3}^{\hat E_k}$$
where the cyclism group $\bar Cyc_{3|3}^{\hat E_k}$, is defined by the condition $b\in[a]\in \bar Cyc_{3|3}^{\hat E_k}$ iff $\partial (a-b)=0$. Since, under the full admissibility hypothesis, one has $\Omega_{3|3}^{\hat E_k}=0$, we get $\bar Bor_{3|3}^{\hat E_k}\cong \bar Cyc_{3|3}^{\hat E_k}$.  Therefore, we can say that inequivalent quantum particles in $\hat E_k$, are identified with $(3|3)$-dimensional bounded compact quantum supermanifolds, Cauchy data in $\hat E_k$, that are in correspondence one-to-one with  bordism classes in $\bar Bor_{3|3}^{\hat E_k}\cong \bar Cyc_{3|3}^{\hat E_k}$. Therefore, we get that the group of such fundamental particles are identified with
$$\bar C_{3|3}(\hat E_k;A)/\bar B_{3|3}(\hat E_k;A)\cong \bar C_{3|3}(\hat E_k;A)/\bar Z_{3|3}(\hat E_k;A).$$
\begin{equation}\label{exact-commutative-diagram-bigraded-bar-quantum-chain-complex-pde}
\scalebox{0.8}{$\xymatrix{&&0\ar[d]&0\ar[d]&&\\
&0\ar[r]&\bar B_{\bullet|\bullet} (\hat E_k;A)\ar[d]\ar[r]&
          \bar Z_{\bullet|\bullet} (\hat E_k;A)\ar[d]\ar[r]&
          \bar H_{\bullet|\bullet} (\hat E_k;A)\ar[r]&0\\
&&\bar C_{\bullet|\bullet} (\hat E_k;A)\ar[d]\ar@{=}[r]&
          \bar C_{\bullet|\bullet} (\hat E_k;A)\ar[d]&&\\
0\ar[r]&{}^A{\underline{\Omega}}_{\bullet|\bullet,s}(\hat
E_k)\ar[r]&\bar Bor_{\bullet|\bullet} (\hat
E_k;A)\ar[d]\ar[r]&\bar Cyc_{\bullet|\bullet} (\hat E_k;A)\ar[d]\ar[r]&0&\\
&&0&0&&}$}\end{equation}
$\bullet$\hskip 2pt Above results mean that all $(3|3)$-dimensional closed integral quantum supermanifolds
$X\subset \hat E_k$ are boundary of some $4|4$-dimensional integral quantum supermanifold $V\subset \hat E_k$, i.e., $X=\partial V$. In other words, fixing two Cauchy data $N_0$ and $N_1$ in $\hat E_k$, we can consider a quantum non-linear propagator $V$ between them, identified by means of a $(3|3)$-dimennsional closed compact integral quantum supermanifold $X\subset \hat E_k$, such that $X=N_0\bigcup P\bigcup N_1$. Therefore a quantum non-linear propagator between $N_0$ and $N_1$ can be identified with a connected $(3|3)$-dimennsional closed compact integral quantum supermanifold $X\subset \hat E_k$
homeomorphic to $N_0\bigcup P\bigcup N_1$. Since $\Omega_{3|3}=0$, it follows that we can represent such quantum non-linear propagators, with a $(3|3)$-dimensional quantum supersphere $\hat S^{3|3}$. Taking into account results contained in \cite{PRAS27}, we can also distinguish between them $(3|3)$-dimensional quantum exotic superspheres $\hat \Sigma^{3|3}$. These last are classified with respect to the equivalence relation induced by quantum diffeomorphisms.

\begin{figure}[h]
\centerline{\includegraphics[width=5cm]{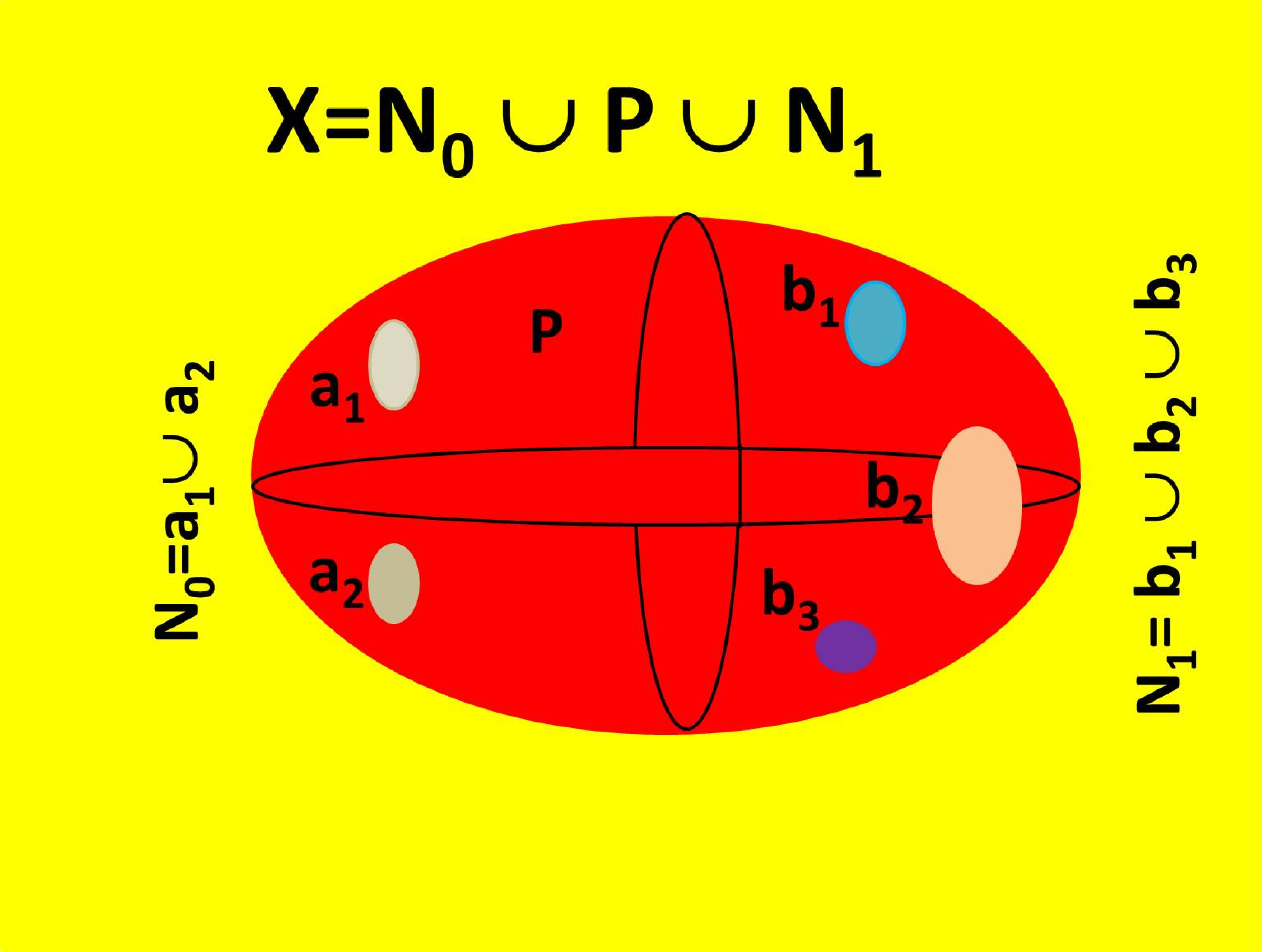}}
\renewcommand{\figurename}{Fig.}
\caption{Representation of an elementary quantum non-linear propagator $V$, identified by means of a $(3|3)$-dimensional quantum supersphere $X=N_0\bigcup P\bigcup N_1$, such that $\partial V=X$.}
\label{figure-quantum-nonlinear-propagator-quantum-supersphere}
\end{figure}
In Fig. \ref{figure-quantum-nonlinear-propagator-quantum-supersphere} is represented the boundary $X=\partial V$, of an elementary quantum non-linear propagator $V$, represented as a $(3|3)$-dimensional quantum supersphere $\hat S^{3|3}$. Therefore, such quantum non-linear propagators are homeomorphic to a $(4|4)$-dimensional quantum superdisk $\hat D^{4|4}$.
Fixing a $(3|3)$-dimensional quantum (exotic) supersphere $X\thickapprox\hat \Sigma^{3|3}\subset \hat E_k$, we can embed there $(3|3)$-dimensional quantum superdisks $a_i$, $i=1,\cdots,p$ and $b_j$, $j=1,\cdots,q$, and consider the $(4|4)$-dimennsional quantum integral submanifold $V\subset \hat E_k$, identified by a fixed boundary $X$. Then we can consider $V$ the quantum non-linear propagator between $N_0\equiv \bigcup_{1\le i\le p}a_i$ and $N_1\equiv \bigcup_{1\le j\le q}b_j$, such that $\partial V=N_0\bigcup P\bigcup N_1$, with $P\equiv \overline{X\setminus(N_0\bigcup N_1)}$.

$\bullet$\hskip 2pt It is important to remark that above quantum non-linear propagators, $V$, cannot be considered regular solutions in the equation $\widehat{(YM)}\subset J\hat{\it D}^2(W)$, but are necessarily singular solutions there. However, by considering the natural embedding $J\hat{\it D}^2(W)\hookrightarrow \hat J^2_{4|4}(W)$ we can consider $V$ as regular solution of  $\widehat{(YM)}\subset\hat J^2_{4|4}(W)$. But, whether $N_0$ is not diffeomorphic to $N_1$, we cannot say that $V$ is diffeomorphic to $N_0\times\hat D^{1|1}$. Hence, the integrable full-quantum vector field $\zeta:V\to\widehat{T}V$, relating $N_0$ to $N_1$, must necessarily be singular one.

$\bullet$\hskip 2pt For any quantum conservation law $\alpha$ of $\widehat{(YM)}$ we have $<\alpha,X>=0$, if $X=\partial V$, where $V$ is a nonlinear quantum propagator between $N_0$ and $N_1$. Therefore we get the relation (\ref{condition-quantum-supernumbers-betwenn-cuachy-data-related-by-quantum-non-linear-propagators}) between quantum numbers induced by the quantum conservation laws:\footnote{
In the particular case that $P=\varnothing$, i.e., $\partial N_0=\partial N_1=\varnothing$, then $<\alpha,N_0>=<\alpha,N_1>$.}

\begin{equation}\label{condition-quantum-supernumbers-betwenn-cuachy-data-related-by-quantum-non-linear-propagators}
<\alpha,N_0>-<\alpha,N_1>=-<\alpha,P>\in B.
\end{equation}
Therefore, identifying nonlinear quantum propagators, for $N_0$ and $N_1$, with $(4|4)$-dimensional quantum integral superdisks $\hat D^{4|4}\subset \hat E_k$, we get that the relation (\ref{condition-quantum-supernumbers-betwenn-cuachy-data-related-by-quantum-non-linear-propagators}) between the corresponding quantum conservation supernumbers, gives $$P\thickapprox\overline{\hat S^{4|4}\setminus(N_0\bigcup N_1)}.$$
\end{proof}

Similar considerations hold for the observed quantum non-linear propagators. (See also \cite{PRAS21,PRAS22}.)
\begin{theorem}
The observed dynamic equation $\widehat{(YM)}[i]$, by means of a quantum relativistic frame, is a quantum extended
crystal super PDE. Moreover, under the full admissibility hypothesis, it becomes a quantum $0$-crystal super PDE.
\end{theorem}
\begin{proof}
The evaluation of $\widehat{(YM)}$ on a macroscopic shell
$i(M_C)\subset M$ is given by the equations reported in Tab. \ref{local-expression-observed-yang-mill-pde-bianchi-identity}.

\begin{table}[h]
\caption{Local expression of $\widehat{(YM)}[i]\subset J\hat D^2(i^*W)$
and Bianchi idenity $(B)[i]\subset J\hat D^2(i^*W)$.}
\label{local-expression-observed-yang-mill-pde-bianchi-identity}
\begin{tabular}{|l|c|}
\hline
{\footnotesize\rm(Observed Field Equations)\hskip 0.5cm$(\partial_{\alpha}.\tilde R^{K\alpha\beta})+[\widehat
C^K_{IJ}\tilde\mu^I_{\alpha},\tilde R^{J\alpha\beta}]_+
=0$}&{\footnotesize $\widehat{(YM)}[i]$} \\
\hline
\hfil{\footnotesize\rm $\bar
R^K_{\alpha_1\alpha_2}=(\partial\xi_{[\alpha_1}.\tilde\mu^K_{\alpha_2]})+
\frac{1}{2}\widehat{C}{}^K_{IJ}\tilde\mu^I_{[\alpha_2}\tilde\mu^J_{\alpha_1]}$}\hfil&{\footnotesize\rm(Observed Fields)}\\
\hline
{\footnotesize\rm(Observed Bianchi Identities)\hskip 2pt$(\partial
\xi_{[\gamma}.\tilde R^K_{\alpha_1\alpha_2]})+\frac{1}{2}
\widehat{C}{}^K_{IJ}\tilde\mu^I_{[\gamma}\tilde
R^J_{\alpha_1\alpha_2]}=0$}&{\footnotesize$(B)[i]$}\\
\hline
\end{tabular}\end{table}

This equation is also formally integrable
and completely integrable. Furthermore, the
$3$-dimensional integral bordism group of $\widehat{(YM)}[i]$ and
its infinity prolongation $\widehat{(YM)}[i]_ +\infty$ are trivial,
under the full admissibility hypothesis:
$$\Omega_3^{\widehat{(YM)}[i]}\cong\Omega_3^{\widehat{(YM)}[i]_
+\infty}\cong 0.$$ So equation $\widehat{(YM)}[i]\subset J\hat
D^2(i^*W)$ becomes a quantum $0$-crystal super PDE and it admits global
(smooth) solutions for any fixed time-like $3$-dimensional (smooth)
boundary conditions. Whether $N_0\not\cong N_1$, has as a consequence that, considering the analog boundary value problem in the observed PDE $\widehat{(YM)}[i]\subset \hat J^2_3{W}$, cannot exist an observed smooth solution $V$, representing an observed quantum non-linear propagator between $N_0$ and $N_1$. The corresponding observed nonlinear quantum propagator must necessarily be singular one. \end{proof}

\section{\bf Some open problems in high energy physics solved}\label{sec-some-open-problems-in-high-energy-physics-solved}

In this section we will use the geometric mathematical architecture of the algebraic topology of quantum super PDEs previously developed to answer to some important open problems in High Energy Physics.

$\bullet$\hskip 2pt Let us, first, consider the following question: {\em ``Does exist a linear dependence between square mass and spin in a solution of the quantum PDE $\widehat{(YM)}$ ?"}.\footnote{This question arises from some semiclassical approaches that suggest to consider the following relation $J\simeq \alpha'\, M^2$, between angular momentum $J$ and mass $M$, for a spherical rotating object. This is just the philosophy of so-called {\em Regge trajectories}. Such linear trajectories are used to interpret scattering amplitudes strong reactions and are usually seen as gluonic strings attached to quarks at the end points. However there are also more exotic points of view where instead are considered pion excitations of light hadrons. (See, e.g., \cite{DIAKONOV-PETROV}.)} We shall see that in general such simple dependence for quantum solutions is not assured, but there is a very more complex relation between quantum mass and quantum spin.

Let us recall our some previous results about quantum torsion.

The quantum vierbein curvature ${}_{\circledR}\hat R$ identifies, by
means of the quantum vierbein $\hat\theta$ a quantum field $\hat
S:M\to Hom_Z(\dot\Lambda^2_0M;TM)$, that we call {\em quantum
torsion}, associated to $\hat\mu$. In quantum coordinates one can
write
\begin{equation}\label{quantum-torsion}
    \hat S=\partial x_C\otimes\hat S^C_{AB}dx^A\triangle dx^B,\quad \hat
    S^C_{AB}=\hat\theta^C_K{}_{\circledR}\hat R^K_{AB}.
\end{equation}

Furthermore, with respect to a quantum relativistic frame $i:N\to
M$, the quantum torsion $\hat S$ identifies a $A$-valued
$(1,2)$-tensor field on $N$, $\widetilde{S}\equiv i^*\hat S:N\to
A\otimes_{\mathbb{R}}\Lambda^0_2N\otimes_{\mathbb{R}}TN$, that we
call {\em quantum torsion of the observed solution}.

We say that an observed solution has a {\em quantum
spin}, if the observed solution has an observed torsion
\begin{equation}\label{observed-quantum-torsion}
\widetilde{S}\equiv i^*\hat S=\partial
x_\gamma\otimes\sum_{0\le\alpha<\beta\le
3}\widetilde{S}^\gamma_{\alpha\beta}dx^\alpha\wedge dx^\beta:N\to
A\otimes_{\mathbb{R}}\mathbf{N}\otimes_{\mathbb{R}}\Lambda^0_2(N)\cong
      A\otimes_{\mathbb{R}}\Lambda^0_2(N)\otimes_{\mathbb{R}}TN
\end{equation}
with
$\widetilde{S}^\gamma_{\alpha\beta}(p)=-\widetilde{S}^\gamma_{\beta\alpha}(p)\in
A$, $p\in N$, that satisfies the following conditions,
{\em(quantum-spin-conditions)}:
\begin{equation}\label{quantum-spin-conditions}
\begin{array}{l}
 \left\{
    \begin{array}{l}
      \widetilde{S}=\widetilde{s}\otimes\dot\psi\\
      \widetilde{s}=\sum_{0\le\alpha<\beta\le
3}\widetilde{s}_{\alpha\beta}dx^\alpha\wedge dx^\beta:N\to
A\otimes_{\mathbb{R}}\Lambda^0_2N,\\
\hskip 1cm\widetilde{s}_{\alpha\beta}(p)=-\widetilde{s}_{\beta\alpha}(p)\in A, \hskip 2pt p\in N,\\
\dot\psi\rfloor \widetilde{S} =0\\
    \end{array}
    \right\}\\
    \Downarrow\\
    \left\{
    \begin{array}{l}
      \widetilde{S}^\lambda_{\alpha\beta}=\widetilde{s}_{\alpha\beta}\dot\psi^\lambda\\
      \widetilde{S}^\lambda_{\alpha\beta}\dot\psi^\alpha =0\\
    \end{array}
    \right\}.\\
 \end{array}
\end{equation}
where $\dot\psi$ is the velocity field on $N$ of the time-like foliation representing the quantum relativistic frame on $N$. When conditions (\ref{quantum-spin-conditions})
are satisfied, we say that the solution considered  admits a {\em quantum spin-structure}, with respect to the quantum relativistic frame. We call $\widetilde{s}$ the {\em quantum $2$-form spin} of the observed solution. Let $\{\xi^\alpha\}_{0\le\alpha\le 3}$ be coordinates on $N$, adapted to the quantum relativistic frame. Then one has the following local representations:
  \begin{equation}\label{local-representations-observed-quantum-spin-observed-quantum-torsion}
 \left\{
    \begin{array}{l}
      \widetilde{s}=\widetilde{s}_{ij}d\xi^i\wedge d\xi^j\\
    \end{array}
    \right\}.
\end{equation}
We define {\em quantum spin-vector-field} of the observed solution
\begin{equation}\label{quantum-spin-vector-field}
\widetilde{\underline{s}}=<\epsilon,\widetilde{S}>=[\epsilon_{\mu\nu\lambda\rho}\dot\psi^\mu\widetilde{s}^{\nu\lambda}]d\xi^\rho
\equiv \widetilde{s}_\rho d\xi^\rho =\widetilde{s}_kdx\xi^k\hskip
2pt\Rightarrow
\widetilde{s}=\partial\xi_i\widetilde{s}_kg^{ki}=\partial\xi_i\widetilde{s}^i
\end{equation}
where
$\epsilon_{\mu\nu\lambda\rho}=\sqrt{|g|}\delta_{\mu\nu\lambda\rho}^{0123}$
is the completely antisymmetric tensor density on $N$. One has
$\widetilde{s}^\rho(p)\in A$, $p\in N$. The classification of the
observed solution on the ground of the spectrum of
$|\widetilde{s}|^2\equiv\widetilde{s}^\rho\widetilde{s}_\rho$, and
its ({\em quantum helicity}), i.e., component $\widetilde{s}_z$, is
reported in Tab. \ref{local-quantum-spectral-spin-classification-observed-yang-mills-pde-solutions}.

\begin{table}[h]
\caption{Local quantum spectral-spin-classification of $\widehat{(YM)}[i]$ solutions.}
\label{local-quantum-spectral-spin-classification-observed-yang-mills-pde-solutions}
\scalebox{0.9}{$\begin{tabular}{|l|l|}
\hline
{\footnotesize{\rm Definition}} &{\footnotesize{\rm Name}} \\
\hline
\hline
{\footnotesize $Sp(|\widetilde{s}(p)|^2)\subset\mathfrak{b}\equiv\{\hbar^2s(s+1)|
s\in\mathbb{N}\equiv\{0,1,2,\dots\}\},\quad Sp(\widetilde{s}_z(p))\subset
 \mathfrak{c}$}&{\footnotesize {\rm bosonic-polarized}}\\
 \hline
{\footnotesize $Sp(|\widetilde{s}(p)|^2)\subset\mathfrak{f}\equiv\{\hbar^2s(s+1)|
s=\frac{2n+1}{2},n\in\mathbb{N}\equiv\{0,1,2,\dots\}\}, \quad
Sp(\widetilde{s}_z(p))\subset
 \mathfrak{c}$}&{\footnotesize{\rm  fermionic-polarized}}\\
 \hline
{\footnotesize $ Sp(|\widetilde{s}(p)|^2)\cap\mathfrak{b}=Sp(|\widetilde{s}(p)|^2)\cap\mathfrak{f}=\varnothing$}&{\footnotesize {\rm  unpolarized}}\\
\hline
{\footnotesize $Sp(|\widetilde{s}(p)|^2)\cap\mathfrak{b}\not=\varnothing, \hbox{\rm   and/or } Sp(|\widetilde{s}(p)|^2)\cap\mathfrak{f}\not=\varnothing$}&{\footnotesize{\rm  mixt-polarized}}\\
\hline
\multicolumn {2}{l}{\rm
\footnotesize $|\widetilde{s}(p)|^2\equiv\widetilde{s}^\rho(p)\widetilde{s}_\rho(p)\in\widehat{A},\hskip
2pt p\in N$. $\mathfrak{c}\equiv\{\hbar m_s|m_s=-s,-s+1,\cdots,s-1,s\}$. $\widetilde{s}_z$ quantum helicity.}\\
\multicolumn {2}{l}{\rm\footnotesize  $s$= spin quantum
number; $ m_s$= spin orientation quantum number.}\\
\end{tabular}$}\end{table}

$\widehat{(YM)}$ is a functional stable quantum super PDE since it is completely integrable and formally integrable.
(See Theorem 2.34 in \cite{PRAS22}). For the same reason it admits
$\widehat{(YM)}_{+\infty}\subset J\hat{\it D}^\infty(W)$ like stable quantum extended crystal PDE. Furthermore,
since its symbol $\hat g_2$ is not trivial, any global solution $V\subset\widehat{(YM)}$ can be unstable, and the
corresponding observed solution, can appear unstable in finite times. However, global smooth solution, result stable
in finite times in $\widehat{(YM)}_{+\infty}$. Finally the asymptotic stability study of global solutions of
$\widehat{(YM)}$, with respect to a quantum relativistic frame, can be performed by means of
Theorem 2.46 in \cite{PRAS22}, since, for any section $s:M\to W$, on the fibers of
$\hat E[s]\to M$ there exists a non-degenerate scalar product.
In fact, $\hat E[s]\cong W$, as $W$ is a vector bundle over $M$. Furthermore, for any section $s$, we can identify on
$M$ a non degenerate metric $\widehat{g}$,\footnote{It should be more precise to denote $\widehat{g}$ with the
symbol $\widehat{g}[s]$, since it is identified by means of the section $s$.} that beside the rigid
metric $\underline{g}$ on $\mathfrak{g}$, identifies a non-degenerate metric on each fiber $\hat E[s]_p\cong W_p=Hom_Z(T_pM;\mathfrak{g})$, $\forall p\in M$. In fact we get $\hat\xi(p)\cdot\hat\xi(p)'=\underline{g}_{KH}\widehat{g}^{AB}(p)\xi^K_A(p)\otimes\xi'{}^H_B(p)\in A$.

Solutions of $\widehat{(YM)}$ that
encode nuclear-charged plasmas, or nuclides, dynamics. These are
described by solutions that, when observed by means of a quantum
relativistic frame have at any $t\in T$, i.e., frame-proper time,
compact sectional support $B_t\subset N$. The {\em global mass} at
the time $t$, i.e. the evaluation
$$m_t=\int_{B_t}m(t,\xi^k)\sqrt{\det(g_{ij})}d\xi^1\wedge
d\xi^2\wedge d\xi^3$$
of such mass on the space-like section $B_t$,
gives the global mass-contents of the nuclear-plasmas or nuclides,
in their ground-eigen-states, at the proper time $t$. Whether such
solutions are asymptotically stable, then they interpret the meaning
of  stable nuclear-plasmas or nuclides.

\begin{lemma}
If a solution admits spin structure, then we get $|\widetilde{S}|^2=|\widetilde{s}|^2$. Furthermore in a frame-adapted coordinates we can write
\begin{equation}\label{torsion-spin-vierbein-fundamental field}
\left\{\begin{array}{ll}
    |\widetilde{s}|^2&=\widetilde{s}_{\alpha\beta}\, \widetilde{s}^{\alpha\beta}\, \dot\psi^\gamma\, \dot\psi_\gamma      \\
         &=\widetilde{s}_{\alpha\beta}\, \widetilde{s}^{\alpha\beta} \\
         &=\widetilde{\mu}^\gamma_K\, {}_{\circledR}\widetilde{R}^K_{\alpha\beta}\, {}_{\circledR}\widetilde{R}_{\bar K}^{\alpha\beta}\, \widetilde{\mu}^{\bar K}_{\gamma}=|\widetilde{S}|^2.\\
       \end{array}
\right.\end{equation}
\end{lemma}
\begin{proof}
Let us first recall that the quantum torsion $\widehat{S}:M\to Hom_Z(\dot\Lambda^2_0M;TM)$ and its dual $\underline{\widehat{S}}:M\to Hom_Z(TM;\dot\Lambda^2_0M)$ are defined, for any $p\in M$, by the compositions reported in (\ref{torsion-spin-vierbein-fundamental-field-a}).

\begin{equation}\label{torsion-spin-vierbein-fundamental-field-a}
\xymatrix@C=1.5cm{\dot\Lambda^2_0(T_pM)\ar@/_2pc/[rr]_{\widehat{S}(p)}\ar[r]^(0.6){{}_{\circledR}\widehat{R}(p)}&
\widehat{\mathfrak{g}}\ar[r]^{\widehat{\theta}^{-1}}&T_pM\\}\hskip 0.5cm
\xymatrix@C=1.5cm{T_pM\ar@/_2pc/[rr]_{\underline{\widehat{S}}(p)}\ar[r]^{\widehat{\theta}}&
\widehat{\mathfrak{g}}\ar[r]^(0.4){{}_{\circledR}\widehat{R}(p)}&\dot\Lambda^2_0(T_pM)\\}
\end{equation}
Then, by taking the pull-back with respect to the embedding $i:N\to M$, identified by the quantum relativistic frame, we get $$|\widetilde{S}|^2=i^*(\widehat{S}\, \underline{\widehat{S}})=\widetilde{S}^\gamma_{\alpha\beta}\, \widetilde{S}_\gamma^{\alpha\beta}=\widetilde{\mu}^\gamma_K\, {}_{\circledR}\widetilde{R}^K_{\alpha\beta}\, {}_{\circledR}\widetilde{R}_{\bar K}^{\alpha\beta}\, \widetilde{\mu}^{\bar K}_{\gamma}.$$
Therefore, taking into account that in the case of solution with spin structure, one has $\widetilde{S}_\gamma^{\alpha\beta}=\dot\psi_\gamma\, \widetilde{s}^{\alpha\beta}$ and $\widetilde{S}^\gamma_{\alpha\beta}=\dot\psi^\gamma\, \widetilde{s}_{\alpha\beta}$, we get (\ref{torsion-spin-vierbein-fundamental field}).
\end{proof}

$\bullet$\hskip 2pt The following shows how thermodynamic functions can be associated to solutions of $\widehat{(YM)}$.

Let $\widetilde{H}$ be the Hamiltonian corresponding to an observed
solution of $\widehat{(YM)}$. Let us recall that $\widetilde{H}$ is a $\widehat{A}$-valued
function on the $4$-dimensional space-time $N$, considered in the quantum relativistic frame. Let us denote by $E\in
Sp(\widetilde{H})$. If $N(E)=\TR\delta(E-\widetilde{H})$ denotes the
degeneracy of $E$, let us define {\em local partition function} of
the observed solution the Laplace transform of the degeneracy
$N(E)$, with respect the spectrum $Sp(\widetilde{H})$ of
$\widetilde{H}$. We get
\begin{equation}\label{partition-function1}
    \left\{
\begin{array}{ll}
  Z(\beta)&= \int_{Sp(\widetilde{H})}e^{-\beta E}N(E) dE\\
  & = \int_{Sp(\widetilde{H})} e^{-\beta E}\TR\delta(E-\widetilde{H}) dE \\
  &=\TR e^{-\beta \widetilde{H} }.\\
\end{array}
    \right.
\end{equation}
So we get the following formula
\begin{equation}\label{partition-function2}
Z(\beta)=\TR e^{-\beta \widetilde{H} }
\end{equation}
where $\beta$ is the Laplace transform variable and it does not
necessitate to be interpreted as the ``inverse temperature", i.e.,
$\beta=\frac{1}{\kappa_B\theta}$, where $\kappa_B$ is the
Boltzmann's constant. If $\beta=\frac{1}{\kappa_B\theta}$
then the system encoded by the observed solution of
$\widehat{(YM)}$, is in equilibrium with a heat bath (canonical
system). Note that all above objects are local functions on the
space-time $N$. The same holds for $\beta$. We can interpret $Z(\beta)$ as a normalization factor
for the local probability density
\begin{equation}\label{probability-density}
P(E)=\frac{1}{Z}N(E)e^{-\beta E}
\end{equation}
that
the system, encoded by the observed solution, should assume the
local energy $E$, with degeneration $N(E)$. In fact we have:
\begin{equation}\label{partition-function-normalization-factor}
1=\int_{Sp(\widetilde{H})} P(E)
dE=\frac{1}{Z}\int_{Sp(\widetilde{H})} N(E)e^{-\beta
E}dE=\frac{Z}{Z}.
\end{equation}
As a by-product we get that the {\em local average energy}
$<E>\equiv e$ can be written, by means of the partition function, in
the following way:

\begin{equation}\label{partition-function-average-energy1}
e=-(\partial\beta\ln Z).
\end{equation}
In fact, one has
\begin{equation}\label{partition-function-average-energy2}
    \left\{
\begin{array}{ll}
  e&= \int_{Sp(\widetilde{H})}E P(E) dE=\frac{1}{Z}\int_{Sp(\widetilde{H})}EN(E)e^{-\beta E}dE\\
  & = \frac{1}{Z}\int_{Sp(\widetilde{H})} E\TR\delta(E-\widetilde{H})e^{-\beta E} dE \\
  &=\frac{1}{Z}\TR(\widetilde{H} e^{-\beta \widetilde{H}})=-\frac{1}{Z}(\partial\beta.Z)=-(\partial\beta.\ln Z).\\
\end{array}
    \right.
\end{equation}

 When we can interpret
$\beta=\frac{1}{\kappa_B\theta}$, then one can write
\begin{equation}\label{partition-function-average-energy3}
e=\kappa_B\theta^2(\partial\theta.\ln Z).
\end{equation}
Then we get also that the local {\em energy fluctuation} is
expressed by means of the variance of $e$:
\begin{equation}\label{partition-function-energy-fluctuation}
<(\triangle E)^2>\equiv <(E-e)^2>=(\partial\beta\partial\beta.\ln
Z).
\end{equation}
Furthermore, we get for the {\em local heat capacity} $C_v$ the following
formula:
\begin{equation}\label{heat-capacity-v}
C_v=(\partial\theta.e)=\frac{1}{\kappa_B\theta^2}<(\triangle E)^2>.
\end{equation}
We can define the {\em local entropy} by means of the following
formula:
\begin{equation}\label{entropy1}
s=-\kappa_B\int_{Sp(\widetilde{H})}P(E) \ln P(E) dE.
\end{equation}
In fact one can prove that one has the usual relation by means of
the energy.  Really we get:

\begin{equation}\label{entropy2}
\left\{
\begin{array}{ll}
  s&  =-\kappa_B\int_{Sp(\widetilde{H})}P(E) \ln P(E) dE\\
  & =\kappa_B(\ln Z+\beta e)=(\partial\theta.(\kappa_B\theta \ln Z)).\\
\end{array}
\right.
\end{equation}
Then from the relation $s=\kappa_B(\ln Z+\beta e)$ we get
$e=\theta s-\kappa_B\ln Z$, hence also
$(\partial s.e)=\theta$. This justifies the definition of entropy
given in (\ref{entropy1}). Furthermore, from {\em(\ref{entropy2})} we get also $\frac{1}{\beta}\ln Z=e-\theta s=f$,
 where $f$is the Helmoltz free
energy. It follows the following
expression of the local Helmoltz free energy, by means of the local
partition function $Z$:
\begin{equation}\label{free-energy-partition-function}
f\equiv e-\theta s=-\kappa_B\theta\ln Z.
\end{equation}
Conversely, from (\ref{free-energy-partition-function}) it follows that the partition function can be expressed by means
of the local Helmoltz free energy
\begin{equation}\label{partition-function-free-energy}
Z=e^{-\beta f}.
\end{equation}

So we see that the local thermodynamic functions, can be expressed as
scalar-valued differential operators on the fiber bundle
$W[i]\times_NT^0_0N\to N$.

The concept of quantum states can be also related to a proof for existence of solutions with mass-gap. In fact, we have proved that equation $\widehat{(YM)}$ admits local and global solutions with mass-gap. These are contained into a sub-equation,
{\em(Higgs-quantum super PDE)},
$\widehat{(Higgs)}\subset\widehat{(YM)}$, that is formally integrable and completely integrable, and also
 a stable quantum super PDE. If $H_3(M;\mathbb{K})=0$, $\widehat{(Higgs)}$ is also a quantum extended crystal super PDE.
 In general
 solutions contained in $\widehat{(Higgs)}$ are not stable in finite times. However there exists an associated stabilized
 quantum super PDE, (resp. quantum extended crystal super PDE), where all global smooth solutions are stable in finite
 times.
 Furthermore, there exists a quantum super partial differential relation, {\em(quantum Goldstone-boundary)},
 $\widehat{(Goldstone)}\subset\widehat{(YM)}$, bounding $\widehat{(Higgs)}$, such that any global solution of
 $\widehat{(YM)}$,
 loses/acquires mass, by crossing $\widehat{(Goldstone)}$.
\begin{example}
Now, let us answer to the question: ``{\em Do pictures in Fig. \ref{figure-proton-proton-chain-reactions-examples} represent possible smooth integral manifolds, i.e., solutions, of $\widehat{(YM)}$, or $\widehat{(YM)}[i]$, with respect to a quantum relativistic frame $i$ ?}"

Taking into account above results on $\widehat{(YM)}$ and $\widehat{(YM)}[i]$, we can answer, ``{\em yes}", accepted that the initial and final Cauchy data have the same quantum numbers, identified by the quantum conservation laws of $\widehat{(YM)}$ and $\widehat{(YM)}[i]$ respectively. Furthermore, we can state that such solutions must necessarily be singular ones, since in order that should be smooth, it is necessary that $V$ should be diffeomorphic to $N_0\times \hat D^{1|1}$ or $N_0\times I$. Then in such a case we should also have $N_0\cong N_1$. Therefore we can conclude that nonlinear quantum propagators representing reactions considered in Example \ref{two-body-high-energy-reactions-examples} cannot be smooth solutions, but singular ones.
\end{example}

\begin{lemma}
The quantum mass of a solution of $\widehat{(YM)}$ can be written, in coordinates adapted to a quantum relativistic frame, in the form reported in {\em(\ref{alternative-form-quantum-mass})}.
\begin{equation}\label{alternative-form-quantum-mass}
 m=\frac{1}{2}|\widetilde{S}|^2-\frac{1}{2}\hat\Delta+m_{\copyright}+m_{\maltese}.
\end{equation}
\end{lemma}
\begin{proof}

The quantum Hamiltonian can be written in the form reported in {\em(\ref{alternative-form-quantum-hamiltonian})}.
\begin{equation}\label{alternative-form-quantum-hamiltonian}
 H={}_{\circledR}H+{}_{\copyright}H+{}_{\maltese}H.
\end{equation}
This follows directly, by considering the expression for the quantum Hamiltonian, by using the splitting induced by the quantum superalgebra $\mathfrak{g}$, and by taking into account the condition that the solution admits a spin structure. The corresponding calculus is reported in (\ref{alternative-form-quantum-hamiltonian-a}) and following ones.
\begin{equation}\label{alternative-form-quantum-hamiltonian-a}
\left\{\begin{array}{ll}
         H&=\frac{1}{2}(\widetilde{R}_{\alpha\beta}^K\, \widetilde{R}^{\alpha\beta}_K-\widetilde{\mu}^K_{\alpha\beta}\, \widetilde{R}^{\alpha\beta}_K) \\
         &=\frac{1}{2}({}_{\circledR}\widetilde{R}^{\alpha\beta}_K\, {}_{\circledR}\widetilde{R}^{\alpha\beta}_K-{}_{\circledR}\widetilde{\mu}^K_{\alpha\beta}\, {}_{\circledR}\widetilde{R}^{\alpha\beta}_K) \\
         &+\frac{1}{2}({}_{\copyright}\widetilde{R}^{\alpha\beta}_K\, {}_{\copyright}\widetilde{R}^{\alpha\beta}_K-{}_{\copyright}\widetilde{\mu}^K_{\alpha\beta}\, {}_{\copyright}\widetilde{R}^{\alpha\beta}_K) \\
         &+\frac{1}{2}({}_{\maltese}\widetilde{R}^{\alpha\beta}_K\, {}_{\maltese}\widetilde{R}^{\alpha\beta}_K-{}_{\maltese}\widetilde{\mu}^K_{\alpha\beta}\, {}_{\maltese}\widetilde{R}^{\alpha\beta}_K) \\
         &={}_{\circledR}H+{}_{\copyright}H+{}_{\maltese}H.
       \end{array}
\right.
\end{equation}
In coordinates adapted to the quantum relativistic frame, we get that $\dot\psi^\lambda=\delta_0^\lambda$, and since
$${}_{\circledR}\widetilde{\mu}_{K}^\lambda\,{}_{\circledR}\widetilde{R}_{\alpha\beta}^K=
\widetilde{S}^\lambda_{\alpha\beta}=\widetilde{s}_{\alpha\beta}\, \dot\psi^\lambda$$
we get ${}_{\circledR}\widetilde{\mu}_{K}^i\,{}_{\circledR}\widetilde{R}_{\alpha\beta}^K=
\widetilde{S}^i_{\alpha\beta}=0$ and ${}_{\circledR}\widetilde{\mu}_{K}^0\,{}_{\circledR}\widetilde{R}_{\alpha\beta}^K=
\widetilde{S}^0_{\alpha\beta}=\widetilde{s}_{\alpha\beta}$. Therefore we have
\begin{equation}\label{alternative-form-quantum-hamiltonian-b}
\left\{\begin{array}{ll}
|\widetilde{S}|^2&=|\widetilde{s}|^2={}_{\circledR}\widetilde{\mu}_{K}^\lambda\,{}_{\circledR}\widetilde{R}_{\alpha\beta}^K\,
{}_{\circledR}\widetilde{R}_{\bar K}^{\alpha\beta}\,{}_{\circledR}\widetilde{\mu}_{\lambda}^{\bar K}\\
&={}_{\circledR}\widetilde{\mu}_{K}^0\,{}_{\circledR}\widetilde{\mu}_{0}^{\bar K}\,{}_{\circledR}\widetilde{R}_{\alpha\beta}^K\, {}_{\circledR}\widetilde{R}_{\bar K}^{\alpha\beta}+\hat\Delta\\
&=\delta_{K}^{\bar K}\,{}_{\circledR}\widetilde{R}_{\alpha\beta}^K\, {}_{\circledR}\widetilde{R}_{\bar K}^{\alpha\beta}+\hat\Delta\\
&={}_{\circledR}\widetilde{R}_{\alpha\beta}^K\, {}_{\circledR}\widetilde{R}_{K}^{\alpha\beta}+\hat\Delta\\
\end{array}
\right.
\end{equation}
with
\begin{equation}\label{alternative-form-quantum-hamiltonian-c}
\hat\Delta\equiv {}_{\circledR}\widetilde{\mu}_{K}^0\,[{}_{\circledR}\widetilde{R}_{\alpha\beta}^K,{}_{\circledR}\widetilde{\mu}_{0}^{\bar K}]{}_{\circledR}\widetilde{R}_{\bar K}^{\alpha\beta}+{}_{\circledR}\widetilde{\mu}^0_K\, {}_{\circledR}\widetilde{R}^{K}_{\alpha\beta}[{}_{\circledR}\widetilde{R}_{\bar K}^{\alpha\beta},{}_{\circledR}\widetilde{\mu}^{\bar K}_0].
\end{equation}

So, in general, we can write
 \begin{equation}\label{alternative-form-quantum-hamiltonian-d}
 {}_{\circledR}\widetilde{R}_{\alpha\beta}^K\, {}_{\circledR}\widetilde{R}_{K}^{\alpha\beta}=|\widetilde{S}|^2-\hat\Delta
 \end{equation}
 and by using equation (\ref{alternative-form-quantum-hamiltonian-a}) we get the expression (\ref{alternative-form-quantum-hamiltonian-e}) for the quantum Hamiltonian in terms of quantum torsion.
  \begin{equation}\label{alternative-form-quantum-hamiltonian-e}
 H=\frac{1}{2}|\widetilde{S}|^2-\frac{1}{2}\hat\Delta+{}_{\copyright}H+{}_{\maltese}H.
 \end{equation}
 As a by-product, we get that the quantum mass of an observed solution by means of a quantum relativistic frame $i:N\to M$, of $\widehat{(YM)}[i]$, the alternative representation (\ref{alternative-form-quantum-mass}) in terms of quantum torsion.
\end{proof}

\begin{definition}\label{reduced-quantum-mass}
Let us define {\em reduced quantum mass} of an observed solution by means of a quantum relativistic frame $i:N\to M$, of $\widehat{(YM)}[i]$: $$M\equiv m+\frac{1}{2}\hat\Delta-m_{\copyright}-m_{\maltese}.$$
 \end{definition}
 \begin{theorem}\label{reduced-quantum-mass-torsion}
There is the following direct relation between reduced quantum mass and quantum torsion: $$M=\frac{1}{2}|\widetilde{S}|^2.$$
\end{theorem}
\begin{definition}\label{phenomenological-quantum-spin}
Let us define {\em phenomenological quantum spin} of an observed solution by means of a quantum relativistic frame $i:N\to M$, of $\widehat{(YM)}[i]$, a $\widehat{A}$-valued quantum scalar $J$ such that:\footnote{The natural meaning of $\sqrt{J}$ is the $\widehat{A}$-valued function on $N$, such that: $\sqrt{J}\, \sqrt{J}=J$.}  $$\sqrt{J}\equiv |\widetilde{S}|^2.$$
 \end{definition}
 \begin{cor}[Quantum Regge-type trajectories]\label{phenomenological-quantum-regge-type-rule}
There is the following direct relation between reduced quantum mass and phenomenological quantum spin:
\begin{equation}\label{quantum-regge-type-relation}
M^2=\frac{1}{4}J.
\end{equation}

We call {\em quantum Regge-type trajectories} the relation {\em(\ref{quantum-regge-type-relation})} between quantum square mass $M^2$ and the phenomenological quantum spin $J$.
 \end{cor}

\begin{figure}[h]
\centerline{\includegraphics[width=10cm]{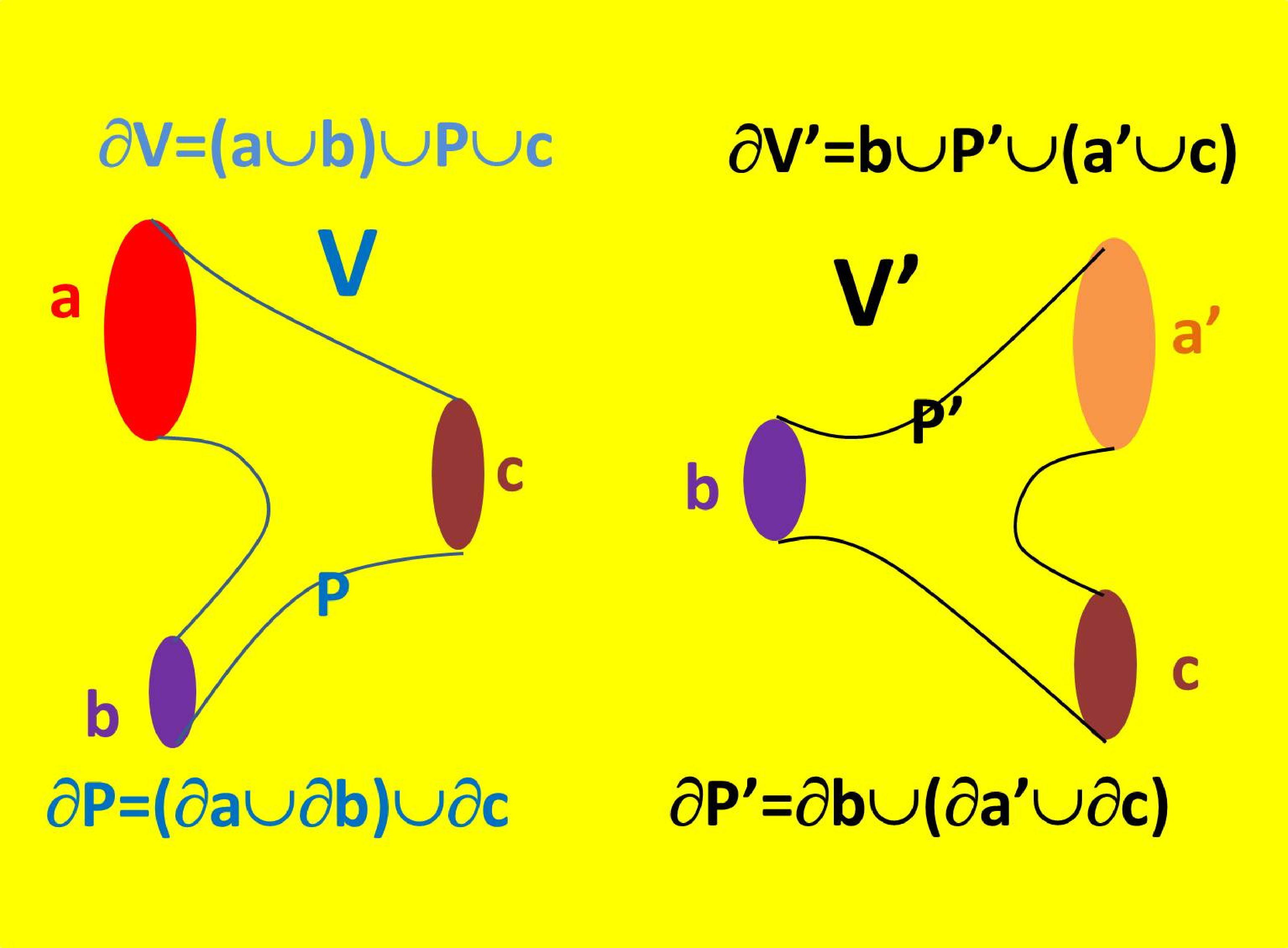}}
\renewcommand{\figurename}{Fig.}
\caption{Representation of a nonlinear quantum propagator $V$, for a reaction $a+b\to c$, and its quantum nonlinear anti-propagator $V'$, corresponding to the reaction $b\to a'+c$, obtained for quantum crossing symmetry.}
\label{figure-quantum-anti-propagator}
\end{figure}

$\bullet$\hskip 2pt Another question, in a sense related to the previous problem, is the dynamical justification of the Gell-Mann-Nishijima formula. This actually has not yet found a dynamical derivation.\footnote{The more interesting effort in this direction was an heuristic semiclassical justification given in \cite{CUI}.} Really dynamical characterizations of hypercharge and $3^{th}$ isospin component, can be obtained by means of electric charge, by giving the dynamical expression of this last quantity with respect to a quantum relativistic frame. This is obtained in the following theorem.

\begin{theorem}[Quantum Gell-Mann-Nishijima formula and electric-charge-gap-solutions]\label{quantum-gell-mann-nishijima-formula}
We define {\em square-fundamental electric charge} of an observed solution $V$ of $\widehat{(YM)}[i]$, the littlest value in the spectrum $Sp(\hat w)$ of the observed quantum electromagnetic energy $\hat w$, and denote it by $q^2$.  We say that $V$ has an {\em electric-charge gap} if $q^2>0$. (In general $q^2\ge 0$.)
The {\em quantum electric charge}, $\hat Q(t)$, of an observed solution $V$ of $\widehat{(YM)}[i]$, at the proper time $t$, is expressed by the formula {\em(\ref{quantum-gell-mann-nishijima-formula-a})}.
\begin{equation}\label{quantum-gell-mann-nishijima-formula-a}
  \scalebox{0.8}{$ \hat Q^2(t)=\int_{\sigma_t}\left[{}_{\copyright}\widetilde{R}^K_{0i}{}_{\copyright}\widetilde{R}^i_{0K}
   +{}_{\maltese}\widetilde{R}^K_{0i}{}_{\maltese}\widetilde{R}^i_{0K}
   +{}_{\circledR}\widetilde{B}^K_{0i}{}_{\circledR}\widetilde{B}^i_{0K}
    +{}_{\copyright}\widetilde{B}^K_{0i}{}_{\copyright}\widetilde{B}^i_{0K}
   +{}_{\maltese}\widetilde{B}^K_{0i}{}_{\maltese}\widetilde{B}^i_{0K}\right]\otimes\eta\in A$}
\end{equation}
where $\eta=\sqrt{g_{ij}}\, d\xi^1\wedge d\xi^2\wedge d\xi^3$ is the canonical space-like volume form on $\sigma_t\subset V$. This last is the $3$-dimensional space-like sub-manifold of $V$ at the time $t$.
Furthermore, we call {\em quantum hypercharge}, $\hat Y(t)$, (resp. {\em quantum $3^{th}$-isospin component}, $\hat I_3(t)$), at the proper time $t$, the elements of the quantum (super)algebra $A$, such that holds the formula {\em(\ref{quantum-gell-mann-nishijima-formula-b})} {\em(quantum Gell-Mann-Nishijima formula)}.
\begin{equation}\label{quantum-gell-mann-nishijima-formula-b}
   \hat Q^2(t)=(\hat I_3(t)+\frac{1}{2}\hat Y(t))^2,\, \forall t \in \triangle
\end{equation}
where $\triangle$ is the definition time-set of the considered solution of $\widehat{(YM)}[i]$.

$\bullet$\hskip 2pt The spectral content of $\hat Q(t)$ is given by $Sp(\pm\sqrt{\hat Q(t)^2})$.

$\bullet$\hskip 2pt There exists an open sub-equation $\widehat{(YM)}[i]_w\subset\widehat{(YM)}[i]$, {\em quantum electromagnetic-Higgs PDE}, that is formally integrable and completely integrable, where live all solutions with electric-charge gap. The boundary $\partial\widehat{(YM)}[i]_w=\overline{\widehat{(YM)}[i]_w}\setminus\widehat{(YM)}[i]_w\subset\widehat{(YM)}[i]$, is a partial differential relation that we call {\em quantum electromagnetic-Goldstone boundary} and denote by $\widehat{(Goldstone)}[i]_w$. An electrically neutral connected, simply connected, $3$-dimensional Cauchy data $N\subset \widehat{(YM)}[i]$, cannot be contained into $\widehat{(YM)}[i]_w$. Let $V$ be a nonlinear quantum propagator, such that $\partial V=N_0\sqcup P\sqcup N_1$, with $N_0\subset \widehat{(YM)}[i]_w$ and $N_1\not\subset \widehat{(YM)}[i]_w$. Let us assume that $N_r$, $r=0,1$, are connected, simply connected particles, hence $N_0$ has an electric-charge gap. The particle $N_1$, instead cannot have electric-charge gap. Thus, $V$, by crossing $\widehat{(Goldstone)}[i]_w$, loses its electric-charge gap, passing from $N_0$ to $N_1$, and vice versa.
\end{theorem}

\begin{proof}
The solution $V$ identifies a quantum-electric-charge field $\hat E=\dot\psi\rfloor\widetilde{R}={}_{\circledR}\hat E+{}_{\copyright}\hat E+{}_{\maltese}\hat E$ and a quantum-magnetic-charge field
$\hat B=\dot\psi\rfloor(\star\widetilde{R})={}_{\circledR}\hat B+{}_{\copyright}\hat B+{}_{\maltese}\hat B$. Assuming that the solution $V$ has a quantum spin, one has ${}_{\circledR}\hat E=0$. (See \cite{PRAS22}.) Then the
quantum electro-magnetic-charge energy of the solution of $\widehat{(YM)}[i]$ is given by equation (\ref{quantum-gell-mann-nishijima-formula-quantum-electromagnetic-energy-density}).
\begin{equation}\label{quantum-gell-mann-nishijima-formula-quantum-electromagnetic-energy-density}
  \hat w=\hat w_e+\hat w_m={}_{\copyright}\widetilde{R}^K_{0i}{}_{\copyright}\widetilde{R}^i_{0K}
   +{}_{\maltese}\widetilde{R}^K_{0i}{}_{\maltese}\widetilde{R}^i_{0K}
   +{}_{\circledR}\widetilde{B}^K_{0i}{}_{\circledR}\widetilde{B}^i_{0K}
    +{}_{\copyright}\widetilde{B}^K_{0i}{}_{\copyright}\widetilde{B}^i_{0K}
   +{}_{\maltese}\widetilde{B}^K_{0i}{}_{\maltese}\widetilde{B}^i_{0k}.
\end{equation}
Therefore the quantum electromagnetic energy of the space-like set $\sigma_t$ is given by the expression on the right in (\ref{quantum-gell-mann-nishijima-formula-a}). We define {\em quantum electric charge} contained, at the proper time $t$, into a space-like set $\sigma_t\subset V$, $\hat Q(t)\in A$, such that $\hat Q^2(t)=\int_{\sigma_t}\hat w\otimes \eta$.\footnote{This is justified taking into account that the electric charge contained into a $3$-dimensional set having the electromagnetic energy $W$, is $q^2=2 C\, W$, with $C$ the capacity of this set. Therefore $Q^2=\frac{q^2}{2C}$. i.e., the formula (\ref{quantum-gell-mann-nishijima-formula-a}) is normalized with respect to a {\em factor form}  $\kappa_{electric-charge}(t)=\frac{1}{2C(t)}$.} Then equation (\ref{quantum-gell-mann-nishijima-formula-b}) gives a dynamical definition for $\hat Y$ and $\hat I_3$.

To prove the last part of the theorem it is enough to consider the continuous mapping $\hat w:\widehat{(YM)}[i]\to \widehat{A}$, defined by means of equation (\ref{quantum-gell-mann-nishijima-formula-quantum-electromagnetic-energy-density}). Set $\widehat{(YM)}[i]_w\equiv (\hat w)^{-1}(G(\widehat{A}))\subset \widehat{(YM)}[i]$. Then $\widehat{(YM)}[i]_w$ is an open quantum PDE of $\widehat{(YM)}[i]$, hence retains the same formal properties of this last equation. The rest of the proof follows the same line of the ones of Theorem 3.28 in \cite{PRAS22} about solutions with mass-gap. Let us only emphasize that a neutral particles $N_0\subset \widehat{(YM)}[i]$, with trivial topology, i.e., connected and simply connected, cannot admit $0\in Sp(\hat Q(t))$, whether $N_0\subset \widehat{(YM)}[i]_w$, since in $\widehat{(YM)}[i]_w$, one has $\lambda>0$, for any $\lambda\in Sp(\hat w)$. Furthermore, if the nonlinear quantum propagator $V$ is such $\partial V=N_0\sqcup P\sqcup N_1$, assuming that $N_1\not\subset \widehat{(YM)}[i]_w$ and that it has trivial topology, then $N_1$ cannot have electric-charge gap, hence $V$ crossing $\widehat{(Goldstone)}[i]_w$ must necessarily lose electric-charge gap.
\end{proof}
\begin{remark}\label{remark-quantum-gell-mann-nishijima-formula}
Theorem \ref{quantum-gell-mann-nishijima-formula} does not necessarily contradict the conservation of the quantum electric charge. In fact, we can have two Cauchy data $a,\, b\subset \widehat{(YM)}[i]_w$ and another one $c\not\subset \widehat{(YM)}[i]_w$, such that $<\hat q,a>=-<\hat q,b>\not=0$, where $<\hat q,a>=\sqrt{\int_a\hat f\otimes \eta}=-\sqrt{\int_b\hat f\otimes\eta}$, with $\hat f$ the function under integral in {\em(\ref{quantum-gell-mann-nishijima-formula-a})}. Furthermore, let $<\hat q,c>=0$. Then in the integral bordism classes $[a\sqcup b=N_0]$ and $[N_1]$, one has $<\hat q,[a\sqcup b]>=0=<\hat q,c>$. Therefore, for a nonlinear quantum propagator $V$, such that $\partial V=N_0\sqcup P\sqcup N_1$ of $\widehat{(YM)}[i]$ one has $<\alpha,a>+<\alpha,b>=-<\alpha,P>$, for $\alpha=\hat q$, since $<\hat q,P> =0$. This is surely down since $P$ is a time-like manifold having $2$-dimensional space-like sections $\sigma_t$, hence $<\hat q,P>=\sqrt{\int_{\sigma_t}\hat f\otimes \sqrt{g_{ij}}\, d\xi^1\wedge d\xi^2\wedge d\xi^3}=0$.
\end{remark}

From Theorem \ref{quantum-gell-mann-nishijima-formula} and Remarks \ref{remark-quantum-gell-mann-nishijima-formula} it follows the following important theorem.

\begin{theorem}[$Q$-exotic nonlinear quantum propagators of ${\widehat{(YM)}[i]}$]\label{non-conservation-of-quantum-electric-charge}
For any observed nonlinear quantum propagator $V$ of $\widehat{(YM)}[i]$, such that $\partial V=N_0\sqcup P\sqcup N_1$, where $N_i$, $i=0,1$, are $3$-dimensional space-like admissible Cauchy data of $\widehat{(YM)}[i]$, and $P$ is a suitable time-like $3$-dimensional integral manifold with $\partial P=\partial N_0\sqcup\partial N_1$, equation {\em(\ref{conserved-quantum-observed-electric-charge-in-quantum observed-nonlinear-propagator})} holds.
\begin{equation}\label{conserved-quantum-observed-electric-charge-in-quantum observed-nonlinear-propagator}
 \hat Q[i|t_0]=\hat Q[i|t_1] \, {\rm mod}\hskip 3pt \mathfrak{Q}[V]\in A
 \end{equation}
 where $\hat Q[i|t_r]\in A$ is the quantum electric charge on $N_r$, $r=0,1$, and $\mathfrak{Q}[V]\in A$, is a term that in general is not zero and that we call {\em defect quantum electric-charge}.
We call {\em $Q$-exotic nonlinear quantum propagators}, nonlinear quantum propagators such that $\mathfrak{Q}[V]\not=0\in A$.\footnote{This agrees with the conservation of the observed quantum Hamiltonian. See Theorem 3.20 in \cite{PRAS28}. In fact, the observed quantum electromagnetic energy is a form of quantum energy, even if, in general, it does not coincide with the observed quantum Hamiltonian.}
\end{theorem}
\begin{proof}
The proof follows the same strategy of one considered to prove Theorem 3.20 in part I \cite{PRAS28}. In fact the gauge invariance of equation $\widehat{(YM)}[i]$ produces a quantum integral characteristic $3$-form, in the sense of Lemma 3.19 in part I, that has a structure similar to $\omega_H$. Let us denote such a conservation law by $$\omega_q=(\omega_q)_0\otimes d\xi^1\wedge d\xi^2\wedge d\xi^3+\sum_{1\le i\le 3}(\omega_q)_i\otimes d\xi^0\wedge d\xi^1\wedge\cdots\wedge \widetilde{d\xi^i}\wedge d\xi^3,$$ with $\widetilde{d\xi^i}$, absent, $i=1,2,3$, and $(\omega_q)_\alpha(p)\in A$, for $p\in V$.
The quantum electric charge, on a space-like section $\sigma_t$ of $V$, is given by means of $\omega_q$, by the following equation
$$\hat Q[i|t]=\int_{\sigma_t}(\omega_q)_0|_{\sigma_t}\otimes d\xi^1\wedge d\xi^2\wedge d\xi^3.$$ Therefore we have
$$\hat Q[i|t_0]-\hat Q[i|t_1]=-\int_{P}(\omega_q)_P= -\mathfrak{Q}[V].$$
In general $\mathfrak{Q}[V]\not=0$, hence we get $\hat Q[i|t_0]=\hat Q[i|t_1]\, {\rm mod}\hskip 3pt  \mathfrak{Q}[V]$. We call $\mathfrak{Q}[V]\in A$ {\em defect quantum electric-charge} of the observed nonlinear quantum propagator $V$ of $\widehat{(YM)}[i]$. In other words, in general one has $\frac{d}{dt}\hat Q(t)\not=0\in A$.
\end{proof}

\begin{remark}
In the Part III we will further characterize $Q$-exotic nonlinear quantum propagators as nonlinear effects of exotic quantum supergravity.
\end{remark}

$\bullet$\hskip 2pt In the following we shall prove that the phenomenological crossing symmetry in particle reactions is well justified in the Pr\'astaro's algebraic topology theory of quantum super PDE.

\begin{figure}[t]
\centerline{\includegraphics[width=6cm]{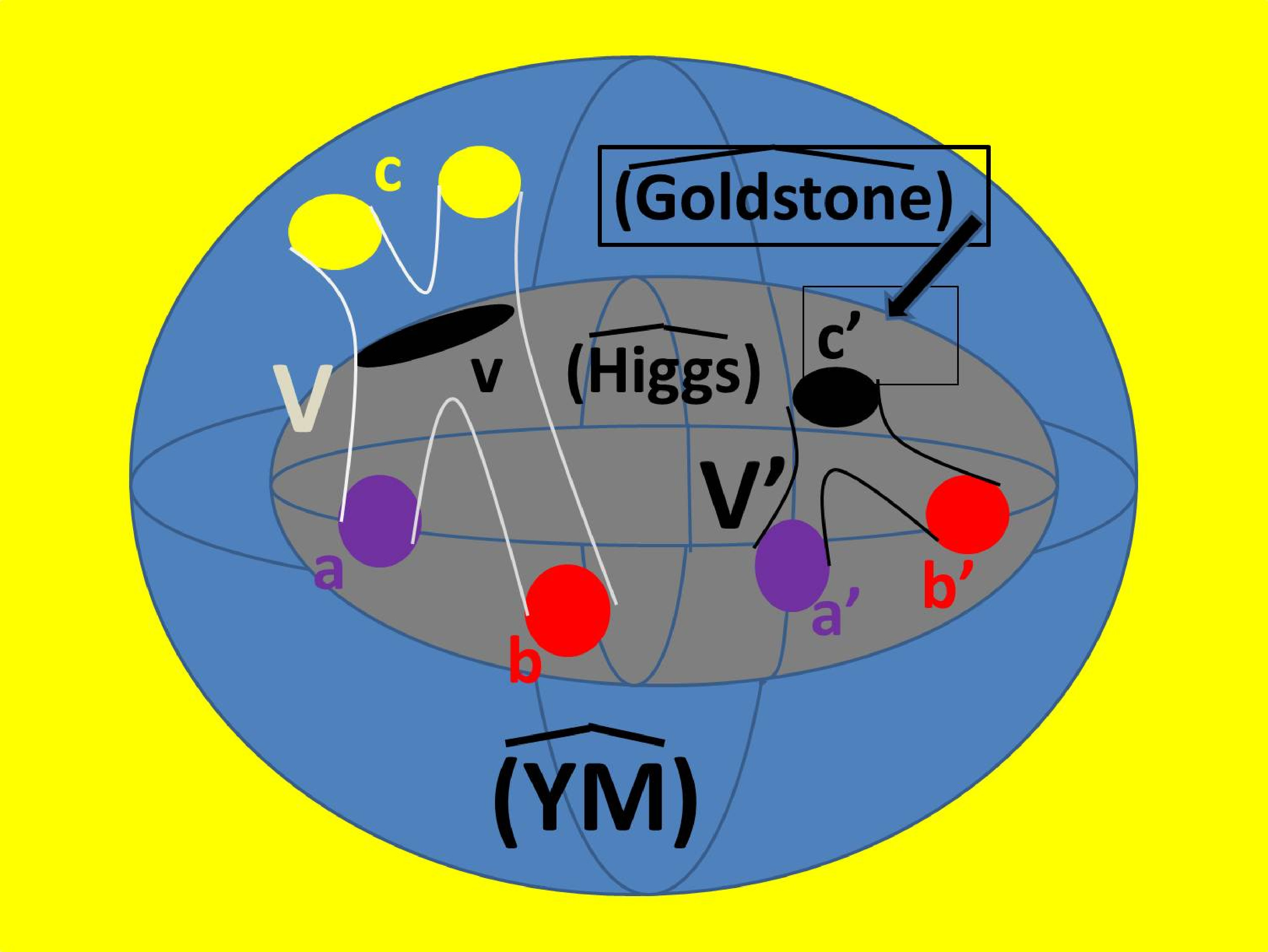}
\includegraphics[width=6.07cm]{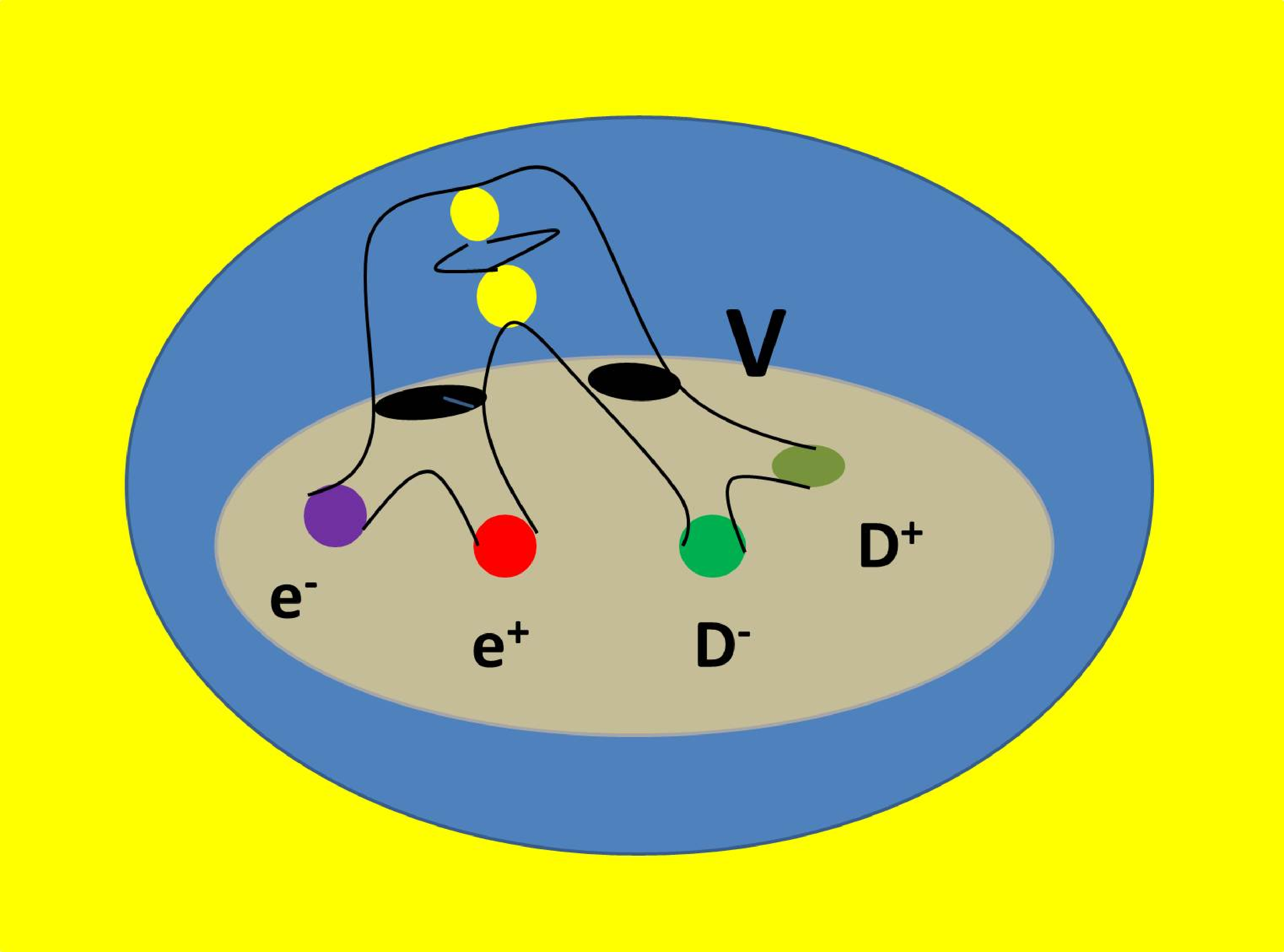}}
\renewcommand{\figurename}{Fig.}
\caption{In the figure on the left it is represented a couple of zero-mass photons ($c\nsubseteq(\widehat{Higgs})$) and massive photons ($c'\subset(\widehat{Higgs})$). The quantum non-linear propagator $V$ bording a couple of massless photons (yellow particles) with electron-positron ($a\sqcup b$) has a {\em Goldstone piece} (black quantum virtual particle). This is identified with a massive photon. In other words the annihilation of electron-positron generates a virtual massive photon before to produce a couple of massless photons. The quantum non-linear propagator $V'$ bording the massive photon (black) with electron-positron ($a'\sqcup b'$) has not a Goldstone piece, since it is completely inside $(\widehat{Higgs})\subset\widehat{(YM)}$. This is a massive photon decay. In the figure on the right, it is represented a nonlinear quantum propagator, $V$, bording massive couple of antiparticles, and having a part outside $\widehat{(Higgs)}$ (the grey region). In this case we can talk of two Goldstone pieces in $V$, i.e., there are two massive photons (black particles) related by $V$.}
\label{figure-zero-mass-and-massive-photons-production}
\end{figure}

\begin{theorem}[Quantum crossing symmetry]\label{quantum-crossing-symmetry}
If $\widehat{(YM)}$ admits a nonlinear quantum propagator $V$ such that $\partial V=N_0\bigcup P\bigcup N_1$, with $N_0=a\sqcup b$ and $N_1=c\sqcup d\sqcup e$, then there exists also a nonlinear quantum propagator $V'$ such that $\partial V'=N'_0\bigcup P'\bigcup N'_1$, with $N'_1=a'\sqcup b'$ and $N'_0=c'\sqcup d'\sqcup e'$, where the primed symbols denote antiparticles. This property is called {\em quantum crossing symmetry}. Similarly there exist nonlinear quantum propagators between $a\sqcup b\sqcup c'$ and $d\sqcup e$ or between $b$ and $a'\sqcup c\sqcup d\sqcup e$.\footnote{The crossing symmetry is a well-known property in the reaction particles phenomenology. These, besides a set of other symmetries and phenomenological conservation laws constitue the Holy Bible for particle-reaction physicists. (See Tab. \ref{figure-zero-mass-and-massive-photons-production}.)  In Tab. \ref{figure-zero-mass-and-massive-photons-production} are reported also two quantum reactions related by means of crossing symmetry. These are the {\em neutron decay} and {\em neutrino detection}.)}
\end{theorem}
\begin{proof}
This follows from the general relation between quantum integral bordism groups and quantum integral characteristic conservation laws. (See \cite{PRAS14}.) In fact if there exists a nonlinear quantum propagator $V$ such that $\partial V=N_0\bigcup P\bigcup N_1$, with $N_0=a\sqcup b$ and $N_1=c\sqcup d\sqcup e$, then this means that $N_1\in[N_0\bigcup P]\in Bor^{\widehat{(YM)}}_{3|3}$, hence $N_0\bigcup P\bigcup N_1\in [0]\in Bor^{\widehat{(YM)}}_{3|3}$. Since antiparticles reverse quantum integral characteristic numbers, it follows that for any quantum integral characteristic conservation law $\alpha$ of $\widehat{(YM)}$, one has $<[\alpha],[N_0\bigcup P\bigcup N_1]>=0=<[\alpha],[N'_1\bigcup P'\bigcup N'_0]>$. Therefore, $(N'_1\bigcup P')\in[N'_0]\in  Bor^{\widehat{(YM)}}_{3|3}$, hence there exists a nonlinear quantum propagator $V'$ such that $\partial V'=N'_0\bigcup P'\bigcup N'_1$. Note also, that as by-product, we get that there exist a $(3|3)$-dimensional quantum integral supermanifold $P'$, such that $\partial P'=\partial N'_1\sqcup \partial N'_0$, similarly to what happens for $P$, namely $\partial P=\partial N_0\sqcup N_1$. This means that one has $\partial N'_1\in[\partial N'_0]\in\Omega^{\widehat{(YM)}}_{2|2}$ and $\partial N_1\in[\partial N_0]\in\Omega^{\widehat{(YM)}}_{2|2}$. Similarly one proves existence of the other reactions obtained from crossing symmetry. For example to the nonlinear quantum propagator $V$, such that $\partial V=(a\sqcup b)\bigcup P\bigcup c$, corresponds for quantum crossing symmetry, the following  nonlinear quantum propagator $V'$, such that $\partial V'=b\bigcup P'\bigcup(a'\sqcup c)$, where $P'$ is the anti-$P$. (See Fig. \ref{figure-quantum-anti-propagator} for a representation of $P$ and its anti-$P$.) There are cases where can exist more than only one anti-$P$.

\end{proof}

$\bullet$\hskip 2pt The following theorem answers to the question: {\em ``Do exist massive photons ?"}.

\begin{theorem}[Existence of quantum massive photons]\label{quantum-massive-photons-theorem}
The quantum super Yang-Mills equation $\widehat{(YM)}$ admits solutions that starting from a Cauchy data $N_0=a\sqcup b$, where $a$ represents a quantum electron and $b$ a quantum positron, bords $N_1$, representing a quantum massive particle. Any annihilation $e^++e^-\to \gamma+\gamma$, must necessarily generate a {\em quantum virtual massive photon}, say $\gamma_{m}$, before to produce a couple of massless photons $\gamma$. Conversely a {\em quantum massive photon} decays giving $\gamma_{m}\to e^++e^-$.
\end{theorem}

\begin{proof}
Let us first identify a quantum massive photon with a quantum massive particle $\gamma_m$ with decay into an electron-positron couple. We know that from Theorem 3.28 in \cite{PRAS22} can exist a quantum propagator $V'$ bording $a'\sqcup b'$ with $c'$, where $a'\sqcup b'$ is a couple electron-positron and $c'$ is a massive particle, hence all contained in $(\widehat{Higgs})\subset \widehat{(YM)}$. On the other hand, from the standard Compton scattering $\gamma+e^-\to e^-+\gamma$ we say, for crossing symmetry, that holds the following reaction ({\em annihilation electron-positron}): $e^++e^-\to\gamma+\gamma$. Theorem 3.28 in \cite{PRAS22} assures the existence of a nonlinear quantum propagator $V$ such that $\partial V=N_0\sqcup P\sqcup N_1$, where $N_0$, is contained in the sub-equation $(\widehat{Higgs})\subset \widehat{(YM)}$, and $N_1\nsubseteq (\widehat{Higgs})$. Therefore $V$, should cross the Goldstone boundary $(\widehat{Goldstone})\subset\widehat{(YM)}$, hence $V$ should contain a Goldstone piece. This proves that the annihilation electron-positron implies to pass across a quantum massive photon:  (Note that whether it is permitted the reaction $e^++e^-\to\gamma_m$, for {\em crossing symmetry} it is also permitted the decay $\gamma_m\to e^-+e^+$.)\footnote{Actually there are attempts to give experimental evidence to the existence of massive photons. (Visit \href{https://www.jlab.org/}{(Thomas Jefferson National Accelerator Facility, Newport News, Virginia, USA)}, and \cite{ARKANI-FINK-SLAT-WEI}.) Massive photons can be identified with neutral massive bosonic particles. In the particle-zoo these could be identified with so-called {\em vector bosons}, like neutral $\rho$-meson and $\varphi$-meson and $Z^0$-boson, all having spin $1$.} In Fig. \ref{figure-zero-mass-and-massive-photons-production} are represented  productions of quantum massless photons and quantum massive photons with respect to the constraints  $(\widehat{Higgs})$ and $(\widehat{Goldstone})$ in $\widehat{(YM)}$. Whether $e^+$ and $e^-$ collide at high energy, their annihilation can be reversed into a massive couple of mesons $D^+\sqcup D^-$. In other words there exists a nonlinear quantum propagator $V$, bording $e^+\sqcup e^-$ with $D^+\sqcup D^-$, such that $\partial V=(e^+\sqcup e^-)\bigcup P\bigcup(D^+\sqcup D^-)$, where $P$ is an integral quantum supermanifold, partially outside $\widehat{(Higgs)}$, but yet contained into $\widehat{(YM)}$. (See Fig. \ref{figure-zero-mass-and-massive-photons-production}.)
\end{proof}

From above considerations, we see that the concept of massive photons can be generalized to more large set of particles.

\begin{definition}[Quantum anti-particles and quantum massive photons]\label{quantum-antiparticles-massive-photons}
We say that two quantum, massive, electric charged particles $a\sqcup a'\subset \widehat{(YM)}$ are {\em quantum anti-particles}, if they have the same mass and opposite quantum numbers.

We call {\em $a$-quantum massive photon} a quantum massive uncharged particle, having a decay into a couple $(a,a')$ of quantum massive, electric charged anti-particles. We denote such a quantum massive photon with the symbol $\gamma^a_m$, or simply $\gamma_m$, if no confusion can arise. Therefore, we can write $\gamma^a_m\to a+a'$.
\end{definition}

\begin{theorem}[Existence of $a$-quantum massive photons]\label{a-quantum-massive-photons}
For any couple $(a,a')$ of quantum massive, electric charged anti-particles, there exists a $a$-quantum massive photon $\gamma^a_m$.
\end{theorem}
\begin{proof}
The proof can be obtained by analogy of the one for Theorem \ref{quantum-massive-photons-theorem}.
\end{proof}
\begin{example}[Existence of $p$-quantum massive photon]\label{existence-of-p-quantum-massive-photon}
Since it is permitted the Compton scattering $p+n\gamma\to p+n\gamma$, for suitable $n\in\mathbb{N}$, it is also permitted for quantum crossing symmetry the following one $p+\bar p\to 2n\gamma$. Therefore, a virtual $p$-quantum massive photon $\gamma_m^p$ necessarily exists.\footnote{From some experimental results it is well-known that the annihilation $p+\bar p$ produces photons $\gamma$, through intermediate reactions and decays coming from massive particles. (For example: $p+\bar p\to3\pi^0$, $p+\bar p\to2\pi^0+\eta$, $p+\bar p\to\pi^0+2\eta$, $\pi^0\to 2\gamma$, $\eta\to 2\gamma$, $Z^0\to e^++e^-\to 2\gamma$.) Therefore, the virtual $p$-quantum massive photon $\gamma_m^p$ has a very complex structure, made by a collection of massive particles, like mesons $\pi$ and $\eta$, bosons like $Z^0$ and leptons $e^{\pm}$, all inside $\widehat{(Higgs)}$. Note that in the virtual massive photon  $\gamma_m^p$ enters also the vector boson $Z^0$ that is more massive than proton. ($m_p=0.938\, GeV/c^2$, $m_Z=91\, GeV/c^2$ and $m_W=80\, GeV/c^2$.) This means that virtual massive photons can have very large masses. (See Remark \ref{quantum-virtual-anomaly-massive-particles}.) Massive intermediate vector bosons $W^{\pm}$ and $Z^0$, were predicted from Steven Weinberg, Abdus Salam and Sheldon Glashow in 1979 (Nobel award 1979) and experimental discovered by Carlo Rubbia in 1983 (Nobel award 1984).}
\end{example}

\begin{definition}[Quantum massive neutrinos]\label{quantum--massive-neutrinos}
We say that a quantum, massive, quasi-particle with zero electric charge is a {\em quantum massive neutrino}, if it admits a decay into a couple (neutrino,anti-neutrino) of the same type.
\end{definition}

\begin{theorem}[Existence of quantum massive neutrinos]\label{existence-of-quantum-massive-neutrino}
The quantum super Yang-Mills equation $\widehat{(YM)}$ admits solutions that represent decays of quantum massive neutrinos.
\end{theorem}
\begin{proof}
In fact, we can find a nonlinear quantum propagator $V\subset\widehat{(YM)}$, such that $\partial V=N_0\bigcup P\bigcup N_1$, where $N_0=a\sqcup b\subset \widehat{(YM)}$, with $(a,b)=(\nu_e,\bar\nu_e)\not\subset\widehat{(Higgs)}$, and $N_1=c\sqcup d\subset \widehat{(Higgs)}$,  with $(c,d)$ a couple of antiparticles in the sense of Definition \ref{quantum-antiparticles-massive-photons}. Then, necessarily $V$ must cross $\widehat{(Goldstone)}$, hence identifies a massive quasi-particle, say $\nu_m\subset \widehat{(Higgs)}$. Then, applying the quantum crossing symmetry to the reaction $\nu_e+\bar\nu_e\to \nu_m\to c+d$, we get also the reaction $d+c\to\nu_m\to \bar\nu_e+\nu_e$, i.e., there exists a nonlinear quantum propagator of $\widehat{(YM)}$ that encodes such a reaction. This proves that the massive, neutral particle $\nu_m$, decays into the couple $(\nu_e,\bar\nu_e)$, therefore $\nu_m$ is a massive neutrino.

Let us emphasize that neutrinos are of three different type: $\nu_e$, $\nu_\mu$ and $\nu_\tau$, called $e$-neutrino, $\mu$-neutrino and $\tau$-neutrino respectively. They differ for the type of decay that produce them. For example see (\ref{different-neutrinos-type-and-decays}).

 \begin{equation}\label{different-neutrinos-type-and-decays}
\scalebox{0.8}{$ \begin{array}{ccc}
\xymatrix{\pi^+\ar[r]&e^++\nu_e\\}&&\\
\xymatrix{&&\mu^+\\
\pi^+\ar[r]&W^+\ar[ru]\ar[rd]&\\
&&\nu_\mu\\}&
\xymatrix{&\nu_\mu&\bar\nu_e\\
\mu^-\ar[ru]\ar[r]&W^-\ar[ru]\ar[rd]&\\
&&e^-\\}&
\xymatrix{&\bar\nu_\tau&\bar\nu_e&\bar\nu_\mu&\bar u\\
\tau^-\ar[ru]\ar[r]&W^-\ar[ru]\ar[rd]\ar@{.>}[rru]\ar@{.>}[rrd]\ar@{-->}[rrru]\ar@{-->}[rrrd]&\\
&&e^-&\mu^-&d\\}\\
e^++e^-\to Z^0\to\nu_e+\bar\nu_e&&\\
\end{array}$}
\end{equation}

In \ref{different-neutrinos-type-and-decays}) the dot-arrows and dash-arrows starting from $W^-$ are alternative to the full-arrows and between them.
In all these cases there is an intermediate production of a quasi particle $Z^0$ or $W^{\pm}$. Such quasi-particles can be usually produced when a massive solution cross $\widehat{(Goldstone)}$. Therefore from this point of view we should have a further argument to consider the usual neutrinos massless particles, according to the Gell-Mann's standard model. However, the quasi-particles $W^{\pm}$ cannot be considered massive neutrinos, since they do not decay into a massless couple (neutrino,anti-neutrino). Instead $Z^0$ has all the properties to be considered a massive neutrino !

Let us also emphasize that any massless particle produced in some reaction from massive particles are encoded by some nonlinear quantum propagators that necessarily must cross the Goldstone boundary to go outside $\widehat{(Higgs)}$, hence they identify massive particles. For example this happens for neutrinos and also for boson gauge particles that should be massless, but in interactions with massive particles appear to have big masses.\footnote{This phenomenon is usually called {\em gauge symmetry breaking}. Inside the Pr\'astaro's Algebraic Topology of quantum super PDEs one can understand that such particles acquire masses by means of geometrodynamic processes related to the structure of quantum Yang-Mills PDEs and nonlinear quantum propagators encoding reactions where they are involved.}
\end{proof}
\begin{example}[What dark matter is ?]\label{massive-virtual-particles-dark-matter}
Theorem \ref{existence-of-quantum-massive-neutrino} agrees with the experimental fact that the annihilation of the couple $(\nu_e,\bar\nu_e)$ produces a $Z^0$ massive boson. This last can again decay in the couple $(\nu_e,\bar\nu_e)$, or in a couple lepton-antilepton, or quark-antiquark, according to energy level considered. For example:
\begin{equation}\label{annihilations-neutrino-antineutrino}
\scalebox{0.8}{$\xymatrix{&&\nu_e+\bar\nu_e\\
\nu_e+\bar\nu_e\ar[r]&Z^0\ar[r]\ar@{.>}[rd]\ar@{-->}[ru]&e^++e^-\\
&&q+\bar q\\}$}
\end{equation}
In {\em(\ref{annihilations-neutrino-antineutrino})} dot-arrow and dash-arrow starting from $Z^0$ are alternative to the full-arrow and each other one.
Then the $Z$ boson can be considered a massive neutrino when decay into the couple $(\nu_e,\bar\nu_e)$. The quantum crossed symmetry reaction $\nu_e+\bar\nu_e\to Z^0$, is an example showing how massless particles can cross $\widehat{(Glodstone)}$, producing massive particles inside $\widehat{(Higgs)}$. The reaction, $e^++e^-\to Z^0\to \nu_e+\bar\nu_e$, obtained, by quantum cross symmetry, from the reaction in the middle line of {\em(\ref{annihilations-neutrino-antineutrino})}, shows that the annihilation of the couple $(e^-,e^+)$, when it is enough energetic, produces a massive quasi-particle $Z^0$, before to destroy its mass, crossing  $\widehat{(Glodstone)}$, and producing the couple of massless particles $(\nu_e,\bar\nu_e)$. This clarifies that the energetic level of an electric-charged, massive couple (particle,antiparticle) decides whether in a scattering will produce a massive photons or a massive neutrino.\footnote{A guess by B. M. Pontecorvo and V. N. Gribov \cite{GRIBOV-PONTECORVO, PONTECORVO}, in order to explain the so-called {\em mystery of the missing solar neutrinos}, assumed that neutrinos have some oscillations that characterize their different types. But such oscillations require massive neutrinos at least for $\tau$-neutrinos and $\mu$-neutrinos. This fact does not agree with the standard model. On the other hand, from Theorem \ref{existence-of-quantum-massive-neutrino} we can understand that massive neutrinos can be identified with some energetic levels of the quasi-particle $Z^0$, and that non-electron neutrinos can retain their usual property to be massless particles. It is important to emphasize that recent experimental results, made by an international collaboration called \href{http://www.nature.com/nature/journal/v512/n7515/full/nature13702.html}{Borexino-Collaboration} \cite{BOREXINO-COLLABORATION}, have been able to detect just the so-called missing solar neutrinos.  With this respect one can understand that nonlinear quantum propagators allow us to reconcile massless neutrinos, according to Standard Model, with the observed massive neutrinos.}\label{dark-matter}

Let us also emphasize that the so-called ``dark matter" can be interpreted as massive virtual particles codified in Theorem \ref{quantum-massive-photons-theorem}, Theorem \ref{a-quantum-massive-photons} and Theorem \ref{existence-of-quantum-massive-neutrino}. In other words ``dark matter" can be considered a generic term to identify virtual massive particles produced when a solution of the quantum super Yang-Mills equation $\widehat{(YM)}$, crosses the Goldstone boundary, coming inside the Higgs quantum super PDE, $\widehat{(Higgs)}\subset \widehat{(YM)}$ contained into $\widehat{(YM)}$. This interpretation for ``dark matter" could theoretically support some recent experimental observations where positron fraction is stably increasing at high energy levels, that are just necessary to produce massive photons, massive neutrinos and $a$-quantum massive photons, according to the above quoted theorems. See \cite{AGUILAR} for recent experiments that should suggest existence of ``dark matter" in the sense here specified.\footnote{Production of virtual $\pi^+$ decaying into positrons and neutrinos, can be obtained by means massive photons identified with virtual states into $\gamma$-nucleons interactions. This is a new point of view looking to long-studied $\pi^+$-photoproduction. (For general information on pions photoproduction see, e.g., \cite{DRECHSEL-TIATOR, DUTTA, DUTTA-GAO-LEE} and references quoted therein.)}\label{dark-matter}
\end{example}

$\bullet$\hskip 2pt The following theorem answers to the question: ``{\em Do quarks are fundamental particles ?}"

\begin{theorem}[Quarks break-down]\label{quarks-break-down}
Quarks are not fundamental particles.
\end{theorem}
\begin{proof}
The proof can be obtained easily by considering Theorem \ref{quantum-massive-photons-theorem}. In fact for any quark $q$, we can associate an anti-quark $\bar q$, so that $(q,\bar q)$ is a couple of massive, electric-charged anti-particles in the sense of Definition \ref{quantum-antiparticles-massive-photons}, hence there exist massive $q$-photons $\gamma_q$. In other words there exist nonlinear quantum propagators allowing reactions $q+\bar q\to n\, \gamma$, $n\in\mathbb{N}$. In these reaction quarks are transformed in particles that do not contain quarks, hence quarks are broken-down ! To similar conclusions one arrives by considering Theorem \ref{existence-of-quantum-massive-neutrino}. In fact, by considering the quantum reaction $\bar q+q\to Z^0\to \bar\nu_e+\nu_e$, obtained by means of quantum cross symmetry of the reaction in the bottom in (\ref{annihilations-neutrino-antineutrino}), we see that the mass of the couple $(q,\bar q)$ completely desappears and appear particles that are not made by quarks. In other words quarks break-down.
\end{proof}

$\bullet$\hskip 2pt In the following we give a precise meaning to the concept of {\em observed quantum time} and to the so-called {\em time-energy uncertainty principle}. This is possible thanks to our geometric theory of quantum (super) PDE's.\footnote{The exact interpretation of the time-energy uncertainty principle was remained an open problem, and not even universally accepted, (see, e.g., L. Landau). Really in order to directly apply the Heisenberg uncertainty principle to the (time,energy) couple, it is necessary to define the meaning of the time as a noncommutative variable, with respect to the energy. But the proper time of a quantum relativistic frame is a commutative variable ! (See \cite{PRAS19}.) (See also some attempts to solve this problem by Dirac \cite{DIRAC} and L. I. Mandelshtam and I. E. Tamm \cite{MANDELS-TAMM}.)  Actually the time-energy uncertainty principle is justified by means adopting the  Mandelshtam-Tamm's point of view. But that approach is not satisfactory, since it refers to a generic auxiliary non-commutative variable, say $B$, that has nothing to do with energy and time. Furthermore, in that description time remains a commutative variable, and one identifies the velocity $\frac{dB}{dt}$ with a finite ratio, i.e., $\frac{dB}{dt}\equiv \frac{\Delta B}{\Delta t}$. In fact one uses the following modified uncertainty relation: $\Delta B\, \frac{\Delta E}{\frac{dB}{dt}}\ge \hbar$. Then by using the approximation $\frac{dB}{dt}\equiv \frac{\Delta B}{\Delta t}$, one obtains $\Delta B\, \frac{\Delta E}{\frac{dB}{dt}}\approx \Delta t\, \Delta E\ge \hbar$. But $\Delta t\, \Delta E\ge \hbar$ is a no sense, since $t$ is a commutative variable !}
\begin{definition}[Quantum relativistic observed time]\label{quantum-relativistic-observed-time}
Let $V\subset\widehat{(YM)}$ be a solution of $\widehat{(YM)}$ and let $\widehat{g}:M\to Hom_Z(\dot T^2_0M;A)$ be the quantum metric induced by the quantum graviton, identified by the solution $V$. Then a quantum relativistic frame $\psi\equiv\{i:N\to M, g\}$, where $(N,g)$ is a pseudo-Riemannian $4$-dimensional manifold, endowed with a time-like flow, identifies on $N$ a $A$-valued metric $i^*\widehat{g}:N\to A\bigotimes S^0_2N$. We define {\em quantum relativistic observed time}, between two points $p_1,\, p_2\in N$, belonging to a same flow line of $\psi$, the $A$-valued quantum length $\hat t[p_1,p_2]\in A$, calculated  with respect to $i^*\widehat{g}$. (For more details see \cite{PRAS19}.) We call {\em spectral content} of the quantum relativistic observed time the spectrum $Sp(\hat t[p_1,p_2])$ of $\hat t[p_1,p_2]$.
\end{definition}

\begin{proposition}
The quantum relativistic observed time cannot, in general, coincide with the proper time of the quantum relativistic frame.
\end{proposition}
\begin{proof}
In fact the length of the arc considered in above proposition differs whether measured with respect to $g$ or to $i^*\widehat{g}$, even if one has the case $A=\mathbb{R}$ ! This can be seen by a direct calculus, by using coordinates $\{\xi^\alpha\}_{0\le\alpha\le 3}$, with $\xi^0=t$ the proper time of the relativistic quantum frame $\psi$. Then the arc length between the points $p_1,\, p_2$ belonging to the same flow line $\phi_{p_1}$, passing for $p_1$ and $p_2$, is given by $s[p_1,p_2]=\int_{[p_1,p_2]}ds=\int_{[t_1,t_2]}\sqrt{g_{\alpha\beta}\dot\xi^\alpha\dot\xi^\beta}\, dt=\int_{[t_1,t_2]}\sqrt{g_{00}}\, dt=t_2-t_1$, assuming $g_{00}=1$. On the other hand we have
$\hat t[p_1,p_2]=\int_{[p_1,p_2]}d\hat s=\int_{[t_1,t_2]}\sqrt{|\gamma^*\widehat{g}_{(s)}}\, dt\in A$, where $\gamma\equiv\phi_{p_1}\circ i$. (For more details see \cite{PRAS19}.) In the particular case that $A=\mathbb{R}$, we have that $\hat t[p_1,p_2]\in \mathbb{R}$, as well as $ s[p_1,p_2]$, but $\gamma^*\widehat{g}_{(s)}$ is different from $\phi_{p_1}^*g$.\footnote{Let us also add that observed solutions, namely solutions of the observed quantum super Yang-Mills PDEs, $\widehat{(YM)}[i]$, encoding reactions, are singular solutions. Therefore, internal times separating two space-like Cauchy data $N_{t_0}$ and $N_{t_1}$, $t_0<t_1$, do not equal the intervals $t_1-t_0$. Therefore, also from this point of view the apparent time of a reaction, (with respect to a quantum relativistic frame), does not give an exact valuation of the internal laps of time in which reaction occurred.}
\end{proof}

\begin{theorem}[Quantum virtual anomaly-massive particles]\label{quantum-virtual-anomaly-massive-particles}
A quantum virtual massive quasi-particle can have an observed mass, in the reaction where it appears, larger than the observed energy of incoming particles entering into the reaction.\footnote{See, e.g., Example \ref{existence-of-p-quantum-massive-photon}.} However this does not really contradict the Einstein's conservation mass-energy equation.
\end{theorem}

\begin{proof}
The justification of this fact, that should appear contradict the Einstein's equation of mass-energy conservation $E=mc^2$, is partially to ascribe to the {\em Heisenberg's uncertainty principle}. (This is, really, not more a ``principle", but a ``theorem". See, \cite{PRAS5}.) In fact, adopting notation used in \cite{PRAS5}, the uncertainty relation for the commutator $[quantum-time,quantum-energy]$, can be written as reported in equation (\ref{apparent-time}).
\begin{equation}\label{apparent-time}
\sigma^2(\hat H)\, \sigma^2(\hat t)\ge\frac{1}{4}|<[\hat H,\hat t]>|^2.
    \end{equation}
Thus in quantum processes, localized in very short ``observed quantum-time", as could happen in the case of virtual particles creation, one could have large masses, independently from the initial energy-mass content. However, this argument does not contradict the Einstein's conservation mass-energy equation ! In fact, the concept of ``quantum-time" is only an apparent time, different from the proper time of the quantum relativistic frame, and existing with respect to the interaction quantum-system-quantum-relativistic-frame. Instead the proper time of the observer, belongs to $\mathbb{R}$, hence commutes with the hamiltonian. In other words with respect to this proper-time, the Heisenberg's uncertainty principle does not apply. However a complete precise mathematical justification of such phenomena, can be obtained in the framework of nonlinear quantum propagators. Really only these structures can support the real interpretation of energy conservation in quantum reactions. (See Theorem 3.20(I), Corollary 3.22(I) and Corollary 3.23(I).)
\end{proof}

Let take this occasion to revise the concept of Heisenberg's uncertainty principle in the framework of Algebraic Topology of quantum super PDEs.

\begin{definition}[Apparent quantum observed time in quantum reaction]\label{apparent-quantum-observed-time-in-quantum-reaction}
We call (observed) {\em apparent quantum time} in a quantum reaction, $\sigma^2(\hat t)$ entering in equation {\em(\ref{apparent-time})}.
\end{definition}
\begin{example}[Quantum time and strengths of the fundamental interactions]\label{quantum-time-and-strengths-fundamental-interactions}
An experimental way to use the concept of quantum time and apparent quantum time, is to estimate the type of interaction present in. In fact, by using the Heisenberg uncertainty relation $\triangle\hat t\, \triangle E\ge \hbar$, applied to the quantum time and the strengths of the interactions, it is possible to guess the type of interaction from the observed interaction time. In Tab. \ref{table-quantum-time-and-strengths-fundamental-interactions} such a correlation is reported. There are also reported Planck's constants and quarks masses for any related convenience.

$\bullet$\hskip 2pt For example, if one observes a decay into an apparent quantum time $\sigma^2(\hat t)\approx  10^{-23}\, s$, one can state that it is a strong decay (or quantum-gravity decay) there.

\begin{table}[h]
\caption{Quantum time and strengths of the fundamental interactions.}
\label{table-quantum-time-and-strengths-fundamental-interactions}
\begin{tabular}{|c|c|}
\hline
{\footnotesize{\rm Apparent interaction time}} &{\footnotesize{\rm Interaction type}}\hfil\\
\hline
\hline
{\footnotesize {\rm $10^{-10}\, s$}}&{\footnotesize {\rm Weak}}\\
\hline
{\footnotesize {\rm $10^{-16}\, s$}}&{\footnotesize {\rm Electromagnetic}}\\
\hline
{\footnotesize {\rm $10^{-23}\, s$}}&{\footnotesize {\rm Strong (quantum-gravity)}}\\
\hline
\multicolumn{2}{l}{\footnotesize {\rm $\hbar=6.584\, \times\, 10^{-25}\, GeV\, s$.}}\\
\end{tabular}
\scalebox{0.65}{$\begin{tabular}{|l|l|c|c|c|c|c|c|c|c|}
\hline
\multicolumn{10}{|c|}{\footnotesize {\rm Comparison with Planck's constants, electron, proton and quark mass-energies}}\\
\hline
\hline
{\footnotesize {\rm Planck length}}&{\footnotesize {\rm $l_{pl}=\sqrt{\frac{\hbar\, G}{c^3}}=1.616 \times\, 10^{-35}\, m$}}&\hskip-0.5cm&{\footnotesize {\rm Quarks}}&{\footnotesize {\rm Mass-energy}}&{\footnotesize {\rm charge}}&{\footnotesize {\rm S'}}&{\footnotesize {\rm C'}}&{\footnotesize {\rm B'}}&{\footnotesize {\rm T'}}\\
\hline
{\footnotesize {\rm Planck mass}}&{\footnotesize {\rm $m_{pl}=\sqrt{\frac{\hbar\, c}{G}}=2.176 \times\, 10^{-8}\, Kg$}}&\hskip-0.5cm&{\footnotesize {\rm up (u)}}&{\footnotesize {\rm $0.3$}}&{\footnotesize {\rm $+\frac{2}{3}$}}&{\footnotesize {\rm $0$}}&{\footnotesize {\rm $0$}}&{\footnotesize {\rm $0$}}&{\footnotesize {\rm $0$}}\\
\hline
{\footnotesize {\rm Planck energy}}&{\footnotesize {\rm $E_{pl}=m_{pl}\, c^2=c^2\, \sqrt{\frac{\hbar\, c}{G}}=1.22 \times\, 10^{19}\, GeV$}}&\hskip-0.5cm&{\footnotesize {\rm down (d)}}&{\footnotesize {\rm $0.3$}}&{\footnotesize {\rm $-\frac{1}{3}$}}&{\footnotesize {\rm $0$}}&{\footnotesize {\rm $0$}}&{\footnotesize {\rm $0$}}&{\footnotesize {\rm $0$}}\\
\hline
{\footnotesize {\rm Planck time}}&{\footnotesize {\rm $t_{pl}=\frac{l_{pl}}{c}=\frac{\hbar}{m_{pl}\, c^2}=\sqrt{\frac{\hbar\, G}{c^5}}=5.391 \times\, 10^{-44}\, s$}}&\hskip-0.5cm&{\footnotesize {\rm charm (c)}}&{\footnotesize {\rm $1.25$}}&{\footnotesize {\rm $+\frac{2}{3}$}}&{\footnotesize {\rm $0$}}&{\footnotesize {\rm $+1$}}&{\footnotesize {\rm $0$}}&{\footnotesize {\rm $0$}}\\
\hline
{\footnotesize {\rm Planck charge}}&{\footnotesize {\rm $q_{pl}=\sqrt{4\pi\epsilon_0\hbar\, c}=1.875 \times\, 10^{-18}\, C$}}&\hskip-0.5cm&{\footnotesize {\rm strange (s)}}&{\footnotesize {\rm $0.5$}}&{\footnotesize {\rm $-\frac{1}{3}$}}&{\footnotesize {\rm $-1$}}&{\footnotesize {\rm $0$}}&{\footnotesize {\rm $0$}}&{\footnotesize {\rm $0$}}\\
\hline
{\footnotesize {\rm Planck temperature}}&{\footnotesize {\rm $\theta_{pl}=\frac{m_{pl}\, c^2}{K_B}=\sqrt{\frac{\hbar\, c^5}{G\, K_B^2}}=1.416 \times\, 10^{38}\, K$}}&\hskip-0.5cm&{\footnotesize {\rm top (t)}}&{\footnotesize {\rm $91.00$}}&{\footnotesize {\rm $+\frac{2}{3}$}}&{\footnotesize {\rm $0$}}&{\footnotesize {\rm $0$}}&{\footnotesize {\rm $0$}}&{\footnotesize {\rm $+1$}}\\
\hline
{\footnotesize {\rm $m_{e}=0.51\, \times\, 10^{-3}\, GeV$}}&{\footnotesize {\rm  $m_{p}=0.938\, GeV$}}&\hskip-0.5cm&{\footnotesize {\rm bottom (b)}}&{\footnotesize {\rm $4.8$}}&{\footnotesize {\rm $-\frac{1}{3}$}}&{\footnotesize {\rm $0$}}&{\footnotesize {\rm $0$}}&{\footnotesize {\rm $-1$}}&{\footnotesize {\rm $0$}}\\
\hline
\multicolumn {10}{l}{\rm \footnotesize Baryon number quarks $B=\frac{1}{3}$. Baryon number antiquarks $B=-\frac{1}{3}$. $\{S',C',B',T'\}=\{$strangeness,charm,bottomness,topness$\}$.}\\
\multicolumn {10}{l}{\rm \footnotesize Mass-energy in GeV. $m_e=$ electron mass-energy. $m_{p}=$ proton mass-energy.}\\
\end{tabular}$}\end{table}

$\bullet$\hskip 2pt An annihilation quark-antiquark, should be observed into an apparent quantum time $\sigma^2(\hat t)\approx 3.61\, \times\, 10^{-27},\, 1.097\, \times\, 10^{-24} \, s$. (Look the quark masses reported in Tab. \ref{table-quantum-time-and-strengths-fundamental-interactions}. The first refers to $t+\bar t$ and the last to $u+\bar u$.)

$\bullet$\hskip 2pt Similarly an annihilation electro-positron, should be observed into an apparent quantum time $\sigma^2(\hat t)\approx 6.4549\, \times\, 10^{-22}\, s$. (Look the electron mass reported in Tab. \ref{table-quantum-time-and-strengths-fundamental-interactions}.)
\end{example}

\begin{table}[h]
\caption{Phenomenological conservation laws in quantum reactions and symmetric reactions.}
\label{phenomenological-conservation-laws-in-quantum-reactions}
\scalebox{0.7}{$\begin{tabular}{|c|c|c|l|}
\hline
{\footnotesize{\rm Name}} &{\footnotesize{\rm Symbol}} &{\footnotesize{\rm Conserved}}&\hfil{\footnotesize{\rm Remark}}\hfil\\
\hline
\hline
{\footnotesize {\rm Baryon number}}&{\footnotesize {\rm $B$}}&{\footnotesize {\rm yes}}&{\footnotesize{\rm $B=\frac{1}{3}(n_q-n_{\bar q})\equiv \frac{1}{3}\triangle n_q$}}\\
&&&{\footnotesize{\rm ($n_q=$ number of quarks, $n_{\bar q}=$ number of antiquarks).}}\\
 \hline
{\footnotesize {\rm Electric charge}}&{\footnotesize{\rm $Q$}}&{\footnotesize {\rm yes}}&{\footnotesize {\rm $\frac{Q}{e}=I_3+\frac{Y}{2}$,  ($e=$ electron's elecetric charge) (See Theorem \ref{non-conservation-of-quantum-electric-charge}.)}}\\
&&&{\footnotesize {\rm (Gell-Mann-Nishijima formula). (See Theorem \ref{quantum-gell-mann-nishijima-formula}.)}}\\
 \hline
{\footnotesize {\rm Color charge}}&&{\footnotesize {\rm  yes}}&\\
\hline
{\footnotesize {\rm   Energy-mass}}&{\footnotesize {\rm $E=m\, c^2$}}&{\footnotesize{\rm  yes}}&{\footnotesize {\rm (See Theorem \ref{quantum-virtual-anomaly-massive-particles}.)}}\\
\hline
{\footnotesize {\rm   Isospin}}&{\footnotesize{\rm $I$}}&{\footnotesize{\rm  yes  in strong int.}}&\\
\hline
{\footnotesize {\rm   $3$-Isospin}}&{\footnotesize{\rm $I_3$}}&{\footnotesize{\rm  yes in strong int.}}&\\
\hline
{\footnotesize {\rm   Lepton number}}&{\footnotesize{\rm $L=L_e+L_\mu+L_\tau$}}&{\footnotesize{\rm  yes}}&{\footnotesize{\rm $L=n_l-n_{\bar l}$}}\\
&&&{\footnotesize{\rm ($n_l=$ number of leptons, $n_{\bar l}=$ number of antileptons).}}\\
{\footnotesize {\rm   Electronic lepton number}}&{\footnotesize{\rm $L_e$}}&{\footnotesize{\rm  yes}}&{\footnotesize{\rm $L_e=n_{e^-+\nu_e}-n_{e^++\bar\nu_e}$}}\\
&&&{\footnotesize{\rm ($n_{e^-+\nu_e}=$ number of electrons+electron-neutrinos)}}\\
&&&{\footnotesize{\rm ($n_{e^++\bar\nu_e}=$ number of positrons+electron-antineutrinos).}}\\
{\footnotesize {\rm   Muonic lepton number}}&{\footnotesize{\rm $L_\mu$}}&{\footnotesize{\rm  yes}}&{\footnotesize{\rm $L_\mu=n_{\mu+\nu_\mu}-n_{\bar\mu+\bar\nu_\mu}$}}\\
&&&{\footnotesize{\rm ($n_{\mu+\nu_\mu}=$ number of muons+muon-neutrinos)}}\\
&&&{\footnotesize{\rm ($n_{\bar\mu+\bar\nu_\mu}=$ number of antimuons+muon-antineutrinos).}}\\
{\footnotesize {\rm   Tauonic lepton number}}&{\footnotesize{\rm $L_\tau$}}&{\footnotesize{\rm  yes}}&{\footnotesize{\rm $L_\tau=n_{\tau+\nu_\tau}-n_{\bar\tau+\bar\nu_\tau}$}}\\
&&&{\footnotesize{\rm ($n_{\tau+\nu_\tau}=$ number of tauons+tau-neutrinos)}}\\
&&&{\footnotesize{\rm ($n_{\bar\tau+\bar\nu_\tau}=$ number of antitauons+tau-antineutrinos).}}\\
\hline
{\footnotesize {\rm   Strangeness number}}&{\footnotesize{\rm $S'$}}&{\footnotesize{\rm  yes in strong int.}}&{\footnotesize {\rm $\triangle S'=1$ in week int.}}\\
\hline
{\footnotesize {\rm   Hypercharge}}&{\footnotesize{\rm $Y$}}&{\footnotesize{\rm  yes in strong int.}}&{\footnotesize {\rm no in week int.}}\\
\hline
{\footnotesize {\rm   Crossing symmetry}}&&{\footnotesize{\rm  yes}}&{\footnotesize {\rm $n\to p+e^-+\bar\nu_e\, \Leftrightarrow\, \bar\nu_e+p\to n+e^+$}}\\
&&&{\footnotesize {\rm (neutron decay)$ \Leftrightarrow $(neutrino dedection).}}\\
\hline
{\footnotesize {\rm  CPT symmetry}}&{\footnotesize{\rm $CPT$}}&{\footnotesize{\rm  yes}}&{\footnotesize {\rm $CPT=$(charge-conjugation)(parity)(time-reversal)}}\\
\hline
\multicolumn {4}{l}{\rm\footnotesize Hypercharge: $Y=B+S'+C'+B'+T'$, flavour $\equiv(I_3,S',C',B',T')=$ ($3$-isospin,strangeness,charm,bottomness,topness).}\\
\multicolumn {4}{l}{\rm \footnotesize $I_3=\frac{1}{2}(\triangle n_u-\triangle n_d)=3$-Isospin is the $3^{th}$-component of isospin $I$. $S'=-\triangle n_s$. $C'=\triangle n_c$. $B'=-\triangle n_b$. $T'=\triangle n_t$.}\\
\multicolumn {4}{l}{\rm \footnotesize $Y=\frac{1}{3}[\triangle n_u+\triangle n_d-2(\triangle n_s+\triangle n_b)+4(\triangle n_c+\triangle n_t)]$, $\triangle n_q\equiv n_q-n_{\bar q}$.}\\
\multicolumn {4}{l}{\rm \footnotesize For complementary information see \href{http://en.wikipedia.org/wiki/List_of_particles}{Wikipedia - List of particles}.}\\
\end{tabular}$}\end{table}

$\bullet$\hskip 2pt In the following we characterize the concept of confinement for quantum systems.

\begin{definition}[(De)confined quantum system]\label{de-confined-quantum-system}

$\bullet$\hskip 2pt A solution $V$ of $(\widehat{YM})$ encodes a {\em confined system} if the spectrum $Sp(H(q))$ of the corresponding quantum Hamiltonian, has non-empty point spectrum: $Sp(H(q))_p\not=\varnothing$, for all $q\in V$. Then we say also that $V$ is a {\em confined quantum solution}.

$\bullet$\hskip 2pt When a solution $V$ of $(\widehat{YM})$ has non-empty point spectrum only in some subsets $N$ of $V$, then we say that $V$ encodes a {\em partially confined system} and that $N$ encodes the {\em deconfined subsystem}. Then we say also that $N$ is the {\em quantum deconfined part} of the solution $V$.

$\bullet$\hskip 2pt Whether a solution $V$ of $(\widehat{YM})$ is such that the following two conditions are satisfied:

{\em(i)}  $Sp(H(q))_p\not=\varnothing$, for all $q\in V$;

{\em(ii)}  $Sp(H(q))_c\not=\varnothing$, for all $q\in V$;

then we say that $V$ encodes a {\em confined system that can be deconfined}.
\end{definition}

\begin{theorem}[(De)confinement criterion ]\label{de-confinement-criterion}
Let us define the following set-mapping:
\begin{equation}\label{ker-delta-mapping}
\left\{
\begin{array}{l}
  \ker_{\triangle}:\widehat{A}\multimap A\\
 \ker_{\triangle}(\hat a)\equiv\bigcup_{\lambda\in Sp(\hat a)}\ker(\hat a-\lambda\, e).\\
\end{array}
\right.
\end{equation}
Then a solution $V$ of $(\widehat{YM})$ encodes a confined system (or a confined system that can be deconfined), iff $\ker_{\triangle}(H(q))\supset\{0\}\in A$, $\forall q\in V$.

Furthermore, there exists a quantum deconfined part $N\subset V$ of $V$, iff there exist points $q\in V$, such that $\ker_{\triangle}(H(q))=\{0\}\in A$.
\end{theorem}
\begin{proof}
Let us note that if $\lambda_1\not= \lambda_2\in Sp(\hat a)_p$, then $(\hat a-\lambda_i\, e)$, $i=1,2$, is not injective and $\ker(\hat a-\lambda_1\, e)\bigcap\ker(\hat a-\lambda_2\, e)=0\in A$. Instead, if $\lambda\in Sp(\hat a)_c\bigcup Sp(\hat a)_r$, then $(\hat a-\lambda\, e)$ is injective, hence $\ker(\hat a-\lambda\, e)=\{0\}\subset A$. Therefore, in the points $q\in V$ where $V$ encodes a confined quantum system, or a confined quantum system that can be deconfined, then necessarily $\ker_{\triangle}(H(q))\supset\{0\}\in A$. Instead, in the points $q\in V$, where $V$ encodes a deconfined quantum system, then necessarily must be $\ker_{\triangle}(H(q))=\{0\}\in A$.
\end{proof}

$\bullet$\hskip 2pt Let us conclude this paper by answering to the following question: {\em ``Does {\em quantum Majorana neutrino} exist ? "}.\footnote{This question arises from the paper {\em Symmetrical theory of electron and positron}, published by E. Majorana \cite{MAJORANA}, about a solution of the Dirac equation, where he first showed existence of a massive, electric-neutral, spin $\frac{1}{2}$, solution, that coincides with its anti-particle. This solution is now usually called {\em Majorana neutrino}. Existence of this particle should allow the so-called {\em neutrinoless double beta decay} of some nuclei. (See, e.g., this beautiful link \href{http://thy.phy.bnl.gov/~vogelsan/GGS/Wilkerson.pdf}{http://thy.phy.bnl.gov/~vogelsan/GGS/Wilkerson.pdf}.)}
\begin{definition}[Quantum Majorana neutrino]\label{quantum-majorana-neutrino}
We define {\em quantum Majorana neutrino} a quantum massive, electrically neutral fermionic particle, its own antiparticle, identified with a Cauchy data of the quantum super Yang-Mills PDE $\widehat{(YM)}$.
\end{definition}

\begin{theorem}[Existence of quantum Majorana neutrino]\label{existence-quantum-majorana-neutrino}
The quantum super Yang-Mills PDE $\widehat{(YM)}$ admits nonlinear quantum propagators $\widetilde{V}\subset\widehat{(YM)}$ such that $\partial\widetilde{V}=\widetilde{N}_0\sqcup \widetilde{P}\sqcup \widetilde{N}_1$, where $\widetilde{N}_1$ can be identified with  Majorana neutrino. Furthermore, the following propositions hold.

{\em(i)} A complex quantum quasi-particle exists that is a quasi neutralino, and that we call {\em quantum Majorana neutralino}. This contains two Majorana neutrinos other than two Higgsinos and two supersymmetric partner of the couple $(\nu_e,\bar\nu_e)$.

{\em(ii)} $\widetilde{V}$ is homeomorphic to a quantum $(4|4)$-superdisk with two attached super-handles.
\end{theorem}

\begin{figure}[h]
\centerline{\includegraphics[width=6cm]{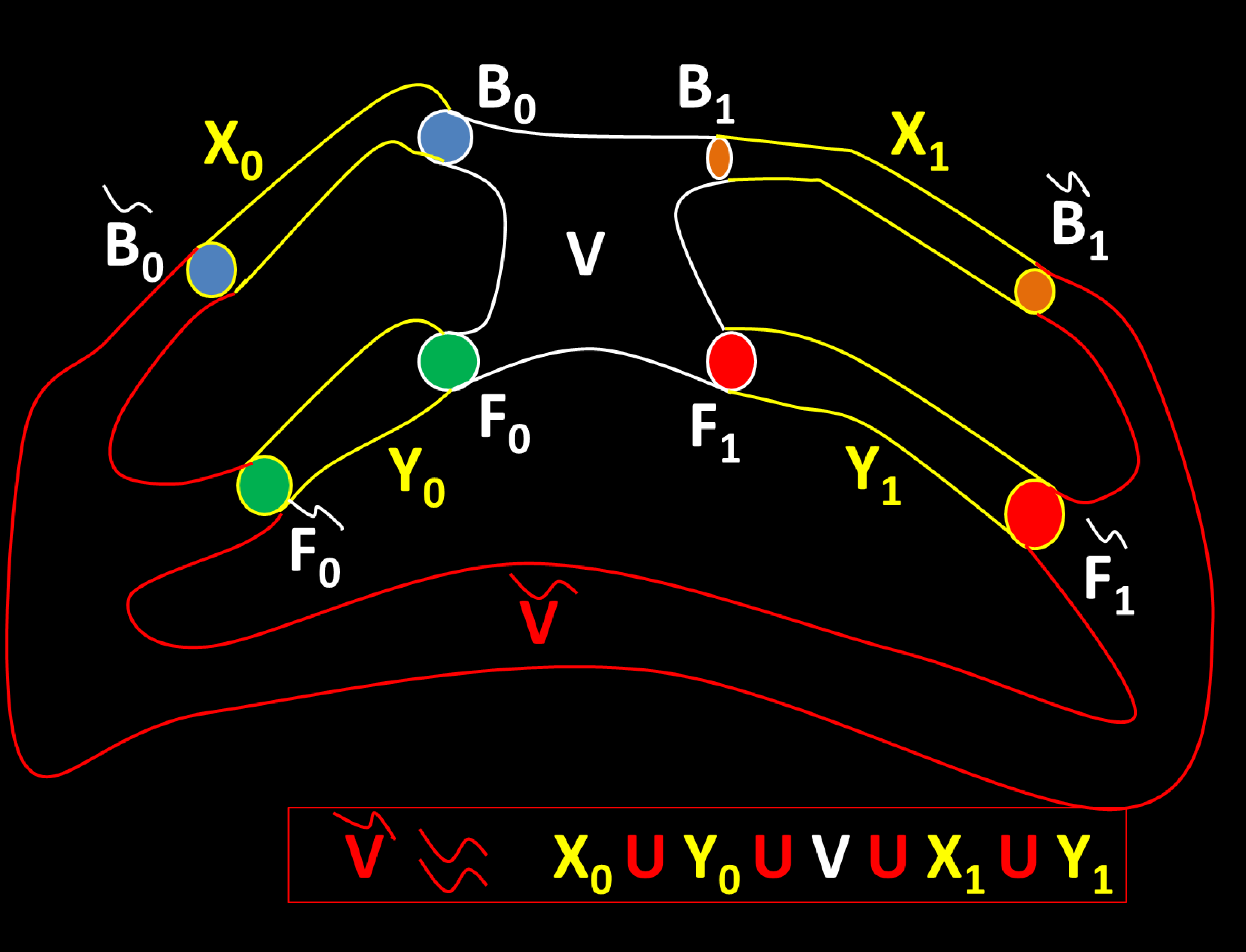}}
\renewcommand{\figurename}{Fig.}
\caption{Rapresentation of the relation between a propagator $V$ and its supersymmetric partner $\widetilde{V}$: $\widetilde{V}\thickapprox X_0\bigcup Y_0\bigcup V\bigcup X_1\bigcup Y_1$, such that $\partial X_i=\widetilde{B_i}\bigcup Q_i\bigcup B_i$, $\partial Y_i=\widetilde{F_i}\bigcup R_i\bigcup F_i$, $i=,0,\,1$, and $\partial V=N_0\bigcup P\bigcup N_1$, with $N_0=(B_0\sqcup F_0)$ and $N_1=(B_1\sqcup F_1)$. Here $B_i$ and $F_i$ are respectively bosonic and fermionic Cauchy data, and $\partial P=(\partial B_0\sqcup \partial F_0)\sqcup (\partial B_1\sqcup \partial F_1)$. Furthermore $\partial Q_i=\partial \widetilde{B_i}\sqcup \partial B_i$ and $\partial R_i=\partial \widetilde{F_i}\sqcup \partial F_i$.}
\label{relation-propagator-supersymmetric-partner}
\end{figure}
\begin{figure}[h]
\centerline{\includegraphics[width=6cm]{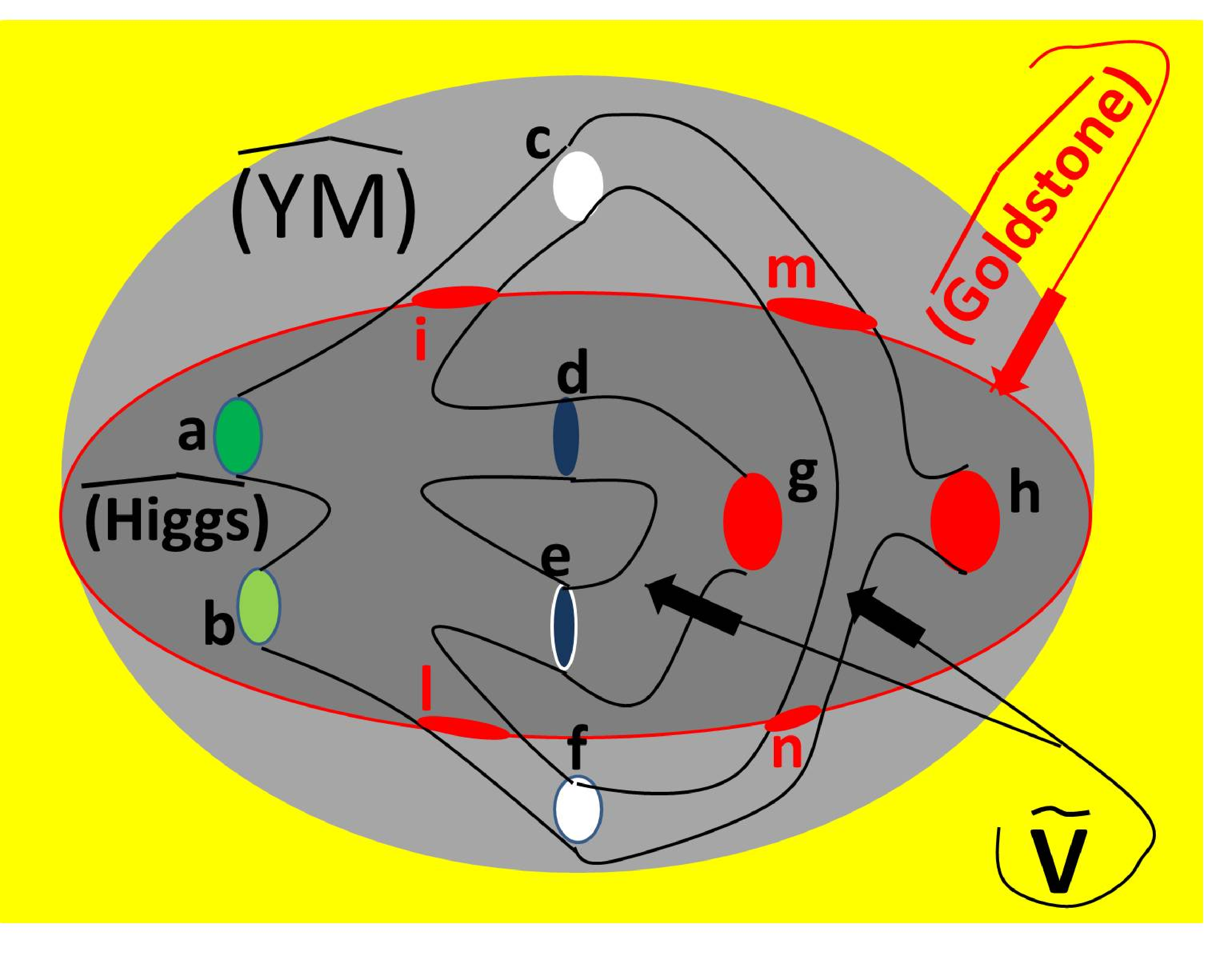}}
\renewcommand{\figurename}{Fig.}
\caption{Nonlinear quantum propagator representing a reaction generating Majorana neutrinos, identified with supersymmetric partners of the vector boson $Z^0$. In the figure one has put: $a=\widetilde{W^+}$, $b=\widetilde{W^-}$, $c=\widetilde{\nu_e}$, $d=\widetilde{e^+}$, $e=\widetilde{e^-}$, $f=\widetilde{\bar\nu_e}$, $h=\widetilde{Z_0}$, $g=\widetilde{Z_0}$. Furthermore, $i,\, l,\, m,\, n$ are the massive quasi-particles identified by means of the intersection of the nonlinear quantum propagator $\widetilde{V}$ with the Goldstone quantum partial differential relation $\widehat{(Golstone)}$, (represented by a red-circle in the picture). The nonlinear quantum propagator $V$, source supersymmetric partner of $\widetilde{V}$, can be similarly pictured.}
\label{figure-majorana-neutrino-reaction}
\end{figure}

\begin{proof}
Let $V\subset\widehat{(YM)}$ be a nonlinear quantum propagator in $\widehat{(YM)}$ such that
\begin{equation}\label{quantum-non-linear-propagator-and-symmetric-one-for-majorana-neutrino}
   \partial V=N_0\sqcup P\sqcup N_1
\end{equation}
with $N_0=W^+\sqcup W^-$ and $N_1=Z^0\sqcup Z^0$. This is possible since $V$ represents the reaction $W^++W^-\to Z^0+ Z^0$ that splits in the intermediate decays $W^+\to \nu_e+e^+$, $W^-\to \bar\nu_e+e^-$ and reactions $e^++e^-\to Z^0$,  $\nu_e+\bar\nu_e\to Z^0$. A representation by means of elementary bordisms is given in Fig. \ref{figure-majorana-neutrino-reaction}. Let us, now, consider the following lemmas.

\begin{lemma}[Boson and Fermion quantum super PDEs]\label{bosonic-fermionic-quantum-super-pds}
In the quantum super Yang-Mills PDE, $\widehat{(YM)}$ there exist two disjoint, open PDEs $\widehat{(Boson)}\subset \widehat{(YM)}$ and $\widehat{(Fermion)}\subset \widehat{(YM)}$, that we call respectively {\em Boson-PDE} and {\em Fermion-PDE}, that are formally integrable and completely integrable. Solutions of $\widehat{(Boson)}$, (resp. $\widehat{(Fermion)}$, are boson-polarized, (resp. fermion-polarized).
\end{lemma}
\begin{proof}
With respect to the definitions in Tab. \ref{local-quantum-spectral-spin-classification-observed-yang-mills-pde-solutions}, we can write $\widehat{A}=\mathfrak{b}\sqcup \mathfrak{f}\sqcup \mathfrak{n}$, where $\mathfrak{n}\equiv \widehat{A}\setminus\mathfrak{b}\sqcup \mathfrak{f}$.\footnote{Let us emphasize that $\mathfrak{b}$ is not empty since $\hat a\equiv \hbar^2\, s(s+1)\, id_A$, $s=n\in\{0,1,2,\cdots\}$, belong to $\mathfrak{b}$. Similarly $\mathfrak{f}$ is not empty since $\hat a\equiv \hbar^2\, s(s+1)\, id_A$, $s=\frac{2n+1}{2}$, $n\in\{0,1,2,\cdots\}$, belong to $\mathfrak{f}$.} We identify such disjoint subsets as equivalence classes of an equivalence relation $\mathcal{R}$ in $\widehat{A}$ and endow the quotient $\widehat{A}/\mathcal{R}$ with the quotient topology, by means of the projection $\pi_{\mathcal{R}}:\widehat{A}\to \widehat{A}/\mathcal{R}$. Then the mapping $[H]\equiv \pi_{\mathcal{R}}\circ  H:\widehat{(YM)}\to \widehat{A}/\mathcal{R}$ is continuous. Let us emphasize that $\widehat{A}/\mathcal{R}=\{[\mathfrak{b}],[\mathfrak{f}],[\mathfrak{n}]\}$ is a discrete topologic space, since $\mathfrak{b}$, $\mathfrak{f}$ and $\mathfrak{n}$ are disjoint subsets of $\widehat{A}$. Therefore $[\mathfrak{b}]$ and $[\mathfrak{f}]$ are open subsets of $\widehat{A}/\mathcal{R}$. As a by-product we get that $\widehat{(Boson)}\equiv[H]^{-1}([\mathfrak{b}])\subset \widehat{(YM)}$ and
$\widehat{(Fermion)}\equiv[H]^{-1}([\mathfrak{f}])\subset \widehat{(YM)}$, are open PDEs contained in $\widehat{(YM)}$, hence retain the formal geometric properties of the quantum super Yang-Mills PDE.
\end{proof}
\begin{lemma}[Supersymmetric partners of particles]\label{supersymmetric-partners-particles}
Cauchy data $N_0\subset \widehat{(Boson)}$, (resp. $N_1\subset \widehat{(Fermion)}$), are called {\em bosonic particles}, (resp. {\em fermionic particles}). All the bosonic particles (resp. fermionic particles), belonging to the same integral bordism class in $\widehat{(Boson)}$, (resp. $\widehat{(Fermion)}$), are considered {\em dynamically equivalent}.
If there exists a quantum nonlinear bordism $V$ of $\widehat{(YM)}$, such that $\partial V=N_0\sqcup P\sqcup N_1$, where $N_0$ is a bosonic particle and $N_1$ is a fermionic particle, and $P$ is an integral quantum supermanifold, such that $\partial P=\partial N_0\sqcup\partial N_1$, then we say that $(N_0,N_1)$ is a {\em couple of supersymmetric partners}.
The following propositions are equivalent.

{\em(i)} $(N_0,N_1)$ is a {\em couple of supersymmetric partners}.

{\em(ii)} $<\alpha,N_0>+<\alpha,N_1> =-<\alpha,P>$, for any quantum conservation law $\alpha$ of $\widehat{(YM)}$.\footnote{In the following we shall denote by $(N,\widetilde{N})$ a couple of supersymmetric partners.}

\end{lemma}
\begin{proof}
The equivalence of propositions (i) and (ii) follows from Theorem 3.6 in part I and Theorem 4.10 in \cite{PRAS19}(II).
\end{proof}

\begin{lemma}[Supersymmetric partners nonlinear quantum propagators]\label{supersymmetric-partners-quantum-nonlinear-propagators}
Let $N_0$ and $N_1$ be respectively initial and final Cauchy data in a quantum reaction in $\widehat{(YM)}$, such that $N_0\equiv b_1\sqcup\cdots\sqcup b_r\sqcup f_1\sqcup\cdots\sqcup f_s$ and $N_1\equiv \bar b_1\sqcup\cdots\sqcup \bar b_{\bar r}\sqcup \bar f_1\sqcup\cdots\sqcup\bar f_{\bar s}$, with $b_i$, $i=1,\cdots,r$, and $\bar b_{\bar i}$, $\bar i =1,\cdots,\bar r$, bosonic, and $f_j$, $j=1,\cdots,s$, $\bar f_{\bar j}$, $\bar j =1,\cdots,\bar s$, fermionic. Let $V$ be a nonlinear quantum propagator of $\widehat{(YM)}$ such that $\partial V=N_0\sqcup P\sqcup N_1$, and $\partial P=\partial N_0\sqcup \partial N_1$. Then there exists a nonlinear quantum propagator $\widetilde{V}$ of $\widehat{(YM)}$ such that $\partial \widetilde{V}=\widetilde{N_0}\sqcup \widetilde{P}\sqcup \widetilde{N_1}$, and $\partial \widetilde{P}=\partial \widetilde{N_0}\sqcup \partial\widetilde{N_1}$, where $(N_0,\widetilde{N_0})$ and $(N_1,\widetilde{N_1})$ are couples of supersymmetric partners. We call $\widetilde{V}$ (resp. $\widetilde{P}$) the {\em quantum supersymmetric partner} of $V$  (resp. $P$). Therefore to any quantum reaction there corresponds a {\em quantum supersymmetric partner reaction}.
\end{lemma}
\begin{proof}
This follows from Lemma \ref{supersymmetric-partners-particles} and from the commutative diagram (\ref{commutative-diagram-supersymmetric-partners-particles}), holding for any quantum conservation law $\alpha$ of $\widehat{(YM)}$.
\begin{equation}\label{commutative-diagram-supersymmetric-partners-particles}
   \xymatrix{<\alpha,N_0>+<\alpha,N_1>\ar@{=}[d]\ar@{=}[r]&-<\alpha,P>\ar@{=}[d]\\
   <\alpha,\widetilde{N_0}>+<\alpha,\widetilde{N_1}>\ar@{=}[r]&-<\alpha,\widetilde{P}>\\}
\end{equation}
In Fig. \ref{relation-propagator-supersymmetric-partner} is rapresented the relation between a propagator $V$ and its supersymmetric partner $\widetilde{V}$: $\widetilde{V}\thickapprox X_0\bigcup Y_0\bigcup V\bigcup X_1\bigcup Y_1$, such that $\partial X_i=\widetilde{B_i}\bigcup Q_i\bigcup B_i$, $\partial Y_i=\widetilde{F_i}\bigcup R_i\bigcup F_i$, $i=,0,\,1$, and $\partial V=N_0\bigcup P\bigcup N_1$, with $N_0=(B_0\sqcup F_0)$ and $N_1=(B_1\sqcup F_1)$. Here $B_i$ and $F_i$ are respectively bosonic and fermionic Cauchy data, and $\partial P=(\partial B_0\sqcup \partial F_0)\bigcup (\partial B_1\sqcup \partial F_1)$.
\end{proof}

Now from Lemma \ref{supersymmetric-partners-quantum-nonlinear-propagators} applied to nonlinear quantum propagator (\ref{quantum-non-linear-propagator-and-symmetric-one-for-majorana-neutrino}), we see that its supersymmetric partner $\widetilde{V}$, represents the following reaction $\widetilde{W}^++\widetilde{W}^-\to \widetilde{Z}^0+ \widetilde{Z}^0$. Here $\widetilde{Z}^0$ is just a massive neutral fermionic particle, its own antiparticle. Thus $\widetilde{Z}^0$ can be identified with a Majorana neutrino.\footnote{In the particle-zoo this is usually called {\em zino}.}  Therefore, we can conclude that Majorana neutrinos can be identified with possible products of quantum reactions in $\widehat{(YM)}$, and the proof of the first part of the theorem is down. In order to prove (i), let us emphasize that the nonlinear quantum propagator $\widetilde{V}$, identifies also a complex quasi-particle, that we call {\em quantum Majorana-neutralino}, $\nu_{Maj-neutralino}\equiv(g,h,i,l,m,n)$. Here $i$ and $m$, are the supersymmetric partners of massive neutrinos, $\nu_m$, $\bar\nu_m$, identified in Theorem \ref{existence-of-quantum-massive-neutrino}, and $m$ and $n$, are supersymmetric partners of Higgs quasi particles, (say {\em Higgsinos}), according to Theorem \ref{quantum-massive-photons-theorem} and Theorem \ref{a-quantum-massive-photons}. In other words the quantum super-Majorana-neutrino is a quasi {\em neutralino}. This last has two photinos (the supersymmetric partner of photons) instead of the couple $(\widetilde{\nu_e},\widetilde{\bar\nu_e})$. So it appears that the name is justified !

In order to complete the proof, i.e., to prove proposition (ii), let us remark that $\partial\widetilde{V}$ belongs to the same singular integral bordism class of a quantum exotic homotopy $(3|3)$-supersphere, $\hat \Sigma^{3|3}$, since $\Omega_{3|3,s}^{\widehat{(YM)}}=0$. However, $\widetilde{V}$ is not homeomorphic to $\hat D^{4|4}$, but to a $(4|4)$-superdisk with two attached super-handles. This means that the corresponding observed solution by means of a quantum relativistic frame $i:N\to M$, is necessarily a singular solution of $\widehat{(YM)}[i]$.\footnote{Let us emphasize that the nonlinear quantum propagators $V$ and $\widetilde{V}$, considered in this proof are not elementary ones in the sense of Theorem \ref{quantum-scattering-processes-in-quantum-yang-mills-d-4-quantum-super-minkowski}. In fact $V$ and $\widetilde{V}$ are not homeomorphic to $\hat D^{4|4}$.}
\end{proof}

